\pdfoutput=1
\documentclass[12pt]{amsart}
\usepackage{lmodern}
\usepackage{ifthen}
\usepackage{amsfonts}
\usepackage{amsxtra}
\usepackage{amssymb}
\usepackage{mathdots}
\usepackage{array}
\usepackage[margin=0.8in]{geometry}
\usepackage{xcolor}
\definecolor{cite}{rgb}{0.30,0.60,1.00}
\definecolor{url}{rgb}{0.00,0.00,0.80}
\definecolor{link}{rgb}{0.40,0.10,0.20}
\usepackage[colorlinks,linkcolor=link,urlcolor=url,citecolor=cite,pagebackref,breaklinks]{hyperref}
\usepackage{bbm}
\usepackage{mathtools}
\usepackage{mathrsfs}
\usepackage{appendix}
\usepackage[all]{xy}
\usepackage[lite,abbrev,msc-links,alphabetic]{amsrefs}

\usepackage{graphicx}
\usepackage{multirow}
\usepackage{pstricks}
\usepackage{pst-pdf}
\usepackage{enumitem}
\usepackage{pifont}
\usepackage{marvosym}
\usepackage{txfonts}

\usepackage[OT2,T1]{fontenc}
\DeclareSymbolFont{cyrletters}{OT2}{wncyr}{m}{n}
\DeclareMathSymbol{\Sha}{\mathalpha}{cyrletters}{"58}

\usepackage{graphicx}

\makeatletter
\providecommand*{\Dashv}{%
  \mathrel{%
    \mathpalette\@Dashv\vDash
  }%
}
\newcommand*{\@Dashv}[2]{%
  \reflectbox{$\m@th#1#2$}%
}
\makeatother

\setitemize[0]{leftmargin=20pt}
\setenumerate[0]{leftmargin=30pt}

\numberwithin{equation}{section}

\theoremstyle{plain}
\newtheorem{proposition}{Proposition}[section]
\newtheorem{conjecture}[proposition]{Conjecture}
\newtheorem{corollary}[proposition]{Corollary}
\newtheorem{lem}[proposition]{Lemma}
\newtheorem{theorem}[proposition]{Theorem}

\theoremstyle{definition}
\newtheorem{definition}[proposition]{Definition}
\newtheorem{construction}[proposition]{Construction}
\newtheorem{notation}[proposition]{Notation}
\newtheorem{assumption}[proposition]{Assumption}

\theoremstyle{remark}
\newtheorem{remark}[proposition]{Remark}

\renewcommand{\b}[1]{\mathbf{#1}}
\renewcommand{\c}[1]{\mathcal{#1}}
\renewcommand{\d}[1]{\mathbb{#1}}
\newcommand{\f}[1]{\mathfrak{#1}}
\renewcommand{\r}[1]{\mathrm{#1}}
\newcommand{\s}[1]{\mathscr{#1}}
\renewcommand{\sf}[1]{\mathsf{#1}}
\renewcommand{\(}{\left(}
\renewcommand{\)}{\right)}
\newcommand{\res}{\mathbin{|}}
\newcommand{\ang}[1]{\langle{#1}\rangle}
\newcommand{\ol}[1]{\overline{#1}{}}

\newcommand{\ul}{\underline}
\renewcommand{\leq}{\leqslant}
\renewcommand{\geq}{\geqslant}

\newcommand{\bA}{\b A}

\newcommand{\bw}{\b w}

\newcommand{\cA}{\c A}
\newcommand{\cB}{\c B}
\newcommand{\cC}{\c C}

\newcommand{\cF}{\c F}

\newcommand{\cV}{\c V}
\newcommand{\cW}{\c W}
\newcommand{\cX}{\c X}

\newcommand{\cZ}{\c Z}

\newcommand{\dC}{\d C}

\newcommand{\dE}{\d E}
\newcommand{\dF}{\d F}
\newcommand{\dG}{\d G}

\newcommand{\dL}{\d L}

\newcommand{\dN}{\d N}

\newcommand{\dQ}{\d Q}
\newcommand{\dR}{\d R}

\newcommand{\dZ}{\d Z}

\newcommand{\fC}{\f C}
\newcommand{\fD}{\f D}
\newcommand{\fE}{\f E}
\newcommand{\fF}{\f F}

\newcommand{\fI}{\f I}

\newcommand{\fK}{\f K}

\newcommand{\fM}{\f M}
\newcommand{\fN}{\f N}

\newcommand{\fS}{\f S}

\newcommand{\fU}{\f U}
\newcommand{\fV}{\f V}
\newcommand{\fW}{\f W}
\newcommand{\fX}{\f X}

\newcommand{\fZ}{\f Z}

\newcommand{\fg}{\f g}

\newcommand{\fm}{\f m}

\newcommand{\fp}{\f p}

\newcommand{\rD}{\r D}
\newcommand{\rE}{\r E}
\newcommand{\rF}{\r F}
\newcommand{\rG}{\r G}
\newcommand{\rH}{\r H}
\newcommand{\rI}{\r I}

\newcommand{\rM}{\r M}

\newcommand{\rR}{\r R}

\newcommand{\rU}{\r U}
\newcommand{\rV}{\r V}
\newcommand{\rW}{\r W}

\newcommand{\rd}{\,\r d}

\newcommand{\sA}{\s A}

\newcommand{\sO}{\s O}

\newcommand{\sS}{\s S}

\newcommand{\sV}{\s V}

\newcommand{\sfC}{\sf C}
\newcommand{\sfD}{\sf D}

\newcommand{\sfF}{\sf F}

\newcommand{\sfP}{\sf P}

\newcommand{\sfh}{\sf h}

\newcommand{\sfv}{\sf v}

\newcommand{\tF}{\mathtt{F}}

\newcommand{\tS}{\mathtt{S}}
\newcommand{\tT}{\mathtt{T}}

\newcommand{\tV}{\mathtt{V}}

\newcommand{\tc}{\mathtt{c}}

\newcommand{\te}{\mathtt{e}}

\renewcommand{\th}{\mathtt{h}}
\newcommand{\ti}{\mathtt{i}}

\newcommand{\ts}{\mathtt{s}}
\renewcommand{\tt}{\mathtt{t}}

\newcommand{\tw}{\mathtt{w}}

\newcommand{\bbu}{\boldsymbol{u}}
\newcommand{\bbv}{\boldsymbol{v}}
\newcommand{\bbz}{\boldsymbol{z}}
\newcommand{\bbZ}{\boldsymbol{Z}}

\newcommand{\llangle}{\left\langle}
\newcommand{\rrangle}{\right\rangle}

\newcommand{\pres}[2]{\prescript{#1}{}{#2}}

\newcommand{\ac}{\r{ac}}

\newcommand{\cusp}{\r{cusp}}

\newcommand{\der}{\r{der}}

\newcommand{\et}{{\acute{\r{e}}\r{t}}}

\newcommand{\fin}{\r{fin}}

\newcommand{\hs}{\r{hs}}
\newcommand{\Herm}{\r{Herm}}

\newcommand{\id}{\r{id}}

\newcommand{\hol}{\r{hol}}

\newcommand{\KM}{\r{FM}}

\newcommand{\num}{\r{num}}
\newcommand{\op}{\r{op}}
\newcommand{\proj}{\r{fl.proj}}

\newcommand{\RE}{\r{Re}\,}

\newcommand{\smp}{\r{sm.proj}}
\newcommand{\smq}{\r{sm.qproj}}

\newcommand{\st}{\r{st}}
\newcommand{\RS}{\r{RS}}
\newcommand{\HA}{\r{HA}}

\newcommand{\ev}{\r{ev}}
\newcommand{\coev}{\r{coev}}

\newcommand{\Sch}{\sf{Sch}}

\newcommand{\AV}{\sf{AV}}
\newcommand{\Cat}{\sf{Cat}}
\newcommand{\Mfd}{\sf{Mfd}}
\newcommand{\Set}{\sf{Set}}
\newcommand{\Hk}{\sf{Hk}}

\newcommand{\Mot}{\sf{Mot}}

\DeclareMathOperator{\Aut}{Aut}
\DeclareMathOperator{\BC}{BC}
\DeclareMathOperator{\CH}{CH}

\DeclareMathOperator{\Diff}{Diff}
\DeclareMathOperator{\Div}{Div}

\DeclareMathOperator{\End}{End}

\DeclareMathOperator{\Fil}{Fil}

\DeclareMathOperator{\Fun}{Fun}
\DeclareMathOperator{\Gal}{Gal}

\DeclareMathOperator{\GL}{GL}

\DeclareMathOperator{\GU}{GU}
\DeclareMathOperator{\Hom}{Hom}

\DeclareMathOperator{\Ker}{ker}

\DeclareMathOperator{\Lie}{Lie}

\DeclareMathOperator{\Mat}{Mat}

\DeclareMathOperator{\Nm}{Nm}

\DeclareMathOperator{\ord}{ord}

\DeclareMathOperator{\rank}{rank}
\DeclareMathOperator{\Res}{Res}

\DeclareMathOperator{\Sh}{Sh}
\DeclareMathOperator{\SC}{SC}

\DeclareMathOperator{\Spec}{Spec}
\DeclareMathOperator{\Spf}{Spf}

\DeclareMathOperator{\tr}{tr}
\DeclareMathOperator{\Tr}{Tr}

\DeclareMathOperator{\vol}{vol}

\begin{document}

\title[Motivic generating functions and $L$-derivatives]
{Motivic generating functions and $L$-derivatives of cohomological automorphic representations}

\author{Yifeng Liu}
\address{Institute for Advanced Study in Mathematics, Zhejiang University, Hangzhou 310058, China}
\email{liuyf0719@zju.edu.cn}

\date{\today}

\begin{abstract}
  Forty years ago, Gross and Zagier proved a remarkable formula computing the central derivative of the Rankin--Selberg $L$-function of a cusp form $f$ of weight $2$ and a character of the class group of an imaginary quadratic field, satisfying the so-called Heegner condition, and deduce that the nonvanishing of the central derivative implies the nonvanishing of the corresponding Mordell--Weil rank. This formula was later generalized in various contexts and degrees, for example, by S.~Zhang for $f$ of higher even weights; and by Yuan--Zhang--Zhang for $f$ a Hilbert cusp form of parallel weight $2$ over an arbitrary totally real field $F$ and $\chi$ a Dirichlet character of an arbitrary CM extension $E/F$, allowing arbitrary ramification. In this article, we prove a Gross--Zagier type formula in the context of doubling method, for $f$ a Hilbert cusp form over $F$ and $\chi$ a Hecke character of $E$ of general weights under a balancing condition, again allowing arbitrary ramification. Indeed, we even formulate a general conjecture which we call the motivic arithmetic inner product formula, for cohomological automorphic representations of $\GL_{2r}$ over $E$ of hermitian type of balanced weights.
\end{abstract}

\maketitle

\tableofcontents

\section{Introduction}

Height is one of the most important concepts in number theory, with earliest prototype dating back to the method of \emph{infinite descent} by Fermat. A century ago, Mordell systemically used the notion of height to prove the famous Mordell theorem on the finite generation of rational points on elliptic curves. The next highlight moment of height was Faltings' proof of the Mordell conjecture and the Gross--Zagier formula toward the Birch and Swinnerton-Dyer conjecture, both in the 1980s.

Recall that the Gross--Zagier formula \cite{GZ86} computes the central derivative of the Rankin--Selberg $L$-function of a cusp newform of weight $2$ and trivial Nebentypus and a ring class group character of an imaginary quadratic field, jointly under the so-called Heegner condition, through the height of the so-called Heegner point. Under the philosophy of this formula, one expects an arithmetic interpretation of the central derivative of $L(s,\sigma,\chi)$, where $\sigma$ is a cuspidal automorphic representation of $\GL_2(\bA_F)$ with $F$ a totally real field such that $\sigma_\infty$ is a discrete series, and $\chi$ is a Hecke character of a CM extension $E/F$, jointly satisfying
\begin{itemize}
  \item \emph{self-duality condition}: the central character of $\sigma$ is inverse to the restriction of $\chi$ to $\bA_F^\times$ (which would imply that $L(s,\sigma,\chi)$ is self-dual);

  \item \emph{odd-parity condition}: the root number of $L(s,\sigma,\chi)$ is $-1$ (so that the central vanishing order is odd);

  \item \emph{weight-balancing condition}: the weight of $\chi$ at every complex embedding of $E$ is at most the weight of $\sigma$ at the underlying real place of $F$ minus $2$.
\end{itemize}

When $\sigma$ has parallel weight $2$ (so that $\chi$ must be a Dirichlet character), such expectation has been achieved by Yuan--Zhang--Zhang \cite{YZZ13} unconditionally. In higher weight case, S.~Zhang \cite{Zha97} obtained a formula, still under the Heegner condition (so that $F=\dQ$), for $\sigma$ of arbitrary even weight and trivial central character and $\chi$ a (ring) class group character, through the Beilinson--Bloch height of the so-called Heegner cycle. However, unlike the weight $2$ case, Zhang's formula crucially replies on the Heegner condition as the construction involves the Kuga--Sato variety, which does not exist for other quaternionic Shimura curves.\footnote{For quaternionic Shimura curves over $\dQ$, one does have an analogue of the Kuga--Sato variety; but it remains the challenge to construct the Heegner cycle.} For this reason, there is almost no further progress in the case of general weights since Zhang's work, except for a recent preprint \cite{LS} where the authors allow the character to have nonzero weights (but still under the Heegner condition).

In this article, we will give an answer to the question of finding an arithmetic interpretation of the central derivative of $L(s,\sigma,\chi)$ in the above generality. Note that, it is possible that the weights of $\sigma$ may have different parities, in which case $\sigma$ is \emph{not} cohomological hence does not appear in the cohomology of quaternionic Shimura curves. For this reason, it is necessary to encode $\sigma$ and $\chi$ simultaneously, which leads naturally to the \emph{unitary} Shimura curves.

\subsection{Special cycles in motivic cohomology}

Let $E/F$ be a CM field. We develop a theory of special cycles and generating functions of general weights for incoherent totally positive definite hermitian spaces over $\bA_E/\bA_F$ (Definition \ref{de:space_coherence}), parallel to the original Kudla--Millson theory of the trivial weight for incoherent totally positive definite quadratic spaces over $\bA_F$ \cite{KM90}.

To reduce unnecessary technical complication, in the main part of the article we will assume that $E$ contains an imaginary quadratic field (see Remark \ref{re:general_cm} for general CM fields). From now on, we assume this and choose an imaginary quadratic subfield $E_0\subseteq E$.

Let $V$ be an incoherent totally positive definite hermitian space over $\bA_E/\bA_F$ of rank $m\geq 1$; and put $H\coloneqq\rU(V)$.

In \S\ref{ss:motive_standard}, we construct a projective system of smooth quasi-projective varieties $\{X_L\}_L$ over $E$ indexed by neat open compact subgroups $L$ of $H(\bA_F^\infty)$, whose base change along any given complex embedding of $E$ recovers usual Shimura varieties for certain unitary groups, together with a relative (Chow) motive $\rM_L$ over $X_L$ with coefficients in $E$ functorial in the level $L$, which we call the \emph{standard motive}. In \S\ref{ss:motive_automorphic}, we use $\rM_L$ to construct, for every hermitian $E$-weight $\lambda$ (Definition \ref{de:weight}), a relative motive $\rM_{L,\lambda}$ over $X_L$ with coefficients in $\dE^\lambda$ -- the field of definition of $\lambda$ as a subfield of $\dC$ -- functorial in the level $L$, which we call the \emph{automorphic motive} of weight $\lambda$. The motive $\rM_{L,\lambda}$ is an explicit motivic object underlying all coefficient systems of weight $\lambda$ in the following sense: for every prime $\ell$, the $\ell$-adic \'{e}tale realization of $\rM_{L,\lambda}$ is naturally the $\ell$-adic local system over $X_L$ of weight $\lambda$; and the similar property holds for Hodge realizations.

Consider a positive integer $r$ and a hermitian $E$-weight $\lambda$ that is $r$-balanced (Definition \ref{de:weight}). In \S\ref{ss:special_class} we define, for every $r$-tuple $x\in(V^\infty)^r$ with $T(x)\in\Herm_r(F)^+$ (Notation \ref{no:cm}), an element
\[
z^\lambda(x)_L\in\tS_{\lambda^+}(\dQ^r)^\vee\otimes\rH^{2r}(X_L,\rM_{L,\lambda}(r)),
\]
called \emph{motivic special class}, for every level $L$, in which $\tS_{\lambda^+}(\dQ^r)$ is the (algebraic) representation of $\GL_r(\dQ)^{\Hom(E,\dC)}$ of weight $\lambda^+$ (Notation \ref{no:weight}) and $\rH^\bullet(X_L,\rM_{L,\lambda}(r))$ denotes the motivic cohomology of the motive $\rM_{L,\lambda}(r)$ over $X_L$. Once constructing special classes, we can formulate generating functions. Denote by $G_r$ the quasi-split unitary group of rank $2r$ over $E/F$. For every Schwartz function $\phi^\infty\in\sS((V^\infty)^r)^L$ for some level $L$, we have the formal series
\[
Z^\lambda(\phi^\infty)_L\coloneqq\sum_{T\in\Herm_r(F)^+}Z^\lambda_T(\phi^\infty)_L
\]
as a function on $G_r(\bA_F)$ (or rather its metaplectic double cover when $m$ is odd) valued in the complex vector space $\tS_{\lambda^+}(\dC^r)^\vee\otimes_{\dE^\lambda}\rH^{2r}(X_L,\rM_{L,\lambda}(r))$, which we call the \emph{motivic generating function} (see Definition \ref{de:generating_function} for the precise definition). Note that when $\lambda=\b{0}$ is the trivial weight, $\tS_{\b{0}}(\dC^r)=\dC$, $\dE^{\b{0}}=\dQ$, and $\rM_{L,\b{0}}$ is the unit motive over $X_L$ so that $\rH^{2r}(X_L,\rM_{L,\b{0}}(r))=\CH^r(X_L)\otimes\dQ$; in this case, $z^{\b{0}}(x)_L$ recovers the usual special cycle and $Z^{\b{0}}(\phi^\infty)_L$ recovers the usual generating function, both in the hermitian case (see \cite{Liu11} or \cite{LL21}). Thus, it is natural to propose the following conjecture.

\begin{conjecture}[Modularity of motivic generating function, Conjecture \ref{co:modularity}]
Let $r$ be a positive integer such that $\lambda$ is $r$-balanced. For every level $L$ and every Schwartz function $\phi^\infty\in\sS((V^\infty)^r)^L$, the motivic generating function $Z^\lambda(\phi^\infty)_L$ is given by a holomorphic vector-valued hermitian Siegel modular form (of genus $r$ with respect to the CM extension $E/F$) for the representation $\tS_{\lambda^+}(\dC^r)^\vee\otimes\delta^{-m/2}$ (see Notation \ref{no:hermitian_siegel}(4) for $\delta$) with coefficients in $\rH^{2r}(X_L,\rM_{L,\lambda}(r))\otimes_{\dE^\lambda}\dC$.
\end{conjecture}

Unlike the case for the trivial weight, even for $r=1$ (or even for $r=1$ and $m=2$), we do not have a proof of the above conjecture. Nevertheless, assuming the above conjecture, one may use $Z^\lambda(\phi^\infty)_L$ as the theta kernel to define what we call \emph{motivic arithmetic theta lift}: For every form $\varphi\in\cA_{\lambda^+}^{m/2}(G_r)_\cusp^\hol$ (Definition \ref{de:hermitian_siegel}) and every $\phi^\infty\in\sS((V^\infty)^r)^L$, we define
\[
\Theta_{\phi^\infty}^\lambda(\varphi)_L\coloneqq\int_{G_r(F)\backslash G_r(\bA_F)}
\tr_\lambda\(Z^\lambda(\phi^\infty)_L(g)\otimes\ol{\varphi(g)}\)\rd g
\in \rH^{2r}(X_L,\rM_{L,\lambda}(r))\otimes_{\dE^\lambda}\dC,
\]
where $\tr_\lambda\colon\tS_\lambda(\dC^r)^\vee\otimes_\dC\ol{\tS_\lambda(\dC^r)^\vee}\to\dC$ is the natural pairing induced from the standard hermitian form on $\dC^r$. One of the main outcomes of \S\ref{ss:4} is to construct, when $m=2r$ and the system $\{X_L\}_L$ is proper, a level-independent height pairing
\[
\llangle\Theta_{\phi^\infty_1}^\lambda(\varphi_1),\Theta_{\phi^\infty_2}^\lambda(\varphi_2)\rrangle_\lambda
\]
between two motivic arithmetic theta lifts (see Remark \ref{re:height_intro} below), so that we may formulate the following conjecture.

\begin{conjecture}[Motivic arithmetic inner product formula, Conjecture \ref{co:maipf}]
Suppose that $m=2r$, that the hermitian $E$-weight $\lambda$ is $r$-balanced, and that the system $\{X_L\}_L$ is proper. Let $(\pi,\cA_\pi)$ be a holomorphic tempered cuspidal automorphic representation of $G_r(\bA_F)$ of weight $\lambda$ (Definition \ref{de:hermitian_representation}). For every pair $\ul\varphi=(\varphi_1,\varphi_2)$ of automorphic forms in $\(\tS_{\lambda^+}(\dC^r)^\vee\otimes_\dC\cA_\pi\)\cap\cA_{\lambda^+}^r(G_r)_\cusp^\hol$, and very pair $\ul\phi^\infty=(\phi^\infty_1,\phi^\infty_2)$ of Schwartz functions in $\sS((V^\infty)^r)$,
\[
\llangle\Theta_{\phi^\infty_1}^\lambda(\varphi_1),\Theta_{\phi^\infty_2}^\lambda(\varphi_2)\rrangle_\lambda
=\frac{L'(\tfrac{1}{2},\pi)}{b_{2r}(0)}\cdot C(\lambda)\cdot
\fZ^\natural_\fin(\ul\varphi,\phi^\infty_1\otimes\ol{\phi^\infty_2})
\]
holds. Here, $C(\lambda)\in\dR^\times$ is the real number of sign $(-1)^{r[F:\dQ]}$ from Proposition \ref{pr:ripf}, and $\fZ^\natural_\fin(\ul\varphi,\phi^\infty_1\otimes\ol{\phi^\infty_2})$ is the normalized doubling zeta integral introduced in Definition \ref{de:doubling_zeta}.
\end{conjecture}

\begin{remark}\label{re:height_intro}
Similar to Beilinson's height pairing on Chow groups \cite{Bei87} (see also \cite{LL21}*{\S4}), the height pairing on the motivic cohomology is actually defined on a certain subspace
\[
\rH^{2r}(X_L,\rM_{L,\lambda}(r))^{\ang{\ell}}\subseteq\rH^{2r}(X_L,\rM_{L,\lambda}(r))
\]
depending on the choice of a (sufficiently large) rational prime $\ell$, and takes values in $\dC\otimes\dQ_\ell$ a priori.

The above (global) height pairing can also be decomposed into a sum of local heights over all places of $E$ for certain representatives of classes, similar to the case of constant coefficients. Gladly, when $r=1$ and $m=2$ (so that $X_L$ is a curve), such local heights coincide with the one defined and studied by Brylinski in \cite{Bry89}. In some sense, our construction fulfills an expectation of Brylinski,\footnote{See the last sentence of \cite{Bry89}*{\S0}.} at least for unitary Shimura curves, that there should be a global theory overseeing these local heights.

Moreover, in the case where $r=1$ and $m=2$, the subspace $\rH^2(X_L,\rM_{L,\lambda}(1))^{\ang{\ell}}$ is independent of $\ell$ (which coincides with $\rH^2(X_L,\rM_{L,\lambda}(1))$ when $\lambda\neq\b{0}$) and the height pairing takes values in $\dC$ that is also independent of $\ell$ (Theorem \ref{th:independence}).
\end{remark}

\begin{remark}
The idea of formulating the motivic arithmetic inner product formula can also be used to formulate the \emph{arithmetic Gan--Gross--Prasad conjecture} \cite{GGP12}*{\S27} for cohomological cuspidal automorphic representations of $\GL_n(\bA_E)\times\GL_{n+1}(\bA_E)$ of general interlacing weights \cite{Liu25}*{Definition~4.1}.
\end{remark}

\subsection{Results for unitary groups of two variables}

In this section, we assume $r=1$ and $m=2$. When $\lambda=\b{0}$, both of the modularity conjecture and the (motivic) arithmetic inner product formula have been confirmed in \cites{Liu11,Liu12}.\footnote{Our argument in this article also simplifies the proof of the arithmetic inner product formula given in \cite{Liu12}; see Remark \ref{re:aipf}.}

Here, we assume $\lambda\neq\b{0}$. Though we cannot prove the modularity conjecture now hence cannot define the motivic arithmetic theta lifts unconditionally, we can prove a numerical analogue of the modularity conjecture in terms of the height pairing of Fourier coefficients of generating functions. Using such modularity, we can produce sufficiently many doubling kernel functions that encode central $L$-derivatives of holomorphic tempered cuspidal automorphic representations of $G_1(\bA_F)$ of weight $\lambda$ and root number $-1$, which would imply the motivic arithmetic inner product formula assuming the modularity conjecture.

Recall that for every $t\in F^+$, we have the $t$-th Fourier coefficient $Z^\lambda_t(\phi^\infty)_L$ as a function on $G_1(\bA_F)$ valued in $\rH^2(X_L,\rM_{L,\lambda}(1))\otimes_{\dE^\lambda}\dC$, among which $Z^\lambda_0(\phi^\infty)_L=0$ as $\lambda\neq\b{0}$. Consider a finite set $\Diamond$ of rational primes containing all those ramified in $E$ and an open compact subgroup $L^\Diamond$ of $H(\bA^{\infty,\Diamond}_F)$ that is the stabilizer of a self-dual lattice of $V^{\infty,\Diamond}$. For every nonempty subset $\Diamondblack$ of places of $F$ above $\Diamond$ containing all those nonsplit in $E$,\footnote{We also require that $\Diamondblack$ has cardinality at least $2$ if every place in $\Diamondblack$ splits in $E$.} we define a $H(F_\Diamond)$-stable subspace
\[
\sS((V^\infty)^2,\dL)^{L^\Diamond}_\Diamondblack\subseteq\sS((V^\infty)^2,\dL)^{L^\Diamond}
\]
in Definition \ref{de:schwartz_special}, and denote by $\sS((V^\infty)^2,\dL)^{L^\Diamond}_{(\Diamondblack)}$ the $G_1(\bA_F^\infty)^2$-stable (via the diagonal Weil representation) subspace generated by $\sS((V^\infty)^2,\dL)^{L^\Diamond}_\Diamondblack$. The subspace $\sS((V^\infty)^2,\dL)^{L^\Diamond}_\Diamondblack$ is large enough in the sense that for every generic irreducible representation $\pi^\infty$ of $G_1(\bA_F^\infty)$, the natural restriction map
\[
\Hom_{G_1(\bA_F^\infty)^2}\(\sS((V^\infty)^2,\dC)^{L^\Diamond},\pi^\infty\boxtimes(\pi^\infty)^\vee\)
\to
\Hom_{G_1(\bA_F^\infty)^2}\(\sS((V^\infty)^2,\dC)^{L^\Diamond}_{(\Diamondblack)},\pi^\infty\boxtimes(\pi^\infty)^\vee\)
\]
is an isomorphism, which has dimension at most one by the multiplicity one property for degenerate principal series.

We assume that the system $\{X_L\}_L$ is proper. As pointed out in Remark \ref{re:height_intro}, we have a \emph{height pairing}
\[
\rH^2(X_L,\rM_{L,\lambda}(1))\times \rH^2(X_L,\rM_{L,\lambda}(1))\to\dC
\]
for every level $L$ as a hermitian form with coefficients in $\dE^\lambda$.

\begin{theorem}[Theorem \ref{th:modularity}]
Assume that the system $\{X_L\}_L$ is proper. There is a unique $G_1(\bA_F^\infty)^2$-equivariant complex linear map
\[
\fI\colon\sS((V^\infty)^2,\dC)^{L^\Diamond}_{(\Diamondblack)}\to
\cA_{\lambda^+}(G_1)_\cusp^\hol\otimes_\dC\ol{\cA_{\lambda^+}(G_1)_\cusp^\hol}
\]
such that for every pair $\ul\phi^\infty=(\phi^\infty_1,\phi^\infty_2)$ of Schwartz functions in $\sS(V^\infty,\dC)^{L^\Diamond}$ satisfying $\phi^\infty_1\otimes\ol{\phi^\infty_2}\in\sS((V^\infty)^2,\dC)^{L^\Diamond}_\Diamondblack$, the identity
\[
\fI(\phi^\infty_1\otimes\ol{\phi^\infty_2})(g_1,g_2)=\vol^\natural(L)\cdot\sum_{t_1\in(F^\times)^+}\sum_{t_2\in(F^\times)^+}
\llangle Z^\lambda_{t_1}(\phi^\infty_1)_L(g_1),Z^\lambda_{t_2}(\phi^\infty_2)_L(g_2)\rrangle_{\lambda,L}
\]
holds with absolute convergence on the right-hand side for every $(g_1,g_2)\in G_1(\bA_F^\Diamondblack)^2$ and every neat open compact subgroup $L$ of $H(\bA_F^\infty))$ fixing both $\phi^\infty_1$ and $\phi^\infty_2$.
\end{theorem}

The theorem below provides a \emph{variant of the motivic arithmetic inner product formula}, showing that the functional $\fI$ is indeed a doubling kernel function for the central $L$-derivative.

\begin{theorem}[Theorem \ref{th:maipf}]\label{th:maipf_intro}
Assume that the system $\{X_L\}_L$ is proper. Let $(\pi,\cA_\pi)$ be a holomorphic tempered cuspidal automorphic representation of $G_1(\bA_F)$ of weight $\lambda$. For every pair $\ul\varphi=(\varphi_1,\varphi_2)$ of automorphic forms in $\cA_\pi\cap\cA_{\lambda^+}(G_1)_\cusp^\hol$ and every $\Phi^\infty\in\sS((V^\infty)^2,\dC)^{L^\Diamond}_{(\Diamondblack)}$,
\[
\iint_{(G_1(F)\backslash G_1(\bA_F))^2}
\fI(\Phi^\infty)(g_1,g_2)\cdot\ol{\varphi_1(g_1)}\varphi_2(g_2)\rd g_1\rd g_2
=\frac{L'(\tfrac{1}{2},\pi)}{b_2(0)}\cdot \(-\frac{\pi}{4}\)^{[F:\dQ]}\cdot
\fZ^\natural_\fin(\ul\varphi,\Phi^\infty)
\]
holds.\footnote{By Remark \ref{re:constant}(2), when $r=1$, $C(\lambda)=\(-\pi/4\)^{[F:\dQ]}$, which is independent of $\lambda$.}
\end{theorem}

By Lemma \ref{le:base_change}, the above theorem indeed fulfills our original objective of finding an arithmetic interpretation of the central derivative of $L(s,\sigma,\chi)$ at the beginning of the introduction by taking $\pi$ such that its base change to $\GL_2(\bA_E)$ is $\BC(\sigma)\otimes\chi$.

\begin{corollary}[Corollary \ref{co:maipf_bis}]
Let the situation be as in the above theorem. Assume the modularity of motivic generating function for $V$. Then the motivic arithmetic inner product holds for $(\pi,\cA_\pi)$.
\end{corollary}

\begin{remark}
Careful readers may notice that we exclude one case in our theory (for $r=1$ and $m=2$), namely, the case where the system $\{X_L\}_L$ is not proper, or equivalently, $F=\dQ$ and $V^\infty$ is split. We expect that our theory works in this case as well with extra work on the boundary extension. Note that, this case is closely related to the one studied in \cite{Zha97} and \cite{LS}, and our theory should be equivalent to the cycle construction in the two articles under their assumptions. We leave the details to the future.
\end{remark}

\subsection{Beilinson--Bloch conjecture for motivic cohomology}

Let $V$ be an incoherent totally positive definite hermitian space over $\bA_E/\bA_F$ of rank $2r$. For every hermitian $E$-weight $\lambda$, put
\[
\rH^{2r}(X,\rM_{\lambda}(r))\coloneqq\varinjlim_L\rH^{2r}(X_L,\rM_{L,\lambda}(r))
\]
as a module over $\dE^\lambda[H(\bA_F^\infty)]$.

The conjecture below can be regarded as an analogue of the Birch--Swinnerton-Dyer conjecture (or rather the Beilinson--Bloch conjecture) for motivic cohomology of automorphic motives.

\begin{conjecture}[Beilinson--Bloch, Conjecture \ref{co:beilinson}]
Suppose that the hermitian $E$-weight $\lambda$ is balanced. Suppose that the system $\{X_L\}_L$ is proper.
\begin{enumerate}
  \item The $\dE^\lambda[H(\bA_F^\infty)]$-module $\rH^{2r}(X,\rM_{\lambda}(r))$ is admissible and semisimple.

  \item For every holomorphic tempered cuspidal automorphic representation $(\pi,\cA_\pi)$ of $G_r(\bA_F)$ of weight $\lambda$, the relation
      \[
      \dim_\dC\Hom_{H(\bA_F^\infty)}\(\tilde\pi^\infty,\rH^{2r}(X,\rM_{\lambda}(r))\otimes_{\dE^\lambda}\dC\)=
      \ord_{s=1/2}L(s,\Pi_{\tilde\pi^\infty})
      \]
      holds for every irreducible admissible representation $\tilde\pi^\infty$ of $H(\bA_F^\infty)$ that is nearly equivalent to $\pi^\infty$. Here, $\Pi_{\tilde\pi^\infty}$ is the unique isobaric factor of the automorphic base change of $\pi$ to $E$ (as an isobaric automorphic representation of $\GL_{2r}(\bA_E)$) determined by $\tilde\pi^\infty$ (see Remark \ref{re:endoscopic}).
\end{enumerate}
\end{conjecture}

For $\pi$ in (2) above, there is, by the local theta dichotomy, a unique up to isomorphism totally positive definite hermitian space $V(\pi)$ over $\bA_E/\bA_F$ such that the theta lift of $\pi^\infty$ to $V(\pi)^\infty$ does not vanish.

Theorem \ref{th:maipf_intro} has a corollary that can be regarded as an evidence toward the above conjecture when $\ord_{s=1/2}L(s,\pi)=1$. Once again we assume $r=1$, and $\lambda$ nontrivial and $1$-balanced. For every intermediate field $\dE^\lambda\subseteq\dL\subseteq\dC$ and every level $L$, we define a subquotient $\ol\SC_\lambda(X_L,\dL)$ of $\rH^2(X_L,\rM_{L,\lambda})\otimes_{\dE^\lambda}\dL$ in Definition \ref{de:special_numerical}, which is, roughly speaking, the subspace generated by motivic special classes modulo the radical with respect to the height pairing.

\begin{corollary}[Corollary \ref{co:derivative}]
Let $(\pi,\cA_\pi)$ be a holomorphic tempered cuspidal automorphic representation of $G_1(\bA_F)$ of weight $\lambda$ such that $V\coloneqq V(\pi)$ is incoherent with the associated system $\{X_L\}_L$ proper. Let $\dE^\lambda\subseteq\dL_\pi\subseteq\dC$ be the field of definition of $\pi^\infty$. If $\ord_{s=1/2}L(s,\pi)=1$, then for every finite set $\Diamond$ of rational primes containing all those ramified in $E$ such that $\pi^{\infty,\Diamond}$ is spherical\footnote{See Corollary \ref{co:derivative} for the exact meaning of \emph{being spherical}.} and every open compact subgroup $L^\Diamond$ of $H(\bA^{\infty,\Diamond}_F)$ that is the stabilizer of a self-dual lattice of $V^{\infty,\Diamond}$,
\[
\varinjlim_{L_\Diamond}\ol\SC_\lambda(X_{L_\Diamond L^\Diamond},\dL_{\pi})_{\Ker\chi_\pi}\neq 0
\]
holds, where $\chi_\pi\colon\dL_\pi[L^\Diamond\backslash H(\bA_F^{\infty,\Diamond})/L^\Diamond]\to\dL_\pi$ is the Satake character of $\pi^{\infty,\Diamond}$.
\end{corollary}

\subsection{Structure of the article}

The article consists of five chapters after the introduction and one appendix.

In Chapter \ref{ss:2}, we collect some notions and facts from representation theory and automorphic forms. In \S\ref{ss:representation_unitary}, we review the construction of algebraic representations of unitary groups of given weights via a way that is symmetric in the standard representation and its dual. In \S\ref{ss:hermitian_space}, we apply the previous construction to hermitian spaces over CM fields and introduce further objects. In \S\ref{ss:cohomological_automorphic}, we review the notion of cohomological automorphic representations of $\GL_m(\bA_E)$ of hermitian type of given weights and their relation with Rankin--Selberg pairs $(\sigma,\chi)$ appearing at the beginning of the introduction when $m=2$. In \S\ref{ss:automorphic_unitary}, we collect the notation for automorphic forms on quasi-split unitary groups of even ranks. At last, we review the doubling method in the theory of theta lifting and deduce some formulae concerning derivatives of $L$-functions in \S\ref{ss:doubling}.

In Chapter \ref{ss:3}, we define several key concepts in our theory, including standard (relative) motives, automorphic (relative) motives, motivic special classes, motivic generating functions, and motivic arithmetic theta lifting. In \S\ref{ss:auxiliary_moduli}, we introduce a CM moduli problem that is auxiliary to our construction of the unitary moduli schemes in the next section. In \S\ref{ss:motive_standard}, we define the unitary moduli schemes over the (abstract) CM field $E$ and construct what we called the standard (relative) motive over those moduli schemes. In \S\ref{ss:speical_morphism}, we review the construction of special morphisms for unitary moduli schemes. In \S\ref{ss:motive_automorphic}, we construct automorphic (relative) motives over unitary moduli schemes of given weights, using the standard motive constructed previously. We compute the \'{e}tale realization and Hodge realization of those automorphic motives in \S\ref{ss:realization_etale} and \S\ref{ss:realization_hodge}, respectively. In \S\ref{ss:special_class}, we define motivic special classes and motivic generating functions, extending Kudla's special cycles and generating functions in the trivial weight case; we also formulate the conjecture on the modularity of motivic generating functions. At the end, we review the construction of canonical smooth integral models of unitary moduli schemes in the hyperspecial case in \S\ref{ss:integral_model}, which will be used later.

In Chapter \ref{ss:4}, we construct and study a height pairing on motivic cohomology of those automorphic motives constructed in the previous chapter. In \S\ref{ss:height_pairing}, we define such height pairing via Beilinson's height pairing on Chow groups (that is, motivic cohomology of the unit motive) with respect to a (nice) auxiliary rational prime $\ell$. In \S\ref{ss:height_decomposition}, we show that such height pairing has a decomposition into local heights for certain representatives of classes, similar to the height pairing on Chow groups. In \S\ref{ss:maipf}, we use the height pairing constructed previously to formulate the motivic arithmetic inner product formula assuming the modularity of generating functions. Finally in \S\ref{ss:independence}, we show that over unitary moduli curves, the height pairing on motivic cohomology is independent of $\ell$.

In Chapter \ref{ss:5}, we compute local heights of motivic special classes in various situations. In \S\ref{ss:height_inert}, we relate such local height at an inert place to the derivative of Whittaker coefficients of certain Eisenstein series in both smooth and semistable cases, by reducing the computation to the trivial weight case, which is already known by Li--Zhang \cite{LZ22}. In \S\ref{ss:height_split}, we show that over a split place, such local height vanishes just like in the trivial weight case. In \S\ref{ss:height_curve}, we state our result on local heights on unitary moduli curves at a general nonsplit nonarchimedean place, and prove it in \S\ref{ss:height_curve_split} and \S\ref{ss:height_curve_nonsplit}, respectively, according to whether the hermitian space at this place is split or not. In \S\ref{ss:height_curve_archimedean}, we collect our results on local heights on unitary moduli curves at an archimedean place from the previous work \cite{Liu26}. In the last section \S\ref{ss:abundance}, we identify a sufficiently large subspace of nonarchimedean Schwartz functions for which we can show the modularity of motivic doubling kernel functions later; we also reprove the arithmetic inner product formula (for the trivial weight) obtained in \cite{Liu12} by a simpler argument.

In Chapter \ref{ss:6}, we prove our main results over unitary moduli curves when the weight is balanced and nontrivial. In \S\ref{ss:modularity}, we show the modularity of motivic doubling kernel functions defined by height pairing of motivic special classes with respect to the Schwartz functions introduced in \S\ref{ss:abundance}. In \S\ref{ss:derivative}, we prove the formula computing derivatives of $L$-functions of cohomological automorphic representations of $\GL_2(\bA_E)$ of hermitian type and deduce its corollaries.

In Appendix \ref{ss:a}, we collect some known facts concerning relative motives. In \S\ref{ss:relative_motive}, we review the notion of relative motives or, more generally, motivic sheaves, together with their \'{e}tale and Hodge realizations. In \S\ref{ss:abelian_motive}, we review the motivic K\"{u}nneth decomposition of abelian schemes and identify its realizations.

\subsection{Notation and conventions}
\label{ss:notation}

\subsubsection*{Generalities:}

\begin{itemize}
  \item The imaginary unit is denoted by $\ti$; and the natural constant is denoted by $\te$.

  \item We denote by $\arg\colon\dC^\times\to\dC^\times$ the argument character.

  \item For $b\in\dN=\{0,1,2,\dots\}$, denote by $\fS_b$ the group of permutations on $\{1,\dots,b\}$. For a finite tuple $b=(b_u)_u$ of elements in $\dN$, put $\fS_b\coloneqq\prod_u\fS_{b_u}$.

  \item We use standard notations from the category theory. The category of sets is denoted by $\Set$. For a category $\sfC$, we denote by $\sfC^\op$ its opposite category. For another category $\sfD$, we denote by $\Fun(\sfC,\sfD)$ the category of functors from $\sfC$ to $\sfD$. In particular, we denote by $\sfP\sfC\coloneqq\Fun(\sfC^\op,\Set)$ the category of presheaves on $\sfC$, which contains $\sfC$ as a full subcategory by the Yoneda embedding.

  \item In a category with endofunctors like ``Tate twist'' $-(n)$ and ``degree shift'' $-[n]$ for $n\in\dZ$, we put $-\ang{n}\coloneqq-(n)[2n]$ as usual.

  \item All rings are commutative and unital; and ring homomorphisms preserve units. For a (topological) ring $L$, a (topological) $L$-ring is a (topological) ring $R$ together with a (continuous) ring homomorphism from $L$ to $R$. However, we use the word \emph{algebra} in the general sense, which is not necessarily commutative or unital.

  \item If a base ring is not specified in the tensor operation $\otimes$, then it is $\dZ$, unless otherwise noted.

  \item All hermitian spaces (over an \'{e}tale quadratic extension $\rE/\rF$ of rings) are assumed to be of finite rank and nondegenerate. For a hermitian space $V$, we usually denote by $(\;,\;)_V$ its hermitian form. When $2\in\rF^\times$, we have the moment matrix
      \[
      T(x)=\(\tfrac{1}{2}(x_i,x_j)_V\)_{1\leq i,j\leq r}
      \]
      for a tuple $x=(x_1,\dots,x_r)\in V^r$; and for an $r$-by-$r$ matrix $T$ with entries in $\rE$, put
      \[
      V^r_T\coloneqq\left\{x\in V^r\res T(x)=T\right\}.
      \]
\end{itemize}

\subsubsection*{Geometry:}

\begin{itemize}
  \item For a scheme $T$, we denote by $\Sch_{/T}$ the category of $T$-schemes, and by
      \[
      \Sch'_{/T},\quad \Sch^\proj_{/T},\quad \Sch^\smp_{/T},\quad  \Sch^\smq_{/T}
      \]
      the full subcategories of locally Noetherian, flat projective, smooth projective, and smooth quasi-projective schemes over $T$, respectively. When $T=\Spec R$ is affine, we also write $\Sch^?_{/R}$ for $\Sch^?_{/T}$.

  \item Denote by $\dG_a$ and $\dG_m$ the additive and multiplicative group schemes over $\Spec\dZ$, respectively.

  \item For every positive integer $r$, denote by $\Mat_r$ the scheme of $r$-by-$r$ matrices, $\GL_r\subseteq\Mat_r$ the subscheme of invertible ones, and $1_r\in\GL_r(\dZ)$ the unit matrix.

  \item Denote by $\Mfd$ the category of complex manifolds. For a field $k$ together with an embedding $u\colon k\to \dC$, we have a functor $-^u\colon\Sch^\smq_{/k}\to\Mfd$ sending $X$ to its underlying complex manifold.
\end{itemize}

\subsubsection*{Number Theory:}

Let $K$ be a number field.

\begin{itemize}
  \item Denote by $\fV_K$ the set of all places of $K$. For a finite set $\Diamond$ of places of $\dQ$, we denote by $\fV_K^{(\Diamond)}$ the subset of $\fV_K$ of places above $\Diamond$; when $\Diamond=\{\diamond_1,\diamond_2,\dots\}$ is written explicitly, we simply write $\fV_K^{(\diamond_1,\diamond_2,\dots)}$ instead of $\fV_K^{(\{\diamond_1,\diamond_2,\dots\})}$. Put $\fV_K^\fin\coloneqq\fV_K\setminus\fV_K^{(\infty)}$.

  \item For a nonarchimedean place $v$ of $K$, we write $\|v\|$ for the cardinality of the residue field of $K_v$ and $p(v)$ for the underlying rational prime of $v$.

  \item Denote by $\bA$ the ring of ad\`eles of $\dQ$, and put $\bA_K\coloneqq\bA\otimes K$.

  \item Denote by $\psi\colon\dQ\backslash\bA\to\dC^\times$ the standard additive automorphic character, that is, the unique character such that $\psi_\infty(x)=\te^{2\pi\ti x}$. Put $\psi_K\coloneqq\psi\circ\Tr_{K/\dQ}\colon K\backslash\bA_K\to\dC^\times$.

  \item We write
      \begin{itemize}
        \item $K_\diamond\coloneqq K\otimes\dQ_\diamond$ for a place $\diamond\in\{\infty,2,3,5,\dots\}$ of $\dQ$,

        \item $O_{K,(p)}\coloneqq O_K\otimes\dZ_{(p)}$ and $O_{K,p}\coloneqq O_K\otimes\dZ_p$ for a rational prime $p\in\{2,3,5,\dots\}$,

        \item more generally, for a finite set $\Diamond$ of rational primes, we put $O_{K,(\Diamond)}\coloneqq\bigcap_{p\in\Diamond}O_{K,(p)}$ (in $K$).
      \end{itemize}

  \item When $K$ is a CM field, we put $\ol{M}\coloneqq M\otimes_{K,\ol{\phantom{a}}}K$ for an object $M$ in a $K$-linear category, where $\ol{\phantom{a}}$ denotes the complex conjugation.
\end{itemize}

\begin{notation}\label{no:contraction}
Consider a number field $\dL$ \emph{contained in} $\dC$. For an $\dL_\infty$-linear object $\rV$ (e.g., vector space, category, ...), we denote by
\[
\kappa\colon \rV\to\rV^\flat
\]
the natural projection map, where $\rV^\flat$ is the $\dC$-linear factor of $\rV$ corresponding to the natural contraction of coefficients $\dL_\infty=\dL\otimes\dR\to\dC$.
\end{notation}

\begin{definition}\label{de:hecke}
Let $K$ be a number field. For a reductive group $H$ over $K$ and a finite set $\Diamondblack$ of places of $K$ including all archimedean ones, we define $\Hk(H(\bA_K^\Diamondblack))$ to be the category
\begin{itemize}
  \item whose objects are neat open compact subgroups $L$ of $H(\bA_K^\Diamondblack)$; and

  \item in which a morphism from $L$ to $L'$ is an element $h\in L\backslash H(\bA_K^\Diamondblack)/L'$ satisfying $h^{-1}L h\subseteq L'$.
\end{itemize}
\end{definition}

\begin{definition}\label{de:dualizable}
Let $\sfC$ be a symmetric monoidal category with product $\otimes$ and unit object $\mathbbm{1}$. Recall that an object $A$ of $\sfC$ is \emph{dualizable} if there exist another object $A^\vee\in\sfC$ together with two maps
\[
\ev_A\colon A\otimes A^\vee\to\mathbbm{1},\quad
\coev_A\colon \mathbbm{1}\to A^\vee\otimes A,
\]
called the \emph{evaluation} and \emph{coevaluation maps}, respectively, satisfying that the composite maps
\[
A \xrightarrow{1_A\otimes\coev_A} A\otimes A^\vee\otimes A \xrightarrow{\ev_A\otimes 1_A} A,\quad
A^\vee \xrightarrow{\coev_A\otimes 1_{A^\vee}} A^\vee\otimes A\otimes A^\vee \xrightarrow{1_{A^\vee}\otimes\ev_A} A^\vee
\]
are both identity maps. The data $A^\vee,\ev_A,\coev_A$, if exist, are unique up to unique isomorphisms. When $\sfC$ is closed, a map $A\otimes A^\vee\to\mathbbm{1}$ is an evaluation map if and only if the induced maps
\[
A\to\ul\Hom(A^\vee,\mathbbm{1}),\quad A^\vee\to\ul\Hom(A,\mathbbm{1})
\]
are both isomorphisms.
\end{definition}

\begin{definition}\label{de:category_section}
Denote by $\Cat$ the (2,1)-category of categories. Consider an (ordinary) category $\sfC$ and a (2,1)-functor $\sfF\colon\sfC^\op\to\Cat$. For every subcategory $\sfC'$ of $\sfC$,
\begin{itemize}
  \item we denote by $\sfF(\sfC')$ the category whose objects are pairs $(c,A_c)$ with $c$ an object in $\sfC'$ and $A_c$ an object in $\sfF(c)$ and morphisms from $(c,A_c)$ to $(c',A_{c'})$ are pairs $(f,A_f)$ with $f\colon c\to c'$ an edge in $\sfC'$ and $A_f\colon A_c\to \sfF(f)A_{c'}$ an edge in $\sfF(c)$;\footnote{The category $\sfF(\sfC')$ is known as the Grothendieck construction (or unstraightening) of the functor $\sfF\res_{\sfC'}$.}

  \item we say that an edge $(f,A_f)$ in $\sfF(\sfC')$ is \emph{Cartesian} if $A_f$ is an isomorphism;

  \item we say that a functor $\sfD\to\sfF(\sfC')$ is \emph{Cartesian} if the image of every edge is Cartesian.
\end{itemize}
\end{definition}

\subsubsection*{Acknowledgements}

The author would like to thank Fangzhou~Jin, Andreas~Mihatsch, and Lei~Wu for helpful discussion. The author is partially supported by the National Key R\&D Program of China No.~2022YFA1005300.

\subsubsection*{AI disclosure}

AI was not used in this work.

\section{Preliminaries on representation theory}
\label{ss:2}

In this chapter, we collect some notions and facts from representation theory and automorphic forms. Let $r$ and $m$ be positive integers.

\subsection{Algebraic representations of unitary groups}
\label{ss:representation_unitary}

We start from the notation of weights.

\begin{notation}\label{no:weight}
Let $m\geq 0$ be an integer.
\begin{enumerate}
     \item Denote by $\fW_m$ the set of tuple $\lambda=(b_1,b_2,\dots,b_m)$ of non-increasing integers. For $\lambda\in\fW_m$, put $b(\lambda)\coloneqq\sum_{i=1}^m b_i$. We denote the zero tuple by $\b{0}$ (when the length $m$ is clear).

     \item Denote by $\ol{\phantom{a}}\colon\fW_m\to\fW_m$ the involution sending $\lambda$ to $\ol\lambda\coloneqq(-b_m,\dots,-b_2,-b_1)$.

     \item For $\lambda\in\fW_m$, denote by $0\leq s=s(\lambda)\leq m$ the largest integer such that $b_s>0$ ($s=0$ if $b_1\leq 0$).

     \item Denote by $\fW_m^+$ the subset of $\fW_m$ consisting of $\lambda$ such that $b_m\geq 0$. We will adopt the convention that an element $\lambda\in\fW_m^+$ is naturally an element in $\fW_r^+$ for every integer $r\geq m$ by adding zeros at the end.

     \item For $\lambda\in\fW_m$, put $\lambda^+\coloneqq(b_1,\dots,b_s)\in\fW_s^+$ with $s=s(\lambda)$.
\end{enumerate}
\end{notation}

Let $\rE/\rF$ be an \'{e}tale quadratic extension of rings with $\ol{\phantom{a}}\in\Aut(\rE/\rF)$ the nontrivial involution. We first recall the construction of algebraic representations of general linear groups.

\begin{construction}\label{co:representation_general}
Consider a free $\rF$-module $U$ of rank $r$ and an element $\lambda\in\fW_r^+$ such that $b!$ is invertible in $\rF$, where $b\coloneqq b(\lambda)$. Let $\pi_\lambda\in\dZ[1/b!][\fS_b]$ be the idempotent corresponding to the standard filling that assigns $1,2,\dots,b$ from left to right starting with the first row then moving to the second row, etc (see, for example, \cite{FH91}*{\S4.2}), which gives rise to a direct summand
\[
\tS_\lambda(U)\coloneqq\pi_\lambda\(\tT^b(U)\)
\]
of $\tT^b(U)\coloneqq U^{\otimes b}$ of $\rF[\GL_\rF(U)]$-modules. In particular, when $\rF$ is a field of characteristic zero, $\tS_\lambda(U)$ is an irreducible representation of $\GL_\rF(U)\simeq\GL_r(\rF)$ of highest weight $\lambda$.

In general, for a finite set $\fU$ and a tuple $\lambda=(\lambda_u)_u\in(\fW_r^+)^\fU$ such that $b_u!$ ($b_u\coloneqq b(\lambda_u)$) is invertible in $\rF$ for every $u\in\fU$, we have
\[
\tT^b(U)\coloneqq\bigotimes_{u\in\fU}\tT^{b_u}(U)
\]
(tensor over $\rF$) and similarly $\tS_\lambda(U)$, together with the projector
\[
\pi_\lambda\coloneqq\bigotimes_{u\in\fU}\pi_{\lambda_u}.
\]
\end{construction}

We then recall the construction of algebraic representations of unitary groups.

\begin{construction}\label{co:representation_unitary}
Consider a (nondegenerate) hermitian space $V$ over $\rE/\rF$ of rank $m$ and an $\rF$-ring $\dL$ such that the set $\fU$ of $\rF$-ring homomorphisms from $\rE$ to $\dL$ is nonempty (hence consists of two elements).

Consider a pair $\lambda\in\fW_m^\fU$ satisfying $\lambda_{\ol{u}}=\ol{\lambda_u}$ for every $u\in\fU$. Put $b=(b_u\coloneqq b(\lambda^+_u))_u\in\dN^\fU$ and
\[
B_\lambda\coloneqq m^{\min_{u\in\fU}\{b_u\}}\cdot\prod_{u\in\fU}b_u!.
\]
Assume that $B_\lambda$ is invertible in $\dL$. Put
\[
\tT_b(V)\coloneqq\bigotimes_{u\in\fU}\tT^{b_u}(V_u)
\]
(tensor over $\dL)$ as an algebraic representation of $\rU(V)$ over $\dL$, where $\rU(V)$ acts on $V_u\coloneqq V\otimes_{E,u}\dL$ via the homomorphism $u$. For every pair $I=(i_u)_u$ with $1\leq i_u\leq b_u$, we have the contraction operator
\[
C_I\colon\tT_b(V)\to\tT_{b-1}(V)
\]
by ``pairing'' the $i_u$-th factor of $\tT_{b_u}(V_u)$ for the two elements $u\in\fU$, and the expansion operator
\[
A_I\colon\tT_{b-1}(V)\to\tT_b(V)
\]
by inserting the coevaluation map $\dL\to\bigotimes_{u\in\fU}V_u$ (Definition \ref{de:dualizable}) to the $i_u$-th spot of $\tT^{b_u}(V_u)$ for $u\in\fU$. If we put
\[
\tT_{[b]}(V)\coloneqq\bigcap_{I}\Ker C_I,
\]
then
\begin{align}\label{eq:representation_unitary}
\tT_b(V)=\tT_{[b]}(V)\oplus\(\bigoplus_{r=1}^{\min_{u\in\fU}\{b_u\}}\tT_{[b]}(V)_{b-r}\),
\end{align}
where
\[
\tT_{[b]}(V)_{b-r}\coloneqq\sum A_{I_1}\circ\cdots\circ A_{I_r}\(\tT_{[b-r]}(V)\).
\]
We call $\tT_{[b]}(V)$ the \emph{harmonic direct summand} of $\tT_b(V)$.

It is clear that the decomposition \eqref{eq:representation_unitary} is stable under the action of both $\rU(V)$ and $\fS_b$. In particular, we may define
\[
\tS_{[\lambda]}(V)\coloneqq(\otimes_{u\in\fU}\pi_{\lambda_u^+})\(\tT_{[b]}(V)\)
\]
(Construction \ref{co:representation_general} with $r$ replaced by $m$), which is a direct summand of the $\dL[\rU(V)]$-module $\tT_{[b]}(V)$ hence $\tT^{b}(V)$.
\end{construction}

\begin{remark}\label{re:representation_unitary}
Suppose that $\dL$ is a field of characteristic zero. Take an element $u\in\fU$. Put $V_u\coloneqq V\otimes_{\rE,u}\dL$ and similarly for $V_{\ol{u}}$.
\begin{enumerate}
  \item The element $u\in\fU$ induces a homomorphism $\rU(V)\to\GL(V_u)\simeq\GL_m(\dL)$; and it is well-known that $\tS_{[\lambda]}(V)$ extends by scalar to an irreducible representation of $\GL_m(\dL)$ of highest weight $\lambda_u$.

  \item The hermitian form induces an isomorphism $V_{\ol{u}}\simeq V_u^\vee$. The map $(\wedge^m V_u)\otimes_\dL V_u^\vee\to\wedge^{m-1}V_u$ given by pairing the last two factors induces an isomorphism $V_u^\vee\simeq (\wedge^{m-1}V_u)\otimes_\dL(\wedge^mV_u)^{\vee}$, which further induces an isomorphism
      \[
      \tT_b(V_u)\simeq \tT^{b_u}(V_u)\otimes_\dL\tT^{b_{\ol{u}}}(\wedge^{m-1}V_u)\otimes_\dL(\wedge^mV_u)^{-b_{\ol{u}}}.
      \]
      In particular, $\tT_b(V_u)$ is canonically a subspace of $\tT^{b_u+(m-1)b_{\ol{u}}}(V_u)\otimes_\dL(\wedge^mV_u)^{-b_{\ol{u}}}$, in which $\fS_{b_u+(m-1)b_{\ol{u}}}$ acts on the latter via the first tensor factor. There exists a unique idempotent $\delta_{\lambda_u,\lambda_{\ol{u}}}\in\dQ[\fS_{b_u+(m-1)b_{\ol{u}}}]$, depending only on the \emph{ordered} pair $(\lambda_u,\lambda_{\ol{u}})$, such that $\tT_{[\lambda]}(V)$ is the kernel of $\delta_{\lambda_u,\lambda_{\ol{u}}}$.
\end{enumerate}

\end{remark}

\begin{construction}\label{co:special_map}
Let the situation be as in Construction \ref{co:representation_unitary} and $r$ a positive integer satisfying $r\geq s(\lambda_u)$ for every $u\in\fU$. For every tuple $x=(x_1,\dots,x_r)\in V^r$, we construct a map
\[
\mu_{[\lambda]}(x)\colon\tS_{[\lambda]}(V)^\vee\to\tS_{\lambda^+}(\dL^r)^\vee
\]
(notation from Construction \ref{co:representation_general}) as follows. Dualizing the $r$ components of $x$ gives an $\rE$-linear map
\[
\mu(x)\colon V^\vee\to(\rE^\vee)^r,
\]
so that we have the induced map $\mu(x)_u\colon V_u^\vee\to(\dL^\vee)^r$ for each $u\in\fU$. Taking tensor products, we have the $\dL$-linear map
\[
\mu_b(x)\colon\tT_b(V)^\vee\to \tT^b((\dL^\vee)^r)=\tT^b(\dL^r)^\vee.
\]
Composing with the dual of the natural map $\tT_b(V)\to\tT_{[b]}(V)$, we obtain the map
\[
\mu_{[b]}(x)\colon\tT_{[b]}(V)^\vee\to \tT^b(\dL^r)^\vee,
\]
which is $\tS_b$-equivariant. The desired map $\mu_{[\lambda]}(x)$ is obtained by applying $\pi_{\lambda^+}$ on both sides.
\end{construction}

\subsection{Hermitian spaces over CM extensions}
\label{ss:hermitian_space}

From now to the end of this article, we fix a CM extension $E/F$ of number fields.

\begin{notation}\label{no:cm}
The following notation will be used throughout the article.
\begin{enumerate}
  \item Denote by $\dE\subset\dC$ the normal closure of $E$ in $\dC$, so that $\Hom(E,\dC)=\Hom(E,\dE)$.

  \item Denote by $D_E\in\dN$ the absolute discriminant of $E$.

  \item Denote by $\chi_{E/F}\colon F^\times\backslash\bA_F^\times\to\{\pm 1\}$ the quadratic character associated with the extension $E/F$. For every place $v$ of $F$, define a local $L$-factor
      \[
      b_{m,v}(s)\coloneqq\prod_{i=1}^m L(2s+i,(\chi_{E/F})_v^i).
      \]

  \item For every place $v$ of $F$, put $E_v\coloneqq E\otimes_FF_v$ and $O_{E_v}\coloneqq O_{F_v}\otimes_{O_F}O_E$.

  \item For every place $v$ of $F$, define
      \[
      d(v)\coloneqq
      \begin{dcases}
      1,&\text{$v$ is split in $E$;}\\
      2,&\text{$v$ is nonsplit in $E$.}
      \end{dcases}
      \]
      For every finite set $\Diamondblack$ of places of $F$, put $d(\Diamondblack)\coloneqq\sum_{v\in\Diamondblack}d(v)$.

  \item For every prime $v$ of $F$, denote by $\fp_v$ the maximal ideal of $O_{F_v}$.

  \item Denote by $\Herm_r\subseteq\Res_{O_E/O_F}\Mat_r$ the $O_F$-scheme of $r$-by-$r$ hermitian matrices with respect to $O_E/O_F$. Put $\Herm^\circ_r\coloneqq\Herm_r\cap\Res_{O_E/O_F}\GL_r$.

  \item Denote by $\Herm_r(F)^+\subseteq\Herm_r(F)$ the subset of totally semi-positive definite hermitian matrices, and put $\Herm_r^\circ(F)^+\coloneqq\Herm_r(F)^+\cap\Herm_r^\circ(F)$.
\end{enumerate}
\end{notation}

\begin{definition}\label{de:weight}
An \emph{$E$-weight} of rank $m$ is a tuple $\lambda=(\lambda_u)_u\in\fW_m^{\Hom(E,\dC)}$.
\begin{enumerate}
  \item We say that an $E$-weight $\lambda$ is \emph{hermitian} (with respect to $E/F$) if $\lambda_{\ol{u}}=\ol{\lambda_u}$ holds for every $u\in\Hom(E,\dC)$.

  \item We say that an $E$-weight $\lambda$ is \emph{$r$-balanced} (for an integer $r$) if $s(\lambda_u)\leq r$ for every $u\in\Hom(E,\dC)$, and is \emph{balanced} if it is $\lfloor m/2\rfloor$-balanced where $m$ is the rank of $\lambda$.

  \item We denote by $\dE^\lambda\subseteq\dC$ field of definition of $\lambda$, that is, the smallest subfield of $\dC$ such that $\Aut(\dC/\dE^\lambda)$ acts trivially on $\lambda$ with respect to the natural action of $\Aut(\dC/\dQ)$ on the set of $E$-weights of rank $m$; $\dE^\lambda$ is indeed contained in $\dE$.
\end{enumerate}
In what follows, we write $\lambda_v$ for the pair $(\lambda_u)_{v(u)=v}$ for every $v\in\fV_F^{(\infty)}$.
\end{definition}

Take a hermitian $E$-weight $\lambda$ of rank $m$ that is $r$-balanced. Put
\begin{align}\label{eq:bad_lambda}
B_\lambda\coloneqq \prod_{v\in\fV_F^{(\infty)}}B_{\lambda_v}
\end{align}
(see Construction \ref{co:representation_unitary}).

\begin{construction}\label{co:harmonic_function}
Consider a $\dZ[(D_E B_\lambda)^{-1}]$-ring $R$ and a hermitian space $V$ over $O_E\otimes R/O_F\otimes R$ of rank $m$. In (2--4) below, we further assume that $[\dE:\dE^\lambda]$ is invertible on $R$.
\begin{enumerate}
  \item Put
      \[
      \tS_{[\lambda]}(V)\coloneqq\bigotimes_{v\in\fV_F^{(\infty)}}\tS_{[\lambda_v]}(V)
      \]
      (tensor over $O_\dE\otimes R$) as an $(O_\dE\otimes R)[\rU(V)(O_F\otimes R)]$-module, where $\tS_{[\lambda_v]}(V)$ is from Construction \ref{co:representation_unitary} with $\rE/\rF=O_E\otimes R/O_F\otimes R$ and $\dL=O_\dE\otimes R$ (as an $O_F\otimes R$-ring via $v$).

  \item For every $\gamma\in\Gal(\dE/\dQ)$, its action on $O_\dE$ induces an isomorphism
      \[
      \gamma_{[\lambda_v]}\colon\tS_{[\lambda_v]}(V)\otimes_{O_\dE,\gamma}O_\dE\xrightarrow\sim\tS_{[(\gamma\lambda)_{\gamma v}]}(V)
      \]
      for every $v\in\fV_F^{(\infty)}$, which further induces an isomorphism
      \[
      \gamma_{[\lambda]}\colon\tS_{[\lambda]}(V)\otimes_{O_\dE,\gamma}O_\dE\xrightarrow\sim\tS_{[\gamma\lambda]}(V)
      \]
      by taking tensor product over all $v\in\fV_F^{(\infty)}$. In particular, if $\gamma\in\Gal(\dE/\dE^\lambda)$, then $\gamma_{[\lambda]}$ is an $O_{\dE^\lambda}$-linear automorphism of $\tS_{[\lambda]}(V)$, from which we obtain an $O_{\dE^\lambda}$-linear action of $\Gal(\dE/\dE^\lambda)$ on $\tS_{[\lambda]}(V)$. Define
      \[
      V_\lambda\coloneqq\tS_{[\lambda]}(V)^{\Gal(\dE/\dE^\lambda)}
      \]
      as an $(O_{\dE^\lambda}\otimes R)[\rU(V)(O_F\otimes R)]$-module, so that $\tS_{[\lambda]}(V)=V_\lambda\otimes_{O_{\dE^\lambda}}O_\dE$.

  \item Take an element $x\in V^r$. We have the map
      \[
      \mu_{[\lambda]}(x)\colon\tS_{[\lambda]}(V)^\vee\to\tS_{\lambda^+}((O_\dE\otimes R)^r)^\vee
      \]
      by taking the tensor product of $\mu_{[\lambda_v]}(x)$ over all $v\in\fV_F^{(\infty)}$, where
      \[
      \mu_{[\lambda_v]}(x)\colon\tS_{[\lambda_v]}(V)^\vee\to\tS_{\lambda_v^+}((O_\dE\otimes R)^r)^\vee
      \]
      is the map from Construction \ref{co:special_map} with $\rE/\rF=O_E\otimes R/O_F\otimes R$ and $\dL=O_\dE\otimes R$. It is straightforward to show the map $\mu_{[\lambda]}(x)$ is $\Gal(\dE/\dE^\lambda)$-equivariant, hence restricts to a map
      \[
      \mu_\lambda(x)\colon(V_\lambda)^\vee\to\tS_{\lambda^+}((O_{\dE^\lambda}\otimes R)^r)^\vee.
      \]
      Finally, denote by
      \[
      \mu^\lambda(x)\colon\tS_{\lambda^+}((O_{\dE^\lambda}\otimes R)^r)\to V_\lambda
      \]
      the dual of $\mu_\lambda(x)$.

  \item We define a function $\Psi_{V^{2r}}^\lambda$ on $V^{2r}$ valued in $\tS_{\lambda^+}((O_{\dE^\lambda}\otimes R)^r)^\vee\otimes_{O_{\dE^\lambda}\otimes R}
      \ol{\tS_{\lambda^+}((O_{\dE^\lambda}\otimes R)^r)^\vee}$ such that for every $\ul{x}=(x_1,x_2)\in V^{2r}=(V^r)^2$,
      \begin{align*}
      \Psi_{V^{2r}}^\lambda(\ul{x})
      =\ol{\mu_\lambda(x_2)}\circ\mu^\lambda(x_1)
      &\in\Hom_{O_{\dE^\lambda}\otimes R}
      \(\tS_{\lambda^+}((O_{\dE^\lambda}\otimes R)^r),\ol{\tS_{\lambda^+}((O_{\dE^\lambda}\otimes R)^r)^\vee}\) \\
      &=\tS_{\lambda^+}((O_{\dE^\lambda}\otimes R)^r)^\vee\otimes_{O_{\dE^\lambda}\otimes R}
      \ol{\tS_{\lambda^+}((O_{\dE^\lambda}\otimes R)^r)^\vee},
      \end{align*}
      where we have identified $V_\lambda$ with $\ol{(V_\lambda)^\vee}$ using the hermitian form on $V$.
\end{enumerate}
\end{construction}

We recall the following definition, which will be used in latter chapters.

\begin{definition}\label{de:lattice}
Let $v$ be a nonarchimedean place of $F$ and $V$ a hermitian space over $E_v/F_v$. For an $O_{E_v}$-lattice $\Lambda$ of $V$, put
\[
\Lambda^\vee\coloneqq\left\{ x\in V\res \Tr_{E_v/F_v}(x,y)_V\in O_{F_v}\text{ for }\forall y\in\Lambda\right\}.
\]
We say that
\begin{enumerate}
  \item $\Lambda$ is an \emph{integral lattice} if $\Lambda\subseteq\Lambda^\vee$;

  \item $\Lambda$ is a \emph{vertex lattice} if $\Lambda^\vee/\Lambda$ is annihilated by a uniformizer of $E_v$ (and we define the index $t(\Lambda)$ of $\Lambda$ to be the length of $\Lambda^\vee/\Lambda$ as an $O_{E_v}$-module);

  \item $\Lambda$ is an \emph{self-dual} or \emph{almost self-dual} lattice if it is a vertex lattice with $t(\Lambda)$ equal to $0$ or $1$.
\end{enumerate}
\end{definition}

\subsection{Cohomological automorphic representations}
\label{ss:cohomological_automorphic}

\begin{definition}\label{de:hermitian_automorphic}
For a hermitian $E$-weight $\lambda$ of rank $m$, we define $\HA(\lambda)$ to be the set of isomorphism classes of automorphic representations $\Pi$ of $\GL_m(\bA_E)$ such that
\begin{enumerate}
  \item $\Pi$ is an isobaric sum of mutually non-isomorphic conjugate self-dual cuspidal automorphic representations;

  \item for every $u\in\Hom(E,\dC)$, $\Pi_{v(u)}$, regarded as an irreducible representation of $\GL_m(\dC)$ via $u$ is isomorphic to the principal series $\boxplus_{i=1}^m\arg^{2b_{ui}+m+1-2i}$, where $\lambda_u=(b_{u1},\dots,b_{um})$.
\end{enumerate}
\end{definition}

We relate the set $\HA(\lambda)$ to the classical Rankin--Selberg setting when $m=2$.

\begin{definition}\label{de:rs_pair}
We introduce the following concepts.
\begin{enumerate}
  \item A \emph{Rankin--Selberg pair} is a pair $(\sigma,\chi)$ in which $\sigma$ is a cuspidal automorphic representation of $\GL_2(\bA_F)$ and $\chi$ is an automorphic character of $\bA_E^\times$, satisfying $\chi_\sigma\cdot\chi\res_{\bA_F^\times}=1$, where $\chi_\sigma$ denotes the central character of $\sigma$.

  \item For a Rankin--Selberg pair $(\sigma,\chi)$, we may associate it with a root number $\epsilon(\sigma,\chi)\in\{\pm 1\}$. We say that $(\sigma,\chi)$ is \emph{coherent/incoherent} if $\epsilon(\sigma,\chi)=+1/{-1}$.

  \item A \emph{Rankin--Selberg weight} is a pair $(k,l)$ in which $k=(k_v)_v\in\dN^{\fV_F^{(\infty)}}$ and $l=(l_u)_u\in\dZ^{\Hom(E,\dC)}$ satisfying $l_u+l_{\ol{u}}=0$ and $l_u\equiv k_{v(u)}\mod 2$ for every $u\in\Hom(E,\dC)$.

  \item We say that a Rankin--Selberg weight $(k,l)$ is \emph{balanced} if $|l_u|\leq k_{v(u)}$ for every $u\in\Hom(E,\dC)$.

  \item For a Rankin--Selberg weight $(k,l)$, denote by $\RS(k,l)$ the set of isomorphism classes of Rankin--Selberg pairs $(\sigma,\chi)$ satisfying that
      \begin{enumerate}
         \item for every $v\in\fV_F^{(\infty)}$, $\sigma_v$ is the (unitary) discrete series representation of $\GL_2(\dR)$ of weight $k_v+2$; and

         \item for every $u\in\Hom(E,\dC)$, $\chi_u$, regarded as a character of $\dC^\times$ via $u$, is $\arg^{l_u}$.
      \end{enumerate}
\end{enumerate}
\end{definition}

For a Rankin--Selberg weight $(k,l)$, define an element $\lambda_{(k,l)}=(b_{u1},b_{u2})_u\in\fW_2^{\Hom(E,\dC)}$ by the formula
\[
b_{u1}=\frac{k_{v(u)}+l_u}{2},\quad b_{u2}=\frac{-k_{v(u)}+l_u}{2}
\]
for $u\in\Hom(E,\dC)$. It is clear that the assignment $(k,l)\mapsto\lambda_{(k,l)}$ defines a bijection between Rankin--Selberg weights and hermitian $E$-weights of rank $2$, under which the two balanced subsets coincide.

\begin{lem}\label{le:base_change}
The assignment $(\sigma,\chi)\mapsto \BC(\sigma)\otimes\chi$ defines a map $\RS(k,l)\to\HA(\lambda_{(k,l)})$, which is surjective. Here, $\BC$ denotes the quadratic automorphic base change from $F$ to $E$.
\end{lem}

\begin{proof}
We start from showing that $\Pi\coloneqq\BC(\sigma)\otimes\chi$ indeed belongs to $\HA(\lambda_{(k,l)})$. Part (2) of Definition \ref{de:hermitian_automorphic} is clear. It remains to check (1). For every automorphic representation $\Sigma$ over $E$, denote by $\Sigma^-$ its conjugation with respect to $E/F$. There are two cases.

First, $\BC(\sigma)$ remains cuspidal so that $\Pi$ is cuspidal. Then
\begin{align*}
\Pi^\vee&\simeq\BC(\sigma^\vee)\otimes\chi^{-1}\simeq\BC(\sigma\otimes\chi_\sigma^{-1})\otimes\chi^{-1}
\simeq\BC(\sigma)\otimes(\chi_\sigma^{-1}\circ\Nm_{E/F})\otimes\chi^{-1} \\
&\simeq\BC(\sigma)\otimes(\chi\circ\Nm_{E/F})\otimes\chi^{-1}
\simeq\BC(\sigma)\otimes\chi^-
\simeq\BC(\sigma)^-\otimes\chi^-\simeq\Pi^-.
\end{align*}

Second, $\BC(\sigma)$ is the isobaric sum of two (distinct) automorphic characters (hence so is $\Pi$). Then we must have $\BC(\sigma)\simeq\tilde\chi\boxplus\tilde\chi^-$ for some automorphic character $\tilde\chi$ of $\bA_E^\times$, hence $\omega_\sigma=\tilde\chi\res_{\bA_F^\times}\cdot\chi_{E/F}$. It follows that $\tilde\chi\cdot\chi\res_{\bA_F^\times}=\tilde\chi^-\cdot\chi\res_{\bA_F^\times}=\chi_{E/F}$, which implies that $(\tilde\chi\cdot\chi)^\vee\simeq(\tilde\chi\cdot\chi)^-$ and $(\tilde\chi^-\cdot\chi)^\vee\simeq(\tilde\chi^-\cdot\chi)^-$.

Then we show the surjectivity. By twisting an automorphic character of $E$ satisfying Definition \ref{de:rs_pair}(5b), we may assume without loss of generality that $l=\b{0}$ (that is, $l_u=0$ for every $u\in\Hom(E,\dC)$). Take an element $\Pi\in\HA(\lambda_{(k,\b{0})})$. Again, there are two cases.

First, $\Pi$ is cuspidal. The central character $\chi_\Pi$ is a Dirichlet character of $E$ satisfying $\chi_\Pi\cdot\chi_\Pi^-=1$. Thus, one may find a Dirichlet character $\chi$ of $E$ such that $\chi_\Pi=\chi/\chi^-$, which implies that $(\Pi\otimes\chi^{-1})^-\simeq\Pi\otimes\chi^{-1}$. Thus, we may find a cuspidal automorphic representation $\sigma$ of $\GL_2(\bA_F)$ such that $\BC(\sigma)\simeq\Pi\otimes\chi^{-1}$, or equivalently, $\Pi\simeq\BC(\sigma)\otimes\chi$. It is clear that $(\sigma,\chi)\in\RS(k,\b{0})$.

Second, $\Pi\simeq\chi_1\boxplus\chi_2$ is an isobaric sum of two distinct conjugate self-dual automorphic characters. Then the quotient $\chi_1/\chi_2^-$ is a Dirichlet character of $E$ satisfying $(\chi_1/\chi_2^-)(\chi_1/\chi_2^-)^-=1$. Thus, one may find a Dirichlet character $\chi$ of $E$ such that $\chi_1/\chi_2^-=\chi/\chi^-$. Put $\tilde\chi\coloneqq\chi_1\chi^{-1}$, so that $\Pi\simeq(\tilde\chi\boxplus\tilde\chi^-)\otimes\chi$. As $\chi_1\neq\chi_2$, one can find a cuspidal automorphic representation $\sigma$ of $\GL_2(\bA_F)$ such that $\BC(\pi)\simeq\tilde\chi\boxplus\tilde\chi^-$, or equivalently, $\Pi\simeq\BC(\sigma)\otimes\chi$. It is clear that $(\sigma,\chi)\in\RS(k,\b{0})$.

The lemma is proved.
\end{proof}

\subsection{Automorphic forms for unitary groups}
\label{ss:automorphic_unitary}

Let $W_r$ be the standard skew-hermitian space over $E/F$ of dimension $2r$, whose skew-hermitian form is given by the matrix \[
\tw_r\coloneqq\begin{pmatrix}  & 1_r \\ -1_r &  \end{pmatrix}.
\]
Put $G_r\coloneqq\rU(W_r)$ as a reductive group over $F$, equipped with homomorphisms
\begin{align*}
m&\colon\Res_{E/F}\GL_r\to G_r,\qquad  a\mapsto
\begin{pmatrix} a  &  \\  &\pres{t}{\ol{a}}^{-1}   \end{pmatrix}, \\
n&\colon\Herm_r\to G_r,\qquad   b\mapsto
\begin{pmatrix} 1_r  & b \\  & 1_r   \end{pmatrix}
\end{align*}
of algebraic groups over $F$, whose image we denote by $M_r$ and $N_r$, respectively.

\begin{notation}\label{no:hermitian_siegel}
We will frequently use the following notation related to $G_r$.
\begin{enumerate}
  \item For every $v\in\fV_F^\fin$, we denote by $K_{r,v}$ the maximal open compact subgroup of $G_r(F_v)$ that is the stabilizer of $O_{E_v}^{2r}$.

  \item For every $v\in\fV_F^{(\infty)}$, the map
      \[
      (U_1,U_2)\mapsto u^{-1}\frac{1}{2}
      \begin{pmatrix}
      U_1+U_2   & -\ti U_1+\ti U_2 \\
      \ti U_1-\ti U_2   & U_1+U_2 \\
      \end{pmatrix},
      \]
      where $u\in\Hom(E,\dC)$ is an arbitrary element above $v$, identifies $\rU(r)\times\rU(r)$ with a maximal compact subgroup $K_{r,v}$ of $G_r(F_v)$ that is independent of the choice of $u$. Moreover, the assignment $(U_1,U_2)\mapsto\det(U_1 U_2^{-1})$ defines a character $\delta_v\colon K_{r,v}\to\dC^\times$, which is also independent of the choice of $u$.

  \item For every finite subset $S\subseteq\fV_F$, denote by $\widetilde{G_r(F_S)}$ the metaplectic double cover of $G_r(F_S)$. Globally, denote by $\widetilde{G_r(\bA_F)}$ the metaplectic double cover of $G_r(\bA_F)$.

  \item Put
      \[
      K_{r,\infty}\coloneqq\prod_{v\in\fV_F^{(\infty)}}K_{r,v},\quad
      \delta\coloneqq\prod_{v\in\fV_F^{(\infty)}}\delta_v\colon K_{r,\infty}\to\dC^\times.
      \]
      Denote by
      \[
      \delta^{1/2}\colon\widetilde{K_{r,\infty}}\to\dC^\times
      \]
      pthe pullback of $\det$ along the map $\dC^\times\to\dC^\times$ sending $z$ to $z^2$, so that $\widetilde{K_{r,\infty}}$ is naturally a maximal compact subgroup of $\widetilde{G_r(F_\infty)}$.

  \item For every tuple $b=(b_u)_u\in\dN^{\Hom(E,\dC)}$, we have the representation $\tT^b(\dC^r)$ (Construction \ref{co:representation_general}) of $\fS_b\times K_{r,\infty}$, where the action of $K_{r,\infty}$ is the unique one such that for every $k\in K_{r,\infty}$ and every $u\in\Hom(E,\dC)$, if we write $u(k)=(U_1,U_2)\in\rU(r)\times\rU(r)$, then $U_1$ acts on $\tT^{b_u}(\dC^r)$ via itself as an element of $\GL_r(\dC)$ and $U_2$ acts on $\tT^{b_{\ol{u}}}(\dC^r)$ via the element $\ol{U_2}$ of $\GL_r(\dC)$.

  \item For every $v\in\fV_F$, we define a Haar measure $\r{d}g_v$ on $G_r(F_v)$ as follows:
    \begin{itemize}
      \item for $v\in\fV_F^\fin$, $\r{d}g_v$ is the Haar measure under which $K_{r,v}$ has volume $1$;

      \item for $v\in\fV_F^{(\infty)}$, $\r{d}g_v$ is the product of the measure on $K_v$ of total volume $1$ and the standard hyperbolic measure on $G(F_v)/K_v$ (see, for example, \cite{EL24}*{\S2.1}).
    \end{itemize}
    Put $\r{d}g\coloneqq\prod_{v}\r{d}g_v$, which is a Haar measure on $G_r(\bA_F)$.

  \item Denote by $\cA(G_r)$ the space of $(\fg_{r,\infty},\widetilde{K_{r,\infty}})$-finite automorphic forms on $\widetilde{G_r(\bA_F)}$, where $\fg_{r,\infty}$ denotes the complexified Lie algebra of $G_{r,\infty}$. Denote by $\cA(G_r)_\cusp$ and $\cA(G_r)^\hol$ the subspaces of $\cA(G_r)$ of cusp and holomorphic forms, respectively. Put $\cA(G_r)_\cusp^\hol\coloneqq\cA(G_r)_\cusp\cap\cA(G_r)^\hol$.
\end{enumerate}
\end{notation}

\begin{definition}\label{de:hermitian_siegel}
Consider an element $\lambda=(\lambda_u)_u\in(\fW_r^+)^{\Hom(E,\dC)}$. Construction \ref{co:representation_general} provide us with a direct summand $\tS_\lambda(\dC^r)$ of $\tT^{b(\lambda)}(\dC^r)$ as an (irreducible) representation of $K_{r,\infty}$, which we denote by $\rho_\lambda$ (also regarded as a representation of $\widetilde{K_{r,\infty}}$ by inflation).

Let $w$ be a half-integer (or an integer).
\begin{enumerate}
  \item Denote by $\cA^w_\lambda(G_r)$ the subspace of
      \[
      \Hom_\dC\(\tS_\lambda(\dC^r),\cA(G_r)\)=\tS_\lambda(\dC^r)^\vee\otimes_\dC\cA(G_r)
      \]
      consisting of elements $\varphi$ (written as a homomorphism) satisfying
      \[
      \varphi(\rho_\lambda(k)v)=\delta^{-w}(k)\cdot R_k(\varphi(v)),\qquad\forall v\in \tS_\lambda(\dC^r),\forall k\in\widetilde{K_{r,\infty}},
      \]
      where $R_k$ denotes the right translation on $\cA(G_r)$ by $k$.\footnote{When $w$ is an integer (resp.\ not an integer), the image of $\varphi$ consists of automorphic forms on $G_r(\bA_F)$ (resp.\ genuine automorphic forms on $\widetilde{G_r(\bA_F)}$).} Similarly, we define $\cA^w_\lambda(G_r)_\cusp$, $\cA^w_\lambda(G_r)^\hol$ and $\cA^w_\lambda(G_r)_\cusp^\hol$. When $w=r$, we simply write $\cA_\lambda(G_r)^?_?$ for $\cA^r_\lambda(G_r)^?_?$.

  \item Define a hermitian pairing
    \[
    (\;,\;)_{G_r}\colon\cA^w_\lambda(G_r)_\cusp\times\cA^w_\lambda(G_r)_\cusp\to\dC
    \]
    by the formula
    \begin{align*}
    (\varphi_1,\varphi_2)_{G_r}\coloneqq\int_{G_r(F)\backslash G_r(\bA_F)}
    \tr_\lambda\(\varphi_1(g)\otimes\ol{\varphi_2(g)}\)\rd g,
    \end{align*}
    where
    \[
    \tr_\lambda\colon\tS_\lambda(\dC^r)^\vee\otimes_\dC\ol{\tS_\lambda(\dC^r)^\vee}\to\dC
    \]
    is the natural pairing induced from the standard hermitian form on $\dC^r$, and $\r{d}g$ is the measure in Notation \ref{no:hermitian_siegel}(6). The integral $(\varphi_1,\varphi_2)_{G_r}$ extends to the case where one of $\varphi_1,\varphi_2$ belongs to $\cA^w_\lambda(G_r)$.

  \item The representation $(\rho_\lambda,\tS_\lambda(\dC^r))$ restricted to $K_{r,\infty}\cap M_r(F_\infty)$ extends naturally to a representation of $M_r(F_\infty)$ in the way that $m(a)\in M_r(F_\infty)$ acts on the factor $\tT^{b_u}(\dC^r)$ via $u(a)\in\GL_r(\dC)$. More generally, we have a \emph{map}
      \[
      \rho_\lambda\colon\widetilde{G_r(F_\infty)}\to\Aut_\dC(\tS_\lambda(\dC^r))
      \]
      sending $g\in\widetilde{G_r(F_\infty)}$ to $\rho_\lambda(m(a))\rho_\lambda(k)$ if we write $g=n(b)m(a)k$ for $b\in\Herm_r(F_\infty)$, $a\in\GL_r(E_\infty)$, and $k\in\widetilde{K_{r,\infty}}$.\footnote{The map $\rho_\lambda$ can be regarded as a matrix-valued automorphy factor.}
\end{enumerate}
\end{definition}

\subsection{Doubling method}
\label{ss:doubling}

Consider a hermitian space $V$ over $\bA_E/\bA_F$ of rank $m$. Put $H\coloneqq\rU(V)$ as a reductive group over $\bA_F$.

\begin{definition}\label{de:space_coherence}
We say that a hermitian space $V$ over $\bA_E/\bA_F$ is \emph{coherent/incoherent} if $\chi_{E/F}(\det V)=+1/{-1}$.
\end{definition}

\begin{definition}\label{de:regular_support}
Let $R$ be a ring and $\Diamondblack\subseteq\fV_F^\fin$ a finite subset.
\begin{enumerate}
  \item We denote by $\sS(V_\Diamondblack^r,R)$ and $\sS((V^{\infty,\Diamondblack})^r,R)$ the $R$-module of $R$-valued Schwartz functions (namely, locally constant functions with compact support) on $V_\Diamondblack^r$ and $(V^{\infty,\Diamondblack})^r$, respectively.

  \item Suppose that $1\leq r\leq m$. We say that a Schwartz function in $\sS(V_\Diamondblack^r,R)$ is \emph{regularly supported} if its support is contained in the set
      \[
      \left\{x\in V_\Diamondblack^r\left| T(x)\in\Herm_r^\circ(F_\Diamondblack) \right.\right\}.
      \]

  \item Suppose that $1\leq r\leq m$ and $\Diamondblack=\{v\}$ in which $v$ is nonsplit in $E$. We say that a Schwartz function in $\sS(V_v^r,R)$ is \emph{deeply regularly supported} if its support is contained in the set
      \[
      \left\{x\in V_v^r \left| T(x)+v^d\Herm_r(O_{F_v})\subseteq\Herm_r^\circ(F_v) \right.\right\},
      \]
      where $d$ is the smallest integer such that $\psi_{F,v}$ is trivial on $v^d O_{F_v}$.
\end{enumerate}
We suppress $R$ from the notation when it is $\dC$.
\end{definition}

We review the notion of Weil representation. For every $v\in\fV_F$, we have the Weil representation $\omega_{r,V_v}$ of $\widetilde{G_r(F_v)}\times H(F_v)$ on the space $\sS(V_v^r)$ with respect to the additive character $\psi_{F,v}$. When $m$ is even, the representation factors through $G_r(F_v)\times H(F_v)$,\footnote{When $m$ is odd, one can choose an admissible character $\chi_v$ of $E_v$ satisfying $\chi_v\res_{F_v^\times}=\chi_{E/F,v}$ to split the ($S^1$-)metaplectic cover, under which the formula of $\omega_{r,V_v}$ can be found in, for example, \cite{Liu11}*{\S2A}.} whose formula we review as follows: Taking $\phi\in\sS(V_v^r)$,
\begin{itemize}
  \item for $a\in\GL_r(E_v)$,
    \[
    \omega_{r,V_v}(m(a))\phi(x)=|\det a|_{E_v}^{m/2}\cdot \phi(x a);
    \]

  \item for $b\in\Herm_r(F_v)$,
    \[
    \omega_{r,V_v}(n(b))\phi(x)=\psi_{F,v}(\tr (T(x)b))\cdot\phi(x);
    \]

  \item for $\tw_r\in G_r(F_v)$,
    \[
    \omega_{r,V_v}\(\tw_r\)\phi(x)=\gamma_{V_v}^r\cdot\widehat\phi(x),
    \]
    where $\gamma_{V_v}$ is the Weil constant of $V_v$ (with respect to $\psi_{F,v}$);

  \item for $h\in H(F_v)$,
    \[
    \omega_{r,V_v}(h)\phi(x)=\phi(h^{-1}x).
    \]
\end{itemize}
When $V$ is clear, we simply write $\omega_{r,v}$ for $\omega_{r,V_v}$. Taking (restricted) tensor product, we have $\omega_{r,\infty}$ for all archimedean places or $\omega_r$ for all places.

\begin{remark}\label{re:schwartz_vector}
Let $\{\tS_v\res v\in\fV_F\}$ be a collection of finite-dimensional complex vector spaces indexed by $\fV_F$ such that $\tS_v=\dC$ for all but finitely many $v$.
\begin{enumerate}
  \item For every $v\in\fV_F$, we regard $\sS(V_v^r)\otimes_\dC\tS_v$ as the space of $\tS_v$-valued Schwartz functions on $V_v^r$, hence receives a similar Weil representation $\omega_{r,v}$. Definition \ref{de:regular_support} extends to this setting as well.

  \item Take $\tS\coloneqq\bigotimes_v\tS_v$. We regard $\sS(V^r)\otimes_\dC\tS$ as
      \[
      \(\bigotimes_{v\in\fV_F}\(\sS(V_v^r)\otimes_\dC\tS_v\)\),
      \]
      so that a decomposable Schwartz function $\phi\in\sS(V^r)\otimes_\dC\tS$ is of the form $\otimes_v\phi_v$ in which $\phi_v\in\sS(V_v^r)\otimes_\dC\tS_v$.
\end{enumerate}
\end{remark}

\begin{notation}\label{no:doubling_space}
We collect some notation concerning the doubling skew-hermitian space
\[
W_r^\Box\coloneqq W_r\oplus(-W_r).
\]
Let $G_r^\Box$ be the unitary group of $W_r^\Box$, which contains $G_r\times G_r$ canonically. We now take a basis $\{e^\Box_1,\dots,e^\Box_{4r}\}$ of $W_r^\Box$ by the formula
\[
e^\Box_i=e_i,\quad e^\Box_{r+i}=-e_i^-,\quad e^\Box_{2r+i}=e_{r+i}, \quad e^\Box_{3r+i}=e_{r+i}^-
\]
(where $\{e_1^-,\cdots,e_{2r}^-\}$ is the standard basis of $-W_r$) for $1\leq i\leq r$, under which we may identify $W_r^\Box$ with $W_{2r}$ and $G_r^\Box$ with $G_{2r}$. Finally, put
\[
\bw_r\coloneqq
\begin{pmatrix}
     &  & 1_r &  \\
     & 1_r &  &  \\
    -1_r & 1_r &  &  \\
     &  & 1_r & 1_r \\
\end{pmatrix}
\in G_r^\Box(F).
\]
Then $P_r^\Box\bw_r(G_r\times G_r)$ is Zariski open in $G_r^\Box$.
\end{notation}

From now to the end of this section, we assume $m=2r$. Take a collection $\{\tS_v\res v\in\fV_F\}$ as in Remark \ref{re:schwartz_vector} and put $\tS\coloneqq\bigotimes_v\tS_v$. For every Schwartz function $\Phi\in\sS(V^{2r})\otimes_\dC\tS$, we have the Siegel--Weil section $f_\Phi^{(s)}$ (see \cite{LL21}*{\S3} for details) which is an ($\tS$-valued) function on $P_{2r}(F)\backslash G_{2r}(\bA_F)$, and further the Eisenstein series
\[
E(s,g,\Phi)\coloneqq E(g,f_{\Phi}^{(s)})=\sum_{\gamma\in P_{2r}(F)\backslash G_{2r}(F)}f_\Phi^{(s)}(\gamma g),
\qquad g\in G_{2r}(\bA_F)
\]
which is meromorphic in $s$ and holomorphic on the imaginary line.

For every place $v$ of $F$, $T^\Box\in\Herm_{2r}^\circ(F_v)$, and $\Phi_v\in\sS(V_v^{2r})\otimes_\dC\tS_v$, we have the local Whittaker function on $G_{2r}(F_v)$ with parameter $s\in\dC$ as
\begin{align}\label{eq:local_whittaker}
W_{T^\Box}(s,g,\Phi_v)\coloneqq\int_{\Herm_{2r}(F_v)}f^{(s)}_{\Phi_v}(\tw_{2r}n(b)g)\psi_{F,v}(\tr bT^\Box)^{-1}\rd b
\end{align}
by meromorphic continuation (which is indeed entire), where $\r{d}b$ is the self-dual measure on $\Herm_{2r}(F_v)$ with respect to $\psi_{F,v}$.

\begin{notation}
Consider an element $T^\Box\in\Herm_{2r}(F)$.
\begin{enumerate}
  \item Put $\partial T^\Box=(T_1,T_2)\in\Herm_r(F)^2$ if
      \[
      T=\begin{pmatrix} T_1 & * \\ * & T_2 \end{pmatrix},
      \]
      which belongs to $(\Herm_r^\circ(F)^+)^2$ if $T^\Box\in\Herm_{2r}^\circ(F)^+$.

  \item If $T^\Box\in\Herm_{2r}^\circ(F)$, then we put
      \[
      \Diff(T^\Box,V)\coloneqq\{v\in\fV_F\res (V_v^{2r})_T=\emptyset\},
      \]
      which is a finite set of even/odd cardinality when $V$ is coherent/incoherent, consisting of places of $F$ that are nonsplit in $E$.
\end{enumerate}
\end{notation}

\begin{proposition}\label{pr:eisenstein_series}
Let $\Diamondblack$ be a nonempty finite subset of $\fV_F^{\fin}$. Let $\Phi=\otimes_v\Phi_v\in\sS(V^{2r})\otimes_\dC\tS$ be a decomposable Schwartz function such that $\Phi_\Diamondblack$ is regularly supported (Definition \ref{de:regular_support}).
\begin{enumerate}
  \item For every $g\in P_{2r}(F_\Diamondblack)G_{2r}(\bA_F^\Diamondblack)$, we have
      \[
      E(0,g,\Phi)=\sum_{T^\Box\in\Herm_{2r}^\circ(F)}\prod_{v\in\fV_F} W_{T^\Box}(0,g_v,\Phi_v).
      \]

  \item Suppose that $\Diamondblack$ has cardinality at least $2$ and that $V$ is incoherent (Definition \ref{de:space_coherence}) so that $E(0,g,\Phi)\equiv 0$. Then for every $g\in P_{2r}(F_\Diamondblack)G_{2r}(\bA_F^\Diamondblack)$, we have
      \[
      E'(0,g,\Phi)=\sum_{v\in\fV_F}\fE(g,\Phi)_v,
      \]
      where\footnote{The term $\fE(g,\Phi)_v$ vanishes if $v$ is split in $E$.}
      \[
      \fE(g,\Phi)_v\coloneqq\sum_{\substack{T^\Box\in\Herm_{2r}^\circ(F)\\\Diff(T^\Box,V)=\{v\}}}W'_{T^\Box}(0,g_v,\Phi_v)
      \prod_{v'\in\fV_F\setminus\{v\}} W_{T^\Box}(0,g_{v'},\Phi_{v'}).
      \]

  \item Suppose that for every $v\in\Diamondblack$, $v$ is nonsplit in $E$ such that $\Phi_v$ is deeply regularly supported (Definition \ref{de:regular_support}), and that $V$ is incoherent so that $E(0,g,\Phi)\equiv 0$. Then for every $g\in G_{2r}(\bA_F^\Diamondblack)$, we have
      \[
      E'(0,g,\Phi)=\sum_{v\in\fV_F}\fE(g,\Phi)_v,
      \]
      where $\fE(g,\Phi)_v$ is the one in (2).
\end{enumerate}
\end{proposition}

\begin{proof}
Parts (1) and (2) are proved in \cite{Liu12}*{Section~2B}. Now we show (3).

Write $E_{T^\Box}(s,g,\Phi)$ for the $T^\Box$-th Fourier coefficient of $E(s,g,\Phi)$, so that
\[
E(s,g,\Phi)=\sum_{T^\Box\in\Herm_{2r}(F)}E_{T^\Box}(s,g,\Phi).
\]
It suffices to show that $E_{T^\Box}(s,g,\Phi)\equiv 0$ for $T^\Box\not\in\Herm_{2r}^\circ(F)$ and $g\in G_{2r}(\bA_F^\Diamondblack)$. Without lost of generality, we may assume $g=1_{4r}$. Take an element $T^\Box\in\Herm_{2r}(F)\setminus\Herm_{2r}^\circ(F)$. By \cite{Liu12}*{Lemma~2.2}, we have
\[
E_{T^\Box}(s,1_{4r},\Phi)=\prod_{v\in\fV_F}W_{T^\Box}(s,1_{4r},\Phi_v)
\]
for $\RE s\gg 0$. By Lemma \ref{le:eisenstein_series}(1) below, we have for $v\in\Diamondblack$ (which is nonempty), $W_{T^\Box}(s,1_{4r},\Phi_v)\equiv 0$, so that $E_{T^\Box}(s,g,\Phi)\equiv 0$.

Part (3) is proved.
\end{proof}

\begin{lem}\label{le:eisenstein_series}
Let $v$ be a nonarchimedean place of $F$ that is nonsplit in $E$. Consider a deeply regularly supported Schwartz function $\Phi_v\in\sS(V_v^{2r})\otimes_\dC\tS_v$ (Definition \ref{de:regular_support}).
\begin{enumerate}
  \item We have
      \[
      W_{T^\Box}(s,1_{4r},\Phi_v)\equiv 0
      \]
      for every $T^\Box\in\(\Herm_{2r}(F_v)\setminus\Herm_{2r}^\circ(F_v)\)+v^d\Herm_{2r}(O_{F_v})$, where $d$ is the integer from Definition \ref{de:regular_support}(3).

  \item Let $V^\star_v$ be the other (up to isomorphism) hermitian space over $E_v/F_v$ of rank $2r$ that is not isomorphic to $V_v$. There exists a regularly supported Schwartz function $\Phi^\star_v\in\sS((V^\star_v)^{2r})\otimes_\dC\tS_v$ satisfying that
      \[
      W_{T^\Box}(0,1_{4r},\Phi^\star_v)=W'_{T^\Box}(0,1_{4r},\Phi_v)
      \]
      holds for every $T^\Box$ in the image in the moment map $T\colon(V^\star_v)^{2r}\to\Herm_{2r}(F_v)$.
\end{enumerate}
\end{lem}

\begin{proof}
By \eqref{eq:local_whittaker}, we have
\begin{align*}
W_{T^\Box}(s,1_{4r},\Phi_v)&=\int_{\Herm_{2r}(F_v)}f^{(s)}_{\Phi_v}(\tw_{2r} n(b))\psi_{F,v}(\tr bT^\Box)^{-1}\rd b \\
&=\int_{\Herm_{2r}(F_v)}\(\int_{V_v^{2r}}\psi_{F,v}(\tr bT(x))\Phi_v(x)\rd x\)\th(\tw_{2r} n(b))^s\psi_{F,v}(-\tr b T^\Box)\rd b,
\end{align*}
where $\th$ is the height function with respect to the standard Iwasawa decomposition. Since $\th(\tw_{2r} n(b))$ is invariant under $\Herm_{2r}(O_{F_v})$, the above expression equals
\begin{align*}
&=\int_{\Herm_{2r}(F_v)/\Herm_{2r}(O_{F_v})}\th(\tw_{2r} n(b'))^s\(\int_{V_v^{2r}}\Phi_v(x)\(\int_{\Herm_{2r}(O_{F_v})}\psi_{F,v}(\tr b(T(x)-T^\Box))\rd b\)\rd x\)\rd b'.
\end{align*}

For (1), suppose that $T^\Box\in\(\Herm_{2r}(F_v)\setminus\Herm_{2r}^\circ(F_v)\)+v^d\Herm_{2r}(O_{F_v})$. For $x$ in the support of $\Phi_v$, the character $b\mapsto\psi_{F,v}(\tr b(T(x)-T^\Box))$ on $\Herm_{2r}(O_{F_v})$ is nontrivial, which implies that
\[
\int_{\Herm_{2r}(O_{F_v})}\psi_{F,v}(\tr b(T(x)-T^\Box))\rd b=0.
\]
Thus, $W_{T^\Box}(s,1_{4r},\Phi_v)\equiv 0$.

For (2), note that the function
\[
(x,T^\Box)\in V_v^{2r}\times\Herm_{2r}^\circ(F_v)
\mapsto\int_{\Herm_{2r}(O_{F_v})}\psi_{F,v}(\tr b(T(x)-T^\Box))\rd b
\]
is locally constant, and compactly supported for every fixed $x$. Since $\Phi_v$ is compactly supported, the assignment $T^\Box\mapsto W_{T^\Box}(s,1_{4r},\Phi_v)$ is uniformly locally constant and compactly supported on $\Herm_{2r}^\circ(F_v)$ with respect to $s$. It follows that the assignment $T^\Box\mapsto W'_{T^\Box}(0,1_{4r},\Phi_v)$ is locally constant and compactly supported on $\Herm_{2r}^\circ(F_v)$. By (1), this assignment extends (by zero) to a locally constant and compactly supported on $\Herm_{2r}(F_v)$. Thus, we can find a regularly supported Schwartz function $\Phi^\star_v\in\sS((V^\star_v)^{2r})\otimes_\dC\tS_v$ whose pushforward along the moment map coincides with the assignment $T^\Box\mapsto W'_{T^\Box}(0,1_{4r},\Phi_v)$ (for $T^\Box$ in the image of $T\colon(V^\star_v)^{2r}\to\Herm_{2r}(F_v)$).

The lemma is proved.
\end{proof}

\begin{definition}\label{de:analytic_kernel}
Consider a pair $(T_1,T_2)$ of elements in $\Herm_r(F)$ and a decomposable Schwartz function $\Phi\in\sS(V^{2r})\otimes_\dC\tS$.
\begin{enumerate}
  \item Suppose that $V$ is coherent. We put
    \[
    E_{T_1,T_2}(g,\Phi)\coloneqq\sum_{\substack{T^\Box\in\Herm_{2r}^\circ(F)\\\partial T^\Box=(T_1,T_2)}}
    \prod_{v\in\fV_F} W_{T^\Box}(0,g_v,\Phi_v).
    \]

  \item Suppose that $V$ is incoherent. For every $v\in\fV_F$, we put
    \[
    \fE_{T_1,T_2}(g,\Phi)_v\coloneqq\sum_{\substack{T^\Box\in\Herm_{2r}^\circ(F)\\ \Diff(T^\Box,V)=\{v\}\\ \partial T^\Box=(T_1,T_2)}}
    W'_{T^\Box}(0,g_v,\Phi_v)\prod_{v'\in\fV_F\setminus\{v\}}W_{T^\Box}(0,g_{v'},\Phi_{v'}).
    \]
\end{enumerate}
\end{definition}

Take a \emph{balanced} hermitian $E$-weight $\lambda$ of rank $2r$ (Definition \ref{de:weight}) and suppose that $V$ is totally positive definite. Apply Remark \ref{re:schwartz_vector} with $\tS_v\coloneqq\tS_{\lambda_v^+}(\dC^r)^\vee\otimes_\dC\ol{\tS_{\lambda_v^+}(\dC^r)^\vee}$ for $v\in\fV_F^{(\infty)}$ and $\tS_v=\dC$ for $v\in\fV_F^\fin$.

\begin{definition}\label{de:doubling_schwartz}
Define a function
\[
\Phi_{V_\infty^{2r}}^\lambda\coloneqq\kappa\Psi_{V_\infty^{2r}}^\lambda\cdot\Phi^0_{V_\infty^{2r}}\colon V_\infty^{2r}\to
\tS_{\lambda^+}(\dC^r)^\vee\otimes_\dC\ol{\tS_{\lambda^+}(\dC^r)^\vee}
\]
(Notation \ref{no:contraction} for $\kappa$ with $\dL=\dE^\lambda$), where $\Psi_{V_\infty^{2r}}^\lambda$ is from Construction \ref{co:harmonic_function}(4) and $\Phi^0_{V_\infty^{2r}}$ is the standard Gaussian function on $V_\infty^{2r}$ given by the formula
\[
\Phi^0_{V_\infty^{2r}}(x)=\prod_{v\in\fV_F^{(\infty)}}\te^{-2\pi\tr T(x_v)}.
\]
In particular, $\Phi^\lambda_{V_\infty^{2r}}$ is an element in $\sS(V_\infty^{2r})\otimes_\dC\tS_{\lambda^+}(\dC^r)^\vee\otimes_\dC\ol{\tS_{\lambda^+}(\dC^r)^\vee}$.
\end{definition}

\begin{remark}\label{re:doubling_schwartz}
The Schwartz function $\Phi_{V_\infty^{2r}}^\lambda$ can be decomposed as the product of local ones over archimedean places of $F$.

Take an element $v\in\fV_F^{(\infty)}$. For every element $x\in V_v^r$, we have the map
\[
\mu_{[\lambda_v]}(x)\colon\tS_{[\lambda_v]}(V_v)^\vee\to\tS_{\lambda_v^+}(\dC^r)^\vee
\]
from Construction \ref{co:special_map} with $\rE/\rF=E_v/F_v$ and $\dL=\dC$ (as the natural $F_v$-ring), and its dual
\[
\mu^{[\lambda_v]}(x)\colon\tS_{\lambda_v^+}(\dC^r)\to\tS_{[\lambda_v]}(V_v).
\]
Define a function
\[
\Phi^{\lambda_v}\colon V_v^{2r}\to
\tS_{\lambda_v^+}(\dC^r)^\vee\otimes_\dC\ol{\tS_{\lambda_v^+}(\dC^r)^\vee}
\]
such that for $(x_1,x_2)\in(V_v^r)^2=V_v^{2r}$,
\[
\Phi^{\lambda_v}(x_1,x_2)=\te^{-2\pi\tr T(x_1,x_2)}\cdot\ol{\mu_{[\lambda_v]}(x_2)}\circ\mu^{[\lambda_v]}(x_1).
\]

Then we have
\[
\Phi_{V_\infty^{2r}}^\lambda=\otimes_{v\in\fV_F^{(\infty)}}\Phi^{\lambda_v}.
\]
\end{remark}

\begin{lem}\label{le:doubling_schwartz}
For every pair $(k_1,k_2)\in K_{r,\infty}\times K_{r,\infty}$, we have
\[
\omega_{2r,\infty}(k_1,k_2)\Phi^\lambda_{V_\infty^{2r}}
=\delta^r(k_1)\ol{\delta^r(k_2)}\cdot
\(\rho_{\lambda^+}(k_1)^\vee\otimes\ol{\rho_{\lambda^+}(k_2)^\vee}\)(\Phi^\lambda_{V_\infty^{2r}}).
\]
\end{lem}

\begin{proof}
This follows from \cite{Liu26}*{Lemma~3.7} with $m=2r$ and $q=0$.
\end{proof}

\begin{definition}\label{de:archimedean_factor}
Take a place $v\in\fV_F^{(\infty)}$.
\begin{enumerate}
  \item Denote by $\pi(\lambda_v)$ the unique up to isomorphism holomorphic discrete series representation of $G_r(F_v)$ such that its base change $\BC(\pi(\lambda_v))$ to $E_v$ (as an irreducible admissible representation of $\GL_{2r}(E_v)$) is the normalized induction of
      \[
      \arg^{2b_{u,1}+2r-1}\boxtimes\arg^{2b_{u,2}+2r-3}\boxtimes\cdots\boxtimes\arg^{2b_{u,2r-1}+3-2r}
      \boxtimes\arg^{2b_{u,2r}+1-2r},
      \qquad \lambda_u=(b_{u,1},\dots,b_{u,2r})
      \]
      if we identify $E_v$ with $\dC$ via $u\in\Hom(E,\dC)$ above $v$, where $\arg\colon\dC^\times\to\dC^\times$ denotes the argument character.

  \item The local $L$-function $L(s,\pi(\lambda_v))$ of $\pi(\lambda_v)$ is defined to be $L(s,\BC(\pi(\lambda_v)))$, so that
      \[
      L(\lambda_v)\coloneqq L(1/2,\pi(\lambda_v))=\prod_{v(u)=v}\prod_{i=1}^r 2(2\pi)^{-(b_{u,i}+r+1-i)}\Gamma(b_{u,i}+r+1-i)\in\dR_{>0}.
      \]

  \item Denote by $\xi^{\lambda_v}\colon G_r(F_v)\to\tS_{\lambda_v^+}(\dC^r)^\vee\otimes_\dC\ol{\tS_{\lambda_v^+}(\dC^r)^\vee}$ the vector-valued matrix coefficient satisfying
      \[
      \xi^{\lambda_v}(g)=(\varphi_1,\pi(\lambda_v)(g)\varphi_2)_{\pi(\lambda_v)}
      \]
      for arbitrary $K_{r,v}$-invariant vectors $\varphi_1\in\tS_{\lambda_v^+}(\dC^r)^\vee\otimes_\dC\pi(\lambda_v)$ and $\varphi_2\in\ol{\tS_{\lambda_v^+}(\dC^r)^\vee}\otimes_\dC\pi(\lambda_v)^\vee$ (each exists uniquely up to a constant) satisfying $\tr_{\lambda_v}(\varphi_1,\varphi_2)_{\pi(\lambda_v)}=1$. Put
      \[
      Z(\lambda_v)\coloneqq\int_{G_r(F_v)}\ol{\xi^{\lambda_v}(g)}\cdot f_{\Phi^{\lambda_v}}^{(0)}(\bw_r(g,1_{2r}))\rd g
      \]
      (see Notation \ref{no:hermitian_siegel}(6) for the measure $\r{d}g$ and Remark \ref{re:doubling_schwartz} for $\Phi^{\lambda_v}$), which is absolutely convergent and depends only on $\lambda_v$. Since the map $g\in G_r(F_v)\mapsto (-1)^r\cdot f_{\Phi^{\lambda_v}}^{(0)}(\bw_r(g,1_{2r}))$ is a matrix coefficient of positive type, $(-1)^rZ(\lambda_v)\in\dR_{\geq 0}$. On the other hand, $Z(\lambda_v)$ is nonzero by \cite{EL24}*{Theorem~1.3.1}.\footnote{The Siegel--Weil section used in \cite{EL24}*{Theorem~1.3.1} is a nonzero multiple of $f_{\Phi^{\lambda_v}}^{(s)}(\bw_r(g,1_{2r}))$ evaluated at a highest weight vector in $\tS_{\lambda_v^+}(\dC^r)\otimes_\dC\ol{\tS_{\lambda_v^+}(\dC^r)}$ of norm one; however, it seems complicated to determine the exact value of this multiple in general cases.} It follows that $Z(\lambda_v)\in(-1)^r\dR_{>0}$.
\end{enumerate}
\end{definition}

\begin{definition}\label{de:hermitian_representation}
We say that a pair $(\pi,\cA_\pi)$, in which $\pi$ is an irreducible representation of $G_r(\bA_F)$ and $\cA_\pi$ is a subspace of $\cA(G_r)_\cusp$ on which $\pi$ realizes,\footnote{In particular, $\cA_\pi$ is contained in the subspace of automorphic forms on $G_r(\bA_F)$, rather than the metaplectic cover.} is a \emph{holomorphic tempered cuspidal automorphic representation of $G_r(\bA_F)$ of weight $\lambda$} if
\begin{enumerate}
  \item $\pi_v$ is tempered for every $v\in\fV_F$;

  \item $\pi_v$ is a discrete series representation for every $v\in\fV_F^{(\infty)}$;

  \item $\(\tS_{\lambda^+}(\dC^r)^\vee\otimes_\dC\cA_\pi\)\cap\cA_{\lambda^+}(G_r)_\cusp^\hol$ (Definition \ref{de:hermitian_siegel}(1)) is nonzero.
\end{enumerate}
Note that (2) plus (3) is equivalent to that $\pi_v\simeq\pi(\lambda_v)$ for every $v\in\fV_F^{(\infty)}$.

Define the (global) $L$-function for $\pi$ to be the Euler product
\[
L(s,\pi)\coloneqq\prod_{v\in\fV_F}L(s,\pi_v),
\]
in which $L(s,\pi_v)$ is the doubling $L$-function defined in \cite{Yam14}*{Theorem~5.2} for $v\in\tV_F^\fin$ and $L(s,\pi_v)=L(s,\pi(\lambda_v))$ for $v\in\tV_F^{(\infty)}$. The $L$-function $L(s,\pi)$ satisfies a self-dual functional equation with the center $s=1/2$ and the epsilon factor $\epsilon(s,\pi)$. Put $\epsilon(\pi)\coloneqq\epsilon(1/2,\pi)\in\{\pm 1\}$.
\end{definition}

Take a holomorphic tempered cuspidal automorphic representation $(\pi,\cA_\pi)$ of $G_r(\bA_F)$ of weight $\lambda$.

\begin{definition}\label{de:doubling_zeta}
Consider a pair $\ul\varphi=(\varphi_1,\varphi_2)$ of automorphic forms in $\(\tS_{\lambda^+}(\dC^r)^\vee\otimes_\dC\cA_\pi\)\cap\cA_{\lambda^+}(G_r)_\cusp^\hol$ and a Schwartz function $\Phi^\infty\in\sS((V^\infty)^{2r})$. For every finite set $\Diamondblack$ of nonarchimedean places of $F$, we have the $\Diamondblack$-local doubling zeta integral
\[
\fZ_\Diamondblack(s,\ul\varphi,\Phi^\infty)\coloneqq\int_{G_r(F_\Diamondblack)}
\ol{(\varphi_1,\pi(g)\varphi_2)_{G_r}}\cdot f_{\Phi^\infty}^{(s)}(\bw_r(g,1_{2r}))\rd g
\]
(Notation \ref{no:hermitian_siegel}(6) and Definition \ref{de:hermitian_siegel}(2)) and its normalized version
\[
\fZ^\natural_\Diamondblack(s,\ul\varphi,\Phi^\infty)\coloneqq
\(\frac{L(s+\tfrac{1}{2},\pi_\Diamondblack)}{b_{2r,\Diamondblack}(s)}\)^{-1}\cdot\fZ_\Diamondblack(s,\ul\varphi,\Phi^\infty),
\]
in which the latter has a holomorphic continuation to the entire complex plane. Here, $b_{2r,\Diamondblack}\coloneqq\prod_{v\in\Diamondblack}b_{2r,v}$ (Notation \ref{no:cm}(3)).

By \cite{Yam14}*{Proposition~7.1 \& (7.2)}, the function $\fZ^\natural_\Diamondblack(s,\ul\varphi,\Phi^\infty)$ stabilizes when $\Diamondblack$ enlarges, which we denote by $\fZ^\natural_\fin(s,\ul\varphi,\Phi^\infty)$. Finally, put $\fZ^\natural_\fin(\ul\varphi,\Phi^\infty)\coloneqq\fZ^\natural_\fin(0,\ul\varphi,\Phi^\infty)$.
\end{definition}

\begin{proposition}\label{pr:ripf}
Let $\lambda$ be a balanced hermitian $E$-weight $\lambda$ of rank $2r$, and $(\pi,\cV_\pi)$ a holomorphic tempered cuspidal automorphic representation of $G_r(\bA_F)$ of weight $\lambda$. Suppose that $V$ is totally positive definite. For every pair $\ul\varphi=(\varphi_1,\varphi_2)$ of automorphic forms in $\(\tS_{\lambda^+}(\dC^r)^\vee\otimes_\dC\cA_\pi\)\cap\cA_{\lambda^+}(G_r)_\cusp^\hol$ and every Schwartz function $\Phi^\infty\in\sS((V^\infty)^{2r})$,
\begin{align*}
&\quad\int_{G_r(F)\backslash G_r(\bA_F)}\tr_\lambda\(\varphi_2(g_2)\otimes
\int_{G_r(F)\backslash G_r(\bA_F)}
\tr_\lambda\(E'(0,(g_1,g_2),\Phi^\lambda_{V_\infty^{2r}}\Phi^\infty)\otimes\ol{\varphi_1(g_1)}\)\rd g_1\)\rd g_2 \\
&=\frac{L'(\tfrac{1}{2},\pi)}{b_{2r}(0)}\cdot C(\lambda)
\cdot\fZ^\natural_\fin(\ul{\varphi},\Phi^\infty)
\end{align*}
holds. Here, $C(\lambda)\coloneqq\prod_{v\in\fV_F^{(\infty)}}C(\lambda_v)$, and
\[
C(\lambda_v)\coloneqq b_{2r,v}(0)\cdot L(\lambda_v)^{-1}\cdot Z(\lambda_v)\in (-1)^r\dR_{>0}
\]
in which the three factors are from Notation \ref{no:cm}(3), Definition \ref{de:archimedean_factor}(2), and Definition \ref{de:archimedean_factor}(3), respectively (in order).
\end{proposition}

\begin{proof}
This follows from the same argument for \cite{LL21}*{Proposition~3.7} together with Remark \ref{re:doubling_schwartz}.
\end{proof}

\begin{remark}\label{re:constant}
Concerning the constant $C(\lambda_v)$ in the above proposition, we have the following remarks.
\begin{enumerate}
  \item When $\lambda_v=\b{0}$ is the trivial weight, $C(\lambda_v)$ recovers the constant $C_r$ in \cite{LL21}*{Proposition~3.7}.

  \item When $r=1$, we have $\lambda_u=(b_u,-b_{\ol{u}})$ for $b_u\geq 0$. If we canonically identify $\tS_{\lambda_v^+}(\dC)$ with $\dC$, then by \cite{Lin}*{Theorem~1.1},\footnote{Our Haar measure on $G_1(F_v)$ is $\pi$ times the one adopted in \cite{Lin}.}
      \[
      Z(\lambda_v)=-\frac{\pi}{1+\sum_{v(u)=v}b_u}\cdot\int_{V_v}\Phi^{\lambda_v}(x,x)\rd x,
      \]
      in which $\r{d}x$ is the self-dual measure on $V_v$ (which is just the standard Euclidean measure). By an exercise in linear algebra, we have
      \[
      \Phi^{\lambda_v}(x,x)=\frac{\prod_{v(u)=v}(b_u!)}{(\sum_{v(u)=v}b_u)!}\cdot (x,x)_{V_v}^{\sum_{v(u)=v}b_u}\te^{-2\pi(x,x)_{V_v}}.
      \]
      It follows that
      \begin{align*}
      Z(\lambda_v)=-\frac{\pi}{1+\sum_{v(u)=v}b_u}
      \cdot\frac{\prod_{v(u)=v}(b_u!)}{(\sum_{v(u)=v}b_u)!}\cdot\frac{(1+\sum_{v(u)=v}b_u)!}{4(2\pi)^{\sum_{v(u)=v}b_u}}
      =-\frac{\pi}{4}\cdot\frac{\prod_{v(u)=v}(b_u!)}{(2\pi)^{\sum_{v(u)=v}b_u}},
      \end{align*}
      hence
      \[
      L(\lambda_v)=-\frac{\pi}{4}.
      \]
\end{enumerate}
\end{remark}

\section{Motivic special cycles and motivic generating functions}
\label{ss:3}

In this chapter, we define several key objects in our theory, including standard (relative) motives, automorphic (relative) motives, motivic special classes, motivic generating functions, and motivic arithmetic theta lifting.

From now to the end of this article except Remark \ref{re:general_cm}, we assume that the CM field $E$ contains an imaginary quadratic field, for which we choose one and denote it by $E_0$. Recall from Notation \ref{no:cm} that $\dE\subseteq\dC$ denotes the normal closure of $E$ in $\dC$, so that $\Hom(E,\dC)=\Hom(E,\dE)$.

\subsection{An auxiliary moduli problem}
\label{ss:auxiliary_moduli}

Put $T\coloneqq\Res_{E_0/\dQ}\dG_m$.

\begin{definition}\label{de:moduli_auxiliary}
We define a functor
\begin{align*}
Y_\bullet\colon\Hk(T(\bA^\infty))\to\sfP\Sch'_{/E_0}
\end{align*}
(Definition \ref{de:hecke}) such that for every objects $U\in\Hk(T(\bA^\infty))$ and $S\in\Sch'_{/E_0}$, $Y_U(S)$ is the set of equivalence classes of triples $(A_0,\lambda_0,\bar\eta_0)$, where
\begin{itemize}
  \item $A_0$ is an object of $\AV(S,E_0)$ (Definition \ref{de:category_av}) such that $\tr_{A_0}(e)=e$ for $e\in E_0$;

  \item $\lambda_0\colon A_0\to A_0^\vee$ is a polarization that turns the complex conjugation of $E_0$ into the Rosati involution;

  \item $\bar\eta_0$ is a $U$-level structure, that is, for a chosen geometric point $s$ on every connected component of $S$, a $\pi_1(S,s)$-invariant $U$-orbit of isomorphisms
      \[
      \eta_0\colon\bA_{E_0}^\infty\xrightarrow\sim\rH_1^{\et}(A_{0s},\bA^\infty)
      \]
      of $\bA_{E_0}^\infty$-modules.
\end{itemize}
Two objects $(A_0,\lambda_0,\bar\eta_0)$ and $(A'_0,\lambda'_0,\bar\eta'_0)$ are equivalent\footnote{In each equivalence class of objects, one always has a representative in which $\lambda_0$ is the canonical principal polarization of the modular curve.} if there exists an isomorphism $\varphi_0\colon A_0\to A'_0$ in the category $\AV(S,E_0)$ satisfying that
\begin{itemize}
  \item $\varphi_0^\vee\circ\lambda'_0\circ\varphi_0$ is a multiple of $\lambda_0$;

  \item the $U$-orbits $\varphi_{0*}\circ\bar\eta_0$ and $\bar\eta'_0$ coincide.
\end{itemize}
A morphism $t\in U\backslash T(\bA^\infty)/U'$ maps $Y_U(S)$ to $Y_{U'}(S)$ by changing $\bar\eta_0$ to $\bar\eta_0\circ t$.
\end{definition}

It is well-known that $Y_\bullet$ is representable and factors through the full subcategory $\Sch^\smp_{/E_0}$ of $\sfP\Sch'_{/E_0}$. Moreover, if we denote by $A_{0,U}\in\AV(Y_U,E_0)$ the universal object, then we obtain a functor
\[
A_{0,\bullet}/Y_\bullet\colon\Hk(T(\bA^\infty))\to\AV(\Sch^\smp_{/E_0},E_0)
\]
over $Y_\bullet$ that is Cartesian (Definition \ref{de:category_section}).

\begin{remark}\label{re:moduli_auxiliary}
For every object $U\in\Hk(T(\bA^\infty))$, the group $T(\bA^\infty)/U$ acts on $Y_U$ via modifying $\bar\eta_0$, which factors through the double quotient $T(\dQ)\backslash T(\bA^\infty)/U$. By the theory of complex multiplication, the induced action of $T(\dQ)\backslash T(\bA^\infty)/U$ on $Y_U$ is free, and the quotient scheme is $\Spec E_0$.
\end{remark}

\subsection{Standard motives over unitary Shimura varieties}
\label{ss:motive_standard}

From now to the end of this chapter, we consider a totally positive definite incoherent hermitian space $V$ over $\bA_E/\bA_F$ of rank $m$ (Definition \ref{de:space_coherence}), and put $H\coloneqq\rU(V)$.

\begin{definition}\label{de:moduli_rational}
We define a functor
\begin{align*}
X_{\bullet,\bullet}\colon\Hk(H(\bA_F^\infty))\times\Hk(T(\bA^\infty))\to\sfP\Sch'_{/E}
\end{align*}
(Definition \ref{de:hecke}) such that for every objects $(L,U)\in\Hk(H(\bA_F^\infty))\times\Hk(T(\bA^\infty))$ and $S\in\Sch'_{/E}$, $X_{L,U}(S)$ is the set of equivalence classes of sextuples $(A_0,\lambda_0,\bar\eta_0;A,\lambda,\bar\eta)$,\footnote{Unfortunately, we have used $\lambda$ here to denote a polarization. But this should not be confused with the weight $\lambda$.} where
\begin{itemize}
  \item $(A_0,\lambda_0,\bar\eta_0)$ is an object of $Y_U(S)$;

  \item $A$ is an object of $\AV(S,E)$ (Definition \ref{de:category_av}) such that $\tr_A(e)=m\Tr_{E/E_0}(e)-e+\ol{e}$ for $e\in E$;

  \item $\lambda\colon A\to A^\vee$ is a polarization that turns the complex conjugation of $E$ into the Rosati involution;

  \item $\bar\eta$ is a \emph{hermitian $L$-level structure}, that is, for a chosen geometric point $s$ on every connected component of $S$, a $\pi_1(S,s)$-invariant $L$-orbit of isometries
      \[
      \eta\colon V^\infty\xrightarrow\sim\Hom_{\bA_E^\infty}^{\lambda_0,\lambda}
      \(\rH_1^{\et}(A_{0s},\bA^\infty)\otimes_{\bA^\infty}\bA_F^\infty,\rH_1^{\et}(A_s,\bA^\infty)\)
      \]
      of hermitian spaces over $\bA_E^\infty/\bA_F^\infty$ (see, for example, \cite{LTXZZ}*{Construction~3.4.4} for the hermitian structure on the right-hand side).
\end{itemize}
Two objects $(A_0,\lambda_0,\bar\eta_0;A,\lambda,\bar\eta)$ and $(A'_0,\lambda'_0,\bar\eta'_0;A',\lambda',\bar\eta')$ are equivalent if there exist isomorphisms $\varphi_0\colon A_0\to A'_0$ in the category $\AV(S,E_0)$ and $\varphi\colon A\to A'$ in the category $\AV(S,E)$ satisfying that
\begin{itemize}
  \item $\varphi_0^\vee\circ\lambda'_0\circ\varphi_0$ and $\varphi^\vee\circ\lambda'\circ\varphi$ are the \emph{same} multiple of $\lambda_0$ and $\lambda$, respectively;

  \item the $U$-orbits $\varphi_{0*}\circ\bar\eta_0$ and $\bar\eta'_0$ coincide;

  \item the $L$-orbit of maps $v\mapsto\varphi_*\circ\bar\eta(v)\circ(\varphi_{0*})^{-1}$ for $v\in V^\infty$ coincides with $\bar\eta'$.
\end{itemize}
A morphism $(h,t)\in L\backslash H(\bA_F^\infty)/L'\times U\backslash T(\bA^\infty)/U'$ maps $X_{L,U}(S)$ to $X_{L',U'}(S)$ by changing $\bar\eta_0$ to $\bar\eta_0\circ t$ and $\bar\eta$ to $\bar\eta\circ h$.
\end{definition}

By \cite{RSZ}*{Theorem~4.4}, $X_{\bullet,\bullet}$ factors through the full subcategory $\Sch^\smq_{/E}$ of $\sfP\Sch'_{/E}$. If we denote by
\[
(A_{0,L,U},A_{L,U})\in\AV(X_{L,U},E_0)\times\AV(X_{L,U},E)
\]
the universal object, then we obtain functors
\begin{align*}
A_{0,\bullet,\bullet}/X_{\bullet,\bullet}&\colon\Hk(H(\bA_F^\infty))\times\Hk(T(\bA^\infty))\to\AV(\Sch^\smq_{/E},E_0),\\
A_{\bullet,\bullet}/X_{\bullet,\bullet}&\colon\Hk(H(\bA_F^\infty))\times\Hk(T(\bA^\infty))\to\AV(\Sch^\smq_{/E},E)
\end{align*}
over $X_{\bullet,\bullet}$ that are Cartesian (Definition \ref{de:category_section}). By definition, we have a natural morphism
\begin{align}\label{eq:moduli_rational}
X_{L,U}\to Y_U\otimes_{E_0}E
\end{align}
for every object $(L,U)\in\Hk(H(\bA_F^\infty))\times\Hk(T(\bA^\infty))$, along which $A_{0,L,U}$ is the pullback of $A_{0,U}\otimes_{E_0}E$ to $X_{L,U}$.

Take an object $L\in\Hk(H(\bA_F^\infty))$. For every $U\in\Hk(T(\bA^\infty))$, the group $T(\bA^\infty)/U$ acts on $X_{L,U}$ via modifying $\bar\eta_0$, which also factors through the quotient
\[
\Xi_U\coloneqq T(\dQ)\backslash T(\bA^\infty)/U.
\]
Indeed, for every element $t\in E_0^\times$, the isomorphisms in $\AV(S,E_0)$ and $\AV(S,E)$ given by the action of $t$ give an equivalence between $(A_0,\lambda_0,\bar\eta_0;A,\lambda,\bar\eta)$ and $(A_0,\lambda_0,\bar\eta_0\circ t;A,\lambda,\bar\eta)$. It is clear that the morphism \eqref{eq:moduli_rational} is equivariant with respect to the action of $\Xi_U$. Thus, by Remark \ref{re:moduli_auxiliary}, the action of $\Xi_U$ on $X_{L,U}$ is free. Since $X_{L,U}$ is smooth and quasi-projective over $\Spec E$, the quotient $X_{L,U}/\Xi_U$ is a scheme. It is clear that for $U'\subseteq U$ in $\Hk(T(\bA^\infty))$, the natural morphism $X_{L,U'}/\Xi_{U'}\to X_{L,U}/\Xi_U$ is an isomorphism.

To summarize, we obtain a functor
\[
X_\bullet\colon\Hk(H(\bA_F^\infty))\to\Sch^\smq_{/E}
\]
such that for every $L$, $X_L=X_{L,U}/\Xi_U$ -- the quotient scheme obtained from the above process, which is canonically independent of the choice of $U\in\Hk(T(\bA^\infty))$. In what follows, we denote by
\begin{align}\label{eq:quotient_morphism}
\pi_{L,U}\colon X_{L,U}\to X_L
\end{align}
the natural quotient morphism.

\begin{remark}\label{re:moduli_shimura}
We discuss the relation between $X_\bullet$ and the usual Shimura variety. Take an element $v\in\fV_F^{(\infty)}$. Denote by $V(v)$ the unique hermitian space over $E/F$ that has signature $(m-1,1)$ at $v$ and such that $V(v)\otimes_F\bA_F^v$ is isometric to $V\otimes_{\bA_F}\bA_F^v$. Put $H(v)\coloneqq\rU(V(v))$, and fix an isometry $V(v)\otimes_F\bA_F^\infty\simeq V^\infty$ under which we identify $H(\bA_F^\infty)$ with $H(v)(\bA_F^\infty)$. For every $u\in\Hom(E,\dC)$ above $v$, we have a Shimura variety functor\footnote{Strictly speaking, the Shimura variety is for the reductive group $\Res_{F^+/\dQ}H(v)$.}
\[
\Sh(H(v))_\bullet^u\colon\Hk(H(\bA_F^\infty))\to\Sch^\smq_{/u(E)}
\]
after we identify $E\otimes_{F,v}\dR$ with $\dC$ using $u$ (see, for example, \cite{LTXZZ}*{Section~3.2}). Then we have a canonical isomorphism
\begin{align}\label{eq:moduli_shimura}
X_\bullet\otimes_{E,u}u(E)\xrightarrow\sim\Sh(H(v))_\bullet^u
\end{align}
in $\Fun(\Hk(H(\bA_F^\infty)),\Sch^\smq_{/u(E)})$. In particular, for every $L\in\Hk(H(\bA_F^\infty))$,
\[
X_L^u=H(v)(F)\backslash\(\fD^u\times H(\bA_F^\infty)/L\)
\]
(as complex manifolds) in which $\fD^u$ denotes the open subset of $V(v)\otimes_{E,u}\dC$ of negative definite lines.
\end{remark}

Take an object $L\in\Hk(H(\bA_F^\infty))$. For every $U\in\Hk(T(\bA^\infty))$, we have the relative motives (see \S\ref{ss:abelian_motive})\footnote{We recommend the readers to go through Appendix \ref{ss:a} first as this point, for a quick summary of the knowledge on motives we use below.}
\begin{align*}
\sfh_1(A_{0,L,U}/X_{L,U},\dQ)&\in\Mot(X_{L,U},E_0),\\
\sfh_1(A_{L,U}/X_{L,U},\dQ)&\in\Mot(X_{L,U},E),
\end{align*}
hence
\[
\rM_{L,U}\coloneqq\ul\Hom_{E_0}\(\sfh_1(A_{0,L,U}/X_{L,U},\dQ),\sfh_1(A_{L,U}/X_{L,U},\dQ)\)\in \Mot(X_{L,U},E).
\]
We would like to descend the above motive to $X_L$. For every $t\in T(\bA^\infty)/U$, denote by
\[
\ang{t}\colon X_{L,U}\to X_{L,U}
\]
the induced (iso)morphism. By definition, the pullback of the universal data $\ang{t}^*(A_0,\lambda_0,\bar\eta_0;A,\lambda,\bar\eta)$ is equivalent to $(A_0,\lambda_0,\bar\eta_0\circ t;A,\lambda,\bar\eta)$. Denote by $\varphi_{t,0}\colon A_{0,L,U}\to\ang{t}^*A_{0,L,U}$ and $\varphi_t\colon A_{L,U}\to\ang{t}^*A_{L,U}$ the unique isomorphisms in $\AV(X_{L,U},E_0)$ and $\AV(X_{L,U},E)$, respectively, that give the unique equivalence from $(A_0,\lambda_0,\bar\eta_0\circ t;A,\lambda,\bar\eta)$ to $\ang{t}^*(A_0,\lambda_0,\bar\eta_0;A,\lambda,\bar\eta)$. Then
\[
\ang{t}^*(\varphi_{t'})\circ\varphi_t=\varphi_{tt'}
\]
for every $t'\in T(\bA^\infty)/U$. Denote by $\ul\Xi_U$ the groupoid given by $\Xi_U$. We have a functor
\[
X_{L,U}/\!\!/\Xi_U\colon\ul\Xi_U\to\Sch^\smq_{/E}
\]
sending the unique object to $X_{L,U}$ and $t\in\Xi_U$ to $\llangle t\rrangle$.

For $t\in T(\bA^\infty)/U$, we define a map
\begin{align*}
\alpha_t\colon\ang{t}^*\rM_{L,U} \to\rM_{L,U}
\end{align*}
as the composition of isomorphisms
\begin{align*}
\ang{t}^*\rM_{L,U}
&\xrightarrow\sim\ul\Hom_{E_0}\(\ang{t}^*\sfh_1(A_{0,L,U}/X_{L,U},\dQ),\ang{t}^*\sfh_1(A_{L,U}/X_{L,U},\dQ)\) \\
&\xrightarrow\sim\ang{t}^*\sfh^1(A_{0,L,U}/X_{L,U},\dQ)\otimes_{E_0}\ang{t}^*\sfh_1(A_{L,U}/X_{L,U},\dQ) \\
&\xrightarrow\sim\sfh^1(\ang{t}^*A_{0,L,U}/X_{L,U},\dQ)\otimes_{E_0}\sfh_1(\ang{t}^*A_{L,U}/X_{L,U},\dQ) \\
&\xrightarrow{\varphi_{t,0}^*\otimes(\varphi_t)_*^{-1}}\sfh^1(A_{0,L,U}/X_{L,U},\dQ)\otimes_{E_0}\sfh_1(A_{L,U}/X_{L,U},\dQ) \\
&\xrightarrow\sim\ul\Hom_{E_0}\(\sfh_1(A_{0,L,U}/X_{L,U},\dQ),\sfh_1(A_{L,U}/X_{L,U},\dQ)\)=\rM_{L,U},
\end{align*}
in which we have used the fact that the canonical maps
\begin{align*}
\varphi_{t,0}^*&\colon\sfh^1(\ang{t}^*A_{0,L,U}/X_{L,U},\dQ)\to\sfh^1(A_{0,L,U}/X_{L,U},\dQ),\\
(\varphi_t)_*&\colon\sfh_1(A_{L,U}/X_{L,U},\dQ)\to\sfh_1(\ang{t}^*A_{L,U}/X_{L,U},\dQ)
\end{align*}
are isomorphisms. Then it is easy to verify that
\[
\alpha_t\circ\ang{t}^*(\alpha_{t'})=\alpha_{tt'}
\]
for every $t'\in T(\bA^\infty)/U$. Moreover, for $t\in E_0^\times$, $\alpha_t$ is the identity isomorphism. Thus, if we denote by $\Mot(X_{L,U}/\!\!/\Xi_U,E)$ the category of Cartesian functors $\ul\Xi_U\to\Mot(\Sch^\smq_{/E},E)$ above $X_{L,U}/\!\!/\Xi_U$, then $\rM_{L,U}$ together with the isomorphisms $\{\alpha_t\}$ define an object
\[
\rM_{L,U}/\!\!/\Xi_U\in\Mot(X_{L,U}/\!\!/\Xi_U,E).
\]
By Lemma \ref{le:motive_descent} below, $\rM_{L,U}/\!\!/\Xi_U$ is an object of $\Mot(X_L,E)$. By the same lemma, for $U'\subseteq U$ in $\Hk(T(\bA^\infty))$, the natural map $\rM_{L,U}/\!\!/\Xi_U\to\rM_{L,U'}/\!\!/\Xi_{U'}$ is an isomorphism in $\Mot(X_L,E)$.

\begin{definition}\label{de:motive_standard}
We define the \emph{standard motive} over $X_L$ to be the object $\rM_{L,U}/\!\!/\Xi_U$ of $\Mot(X_L,E)$ defined above, which is independent of $U$, and we denote it by $\rM_L$.
\end{definition}

\begin{remark}
In general, neither $\sfh_1(A_{0,L,U}/X_{L,U},\dQ)$ nor $\sfh_1(A_{L,U}/X_{L,U},\dQ)$ descends to a relative motive over $X_L$.
\end{remark}

It is straightforward to verify that
\[
\rM_\bullet\colon\Hk(H(\bA_F^\infty))\to\Mot(\Sch^\smq_{/E},E)
\]
is a Cartesian functor above $X_\bullet$.

\begin{lem}\label{le:motive_descent}
Let the notation be as above. The natural functor $\Mot(X_L,E)\to\Mot(X_{L,U}/\!\!/\Xi_U,E)$ induced by pullback is an equivalence of categories.
\end{lem}

\begin{proof}
By construction, the natural morphism $X_{L,U}\to X_L\times_{E_0}Y_U$ is an isomorphism. We would like to apply \cite{Kah}*{Theorem~2.17} with $\Pi=T(\dQ)\backslash T(\bA^\infty)$. Let $B\Pi$ be the classifying stack of the profinite group $\Pi$, that is, the category of finite continuous $\Pi$-sets. We have a functor $B\Pi\to\Sch_{/E_0}$ that sends finite coproducts to finite products and such that for every $U\in\Hk(T(\bA^\infty))$, $Y(\Pi/U)=Y_U$. Now let $\cA$ be the contravariant pseudo-functor from $B\Pi$ to the 2-category of pseudo-abelian categories that sends a finite continuous $\Pi$-set $S$ to $\Mot(X_L\times_{E_0}Y(S),\dQ)$. By \cite{Kah}*{Theorem~2.17(b)}, we know that $\cA$ is a stack over $B\Pi$, which implies the lemma.
\end{proof}

\begin{remark}
It follows from Remark \ref{re:moduli_shimura} that $X_\bullet$ does not depend on the choice of $E_0$ (an imaginary quadratic field contained in $E$). However, we expect but do not know how to prove that $\rM_\bullet$ also does not depend on such a choice. This is why our theory relies on the choice of $E_0$.
\end{remark}

In the rest of this section, we construct a map
\begin{align}\label{eq:evaluation_motive_standard}
\rM_L\otimes_E\ol{\rM_L}\to\mathbbm{1}_{X_L,E},
\end{align}
where $\mathbbm{1}_{X_L,E}$ denotes the unit object of $\Mot(X_L,E)$, which turns out to be an evaluation map (Definition \ref{de:dualizable}).

Take an arbitrary object $U\in\Hk(T(\bA^\infty))$. Using the isomorphism in Lemma \ref{le:poincare}, we have successive isomorphisms in the category $\Mot(X_{L,U},E)$:
\begin{align}\label{eq:evaluation_motive_standard_1}
\rM_{L,U}&\xrightarrow\sim\sfh^1(A_{0,L,U}/X_{L,U},\dQ)\otimes_{E_0}\sfh_1(A_{L,U}/X_{L,U},\dQ)\\
&\xrightarrow\sim\sfh^1(A_{0,L,U}/X_{L,U},\dQ)\otimes_{E_0}\sfh_1(A_{L,U}/X_{L,U},\dQ)  \notag\\
&\xrightarrow\sim\sfh^1(A_{0,L,U}/X_{L,U},\dQ)\otimes_{E_0}\sfh^1(A_{L,U}^\vee /X_{L,U},\dQ)\ang{1}. \notag
\end{align}
Using Lemma \ref{le:poincare} again, we have
\begin{align*}
\rM_{L,U}^\vee&\simeq\sfh^1(A_{0,L,U}/X_{L,U},\dQ)^\vee\otimes_{E_0}\sfh^1(A_{L,U}^\vee /X_{L,U},\dQ)^\vee\ang{-1} \\
&\xrightarrow\sim\sfh^1(A_{0,L,U}^\vee/X_{L,U},\dQ)\ang{1}\otimes_{E_0}\sfh^1(A_{L,U}^{\vee\vee}/X_{L,U},\dQ) \\
&\xrightarrow\sim\sfh^1(A_{0,L,U}^\vee/X_{L,U},\dQ)\otimes_{E_0}\sfh^1(A_{L,U}/X_{L,U},\dQ)\ang{1}.
\end{align*}
Now the polarizations $\lambda_0$ and $\lambda$ induce isomorphisms
\begin{align*}
\lambda_0^*&\colon\sfh^1(A_{0,L,U}^\vee/X_{L,U},\dQ)\xrightarrow\sim\sfh^1(A_{0,L,U}/X_{L,U},\dQ)\otimes_{E_0,\ol{\phantom{a}}}E_0, \\
(\lambda^{-1})^*&\colon\sfh^1(A_{L,U}/X_{L,U},\dQ)\xrightarrow\sim\sfh^1(A_{L,U}^\vee/X_{L,U},\dQ)\otimes_{E,\ol{\phantom{a}}}E
\end{align*}
in the categories $\Mot(X_{L,U},E_0)$ and $\Mot(X_{L,U},E)$, respectively. Thus, we obtain an isomorphism
\begin{align}\label{eq:evaluation_motive_standard_2}
\rM_{L,U}^\vee\xrightarrow\sim
\(\sfh^1(A_{0,L,U}/X_{L,U},\dQ)\otimes_{E_0}\sfh^1(A_{L,U}^\vee /X_{L,U},\dQ)\ang{1}\)\otimes_{E,\ol{\phantom{a}}}E\xrightarrow\sim\rM_{L,U}\otimes_{E,\ol{\phantom{a}}}E
\end{align}
in the category $\Mot(X_{L,U},E)$. It is straightforward to check that the isomorphism \eqref{eq:evaluation_motive_standard_2} is indeed an isomorphism in $\Mot(X_{L,U}/\!\!/\Xi_U,E)$, hence descends to an isomorphism $\rM_L^\vee\xrightarrow\sim\ol{\rM_L}$, which is independent of the choice of $U$. In particular, we obtain the evaluation map \eqref{eq:evaluation_motive_standard}.

\subsection{Special morphisms of unitary Shimura varieties}
\label{ss:speical_morphism}

Consider a tuple $x=(x_1,\dots,x_r)\in(V^\infty)^r$.

\begin{definition}\label{de:special_morphism}
For every object $L\in\Hk(H(\bA_F^\infty))$, we define a functor
\begin{align*}
Z(x)_{L,\bullet}\colon\Hk(T(\bA^\infty))\to\sfP\Sch'_{/E}
\end{align*}
such that for every objects $U\in\Hk(T(\bA^\infty))$ and $S\in\Sch'_{/E}$, $Z(x)_{L,U}(S)$ is the set of equivalence classes of $(6+r)$-tuples $(A_0,\lambda_0,\bar\eta_0;A,\lambda,\bar\eta;\nu_1,\dots,\nu_r)$, where
\begin{itemize}
  \item $(A_0,\lambda_0,\bar\eta_0;A,\lambda,\bar\eta)$ is an object of $X_{L,U}(S)$;

  \item $\nu_1,\dots,\nu_r$ are morphisms from $A_0$ to $A$ in the category $\AV(S,E_0)$ satisfying that (on every connected component of $S$) there exists an element $\eta\in\bar\eta$ such that $\eta(x_i)$ coincides with the induced map $\nu_{i*}\colon\rH_1^{\et}(A_{0s},\bA^\infty)\otimes_{\bA^\infty}\bA_F^\infty\to\rH_1^{\et}(A_s,\bA^\infty)$ for $1\leq i\leq r$.
\end{itemize}
Two objects $(A_0,\lambda_0,\bar\eta_0;A,\lambda,\bar\eta;\nu_1,\dots,\nu_r)$ and $(A'_0,\lambda'_0,\bar\eta'_0;A',\lambda',\bar\eta';\nu'_1,\dots,\nu'_r)$ are equivalent if there exist isomorphisms $\varphi_0\colon A_0\to A'_0$ in the category $\AV(S,E_0)$ and $\varphi\colon A\to A'$ in the category $\AV(S,E)$ that induce equivalence between $(A_0,\lambda_0,\bar\eta_0;A,\lambda)$ and $(A'_0,\lambda'_0,\bar\eta'_0;A',\lambda',\bar\eta')$ and such that $\varphi\circ\nu_i=\nu'_i\circ\varphi_0$ for $1\leq i\leq r$.

A morphism $t\in U\backslash T(\bA^\infty)/U'$ maps $Z(x)_{L,U}(S)$ to $Z(x)_{L,U'}(S)$ by changing $\bar\eta_0$ to $\bar\eta_0\circ t$.
\end{definition}

It is clear that $Z(x)_{L,\bullet}$ is representable and factors through the full subcategory $\Sch^\smq_{/E}$ of $\sfP\Sch'_{/E}$.

For every $U\in\Hk(T(\bA^\infty))$, we have the natural morphism
\[
\varsigma(x)_{L,U}\colon Z(x)_{L,U}\to X_{L,U}
\]
by forgetting $(\nu_1,\dots,\nu_r)$, which is finite and unramified. For a morphism $(h,t)\in L\backslash H(\bA_F^\infty)/L'\times U\backslash T(\bA^\infty)/U'$, we have the induced morphism
\[
\ang{h,t}\colon Z(x)_{L,U}\to Z(h^{-1}x)_{L',U'}
\]
sending $(A_0,\lambda_0,\bar\eta_0;A,\lambda,\bar\eta;\nu_1,\dots,\nu_r)$ to $(A_0,\lambda_0,\bar\eta_0\circ t;A,\lambda,\bar\eta\circ h;\nu_1,\dots,\nu_r)$, which fits into the commutative diagram
\[
\xymatrix{
Z(x)_{L,U} \ar[r]^-{\ang{h,t}}\ar[d]_-{\varsigma(x)_{L,U}} & Z(h^{-1}x)_{L',U'} \ar[d]^-{\varsigma(x)_{L',U'}} \\
X_{L,U} \ar[r]^-{\ang{h,t}} & X_{L',U'}
}
\]
in which the bottom morphism is from the functor $X_{\bullet,\bullet}$.

The natural action of $T(\bA^\infty)/U$ on $Z(x)_{L,U}$ also factors through $\Xi_U$ under which the morphism $\varsigma(x)_{L,U}$ is equivariant. Similar to $X_L$, we have the quotient scheme $Z(x)_L\coloneqq Z(x)_{L,U}/\Xi_U$, which is independent of the choice of $U$, together with the induced morphism
\begin{align}\label{eq:special_morphism_1}
\varsigma(x)_L\colon Z(x)_L\to X_L
\end{align}
fitting into the commutative diagrams
\begin{align}\label{eq:special_morphism_2}
\xymatrix{
Z(x)_{L,U} \ar[r]^-{\varsigma(x)_{L,U}}\ar[d]_-{\pi(x)_{L,U}} & X_{L,U} \ar[d]^-{\pi_{L,U}} \\
Z(x)_L \ar[r]^-{\varsigma(x)_L} & X_L
}
\end{align}
and
\[
\xymatrix{
Z(x)_L \ar[r]^-{\ang{h}}\ar[d]_-{\varsigma(x)_L} & Z(h^{-1}x)_{L'} \ar[d]^-{\varsigma(x)_{L'}} \\
X_L \ar[r]^-{\ang{h}} & X_{L'}
}
\]
in which the bottom morphism is from the functor $X_\bullet$. The morphism $\varsigma(x)_L$ is finite and unramified, and depends only on the coset $Lx$.

\begin{remark}\label{re:special_morphism}
The scheme $Z(x)_L$ is nonempty only when $T(x)$ belongs to $\Herm_r(F)^+$ whose rank equals $\rank x$. In this case, $\dim X_L-\dim Z(x)_L=\rank x$.
\end{remark}

\begin{remark}\label{re:special_morphism_shimura}
We discuss the relation between $\varsigma(x)_L$ and the special morphisms of Shimura varieties when $T(x)\in\Herm_r^\circ(F)^+$. Take an element $u\in\Hom(E,\dC)$ with $v\coloneqq v(u)\in\fV_F^{(\infty)}$, and fix the data as in Remark \ref{re:moduli_shimura}. We may choose a tuple $y=(y_1,\dots,y_r)\in V(v)^r$ and an element $h\in H(\bA_F^\infty)$ such that $y_i=hx_i$ holds in $V^\infty$ for $1\leq i\leq r$. Denote by $V(v)^y$ the orthogonal complement of $y$ in $V(v)$ and put $H(v)^y\coloneqq\rU(V(v)^y)$. Then the composite morphism
\[
Z(x)_L\otimes_{E,u}u(E)\xrightarrow{\ang{h}^{-1}}Z(y)_{hLh^{-1}}\otimes_{E,u}u(E)
\xrightarrow{\varsigma(y)_{hLh^{-1}}\otimes_{E,u}u(E)}X_{hLh^{-1}}\otimes_{E,u}u(E)
\xrightarrow{\eqref{eq:moduli_shimura}}\Sh(H(v))_{hLh^{-1}}^u
\]
factors through the morphism
\[
\Sh(H(v)^y)^u_{(hLh^{-1})^y} \to \Sh(H(v))_{hLh^{-1}}^u
\]
functorial in $L$ (uniquely), where $(hLh^{-1})^y\coloneqq hLh^{-1}\cap H(v)^y(\bA^\infty)$, and induces an isomorphism
\[
Z(x)_L\otimes_{E,u}u(E)\xrightarrow\sim\Sh(H(v)^y)^u_{(hLh^{-1})^y}.
\]
As a consequence, we have the following commutative diagram
\[
\xymatrix{
Z(x)_L\otimes_{E,u}u(E) \ar[rr]^-{\varsigma(x)_L\otimes_{E,u}u(E)}\ar[d]_-\simeq && X_L\otimes_{E,u}u(E) \ar[d]_-\simeq^-{\eqref{eq:moduli_shimura}} \\
\Sh(H(v)^y)^u_{(hLh^{-1})^y} \ar[r]& \Sh(H(v))_{hLh^{-1}}^u  \ar[r]^-{\ang{h}} & \Sh(H(v))_L^u
}
\]
in $\Sch^\smq_{/u(E)}$.
\end{remark}

We now construct, when $T(x)\in\Herm_r(F)^+$, a canonical map
\begin{align}\label{eq:special_morphism_3}
\mu(x)_L\colon\varsigma(x)_L^*\rM_L^\vee\to\mathbbm{1}_{Z(x)_L,E}^{\oplus r}
\end{align}
in $\Mot(Z(x)_L,E)$. By \eqref{eq:special_morphism_2}, it suffices to construct a map $\mu(x)_{L,U}\colon\varsigma(x)_{L,U}^*\pi_{L,U}^*\rM_L^\vee\to\mathbbm{1}_{Z(x)_{L,U},E}^{\oplus r}$ compatible with the descent data, that is, a map in $\Mot(Z(x)_{L,U}/\!\!/\Xi_U,E)$. Since $\pi_{L,U}$ is \'{e}tale, the natural map $\pi_{L,U}^*\rM_L^\vee\to(\pi_{L,U}^*\rM_L)^\vee$ is an isomorphism. Thus, we have
\begin{align*}
\varsigma(x)_{L,U}^*\pi_{L,U}^*\rM_L^\vee&\simeq
\varsigma(x)_{L,U}^*(\pi_{L,U}^*\rM_L)^\vee \\
&\simeq\varsigma(x)_{L,U}^*\rM_{L,U}^\vee \\
&\simeq\varsigma(x)_{L,U}^*\(\sfh_1(A_{0,L,U}/X_{L,U},\dQ)\otimes_{E_0}\sfh^1(A_{L,U}/X_{L,U},\dQ)\) \\
&\simeq\varsigma(x)_{L,U}^*\sfh_1(A_{0,L,U}/X_{L,U},\dQ)\otimes_{E_0}\varsigma(x)_{L,U}^*\sfh^1(A_{L,U}/X_{L,U},\dQ) \\
&\simeq\varsigma(x)_{L,U}^*\sfh_1(A_{0,L,U}/X_{L,U},\dQ)\otimes_{E_0}\sfh^1(\varsigma(x)_{L,U}^*A_{L,U}/Z(x)_{L,U},\dQ).
\end{align*}
Note that over $Z(x)_{L,U}$, we have homomorphisms $\nu_i\colon\varsigma(x)_{L,U}^*A_{0,L,U}\to\varsigma(x)_{L,U}^*A_{L,U}$ for $1\leq i\leq r$, which induce maps
\begin{align*}
\nu_i^*&\colon\sfh^1(\varsigma(x)_{L,U}^*A_{L,U}/Z(x)_{L,U},\dQ)
\to\sfh^1(\varsigma(x)_{L,U}^*A_{0,L,U}/Z(x)_{L,U},\dQ)\otimes_\dQ F
\end{align*}
in $\Mot(Z(x)_{L,U},E)$. Taking sum, we obtain a map
\begin{align*}
\varsigma(x)_{L,U}^*\pi_{L,U}^*\rM_L^\vee&\to
\varsigma(x)_{L,U}^*\sfh_1(A_{0,L,U}/X_{L,U},\dQ)\otimes_{E_0}\(\sfh^1(\varsigma(x)_{L,U}^*A_{0,L,U}/Z(x)_{L,U},\dQ)\otimes_\dQ F\)^{\oplus r} \\
&\xrightarrow\sim\varsigma(x)_{L,U}^*\(\sfh_1(A_{0,L,U}/X_{L,U},\dQ)\otimes_{E_0}\sfh^1(A_{0,L,U}/Z(x)_{L,U},\dQ)\otimes_\dQ F\)^{\oplus r}
\end{align*}
Composing with the evaluation map
\[
\sfh_1(A_{0,L,U}/X_{L,U},\dQ)\otimes_{E_0}\sfh^1(A_{0,L,U}/Z(x)_{L,U},\dQ)\to\mathbbm{1}_{X_{L,U},E_0},
\]
we obtain
\[
\mu(x)_{L,U}\colon\varsigma(x)_{L,U}^*\pi_{L,U}^*\rM_L^\vee\to\varsigma(x)_{L,U}^*\mathbbm{1}_{X_{L,U},E}^{\oplus r}
=\mathbbm{1}_{Z(x)_{L,U},E}^{\oplus r}.
\]
It is straightforward to check that $\mu(x)_{L,U}$ is compatible with descent data for $\Xi_U$, hence descends to a map $\mu(x)_L$ \eqref{eq:special_morphism_3} we look for.

\begin{remark}
The special morphism $\varsigma(x)_L$ depends only on the span of components of $x$. However, the map $\mu(x)$ indeed depends on $x$.
\end{remark}

\subsection{Automorphic motives over unitary Shimura varieties}
\label{ss:motive_automorphic}

Take a hermitian $E$-weight $\lambda$ of rank $m$ (Definition \ref{de:weight}), and recall from Definition \ref{de:weight}(3) that $\dE^\lambda\subseteq\dE$ is the field of definition of $\lambda$. In this section, we define a Cartesian functor $\rM_{\bullet,\lambda}\colon\Hk(H(\bA_F^\infty))\to\Mot(\Sch^\smq_{/E},\dE^\lambda)$ over $X_\bullet$.

\begin{notation}\label{no:place_decompostion}
For an $E$-linear object $M$ in a pseudo-abelian category $\sfC$, we have a unique decomposition
\[
M\otimes_\dQ\dE=\bigoplus_{u\in\Hom(E,\dE)}M[u]
\]
in the category of $\dE$-linear objects in $\sfC$, in which $M[u]$ is the direct summand on which $E$ acts via the embedding $u\colon E\to\dE$.
\end{notation}

Put $b=(b_u\coloneqq b(\lambda_u^+))_u\in\dN^{\Hom(E,\dC)}$ and $b_v\coloneqq(b_u)_{v(u)=v}$ for every $v\in\fV_F^{(\infty)}$. Similarly, we have $\ol{b}$ and $\ol{b}_v$ for $\ol\lambda$ instead of $\lambda$.

Take an object $L\in\Hk(H(\bA_F^\infty))$. Take an element $v\in\fV_F^{(\infty)}$. Put
\[
\tT_{b_v}(\rM_L)\coloneqq\bigotimes_{v(u)=v}\tT^{b_u}(\rM_L[u])
\]
(tensor over $\dE$) as an object in $\Mot(X_L,\dE)$, which admits a natural action by $\fS_{b_v}$ by switching factors. The evaluation map \eqref{eq:evaluation_motive_standard} induces an evaluation map
\[
\bigotimes_{v(u)=v}\rM_L[u]\to\mathbbm{1}_{X_L,\dE},
\]
(tensor over $\dE$) which further induces
\begin{itemize}
  \item an evaluation map
      \begin{align}\label{eq:evaluation_motive_automorphic_0}
      \tT_{b_v}(\rM_L)\otimes_\dE\tT_{\ol{b}_v}(\rM_L)\to\mathbbm{1}_{X_L,\dE};
      \end{align}

  \item for every pair $I=(i_u)_{v(u)=v}\in\prod_{v(u)=v}\{1,\dots,b_u\}$, a contraction map
      \[
      C_I\colon \tT_{b_v}(\rM_L)\to\tT_{b_v-1}(\rM_L)
      \]
      by evaluating $i_u$-th factors of $\tT^{b_u}(\rM_L[u])$ for $v(u)=v$ (set $\tT_{b_v-1}(\rM_L)=0$ if $b_u-1<0$ for some $u$).
\end{itemize}

\begin{lem}\label{le:harmonic_motive}
The following statements hold.
\begin{enumerate}
  \item The map
      \[
      \sum_{I}C_I\colon \tT_{b_v}(\rM_L)\to\bigoplus_{\prod_{v(u)=v}\{1,\dots,b_u\}}\tT_{b_v-1}(\rM_L)
      \]
      admits kernel in the category $\Mot(X_L,\dE)$, which we denote by $\tT_{[b_v]}(\rM_L)$.

  \item The induced map
      \[
      \tT_{b_v}(\rM_L)\xrightarrow{\sim}\tT_{\ol{b}_v}(\rM_L)^\vee\to\tT_{[\ol{b}_v]}(\rM_L)^\vee
      \]
      admits kernel as well, which we denote by $\tT_{[b_v]}^\perp(\rM_L)$.

  \item The induced map
      \[
      \tT_{[b_v]}(\rM_L)\oplus\tT_{[b_v]}^\perp(\rM_L)\to\tT_{b_v}(\rM_L)
      \]
      is an isomorphism.
\end{enumerate}
\end{lem}

\begin{proof}
For the proof of this lemma, we need to use \'{e}tale realization in the next section until Definition \ref{de:harmonic_etale} and the conservativity property for abelian motives in \cite{Wil18}.

Choose a rational prime $\ell$ and take an element $u\in\Hom(E,\dE)$ with $v(u)=v$. Similar to Remark \ref{re:representation_unitary}(2) (with $\dL=\dE_\ell$), the map
\begin{align}\label{eq:harmonic_motive_1}
(\wedge^m \rM_{L,\ell}[u])\otimes_{\dE_\ell} \rM_{L,\ell}[u]^\vee\to\wedge^{m-1}\rM_{L,\ell}[u]
\end{align}
(see \S\ref{ss:realization_etale} for the notation) given by pairing the last two factors induces an isomorphism
\begin{align}\label{eq:harmonic_motive_2}
\rM_{L,\ell}[u]^\vee\xrightarrow\sim (\wedge^{m-1}\rM_{L,\ell}[u])\otimes_{\dE_\ell}(\wedge^m\rM_{L,\ell}[u])^{\vee},
\end{align}
which further induces an isomorphism
\begin{align}\label{eq:harmonic_motive_3}
\tT_{b_v}(\rM_{L,\ell})\xrightarrow\sim \tT^{b_u}(\rM_{L,\ell}[u])\otimes_{\dE_\ell}\tT^{b_{\ol{u}}}(\wedge^{m-1}\rM_{L,\ell}[u])
\otimes_{\dE_\ell}(\wedge^m\rM_{L,\ell}[u])^{-b_{\ol{u}}},
\end{align}
so that $\tT_{[b_v]}(\rM_{L,\ell})$ (Definition \ref{de:harmonic_etale}) is the kernel of the idempotent $\delta_{\lambda_u,\lambda_{\ol{u}}}\in\dQ[\fS_{b_u+(m-1)b_{\ol{u}}}]$ (as in Remark \ref{re:representation_unitary}(2)) via the action of $\fS_{b_u+(m-1)b_{\ol{u}}}$ on $\tT^{b_u+(m-1)b_{\ol{u}}}(\rM_{L,\ell}[u])\otimes_{\dE_\ell}(\wedge^m\rM_{L,\ell}[u])^{-b_{\ol{u}}}$. In particular, $\tT_{[b_v]}(\rM_{L,\ell})$ is a direct summand of $\tT_{b_v}(\rM_{L,\ell})$.

Now we repeat the same process for motives. We have the similar map
\begin{align}\label{eq:harmonic_motive_4}
(\wedge^m \rM_L[u])\otimes_\dE \rM_L[u]^\vee\to\wedge^{m-1}\rM_L[u]
\end{align}
in the category $\Mot(X_L,\dE)$, which induces a map
\begin{align}\label{eq:harmonic_motive_5}
\rM_L[u]^\vee\to(\wedge^{m-1}\rM_L[u])\otimes_\dE(\wedge^m\rM_L[u])^{\vee},
\end{align}
which further induces a map
\begin{align}\label{eq:harmonic_motive_6}
\tT_{b_v}(\rM_L)\to\tT^{b_u}(\rM_L[u])\otimes_\dE\tT^{b_{\ol{u}}}(\wedge^{m-1}\rM_L[u])\otimes_\dE(\wedge^m\rM_L[u])^{-b_{\ol{u}}}
\end{align}
(here, $(\wedge^m\rM_L[u])^{-b_{\ol{u}}}$ really means $\tT^{b_{\ol{u}}}(\wedge^m\rM_L[u])^\vee$). By Lemma \ref{le:realization_etale}, the ($\ell$-adic) \'etale realizations of \eqref{eq:harmonic_motive_4}, \eqref{eq:harmonic_motive_5} and \eqref{eq:harmonic_motive_6} are canonically \eqref{eq:harmonic_motive_1}, \eqref{eq:harmonic_motive_2} and \eqref{eq:harmonic_motive_3}, respectively. By \cite{Wil18}*{Theorem~4.4}, \eqref{eq:harmonic_motive_6} is an isomorphism since \eqref{eq:harmonic_motive_3} is. Note that the target of \eqref{eq:harmonic_motive_6} is naturally a direct summand of $\tT^{b_u+(m-1)b_{\ol{u}}}(\rM_L[u])\otimes_\dE(\wedge^m\rM_L[u])^{-b_{\ol{u}}}$. Now we define $\tT_{[b_v]}(\rM_L)$ to be the kernel of $\delta_{\lambda_u,\lambda_{\ol{u}}}\in\dQ[\fS_{b_u+(m-1)b_{\ol{u}}}]$. Using \cite{Wil18}*{Theorem~4.4} and the situation for ($\ell$-adic) \'etale realizations, we know that the natural map $\tT_{[b_v]}(\rM_L)\to\tT^{b_u+(m-1)b_{\ol{u}}}(\rM_L[u])\otimes_\dE(\wedge^m\rM_L[u])^{-b_{\ol{u}}}$ factors through $\tT_{b_v}(\rM_L)$; and that the induced map $\tT_{[b_v]}(\rM_L)\to\tT_{b_v}(\rM_L)$ is the kernel of $\sum_{I}C_I$. Part (1) is proved. Part (2) and (3) follows from a similar argument using the obvious case for ($\ell$-adic) \'etale realizations and \cite{Wil18}*{Theorem~4.4}.
\end{proof}

It is clear that the action of $\fS_{b_v}$ on $\tT_{b_v}(\rM_L)$ induces an action on $\tT_{[b_v]}(\rM_L)$. Applying the idempotent $\pi_{\lambda_v}\coloneqq\otimes_{v(u)=v}\pi_{\lambda_u}\in\dQ[\fS_{b_v}]$, we obtain a direct summand
\[
\tS_{[\lambda_v]}(\rM_L)\coloneqq\pi_{\lambda_v}\(\tT_{[b_v]}(\rM_L)\)
\]
of $\tT_{[b_v]}(\rM_L)$.

The remaining steps are motivic analogue of Construction \ref{co:harmonic_function}(1,2). Taking into account all $v\in\fV_F^{(\infty)}$, we define
\[
\tS_{[\lambda]}(\rM_L)\coloneqq\bigotimes_{v\in\fV_F^{(\infty)}}\tS_{[\lambda_v]}(\rM_L)
\]
(tensor over $\dE$) as an object of $\Mot(X_L,\dE)$. For every $\gamma\in\Gal(\dE/\dQ)$, its action on $\dE$ induces an isomorphism
\[
\gamma_{[\lambda_v]}\colon\tS_{[\lambda_v]}(\rM_L)\otimes_{\dE,\gamma}\dE\xrightarrow\sim\tS_{[(\gamma\lambda)_{\gamma v}]}(\rM_L)
\]
for every $v\in\fV_F^{(\infty)}$, which further induces an isomorphism
\[
\gamma_{[\lambda]}\colon\tS_{[\lambda]}(\rM_L)\otimes_{\dE,\gamma}\dE\xrightarrow\sim\tS_{[\gamma\lambda]}(\rM_L)
\]
by taking the tensor products over all $v\in\fV_F^{(\infty)}$. In particular, for $\gamma\in\Gal(\dE/\dE^\lambda)$, $\gamma_{[\lambda]}$ is an $\dE^\lambda$-linear automorphism of $\tS_{[\lambda]}(\rM_L)$. It is straightforward to check that this defines a homomorphism from $\Gal(\dE/\dE^\lambda)$ to $\Aut_{\dE^\lambda}(\tS_{[\lambda]}(\rM_L))$. Thus, we may define
\[
\rM_{L,\lambda}\coloneqq\tS_{[\lambda]}(\rM_L)^{\Gal(\dE/\dE^\lambda)}\in\Mot(X_L,\dE^\lambda),
\]
so that the natural map $\rM_{L,\lambda}\otimes_{\dE^\lambda}\dE\to\tS_{[\lambda]}(\rM_L)$ is an isomorphism. For general $\gamma\in\Gal(\dE/\dQ)$, the isomorphism $\gamma_{[\lambda]}$ induces an isomorphism
\[
\gamma_\lambda\colon\rM_{L,\lambda}\otimes_{\dE^\lambda,\gamma}\dE^{\gamma\lambda}\xrightarrow\sim\rM_{L,\gamma\lambda}
\]
in the category $\Mot(X_L,\dE^{\gamma\lambda})$.

\begin{definition}\label{de:motive_automorphic}
We define the automorphic (relative) motive of weight $\lambda$ over $X_L$ to be $\rM_{L,\lambda}$. It is straightforward to check that
\begin{align*}
\rM_{\bullet,\lambda}\colon\Hk(H(\bA_F^\infty))\to\Mot(\Sch^\smq_{/E},\dE^\lambda)
\end{align*}
is a Cartesian functor over $X_\bullet$.
\end{definition}

\begin{remark}
The construction of the functor $\rM_{\bullet,\lambda}$ a priori depends on the choice of $E_0$. We expect that it does not; however, we do not know how to prove yet.
\end{remark}

Lemma \ref{le:harmonic_motive} implies that the map
\[
\tT_{[b_v]}(\rM_L)\otimes_\dE\tT_{[\ol{b}_v]}(\rM_L)\to\mathbbm{1}_{X_L,\dE}
\]
induced from \eqref{eq:evaluation_motive_automorphic_0} is also an evaluation map. Since the above map is $\fS_{b_v}$-invariant, the induced map
\[
\tS_{[\lambda_v]}(\rM_L)\otimes_\dE\tS_{[\ol\lambda_v]}(\rM_L)\to\mathbbm{1}_{X_L,\dE}
\]
is an evaluation map. Taking products over all $v\in\fV_F^{(\infty)}$, the induced map
\[
\tS_{[\lambda]}(\rM_L)\otimes_\dE\tS_{[\ol\lambda]}(\rM_L)\to\mathbbm{1}_{X_L,\dE}
\]
is an evaluation map. Since the above map is $\Gal(\dE/\dE^\lambda)$-equivariant, the induced map
\begin{align}\label{eq:evaluation_motive_automorphic}
e_{L,\lambda}\colon \rM_{L,\lambda}\otimes_{\dE^\lambda}\rM_{L,\ol\lambda}\to\mathbbm{1}_{X_L,\dE^\lambda}
\end{align}
is an evaluation map.

In the rest of this section, we will relate \eqref{eq:evaluation_motive_automorphic} to the Poincar\'{e} duality pairing on the (product of) universal abelian scheme, which will be used in \S\ref{ss:height_pairing} for the definition of height pairing.

Choose an object $U\in\Hk(T(\bA^\infty))$. Put $\rM_{L,U,\lambda}\coloneqq\pi_{L,U}^*\rM_{L,\lambda}$. Then $\rM_{L,U,\lambda}$ can be obtained from $\rM_{L,U}$ via the same process as $\rM_L$. In particular, $\rM_{L,U,\lambda}$ is naturally an $\dE^\lambda$-linear direct summand of
\[
\bigotimes_{u\in\Hom(E,\dE)}\tT^{b_u}(\rM_{L,U}[u])
\]
(tensor over $\dE$), and we have an evaluation map
\[
e_{L,U,\lambda}\colon\rM_{L,U,\lambda}\otimes_{\dE^\lambda}\rM_{L,U,\ol\lambda}\to\mathbbm{1}_{X_{L,U},\dE^\lambda}
\]
as the restriction of the evaluation map
\[
\bigotimes_{u\in\Hom(E,\dE)}\tT^{b_u}(\rM_{L,U}[u])\otimes\bigotimes_{u\in\Hom(E,\dE)}\tT^{b_u}(\rM_{L,U}^\vee[u])
\to\mathbbm{1}_{X_{L,U},\dE}
\]
(tensor over $\dE$) induced from \eqref{eq:evaluation_motive_standard_2}. The latter is in turn the restriction of the evaluation map
\[
\rM_{L,U}^{\otimes|b|}\otimes_\dQ\rM_{L,U}^{\vee\otimes|b|}\to\mathbbm{1}_{X_{L,U},\dQ}
\]
after base change to $\dE$, where $|b|\coloneqq\sum_{u\in\Hom(E,\dE)}b_u$. By \eqref{eq:evaluation_motive_standard_1}, $\rM_{L,U}$ is a natural direct summand of
\[
\sfh^1(A_{0,L,U}/X_{L,U},\dQ)\otimes_\dQ\sfh^1(A_{L,U}^\vee/X_{L,U},\dQ)\ang{1}\in\Mot(X_{L,U},\dQ),
\]
so that $\rM_{L,U}^{\otimes|b|}$ is a natural direct summand of
\[
\sfh^1(A_{0,L,U}/X_{L,U},\dQ)^{\otimes|b|}\otimes_\dQ\sfh^1(A_{L,U}^\vee/X_{L,U},\dQ)^{\otimes|b|}\ang{|b|}.
\]
The latter is a natural direct summand of
\[
\sfh^{2|b|}(B_{L,U}^{|b|}/X_{L,U},\dQ)\ang{|b|},
\]
where
\[
B_{L,U}\coloneqq A_{0,L,U}\times_{X_{L,U}}A_{L,U}^\vee
\]
and
\begin{align}\label{eq:kuga_sato}
B_{L,U}^n\coloneqq\underbrace{B_{L,U}\times_{X_{L,U}}\cdots\times_{X_{L,U}}B_{L,U}}_n
\end{align}
for $n\geq 0$. Note that the relative dimension of $B_{L,U}/X_{L,U}$ is $g\coloneqq m[F:\dQ]+1$.

\begin{definition}\label{de:straightening}
Define two maps in the category $\Mot(X_{L,U},\dE^\lambda)$:
\begin{enumerate}
  \item the \emph{straightening map}
      \[
      \st_{L,U}^{\lambda,+}\colon\rM_{L,U,\lambda}\to\sfh^{2|b|}(B_{L,U}^{|b|}/X_{L,U},\dQ)\ang{|b|}\otimes_\dQ\dE^\lambda
      \]
      that is induced from the above discussion, and

  \item the \emph{complimentary straightening map}
      \[
      \st_{L,U}^{\lambda,-}\coloneqq
      ((\lambda_0\times\lambda^{-1})^{|b|})^*\circ\tF^{2(g-1)|b|}_{B_{L,U}^{|b|}/X_{L,U}}\ang{|b|}\circ\st_{L,U}^{\lambda,+}
      \colon\rM_{L,U,\lambda}\to\sfh^{2(g-1)|b|}(B_{L,U}^{|b|}/X_{L,U},\dQ)\ang{(g-1)|b|}\otimes_\dQ\dE^\lambda,
      \]
      where $\tF^{2(g-1)|b|}_{B_{L,U}^{|b|}/X_{L,U}}$ stands for the relative Fourier transform \eqref{eq:fourier}.
\end{enumerate}
\end{definition}

\begin{lem}\label{le:straightening}
Let the notation be as above. We have
\[
\pi_{L,U}^*e_{L,\lambda}=\frac{1}{|\Xi_U|}\cdot (e_B^{|b|}\otimes_\dQ\dE^\lambda) \circ(\st_{L,U}^{\lambda,+}\otimes\st_{L,U}^{\ol\lambda,-})\circ(\pi_{L,U}\otimes\pi_{L,U})^*,
\]
where
\[
e_B^{|b|}\colon\sfh^{2|b|}(B_{L,U}^{|b|}/X_{L,U},\dQ)\ang{|b|}\otimes_\dQ
\sfh^{2(g-1)|b|}(B_{L,U}^{|b|}/X_{L,U},\dQ)\ang{(g-1)|b|}\to\mathbbm{1}_{X_{L,U},\dQ}
\]
denotes the Poincar\'{e} duality map.
\end{lem}

\begin{proof}
This follows from the construction and the simple observation that $e_{L,U,\lambda}\circ(\pi_{L,U}\otimes\pi_{L,U})^*=|\Xi_U|\cdot\pi_{L,U}^*e_{L,\lambda}$.
\end{proof}

\begin{remark}\label{re:abstract_motive}
If one forgets about the descent from $X_{L,U}$ to $X_L$, then the construction of the standard and the automorphic (relative) motives can be carried over to a much more general context: Let $S$ be a regular Noetherian scheme on which at least one rational prime is invertible, and consider two pairs $(A_0,\lambda_0)$ and $(A,\lambda)$ in which
\begin{itemize}
  \item $A_0$ is an object of $\AV(S,E_0)$ of relative dimension one, and $\lambda_0\colon A_0\to A_0^\vee$ is a polarization that turns the complex conjugation of $E_0$ into the Rosati involution;

  \item $A$ is an object of $\AV(S,E)$ of relative dimension $m[F:\dQ]$, and $\lambda\colon A\to A^\vee$ is a polarization that turns the complex conjugation of $E$ into the Rosati involution.
\end{itemize}
Put
\[
\rM\coloneqq\ul\Hom_{E_0}\(\sfh_1(A_0/S,\dQ),\sfh_1(A/S,\dQ)\)\in \Mot(S,E),
\]
together with the canonical evaluation map $\rM\otimes_E\ol\rM\to\mathbbm{1}_{S,E}$ similar to \eqref{eq:evaluation_motive_standard_2}. Applying the process of obtaining $\rM_{L,\lambda}$ from $\rM_L$ in \S\ref{ss:motive_automorphic} to $\rM$, we obtain an object
\[
\rM_\lambda\in\Mot(S,\dE^\lambda),
\]
which we call the \emph{abstract automorphic (relative) motive} over $S$ of weight $\lambda$ associated with $(A_0,\lambda_0)$ and $(A,\lambda)$.

Put $B\coloneqq A_0\times_SA^\vee$ and
\[
B^n\coloneqq\underbrace{B\times_S\cdots\times_S B}_n
\]
for $n\geq 0$. Note that the relative dimension of $B/S$ is $g\coloneqq m[F:\dQ]+1$. Similar to Definition \ref{de:straightening}, we have two maps
\begin{align*}
\st^{\lambda,+}&\colon\rM_\lambda\to\sfh^{2|b|}(B^{|b|}/S,\dQ)\ang{|b|}\otimes_\dQ\dE^\lambda, \\
\st^{\lambda,-}&\colon\rM_\lambda\to\sfh^{2(g-1)|b|}(B^{|b|}/S,\dQ)\ang{(g-1)|b|}\otimes_\dQ\dE^\lambda
\end{align*}
in the category $\Mot(S,\dE^\lambda)$, which further induce two maps
\begin{align*}
\st^{\lambda,+}&\colon\rH^{2r}(S,\rM_\lambda(r))\to\CH^{r+|b|}(B^{|b|})\otimes\dE^\lambda, \\
\st^{\lambda,-}&\colon\rH^{2r}(S,\rM_\lambda(r))\to\CH^{r+(g-1)|b|}(B^{|b|})\otimes\dE^\lambda
\end{align*}
on the level of motivic cohomology for every $r\geq 0$ by Remark \ref{re:relative_motivic}(1).
\end{remark}

\subsection{\'{E}tale realizations of automorphic relative motives}
\label{ss:realization_etale}

Fix a rational prime $\ell$. In this section, we compute the image of $\rM_{\bullet,\lambda}$ under the $\ell$-adic \'etale realization functor $r_\ell\colon\sfD(X_\bullet,\dE^\lambda)\to\sfD_{\et}(X_\bullet,\dE^\lambda_\ell)$ (Notation \ref{no:relative_motive_etale}).

Take an object $L\in\Hk(H(\bA_F^\infty))$. We first consider the \'etale realization of the standard motive $\rM_L$.

Choose an object $U\in\Hk(T(\bA^\infty))$. We have $\ell$-adic complexes
\begin{align*}
\sfh_1^{\et}(A_{0,L,U}/X_{L,U},\dQ_\ell)&\in\sfD_{\et}(X_{L,U},E_{0,\ell}),\\
\sfh_1^{\et}(A_{L,U}/X_{L,U},\dQ_\ell)&\in\sfD_{\et}(X_{L,U},E_\ell)
\end{align*}
from Remark \ref{re:abelian_realization}(1). Put
\[
\rM_{L,U,\ell}\coloneqq\ul\Hom_{E_{0,\ell}}\(\sfh_1^{\et}(A_{0,L,U}/X_{L,U},\dQ_\ell),\sfh_1^{\et}(A_{L,U}/X_{L,U},\dQ_\ell)\)
\in\sfD_{\et}(X_{L,U},E_\ell),
\]
which is indeed an object in $\sfD_{\et}^\heartsuit(X_{L,U},E_\ell)$, that is, an $\ell$-adic sheaf (as the above internal hom is concentrated in degree zero). Similar to the case of $\rM_{L,U}$, we have an object
\[
\rM_{L,U,\ell}/\!\!/\Xi_U\in\sfD_{\et}^\heartsuit(X_{L,U}/\!\!/\Xi_U,E_\ell),
\]
where the target denotes the category of Cartesian functors $\ul\Xi_U\to\sfD_{\et}^\heartsuit(\Sch^\smq_{/E},E_\ell)$ above $X_{L,U}/\!\!/\Xi_U$. Then by the usual descent in abelian category, we obtain an object
\[
\rM_{L,\ell}\coloneqq\rM_{L,U,\ell}/\!\!/\Xi_U\in\sfD_{\et}^\heartsuit(X_{L,U}/\!\!/\Xi_U,E_\ell)
=\sfD_{\et}^\heartsuit(X_L,E_\ell),
\]
which is independent of the choice of $U$. Indeed, $\rM_{L,\ell}$ is a local system over $X_L$ of rank $m$ with coefficients in $E_\ell$. Similar to \eqref{eq:evaluation_motive_standard}, we have an evaluation map
\begin{align}\label{eq:evaluation_etale}
\rM_{L,\ell}\otimes_{E_\ell}\ol{\rM_{L,\ell}}\to\mathbbm{1}^\et_{X_L,E_\ell},
\end{align}
where $\mathbbm{1}^\et_{X_L,E_\ell}$ denotes the unit object of $\sfD_{\et}(X_L,E_\ell)$.

\begin{lem}\label{le:realization_etale}
The \'etale realization $r_\ell\rM_L$ is canonically isomorphic to $\rM_{L,\ell}$.
\end{lem}

\begin{proof}
It suffices to construct, for $U\in\Hk(T(\bA^\infty))$, a canonical isomorphism $r_\ell\rM_{L,U}\simeq\rM_{L,U,\ell}$ that is $\Xi_U$-equivariant. From Remark \ref{re:abelian_realization}(1), we have
\[
r_\ell\sfh_1(A_{0,L,U}/X_{L,U},\dQ)\simeq\sfh_1^{\et}(A_{0,L,U}/X_{L,U},\dQ_\ell),\qquad
r_\ell\sfh_1(A_{L,U}/X_{L,U},\dQ)\simeq\sfh_1^{\et}(A_{L,U}/X_{L,U},\dQ_\ell)
\]
canonically. Since $\sfh_1(A_{0,L,U}/X_{L,U},\dQ)$ is compact, we have a canonical isomorphism
\begin{align*}
r_\ell\rM_{L,U}&=r_\ell\ul\Hom_{E_0}\(\sfh_1(A_{0,L,U}/X_{L,U},\dQ),\sfh_1(A_{L,U}/X_{L,U},\dQ)\) \\
&\simeq\ul\Hom_{E_{0,\ell}}\(\sfh_1^{\et}(A_{0,L,U}/X_{L,U},\dQ_\ell),\sfh_1^{\et}(A_{L,U}/X_{L,U},\dQ_\ell)\)=\rM_{L,U,\ell}
\end{align*}
by \cite{Ayo14}*{Th\'eor\`eme~6.6}. It is clear that the above isomorphism is $\Xi_U$-equivariant, which implies the lemma.
\end{proof}

Take an element $v\in\fV_F^{(\infty)}$. Put
\[
\tT_{b_v}(\rM_{L,\ell})\coloneqq\bigotimes_{v(u)=v}\tT^{b_u}(\rM_{L,\ell}[u])
\]
(tensor over $\dE_\ell$), which is an object of $\sfD_{\et}^\heartsuit(X_L,\dE_\ell)$ and admits a natural action by $\fS_{b_v}$ by switching factors. The evaluation map \eqref{eq:evaluation_etale} induces an evaluation map
\[
\bigotimes_{v(u)=v}\rM_{L,\ell}[u]\to\mathbbm{1}^\et_{X_L,\dE_\ell},
\]
(tensor over $\dE_\ell$) which further induces
\begin{itemize}
  \item an evaluation map
      \begin{align*}
      \tT_{b_v}(\rM_{L,\ell})\otimes_{\dE_\ell}\tT_{\ol{b}_v}(\rM_{L,\ell})\to\mathbbm{1}^\et_{X_L,\dE_\ell};
      \end{align*}

  \item for every pair $I=(i_u)_{v(u)=v}\in\prod_{v(u)=v}\{1,\dots,b_u\}$, a contraction map
      \[
      C_{I,\ell}\colon \tT_{b_v}(\rM_{L,\ell})\to\tT_{b_v-1}(\rM_{L,\ell})
      \]
      by evaluating $i_u$-th factors of $\tT^{b_u}(\rM_{L,\ell}[u])$ for $v(u)=v$ (set $\tT_{b_v-1}(\rM_{L,\ell})=0$ if $b_u-1<0$ for some $u$).
\end{itemize}
It follows from Lemma \ref{le:realization_etale} that we have a canonical isomorphism
\begin{align*}
r_\ell C_I\simeq C_{I,\ell}
\end{align*}
for every pair $I\in\prod_{v(u)=v}\{1,\dots,b_u\}$.

\begin{definition}\label{de:harmonic_etale}
Define $\tT_{[b_v]}(\rM_{L,\ell})$ to be the kernel of the natural map
\[
\sum_{I}C_{I,\ell}\colon \tT_{b_v}(\rM_{L,\ell})\to\bigoplus_{\prod_{v(u)=v}\{1,\dots,b_u\}}\tT_{b_v-1}(\rM_{L,\ell})
\]
in the abelian category $\sfD_{\et}^\heartsuit(X_L,\dE_\ell)$, which admits a natural action by $\fS_{b_v}$. Similar to Remark \ref{re:representation_unitary}(2), $\tT_{[b_v]}(\rM_{L,\ell})$ is a direct summand of $\tT_{b_v}(\rM_{L,\ell})$, called the \emph{harmonic direct summand}.
\end{definition}

By Lemma \ref{le:harmonic_motive}, we have a canonical isomorphism
\[
r_\ell\tT_{[b_v]}(\rM_L)\simeq\tT_{[b_v]}(\rM_{L,\ell}).
\]
By the similar process of constructing $\tS_{[\lambda_v]}(\rM_L)$, $\tS_{[\lambda]}(\rM_L)$, and $\rM_{L,\lambda}$ from $\tT_{[b_v]}(\rM_L)$, we obtain the $\ell$-adic analogues $\tS_{[\lambda_v]}(\rM_{L,\ell}),\tS_{[\lambda]}(\rM_{L,\ell})\in\sfD_{\et}^\heartsuit(X_L,\dE_\ell)$, and $\rM_{L,\ell,\lambda}\in\sfD_{\et}^\heartsuit(X_L,\dE^\lambda_\ell)$ from $\tT_{[b_v]}(\rM_{L,\ell})$, with canonical isomorphisms:
\begin{align}\label{eq:realization_etale}
r_\ell\tS_{[\lambda_v]}(\rM_L)\simeq\tS_{[\lambda_v]}(\rM_{L,\ell}),\qquad
r_\ell\tS_{[\lambda]}(\rM_L)\simeq\tS_{[\lambda]}(\rM_{L,\ell}),\qquad
r_\ell\rM_{L,\lambda}\simeq\rM_{L,\ell,\lambda}.
\end{align}

\begin{remark}
The evaluation map $e_{L,\lambda}$ \eqref{eq:evaluation_motive_automorphic} induces an evaluation map
\begin{align}\label{eq:evaluation_etale_automorphic}
e_{L,\ell,\lambda}\colon\rM_{L,\ell,\lambda}\otimes_{\dE^\lambda_\ell}\rM_{L,\ell,\ol\lambda}
\to\mathbbm{1}^{\et}_{X_L,\dE^\lambda_\ell}
\end{align}
in $\sfD_{\et}^\heartsuit(\Sch^\smq_{/E},\dE^\lambda_\ell)$ for every $L\in\Hk(H(\bA_F^\infty))$.
\end{remark}

To end this section, we give an alternative description of $\rM_{L,\ell,\lambda}$. Choose an $O_{E,\ell}$-lattice $\Lambda_\ell$ of $V_\ell$ such that $L_\ell$ is contained in the stabilizer of $\Lambda_\ell$ in $H(F_\ell)$. We review below the well-known construction of constructing a $O_{E,\ell}$-linear local system over $X_L$ from $\Lambda_\ell$.

\begin{construction}\label{co:local_system}
For every integer $n\geq 1$, denote by $L_n\subseteq L$ the subgroup that induces the trivial action on $\Lambda_\ell/\ell^n\Lambda_\ell$ (via $\ell$-adic components), so that we have an $O_{E,\ell}$-linear sheaf
\[
\Lambda_{\ell,n,L}\coloneqq\(X_{L_n}\times\Lambda_\ell/\ell^n\Lambda_\ell\)/(L/L_n)
\]
over $X_L$, in which the action is given by $(x,v).h=(xh,h^{-1}v)$. Then we obtain an obvious projective system $(\Lambda_{\ell,n,L})_n$ of $O_{E,\ell}$-linear sheaves over $X_L$; and we put
\[
\Lambda_{\ell,L}\coloneqq\varprojlim_{n\geq 1}\Lambda_{\ell,n,L},
\]
which is an $O_{E,\ell}$-linear local system over $X_L$.
\end{construction}

It is well-known that the $E_\ell$-linear local $V_{\ell,L}\coloneqq\Lambda_{\ell,L}\otimes_{\dZ_\ell}\dQ_\ell$ on $X_L$ does not depend on the choice of $\Lambda_\ell$, and that the hermitian level structures at $\ell$-adic places in Definition \ref{de:moduli_rational} induce a canonical isomorphism $\rM_{L,\ell}\simeq V_{\ell,L}$. Repeating the same process for $\rM_{L,\ell}$, we obtain $\dE_\ell$-linear local systems $\tS_{[\lambda_v]}(V_{\ell,L})$ for $v\in\fV_F^{(\infty)}$ and $\tS_{[\lambda]}(V_{\ell,L})$, and finally the $\dE^\lambda_\ell$-linear local system $V_{\ell,L,\lambda}$ on $X_L$, with canonical isomorphisms:
\begin{align}\label{eq:realization_local_system}
\tS_{[\lambda_v]}(\rM_{L,\ell})\simeq\tS_{[\lambda_v]}(V_{\ell,L}),\qquad
\tS_{[\lambda]}(\rM_{L,\ell})\simeq\tS_{[\lambda]}(V_{\ell,L}),\qquad
\rM_{L,\ell,\lambda}\simeq V_{\ell,L,\lambda}.
\end{align}

\begin{remark}\label{re:local_system_integral}
When $\ell$ avoids certain divisibility, the $\dE^\lambda_\ell$-linear local system $V_{\ell,L,\lambda}$ has an integral model that can be constructed from an $O_{\dE^\lambda,\ell}$-lattice via the similar process in Construction \ref{co:local_system}.

Indeed, suppose that $\ell\nmid D_EB_\lambda[\dE\colon\dE^\lambda]$ (see \eqref{eq:bad_lambda} for $B_\lambda$). Furthermore, suppose that $L_\ell$ is contained in the stabilizer of a self-dual lattice (Definition \ref{de:lattice}) $\Lambda_\ell$ of $V_\ell$. Applying Construction \ref{co:harmonic_function}(2) to $R=\dZ_\ell$ and $V=\Lambda_\ell$, we obtain an $O_{\dE^\lambda,\ell}$-lattice $\Lambda_{\ell,\lambda}\coloneqq(\Lambda_\ell)_\lambda$ of $V_{\ell,\lambda}\coloneqq(V_\ell)_\lambda$. Applying the similar process in Construction \ref{co:local_system}, we obtain an $O_{\dE^\lambda,\ell}$-linear local system $\Lambda_{\ell,\lambda,L}$ from $\Lambda_{\ell,\lambda}$, so that we have a canonical isomorphism
\begin{align}\label{eq:local_system_integral}
\Lambda_{\ell,\lambda,L}\otimes_{\dZ_\ell}\dQ_\ell\simeq V_{\ell,L,\lambda}
\end{align}
of $\dE^\lambda_\ell$-linear local systems over $X_L$.
\end{remark}

\subsection{Hodge realizations of automorphic relative motives}
\label{ss:realization_hodge}

Fix an element $\bbu\in\Hom(E,\dC)$. In this section, we compute the image of $\rM_{\bullet,\lambda}$ under the $\bbu$-Hodge realization functor $r^{\bbu}\colon\sfD(X_\bullet,\dE^\lambda)\to\sfD_{\hs}(X_\bullet^{\bbu},\dE^\lambda_\infty)$ (Notation \ref{no:relative_motive_hodge}).

Take an object $L\in\Hk(H(\bA_F^\infty))$. We first consider the $\bbu$-Hodge realization of the standard motive $\rM_L$.

Choose an object $U\in\Hk(T(\bA^\infty))$. We have complexes of relative Hodge structures
\begin{align*}
\sfh_1^{\hs}(A_{0,L,U}^{\bbu}/X_{L,U}^{\bbu},\dR)&\in\sfD_{\hs}(X_{L,U}^{\bbu},E_{0,\infty}),\\
\sfh_1^{\hs}(A_{L,U}^{\bbu}/X_{L,U}^{\bbu},\dR)&\in\sfD_{\hs}(X_{L,U}^{\bbu},E_\infty)
\end{align*}
from Remark \ref{re:abelian_realization}(2). Put
\[
\rM_{L,U}^{\bbu}\coloneqq\ul\Hom_{E_{0,\infty}}
\(\sfh_1^{\hs}(A_{0,L,U}^{\bbu}/X_{L,U}^{\bbu},\dR),\sfh_1^{\hs}(A_{L,U}^{\bbu}/X_{L,U}^{\bbu},\dR)\)
\in\sfD_{\hs}(X_{L,U}^{\bbu},E_\infty),
\]
which is indeed an object in $\sfD_{\hs}^\heartsuit(X_{L,U}^{\bbu},E_\infty)$, that is, a relative Hodge structure with coefficients in $E_\infty$ (as the above $\ul\Hom_{E_{0,\infty}}$ is concentrated in degree zero). Similar to the case of $\rM_{L,U}$, we have an object
\[
\rM_{L,U}^{\bbu}/\!\!/\Xi_U\in\sfD_{\hs}^\heartsuit(X_{L,U}^{\bbu}/\!\!/\Xi_U,E_\infty),
\]
where the target denotes the category of Cartesian functors $\ul\Xi_U\to\sfD_{\hs}^\heartsuit(\Mfd,E_\infty)$ above $X_{L,U}/\!\!/\Xi_U$. Then by the usual descent in abelian category, we obtain an object
\[
\rM_L^{\bbu}\coloneqq\rM_{L,U}^{\bbu}/\!\!/\Xi_U\in
\sfD_{\hs}^\heartsuit(X_{L,U}^{\bbu}/\!\!/\Xi_U,E_\infty)=\sfD_{\hs}^\heartsuit(X_L^{\bbu},E_\infty),
\]
which is independent of the choice of $U$. Indeed, $\rM_L^{\bbu}$ is a relative Hodge structure over $X_L^{\bbu}$ of rank $m$ with coefficients in $E_\infty$. Similar to \eqref{eq:evaluation_motive_standard}, we have an evaluation map
\begin{align}\label{eq:evaluation_hodge}
\rM_L^{\bbu}\otimes_{E_\infty}\ol{\rM_L^{\bbu}}\to\mathbbm{1}^\hs_{X_L^{\bbu},E_\infty},
\end{align}
where $\mathbbm{1}^\hs_{X_L,E_\infty}$ denotes the unit object of $\sfD_{\hs}(X_L^{\bbu},E_\infty)$.

Using \cite{Tub25}, we obtain the following lemma similar to Lemma \ref{le:realization_etale}.

\begin{lem}\label{le:realization_hodge}
The $\bbu$-Hodge realization $r^{\bbu}\rM_L$ is canonically isomorphic to $\rM_L^{\bbu}$.
\end{lem}

Take an element $v\in\fV_F^{(\infty)}$. Put
\[
\tT_{b_v}(\rM_L^{\bbu})\coloneqq\bigotimes_{v(u)=v}\tT^{b_u}(\rM_L^{\bbu}[u])
\]
(tensor over $\dE_\infty$), which is an object of $\sfD_{\hs}^\heartsuit(X_L^{\bbu},\dE_\infty)$ and admits a natural action by $\fS_{b_v}$ by switching factors. The evaluation map \eqref{eq:evaluation_hodge} induces an evaluation map
\[
\bigotimes_{v(u)=v}\rM_L^{\bbu}[u]\to\mathbbm{1}^\hs_{X_L^{\bbu},\dE_\infty},
\]
(tensor over $\dE_\infty$) which further induces
\begin{itemize}
  \item an evaluation map
      \begin{align*}
      \tT_{b_v}(\rM_L^{\bbu})\otimes_{\dE_\infty}\tT_{\ol{b}_v}(\rM_L^{\bbu})\to\mathbbm{1}^\hs_{X_L^{\bbu},\dE_\infty};
      \end{align*}

  \item for every pair $I=(i_u)_{v(u)=v}\in\prod_{v(u)=v}\{1,\dots,b_u\}$, a contraction map
      \[
      C_I^{\bbu}\colon \tT_{b_v}(\rM_L^{\bbu})\to\tT_{b_v-1}(\rM_L^{\bbu})
      \]
      by evaluating $i_u$-th factors of $\tT^{b_u}(\rM_L^{\bbu}[u])$ for $v(u)=v$ (set $\tT_{b_v-1}(\rM_L^{\bbu})=0$ if $b_u-1<0$ for some $u$).
\end{itemize}
It follows from Lemma \ref{le:realization_hodge} that we have a canonical isomorphism
\begin{align*}
r^{\bbu} C_I\simeq C_I^{\bbu}
\end{align*}
for every pair $I\in\prod_{v(u)=v}\{1,\dots,b_u\}$.

\begin{definition}\label{de:harmonic_hodge}
Define $\tT_{[b_v]}(\rM_L^{\bbu})$ to be the kernel of the natural map
\[
\sum_{I}C_I^{\bbu}\colon \tT_{b_v}(\rM_L^{\bbu})\to\bigoplus_{\prod_{v(u)=v}\{1,\dots,b_u\}}\tT_{b_v-1}(\rM_L^{\bbu})
\]
in the abelian category $\sfD_{\hs}^\heartsuit(X_L^{\bbu},\dE_\infty)$, which admits a natural action by $\fS_{b_v}$.
\end{definition}

Similar to Remark \ref{re:representation_unitary}(2), $\tT_{[b_v]}(\rM_L^{\bbu})$ is a direct summand of $\tT_{b_v}(\rM_L^{\bbu})$. Thus, by Lemma \ref{le:harmonic_motive}, we have a canonical isomorphism
\[
r^{\bbu}\tT_{[b_v]}(\rM_L)\simeq\tT_{[b_v]}(\rM_L^{\bbu}).
\]
By the similar process of constructing $\tS_{[\lambda_v]}(\rM_L)$, $\tS_{[\lambda]}(\rM_L)$, and $\rM_{L,\lambda}$ from $\tT_{[b_v]}(\rM_L)$, we obtain the $\bbu$-Hodge analogues $\tS_{[\lambda_v]}(\rM_L^{\bbu}),\tS_{[\lambda]}(\rM_L^{\bbu})\in\sfD_{\hs}^\heartsuit(X_L^{\bbu},\dE_\infty)$, and $\rM_{L,\lambda}^{\bbu}\in\sfD_{\hs}^\heartsuit(X_L^{\bbu},\dE^\lambda_\infty)$ from $\tT_{[b_v]}(\rM_L^{\bbu})$, with canonical isomorphisms:
\begin{align}\label{eq:realization_hodge}
r^{\bbu}\tS_{[\lambda_v]}(\rM_L)\simeq\tS_{[\lambda_v]}(\rM_L^{\bbu}),\qquad
r^{\bbu}\tS_{[\lambda]}(\rM_L)\simeq\tS_{[\lambda]}(\rM_L^{\bbu}),\qquad
r^{\bbu}\rM_{L,\lambda}\simeq\rM_{L,\lambda}^{\bbu}.
\end{align}

\begin{remark}
The evaluation map $e_{L,\lambda}$ \eqref{eq:evaluation_motive_automorphic} induces an evaluation map
\begin{align}\label{eq:realization_hodge_automorphic}
e_{L,\lambda}^{\bbu}\colon\rM_{L,\lambda}^{\bbu}\otimes_{\dE^\lambda_\infty}\rM_{L,\ol\lambda}^{\bbu}
\to\mathbbm{1}^{\hs}_{X_L^{\bbu},\dE^\lambda_\infty}
\end{align}
in $\sfD_{\hs}^\heartsuit(\Mfd,\dE^\lambda_\infty)$ for every $L\in\Hk(H(\bA_F^\infty))$.
\end{remark}

To end this section, we give an explicit description of the complex relative Hodge structure \cite{Liu26}*{Definition~4.3} $\kappa\rM_{L,\lambda}^{\bbu}\in\sfD_{\hs}^\heartsuit(X_L^{\bbu},\dC)$ (see Notation \ref{no:contraction} for $\kappa$) via the complex uniformization
\[
X_L^{\bbu}=H(\bbv)(F)\backslash\fD^{\bbu}\times H(\bbv)(\bA_F^\infty)/L
\]
in Remark \ref{re:moduli_shimura}, where $\bbv\coloneqq v(\bbu)$.

In \cite{Liu26}*{\S3}, we have constructed the complex representation
\[
\tS_{[\lambda]}(V(\bbv)_\infty)\coloneqq\bigotimes_{v\in\fV_F^{(\infty)}}\tS_{[\lambda_v]}(V(\bbv)_v)
\]
(tensor over $\dC$) as a representation of $H(\bbv)(F_\infty)$. Recall that
\[
\tS_{[\lambda_v]}(V(\bbv)_v)=(\otimes_{v(u)=v}\pi_{\lambda_u^+})\(\tT_{[b_v]}(V(\bbv)_v)\),
\]
where $\tT_{[b_v]}(V(\bbv)_v)$ is the \emph{harmonic direct summand} of
\[
\tT_{b_v}(V(\bbv)_v)\coloneqq\bigotimes_{v(u)=v}\tT^{b_u}(V(\bbv)_u)
\]
similar to \eqref{eq:representation_unitary}, where $V(\bbv)_u\coloneqq V(\bbv)\otimes_{E,u}\dC$.

\begin{definition}\label{de:hodge_unitary}
We define complex Hodge structures on the local system $V(\bbv)_u$ over $\fD^{\bbu}$. Put $\sV^{\bbu}_u\coloneqq V(\bbv)_u\otimes_\dC\sO_{\fD^{\bbu}}$.
\begin{itemize}
  \item For $u\in\Hom(E,\dC)\setminus\{\bbu,\ol{\bbu}\}$, define $\cV^{\bbu}_u$ to be the unique complex Hodge structure on $V(\bbv)_u$ satisfying $(\cV^{\bbu}_u)^{(0,0)}=(\sV^{\bbu}_u)^\infty$.

  \item For $u=\bbu$, denote by $V_{\fD^{\bbu}}^-$ the tautological complex subbundle of negative lines of $(\sV^{\bbu}_u)^\infty=V(\bbv)_u\otimes_\dC\sA^0_{\fD^{\bbu}}$ over (the underlying smooth manifold of) $\fD^{\bbu}$, which is a holomorphic subbundle of $(\sV^{\bbu}_u)^\infty$, and by $V_{\fD^{\bbu}}^+$ its orthogonal complement. Define $\cV^{\bbu}_u$ to be the unique complex Hodge structure on $V(\bbv)_u$ satisfying that $(\cV^{\bbu}_u)^{(0,0)}=V_{\fD^{\bbu}}^+$ and $(\cV^{\bbu}_u)^{(1,-1)}=V_{\fD^{\bbu}}^-$.\footnote{Indeed, $\cV^{\bbu}_u$ is the unique complex Hodge structure of weight zero on $V(\bbv)_u$ of which the hermitian form on $V(\bbv)_u$ induces a polarization \cite{Liu26}*{Definition~4.5}, satisfying $(\cV^{\bbu}_u)^{(p,-p)}=0$ for $|p|\geq 2$.}

  \item For $u=\ol{\bbu}$, define $\cV^{\bbu}_u$ to be the complex conjugation of $\cV^{\bbu}_{\bbu}$, which is a complex Hodge structure on $V(\bbv)_u=\ol{V(\bbv)_{\bbu}}$. Note that $\cV^{\bbu}_u$ is also canonically isomorphic to the dual of $\cV^{\bbu}_{\bbu}$ induced by the hermitian form on $V(\bbv)_u$.
\end{itemize}

Taking tensor products of $\cV^{\bbu}_u$ for $u\in\Hom(E,\dC)$ induces a complex Hodge structure $\cV^{\bbu}_\lambda$ on the local system $\tS_{[\lambda]}(V(\bbv)_\infty)$ over $\fD^{\bbu}$. The complex Hodge structures $\cV^{\bbu}_u$ and $\cV^{\bbu}_\lambda$ are all $H(\bbv)(F_\infty)$-invariant and descend to the complex unitary Shimura variety
\[
X_L^{\bbu}=H(\bbv)(F)\backslash\fD^{\bbu}\times H(\bbv)(\bA_F^\infty)/L,
\]
for which we denote by $\cV^{\bbu}_{u,L}$ and $\cV^{\bbu}_{\lambda,L}$, respectively.
\end{definition}

We now construct a canonical isomorphism
\begin{align}\label{eq:hodge_comparison_0}
\kappa\(\rM_L^{\bbu}[u]\)\xrightarrow\sim\cV^{\bbu}_{u,L}
\end{align}
of complex Hodge structures over $X_L^{\bbu}$ for every $u\in\Hom(E,\dC)$. To do this, it suffices to construct such an isomorphism over $X_{L,U}^{\bbu}$ for an arbitrarily chosen object $U\in\Hk(T(\bA^\infty))$ that is $\Xi_U$-equivariant. Recall that $\rM_{L,U}=\pi_{L,U}^*\rM_L$ and put $\cV^{\bbu}_{u,L,U}\coloneqq(\pi_{L,U}^{\bbu})^*\cV^{\bbu}_{u,L}$. Note that
\[
X_{L,U}^{\bbu}=\(H(\bbv)(F)\backslash\fD^{\bbu}\times H(\bbv)(\bA_F^\infty)/L\)\times\(T(\dQ)\backslash T(\bA^\infty)/U\).
\]
Choose a purely imaginary element $\jmath\in E_0^\times$ and denote by $W_0$ the space $E_0$ equipped with the rational skew-hermitian form $\llangle x,y\rrangle\mapsto\Tr_{E_0/\dQ} \jmath x\overline{y}$. Then $W\coloneqq W_0\otimes_{E_0}V$ has an induced rational skew-hermitian form using the hermitian form on $V$. Over the complex manifold $\fD^{\bbu}_{L,U}\coloneqq\fD^{\bbu}\times H(\bbv)(\bA_F^\infty)/L\times T(\bA^\infty)/U$, the relative homology groups $\rH_1(A_{0,L,U}^{\bbu}/\fD_{L,U}^{\bbu},\dQ)$ and $\rH_1(A_{L,U}^{\bbu}/\fD_{L,U}^{\bbu},\dQ)$ can be trivialized as
\begin{align}\label{eq:hodge_comparison_1}
W_0\simeq\rH_1(A_{0,L,U}^{\bbu}/\fD_{L,U}^{\bbu},\dQ),\quad
W\simeq\rH_1(A_{L,U}^{\bbu}/\fD_{L,U}^{\bbu},\dQ)
\end{align}
as $E_0$-linear and $E$-linear spaces equipped with rational skew-hermitian forms, respectively. Then we obtain an isomorphism
\begin{align}\label{eq:hodge_comparison_2}
\Hom_{E_0}\(\rH_1(A_{0,L,U}^{\bbu}/\fD_{L,U}^{\bbu},\dQ),\rH_1(A_{L,U}^{\bbu}/\fD_{L,U}^{\bbu},\dQ)\)\simeq V(\bbv)
\end{align}
of $E$-linear local systems over $\fD^{\bbu}_{L,U}$. Take an element $u\in\Hom(E,\dC)$. By Lemma \ref{le:hodge_comparison} below, the above isomorphism induces an isomorphism $\kappa\(\rM_{L,U}^{\bbu}[u]\)\to\cV^{\bbu}_{u,L,U}$ of complex Hodge structures over $\fD^{\bbu}_{L,U}$, which is clearly $H(\bbv)(F)\times T(\dQ)$-invariant hence descends to an isomorphism over $X_{L,U}^{\bbu}$. The resulting isomorphism $\kappa\(\rM_{L,U}^{\bbu}[u]\)\to\cV^{\bbu}_{u,L,U}$ over $X_{L,U}^{\bbu}$ does not depend on the choice of the trivialization \eqref{eq:hodge_comparison_1} and is $\Xi_U$-equivariant, so that we obtain a canonical isomorphism \eqref{eq:hodge_comparison_0}.

Finally, the isomorphisms \eqref{eq:hodge_comparison_0} for all $u\in\Hom(E,\dC)$ induce an isomorphism
\begin{align}\label{eq:hodge_comparison}
\kappa\rM_{L,\lambda}^{\bbu}\xrightarrow\sim\cV^{\bbu}_{\lambda,L}
\end{align}
in the category $\sfD_{\hs}^\heartsuit(X_L^{\bbu},\dC)$.

\begin{lem}\label{le:hodge_comparison}
For every $u\in\Hom(E,\dC)$, the isomorphism
\[
\Hom_{E_0}\(\rH_1(A_{0,L,U}^{\bbu}/\fD_{L,U}^{\bbu},\dQ),\rH_1(A_{L,U}^{\bbu}/\fD_{L,U}^{\bbu},\dQ)\)\otimes_{E,u}\dC
\xrightarrow\sim V(\bbv)_u
\]
induced from \eqref{eq:hodge_comparison_2} is an isomorphism of complex Hodge structures $\kappa\(\rM_{L,U}^{\bbu}[u]\)\to\cV^{\bbu}_{u,L,U}$ over $\fD^{\bbu}_{L,U}$.
\end{lem}

\begin{proof}
First, we realize that the underlying complex local system of $\kappa\(\rM_{L,U}^{\bbu}[u]\)$ is indeed the left-hand side of the isomorphism in question.

Recall that we have the Hodge sequence
\[
0 \to \omega_{(A_{L,U}^{\bbu})^\vee/\fD^{\bbu}_{L,U}}\to\rH_1(A_{L,U}^{\bbu}/\fD_{L,U}^{\bbu},\dC)
\to \Lie_{A_{L,U}^{\bbu}/\fD^{\bbu}_{L,U}}\to 0
\]
of complex local systems, which is the direct sum of $u$-isotypic components
\begin{align}\label{eq:hodge_comparison_3}
0 \to \omega_{(A_{L,U}^{\bbu})^\vee/\fD^{\bbu}_{L,U}}[u] \to \rH_1(A_{L,U}^{\bbu}/\fD_{L,U}^{\bbu},\dQ)\otimes_{E,u}\dC
\to \Lie_{A_{L,U}^{\bbu}/\fD^{\bbu}_{L,U}}[u] \to 0
\end{align}
for $u\in\Hom(E,\dC)$. The sequence \eqref{eq:hodge_comparison_3} gives a complex Hodge structure $\cW$ on $W_u\coloneqq W\otimes_{E,u}\dC$ over $\fD^{\bbu}_{L,U}$\footnote{Note that \eqref{eq:hodge_comparison_1} induces an isomorphism $\rH_1(A_{L,U}^{\bbu}/\fD_{L,U}^{\bbu},\dQ)\otimes_{E,u}\dC\simeq W_u\otimes_\dC\sO_{\fD^{\bbu}_{L,U}}$ of complex local systems.} of weight $-1$ satisfying
\[
\cW^{(0,-1)}=\omega_{(A_{L,U}^{\bbu})^\vee/\fD^{\bbu}_{L,U}}[u],\quad \cW^{(-1,0)}=\Lie_{A_{L,U}^{\bbu}/\fD^{\bbu}_{L,U}}[u],
\]
which is nothing but $\kappa\(\sfh_1^{\hs}(A_{L,U}^{\bbu}/X_{L,U}^{\bbu},\dR)[u]\)$.

Similarly, from $A_{0,L,U}$, we obtain a complex Hodge structure $\cW_0$ on $W_{0,u}\coloneqq W_0\otimes_{E_0,u}\dC$ over $\fD^{\bbu}_{L,U}$ of weight $-1$ satisfying
\[
\cW_0^{(0,-1)}=\omega_{(A_{0,L,U}^{\bbu})^\vee/\fD^{\bbu}_{L,U}}[u],\quad \cW_0^{(-1,0)}=\Lie_{A_{0,L,U}^{\bbu}/\fD^{\bbu}_{L,U}}[u],
\]
which is nothing but $\kappa\(\sfh_1^{\hs}(A_{0,L,U}^{\bbu}/X_{L,U}^{\bbu},\dR)[u]\)$.

Together, we obtain that $\kappa\(\rM_{L,U}^{\bbu}[u]\)=\ul\Hom_{E_0}(\cW_0,\cW)$. It is easy to see that under the canonical isomorphism $\Hom_\dC(W_{0,u},W_u)=V(\bbv)_u$, we have $\ul\Hom_{E_0}(\cW_0,\cW)=\cV^{\bbu}_{u,L,U}$ as complex Hodge structures over $\fD^{\bbu}_{L,U}$.

The lemma is proved.
\end{proof}

\begin{remark}
The various realizations of the automorphic motive $\rM_{\bullet,\lambda}$ (Definition \ref{de:motive_automorphic}) studied in this and the previous sections can be summarized in the following picture:
\[
\xymatrix{
\cV_\lambda^{\bbu} \ar[d]^-{\text{globalization}} &&&&&&&&&&
\Lambda_{\ell,\lambda} \ar@{-->}[d]_-{\text{globalization}}\\
\cV_{\lambda,\bullet}^{\bbu} && \rM_{\bullet,\lambda}^{\bbu}\ar[ll]_-{\kappa}^-{\eqref{eq:hodge_comparison}} &&
\rM_{\bullet,\lambda} \ar[ll]_-{r^{\bbu}}^-{\eqref{eq:realization_hodge}} \ar[rr]^-{r_\ell}_-{\eqref{eq:realization_etale}} &&
\rM_{\bullet,\ell,\lambda} \ar[rr]^-\simeq_-{\eqref{eq:realization_local_system}}
&&  V_{\ell,\bullet,\lambda}
&&  \Lambda_{\ell,\lambda,\bullet} \ar@{-->}[ll]_-{-\otimes_{\dZ_\ell}\dQ_\ell}^-{\eqref{eq:local_system_integral}}
}
\]
in which the dashed arrows require $\ell\nmid D_EB_\lambda[\dE\colon\dE^\lambda]$ and that $L_\ell$ is contained in the stabilizer of the (self-dual) lattice $\Lambda_\ell$ (see Remark \ref{re:local_system_integral}).
\end{remark}

\subsection{Motivic special classes and motivic generating functions}
\label{ss:special_class}

Take a hermitian $E$-weight $\lambda$ of rank $m$ (Definition \ref{de:weight}) and consider a positive integer $r$ such that $\lambda$ is $r$-balanced (Definition \ref{de:weight}(2)).

We construct, for every tuple $x=(x_1,\dots,x_r)\in(V^\infty)^r$ satisfying $T(x)\in\Herm_r(F)^+$ a \emph{motivic special class}, which is an $\dE^\lambda$-linear map
\[
z^\lambda(x)_L\colon\tS_{\lambda^+}((\dE^\lambda)^r)\to \rH^{2r}(X_L,\rM_{L,\lambda}(r))
\]
(see Construction \ref{co:representation_general} for the source) for every $L\in\Hk(H(\bA_F^\infty))$.

Recall from \eqref{eq:special_morphism_3} that we have the canonical map
\[
\mu(x)_L\colon\varsigma(x)_L^*\rM_L^\vee\to\mathbbm{1}_{Z(x)_L,E}^{\oplus r}
\]
in $\Mot(Z(x)_L,E)$. Replace $\mathbbm{1}_{Z(x)_L,E}$ by $\mathbbm{1}_{Z(x)_L,E}^\vee$ by the tautological duality of the unit object, so that we have the induced map
\[
\mu(x)_L[u]\colon\varsigma(x)_L^*(\rM_L[u])^\vee\to(\mathbbm{1}_{Z(x)_L,E}^\vee)^{\oplus r}[u]=(\mathbbm{1}_{Z(x)_L,\dE}^\vee)^{\oplus r}
\]
(Notation \ref{no:place_decompostion}) for every $u\in\Hom(E,\dC)$. Taking tensor products, we have the map
\[
\mu_b(x)_L\colon\varsigma(x)_L^*\tT_b(\rM_L)^\vee\to\tT^b((\mathbbm{1}_{Z(x)_L,\dE}^\vee)^{\oplus r})
=\tT^b(\mathbbm{1}_{Z(x)_L,\dE}^{\oplus r})^\vee
\]
in $\Mot(Z(x)_L,\dE)$. Composing with the dual of the natural map $\tT_b(\rM_L)\to\tT_{[b]}(\rM_L)$, we obtain the map
\[
\mu_{[b]}(x)_L\colon\varsigma(x)_L^*\tT_{[b]}(\rM_L)^\vee\to\tT^b(\mathbbm{1}_{Z(x)_L,\dE}^{\oplus r})^\vee,
\]
which is $\tS_b$-equivariant. Applying $\pi_{\lambda^+}$, we obtain the map
\begin{align*}
\mu_{[\lambda]}(x)_L\colon\varsigma(x)_L^*\tS_{[\lambda]}(\rM_L)^\vee\to\tS_{\lambda^+}(\mathbbm{1}_{Z(x)_L,\dE}^{\oplus r})^\vee.
\end{align*}
It is straightforward to verify that $\mu_{[\lambda]}(x)_L$ is $\Gal(\dE/\dE^\lambda)$-equivariant, hence induces a map
\begin{align*}
\mu_\lambda(x)_L\colon\varsigma(x)_L^*\rM_{L,\lambda}^\vee\to\tS_{\lambda^+}(\mathbbm{1}_{Z(x)_L,\dE^\lambda}^{\oplus r})^\vee
\end{align*}
in the category $\Mot(Z(x)_L,\dE^\lambda)$.\footnote{The construction of $\mu_\lambda(x)_L$ is the motivic analogue of the map $\mu_\lambda(x)$ in Construction \ref{co:harmonic_function}(3).}

Denote by $\rD_-$ the Verdier duality functor in $\sfD(-,\dE^\lambda)$. We have the induced map
\[
\rD_{Z(x)_L}\mu_\lambda(x)_L\colon\rD_{Z(x)_L}\(\tS_{\lambda^+}(\mathbbm{1}_{Z(x)_L,\dE^\lambda}^{\oplus r})^\vee\)
\to\rD_{Z(x)_L}\(\varsigma(x)_L^*\rM_{L,\lambda}^\vee\)
\]
in the category $\sfD(Z(x)_L,\dE^\lambda)$. By the bi-duality property \cite{Ayo}*{Theorem~3.10}, we have canonical isomorphisms
\begin{align*}
\rD_{Z(x)_L}\(\tS_{\lambda^+}(\mathbbm{1}_{Z(x)_L,\dE^\lambda}^{\oplus r})^\vee\)
&\simeq\varsigma(x)_L^!\(\rD_{X_L}\tS_{\lambda^+}(\mathbbm{1}_{X_L,\dE^\lambda}^{\oplus r})^\vee\), \\
\rD_{Z(x)_L}\(\varsigma(x)_L^*\rM_{L,\lambda}^\vee\)
&\simeq\varsigma(x)_L^!\(\rD_{X_L}\rM_{L,\lambda}^\vee\).
\end{align*}
Since $X_L$ is regular, we have canonical isomorphisms
\[
\rD_{X_L}\tS_{\lambda^+}(\mathbbm{1}_{X_L,\dE^\lambda}^{\oplus r})^\vee\simeq
\tS_{\lambda^+}(\mathbbm{1}_{X_L,\dE^\lambda}^{\oplus r})\llangle\dim X_L\rrangle,\qquad
\rD_{X_L}\rM_{L,\lambda}^\vee \simeq\rM_{L,\lambda}\llangle\dim X_L\rrangle.
\]
Combining together, $\rD_{Z(x)_L}\mu_\lambda(x)_L$ is indeed a map
\begin{align}\label{eq:special_map_motive_0}
\mu^\lambda(x)_L\colon\tS_{\lambda^+}((\dE^\lambda)^r)\otimes\varsigma(x)_L^!\mathbbm{1}_{X_L,\dQ}=
\varsigma(x)_L^!\tS_{\lambda^+}(\mathbbm{1}_{X_L,\dE^\lambda}^{\oplus r})
\to\varsigma(x)_L^!\rM_{L,\lambda}
\end{align}
in the category $\sfD(Z(x)_L,\dE^\lambda)$.

Apply the motivic cohomology $\rH^{2\rank x}(Z(x)_L,-(\rank x))$ to $\mu^\lambda(x)$, we have the induced map
\[
\tS_{\lambda^+}((\dE^\lambda)^r)\otimes\rH^{2\rank x}(Z(x)_L,\varsigma(x)_L^!\mathbbm{1}_{X_L,\dQ}(\rank x))
\to\rH^{2\rank x}(Z(x)_L,\varsigma(x)_L^!\rM_{L,\lambda}(\rank x)).
\]
For the source, since $Z(x)_L$ is regular, we have canonically
\[
\rH^{2\rank x}(Z(x)_L,\varsigma(x)_L^!\mathbbm{1}_{X_L,\dQ}(\rank x))=\rH^0(Z(x)_L,\mathbbm{1}_{Z(x)_L,\dQ}).
\]
For the target, we have the natural map
\[
\rH^{2\rank x}(Z(x)_L,\varsigma(x)_L^!\rM_{L,\lambda}(\rank x))
\to\rH^{2\rank x}_{Z[x]_L}(X_L,\rM_{L,\lambda}(\rank x)),
\]
where $Z[x]_L$ denotes the image of the finite unramified morphism $\varsigma(x)_L$ \eqref{eq:special_morphism_1} which is a Zariski closed subset of $X_L$. Together, we obtain the map
\begin{align}\label{eq:special_class}
\tS_{\lambda^+}((\dE^\lambda)^r)\otimes\rH^0(Z(x)_L,\mathbbm{1}_{Z(x)_L,\dQ})\to
\rH^{2\rank x}_{Z[x]_L}(X_L,\rM_{L,\lambda}(\rank x))
\end{align}
induced from $\mu^\lambda(x)_L$.

\begin{definition}\label{de:special_class}
Take a tuple $x=(x_1,\dots,x_r)\in(V^\infty)^r$ satisfying $T(x)\in\Herm_r(F)^+$.
\begin{enumerate}
  \item Define the \emph{refined motivic special class} $\bbz^\lambda(x)_L$ to be the composite map
      \[
      \tS_{\lambda^+}((\dE^\lambda)^r)\to\tS_{\lambda^+}((\dE^\lambda)^r)\otimes\rH^0(Z(x)_L,\mathbbm{1}_{Z(x)_L,\dQ})
      \xrightarrow{\eqref{eq:special_class}}\rH^{2\rank x}_{Z[x]_L}(X_L,\rM_{L,\lambda}(\rank x)),
      \]
      in which the first one is induced by tensoring with the fundamental class of $Z(x)_L$.

  \item Define the \emph{motivic special class} $z^\lambda(x)_L$ to be the composite map
      \[
      \tS_{\lambda^+}((\dE^\lambda)^r)\xrightarrow{\bbz^\lambda(x)_L}\rH^{2\rank x}_{Z[x]_L}(X_L,\rM_{L,\lambda}(\rank x))
      \to\rH^{2\rank x}(X_L,\rM_{L,\lambda}(\rank x))\xrightarrow{\wedge e^{r-\rank x}}
      \rH^{2r}(X_L,\rM_{L,\lambda}(r)),
      \]
      in which $e\in\rH^2(X_L,\mathbbm{1}_{X_L,\dQ}(1))$ denotes the class given by the dual Hodge line bundle over $X_L$.
\end{enumerate}

\end{definition}

\begin{remark}
When $\lambda=\b{0}$, we have $\dE^\lambda=\dQ$, $\rM_{L,\lambda}=\mathbbm{1}_{X_L,\dQ}$ hence $\rH^{2r}(X_L,\rM_{L,\lambda}(r))=\CH^r(X_L)\otimes\dQ$ by Remark \ref{re:relative_motivic}(1). In this case, $z^{\b{0}}(x)_L$ recovers Kudla's special cycle (class) in the hermitian case (see, for example, \cite{Liu11}*{\S3A}).
\end{remark}

Once motivic special classes are constructed, we may define generating functions using those classes. Take an object $L\in\Hk(H(\bA_F^\infty))$ and a Schwartz function $\phi^\infty\in\sS((V^\infty)^r)^L$.

\begin{definition}\label{de:special_cycle}
Take an element $T\in\Herm_r(F)^+$. Put
\[
Z_T[\phi^\infty]_L\coloneqq\bigcup_{\substack{x\in (V^\infty)^r_T\\ \phi^\infty(x)\neq 0}}Z[x]_L,
\]
which is a finite union, as a Zariski closed subset of $X_L$.
\begin{enumerate}
  \item When $T\in\Herm_r^\circ(F)^+$, put
      \[
      \bbZ^\lambda_T(\phi^\infty)_L\coloneqq\sum_{x\in L\backslash (V^\infty)^r_T} \phi^\infty(x)\cdot \bbz^\lambda(x)_L\in\tS_{\lambda^+}(\dC^r)^\vee\otimes_{\dE^\lambda}\rH^{2r}_{Z_T[\phi^\infty]_L}(X_L,\rM_{L,\lambda}(r))
      \]
      where the tensor in the target is taken with respect to the inclusion $\dE^\lambda\subseteq\dC$. We extend $\bbZ^\lambda_T(\phi^\infty)_L$ to a function on $\widetilde{G_r(\bA_F)}$ via the formula
      \begin{align*}
      g\in\widetilde{G_r(\bA_F)}\mapsto
      \bbZ^\lambda_T(\phi^\infty)_L(g)&\coloneqq\omega_{r,\infty}(g_\infty)\phi^0_{V_\infty^r}(T) \cdot\rho_{\lambda^+}(g_\infty)^\vee\(\bbZ^\lambda_T(\omega_r^\infty(g^\infty)\phi^\infty)_L\)\\
      &\in\tS_{\lambda^+}(\dC^r)^\vee\otimes_{\dE^\lambda}
      \rH^{2r}_{Z_T[\omega_r^\infty(g^\infty)\phi^\infty]_L}(X_L,\rM_{L,\lambda}(r)),
      \end{align*}
      where $\rho_{\lambda^+}$ is the map from Definition \ref{de:hermitian_siegel}(3) and $\phi^0_{V_\infty^r}$ is standard Gaussian function on $V_\infty^r$.

  \item For $T$ general, put
      \[
      Z^\lambda_T(\phi^\infty)_L\coloneqq\sum_{x\in L\backslash (V^\infty)^r_T} \phi^\infty(x)\cdot z^\lambda(x)_L\in\tS_{\lambda^+}(\dC^r)^\vee\otimes_{\dE^\lambda}\rH^{2r}(X_L,\rM_{L,\lambda}(r)).
      \]
      Similar to (1), we extend $Z^\lambda_T(\phi^\infty)_L$ to a function on $\widetilde{G_r(\bA_F)}$ valued in $\tS_{\lambda^+}(\dC^r)^\vee\otimes_{\dE^\lambda}\rH^{2r}(X_L,\rM_{L,\lambda}(r))$ via the formula
      \[
      g\in\widetilde{G_r(\bA_F)}\mapsto
      Z^\lambda_T(\phi^\infty)_L(g)\coloneqq\omega_{r,\infty}(g_\infty)\phi^0_{V_\infty^r}(T) \cdot\rho_{\lambda^+}(g_\infty)^\vee\(Z^\lambda_T(\omega_r^\infty(g^\infty)\phi^\infty)_L\).
      \]
\end{enumerate}
\end{definition}

\begin{lem}\label{le:hecke_translate}
For every element $\tt\in\dC[L\backslash H(\bA_F^\infty)/L]$, we have $\tt^*Z^\lambda_T(\phi^\infty)_L=Z^\lambda_T(\tt\phi^\infty)_L$.
\end{lem}

\begin{proof}
This follows from the same proof of \cite{LL21}*{Lemma~4.4} (see also \cite{DL24}*{Remark~4.51}).
\end{proof}

\begin{definition}\label{de:generating_function}
Define the \emph{motivic generating function} $Z^\lambda(\phi^\infty)_L$ to be the formal sum
\[
Z^\lambda(\phi^\infty)_L\coloneqq\sum_{T\in\Herm_r(F)^+}Z^\lambda_T(\phi^\infty)_L
\]
as a function on $\widetilde{G_r(\bA_F)}$ valued in $\tS_{\lambda^+}(\dC^r)^\vee\otimes_{\dE^\lambda}\rH^{2r}(X_L,\rM_{L,\lambda}(r))$.
\end{definition}

\begin{conjecture}[Modularity of motivic generating function]\label{co:modularity}
Let $r$ be a positive integer such that $\lambda$ is $r$-balanced. For every element $\phi^\infty\in\sS((V^\infty)^r)^L$, the motivic generating function $Z^\lambda(\phi^\infty)_L$ is an absolutely convergent sum of terms in a finite-dimensional subspace, which further belongs to $\cA^{m/2}_{\lambda^+}(G_r)^\hol\otimes_{\dE^\lambda}\rH^{2r}(X_L,\rM_{L,\lambda}(r))$ (Definition \ref{de:hermitian_siegel}(1)).
\end{conjecture}

\begin{remark}
We have several remarks concerning Conjecture \ref{co:modularity}.
\begin{enumerate}
  \item When $\lambda=\b{0}$ so that $\rH^{2r}(X_L,\rM_{L,\lambda}(r))=\CH^r(X_L)\otimes\dQ$, the function $Z^{\b{0}}(\phi^\infty)_L$ recovers Kudla's generating function in the hermitian case and Conjecture \ref{co:modularity} recovers Kudla's modularity conjecture in the hermitian case.\footnote{A proof of Conjecture \ref{co:modularity} for $\lambda=\b{0}$ for a general CM extension $E/F$ has been recently claimed by Raum \cite{Rau26}.}

  \item For a general weight $\lambda$, Conjecture \ref{co:modularity} is not known even for $r=1$ and $m=2$.

  \item Assuming the conjecture, $Z^\lambda(\phi^\infty)_L$ further belongs to $\cA^{m/2}_{\lambda^+}(G_r)^\hol_\cusp\otimes_{\dE^\lambda}\rH^{2r}(X_L,\rM_{L,\lambda}(r))$ if $r=\max_u\{s(\lambda_u)\}$.
\end{enumerate}
\end{remark}

The definition below generalizes \cite{Liu11}*{Definition~3.8} or \cite{LL21}*{Definition~4.8}.

\begin{definition}\label{de:theta_lift}
Assume Conjecture \ref{co:modularity} for $r$. For every form $\varphi\in\cA^{m/2}_{\lambda^+}(G_r)_\cusp$, we define the \emph{motivic arithmetic theta lift of weight $\lambda$} of $\varphi$, with respect to a Schwartz function $\phi^\infty\in\sS((V^\infty)^r)^L$, to be
\[
\Theta_{\phi^\infty}^\lambda(\varphi)_L\coloneqq(Z^\lambda(\phi^\infty)_L,\varphi)_{G_r}
\in \rH^{2r}(X_L,\rM_{L,\lambda}(r))\otimes_{\dE^\lambda}\dC
\]
(Definition \ref{de:hermitian_siegel}(2)).
\end{definition}

\begin{remark}\label{re:theta_lift}
In the above definition, it is straightforward to verify that for another object $L'\subseteq L$ of $\Hk(H(\bA_F^\infty))$, the pullback of $\Theta_{\phi^\infty}^\lambda(\varphi)_L$ along the transition map $\tau_{L',L}\colon X_{L'}\to X_L$ coincides with $\Theta_{\phi^\infty}^\lambda(\varphi)_{L'}$. Thus, the collection $\{\Theta_{\phi^\infty}^\lambda(\varphi)_L\}_L$ defines an element
\[
\Theta_{\phi^\infty}^\lambda(\varphi)\in\rH^{2r}(X,\rM_{\lambda}(r))\otimes_{\dE^\lambda}\dC,
\]
where
\[
\rH^{2r}(X,\rM_{\lambda}(r))\coloneqq\varinjlim_L\rH^{2r}(X_L,\rM_{L,\lambda}(r)).
\]
\end{remark}

\begin{remark}
The kernel of the geometric $\ell$-adic realization map
\[
\rH^{2r}(X_L,\rM_{L,\lambda}(r))\to\rH^{2r}(X_L\otimes_EE^\ac,\rM_{L,\ell,\lambda}(r))
\]
on the level of cohomology induced from \eqref{eq:realization_etale} does not depend on the prime $\ell$. Here, $E^\ac$ denotes an algebraic closure of $E$.

Indeed, choose an object $U\in \Hk(T(\bA^\infty))$ so that we have the natural quotient map $\pi_{L,U}\colon X_{L,U}\to X_L$ \eqref{eq:quotient_morphism}. Put $\rM_{L,U,\lambda}\coloneqq\pi_{L,U}^*\rM_{L,\lambda}$ as in \S\ref{ss:motive_automorphic}, which is a direct summand of $\sfh^{2|b|}(B_{L,U}^{|b|}/X_{L,U},\dQ)\ang{|b|}\otimes_\dQ\dE^\lambda$ via the straightening map $\st_{L,U}^{\lambda,+}$ (Definition \ref{de:straightening}). In particular, we obtain a commutative diagram
\[
\xymatrix{
\rH^{2r}(X_L,\rM_{L,\lambda}(r)) \ar[r]\ar[d] & \rH^{2r}_{\et}(X_L\otimes_EE^\ac,\rM_{L,\ell,\lambda}(r)) \ar[d] \\
\rH^{2r}(X_{L,U},\rM_{L,U,\lambda}(r))
\ar[r]\ar[d] & \rH^{2r}_{\et}(X_{L,U}\otimes_EE^\ac,r_\ell\rM_{L,U,\lambda}(r)) \ar[d] \\
\rH^{2r+2|b|}(B_{L,U}^{|b|},\mathbbm{1}_{B_{L,U}^{|b|},\dE^\lambda}(r+|b|)) \ar[r]& \rH^{2r+2|b|}_{\et}(B_{L,U}^{|b|}\otimes_EE^\ac,\dE^\lambda_\ell(r+|b|))
}
\]
in which all vertical maps are injective. Thus, it suffices to show that the kernel of the bottom map does not depend on the choice of $\ell$, which follows from the comparison theorem between geometric $\ell$-adic cohomology groups (together with cycle class maps).
\end{remark}

\subsection{Canonical smooth integral models}
\label{ss:integral_model}

In this section, we recall the construction of canonical smooth integral models at good primes. Take a finite set $\Diamond$ of rational primes.

\begin{definition}\label{de:moduli_auxiliary_integral}
We define a functor
\begin{align*}
Y_{\bullet,(\Diamond)}\colon\Hk(T(\bA^{\infty,\Diamond}))\to\sfP\Sch'_{/O_{E_0,(\Diamond)}}
\end{align*}
such that for every objects $U^\Diamond\in\Hk(T(\bA^{\infty,\Diamond}))$ and $S\in\Sch'_{/O_{E_0,(\Diamond)}}$, $Y_{U^\Diamond,(\Diamond)}(S)$ is the set of equivalence classes of triples $(A_0,\lambda_0,\bar\eta_0^\Diamond)$, where
\begin{itemize}
  \item $A_0$ is an object of $\AV(S,O_{E_0,(\Diamond)})$ such that $\tr_{A_0}(e)=e$ for $e\in O_{E_0,(\Diamond)}$;

  \item $\lambda_0\colon A_0\to A_0^\vee$ is a $\Diamond$-principal polarization that turns the complex conjugation of $O_{E_0,(\Diamond)}$ into the Rosati involution;

  \item $\bar\eta_0^\Diamond$ is a $U^\Diamond$-level structure, that is, for a chosen geometric point $s$ on every connected component of $S$, a $\pi_1(S,s)$-invariant $U^\Diamond$-orbit of isomorphisms
      \[
      \eta_0^\Diamond\colon\bA_{E_0}^{\infty,\Diamond}\xrightarrow\sim\rH_1^{\et}(A_{0s},\bA^{\infty,\Diamond})
      \]
      of $\bA_{E_0}^{\infty,\Diamond}$-modules.
\end{itemize}
Two objects $(A_0,\lambda_0,\bar\eta_0^\Diamond)$ and $(A'_0,\lambda'_0,\bar\eta_0^{p\prime})$ are equivalent if there exists an isomorphism $\varphi_0\colon A_0\to A'_0$ in the category $\AV(S,O_{E_0,(\Diamond)})$ satisfying that
\begin{itemize}
  \item $\varphi_0^\vee\circ\lambda'_0\circ\varphi_0$ is a prime-to-$\Diamond$ multiple of $\lambda_0$;

  \item the $U^\Diamond$-orbits $\varphi_{0*}\circ\bar\eta_0^\Diamond$ and $\bar\eta_0^{\Diamond\prime}$ coincide.
\end{itemize}
A morphism $t\in U^\Diamond\backslash T(\bA^{\infty,\Diamond})/U^{\prime\Diamond}$ maps $Y_{U^\Diamond,(\Diamond)}(S)$ to $Y_{U^{\prime\Diamond},(\Diamond)}(S)$ by changing $\bar\eta_0^\Diamond$ to $\bar\eta_0^\Diamond\circ t$.
\end{definition}

It is well-known that $Y_{\bullet,(\Diamond)}$ is representable and factors through the full subcategory $\Sch^\smp_{/O_{E_0,(\Diamond)}}$ of $\sfP\Sch'_{/O_{E_0,(\Diamond)}}$. Moreover, if we denote by
\[
A_{0,U^\Diamond,(\Diamond)}\in\AV(Y_{U^\Diamond,(\Diamond)},O_{E_0,(\Diamond)})
\]
the universal object, then we obtain a functor
\[
A_{0,\bullet,(\Diamond)}/Y_{\bullet,(\Diamond)}\colon
\Hk(T(\bA^{\infty,\Diamond}))\to\AV(\Sch^\smp_{/O_{E_0,(\Diamond)}},O_{E_0,(\Diamond)})
\]
over $Y_{\bullet,(\Diamond)}$ that is Cartesian.

\begin{definition}\label{de:moduli_integral}
We define a functor
\begin{align*}
X_{\bullet,\bullet,(\Diamond)}\colon
\Hk(H(\bA_F^{\infty,\Diamond}))\times\Hk(T(\bA^{\infty,\Diamond}))\to\sfP\Sch'_{/O_{E,(\Diamond)}}
\end{align*}
such that for every objects $(L^\Diamond,U^\Diamond)\in\Hk(H(\bA_F^{\infty,\Diamond}))\times\Hk(T(\bA^{\infty,\Diamond}))$ and $S\in\Sch'_{/O_{E,(\Diamond)}}$, $X_{L^\Diamond,U^\Diamond,(\Diamond)}(S)$ is the set of equivalence classes of sextuples $(A_0,\lambda_0,\bar\eta_0^\Diamond;A,\lambda,\bar\eta^\Diamond)$, where
\begin{itemize}
  \item $(A_0,\lambda_0,\bar\eta_0^\Diamond)$ is an object of $Y_{U^\Diamond,(\Diamond)}(S)$;

  \item $A$ is an object of $\AV(S,O_{E,(\Diamond)})$ such that $\tr_A(e)=m\Tr_{E/E_0}(e)-e+\ol{e}$ for $e\in O_{E,(\Diamond)}$;

  \item $\lambda\colon A\to A^\vee$ is a $\Diamond$-principal polarization that turns the complex conjugation of $O_{E,(\Diamond)}$ into the Rosati involution;

  \item $\bar\eta^\Diamond$ is a \emph{hermitian $L^\Diamond$-level structure}, that is, for a chosen geometric point $s$ on every connected component of $S$, a $\pi_1(S,s)$-invariant $L^\Diamond$-orbit of isometries
      \[
      \eta^\Diamond\colon V^{\infty,\Diamond}\xrightarrow\sim\Hom_{\bA_E^{\infty,\Diamond}}^{\lambda_0,\lambda}
      \(\rH_1^{\et}(A_{0s},\bA^{\infty,\Diamond})\otimes_{\bA^{\infty,\Diamond}}\bA_F^{\infty,\Diamond},
      \rH_1^{\et}(A_s,\bA^{\infty,\Diamond})\)
      \]
      of hermitian spaces over $\bA_E^{\infty,\Diamond}/\bA_F^{\infty,\Diamond}$.
\end{itemize}
Two objects $(A_0,\lambda_0,\bar\eta_0^\Diamond;A,\lambda,\bar\eta^\Diamond)$ and $(A'_0,\lambda'_0,\bar\eta_0^{\Diamond\prime};A',\lambda',\bar\eta^{\Diamond\prime})$ are equivalent if there exist isomorphisms $\varphi_0\colon A_0\to A'_0$ in the category $\AV(S,O_{E_0,(\Diamond)})$ and $\varphi\colon A\to A'$ in the category $\AV(S,O_{E,(\Diamond)})$ satisfying that
\begin{itemize}
  \item $\varphi_0^\vee\circ\lambda'_0\circ\varphi_0$ and $\varphi^\vee\circ\lambda'\circ\varphi$ are the \emph{same} prime-to-$\Diamond$ multiple of $\lambda_0$ and $\lambda$, respectively;

  \item the $U^\Diamond$-orbits $\varphi_{0*}\circ\bar\eta_0^\Diamond$ and $\bar\eta_0^{\Diamond\prime}$ coincide;

  \item the $L^\Diamond$-orbit of maps $v\mapsto\varphi_*\circ\bar\eta^\Diamond(v)\circ(\varphi_{0*})^{-1}$ for $v\in V^{\infty,\Diamond}$ coincides with $\bar\eta^{\Diamond\prime}$.
\end{itemize}
A morphism $(h,t)\in L^\Diamond\backslash H(\bA_F^{\infty,\Diamond})/L^{\prime\Diamond}\times U^\Diamond\backslash T(\bA^{\infty,\Diamond})/U^{\prime\Diamond}$ maps $X_{L^\Diamond,U^\Diamond,(\Diamond)}(S)$ to $X_{L^{\prime\Diamond},U^{\prime\Diamond},(\Diamond)}(S)$ by changing $\bar\eta_0^\Diamond$ to $\bar\eta_0^\Diamond\circ t$ and $\bar\eta^\Diamond$ to $\bar\eta^\Diamond\circ h$.
\end{definition}

Now we suppose that every prime in $\Diamond$ is unramified in $E$ and $V_v$ is split for every $v\in\fV_F^{(\Diamond)}$.

By \cite{RSZ}*{Theorem~4.4}, $X_{\bullet,\bullet,(\Diamond)}$ factors through the full subcategory $\Sch^\smq_{/O_{E,(\Diamond)}}$ of $\sfP\Sch'_{/O_{E,(\Diamond)}}$. Moreover, if we denote by
\[
(A_{0,L^\Diamond,U^\Diamond,(\Diamond)},A_{L^\Diamond,U^\Diamond,(\Diamond)})\in \AV(X_{L^\Diamond,U^\Diamond,(\Diamond)},O_{E_0,(\Diamond)})\times\AV(X_{L^\Diamond,U^\Diamond},O_{E,(\Diamond)})
\]
the universal object, then we obtain functors
\begin{align*}
A_{0,\bullet,\bullet,(\Diamond)}/X_{\bullet,\bullet,(\Diamond)}&
\colon\Hk(H(\bA_F^{\infty,\Diamond}))\times\Hk(T(\bA^{\infty,\Diamond}))
\to\AV(\Sch^\smq_{/O_{E,(\Diamond)}},O_{E_0,(\Diamond)}),\\
A_{\bullet,\bullet,(\Diamond)}/X_{\bullet,\bullet,(\Diamond)}&
\colon\Hk(H(\bA_F^{\infty,\Diamond}))\times\Hk(T(\bA^{\infty,\Diamond}))
\to\AV(\Sch^\smq_{/O_{E,(\Diamond)}},O_{E,(\Diamond)})
\end{align*}
over $X_{\bullet,\bullet,(\Diamond)}$ that are Cartesian. By definition, we have a natural morphism
\begin{align*}
X_{L^\Diamond,U^\Diamond,(\Diamond)}\to Y_{U^\Diamond,(\Diamond)}\otimes_{O_{E_0}}O_E
\end{align*}
for every object $(L^\Diamond,U^\Diamond)\in\Hk(H(\bA_F^{\infty,\Diamond}))\times\Hk(T(\bA^{\infty,\Diamond}))$, along which $A_{0,L^\Diamond,U^\Diamond,(\Diamond)}$ is the pullback of $A_{0,U^\Diamond,(\Diamond)}\otimes_{O_{E_0}}O_E$ to $X_{L^\Diamond,U^\Diamond,(\Diamond)}$.

If we choose an open compact subgroup $L_\Diamond$ of $H(F_\Diamond)$ that is the stabilizer of a self-dual lattice of $V_\Diamond$ and denote by $U_\Diamond$ the maximal open compact subgroup of $T(\dQ_\Diamond)$, then it is clear from the construction that
\[
X_{L^\Diamond,U^\Diamond,(\Diamond)}\otimes_{O_{E,(\Diamond)}}E=X_{L_\Diamond L^\Diamond,U_\Diamond U^\Diamond}
\]
functorial in $(L^\Diamond,U^\Diamond)$.

\begin{remark}
Take an object $L^\Diamond\in\Hk(H(\bA_F^{\infty,\Diamond}))$. For every $U^\Diamond\in\Hk(T(\bA^{\infty,\Diamond}))$, the group $T(\bA^{\infty,\Diamond})/U^\Diamond$ acts on $X_{L^\Diamond,U^\Diamond,(\Diamond)}$ via modifying $\bar\eta_0^\Diamond$, which also factors through the quotient
\[
\Xi_{U^\Diamond}\coloneqq O_{E_0,(\Diamond)}^\times\backslash T(\bA^{\infty,\Diamond})/U^\Diamond
\]
by a way that is similar to the case for $X_{L,U}$,\footnote{Note that the natural map $\Xi_{U^\Diamond}\to\Xi_{U_\Diamond U^\Diamond}$ is an isomorphism, where $U_\Diamond$ denotes the maximal open compact subgroup of $T(\dQ_\Diamond)$.} so that the action of $\Xi_{U^\Diamond}$ on $X_{L^\Diamond,U^\Diamond,(\Diamond)}$ is free. Since $X_{L^\Diamond,U^\Diamond,(\Diamond)}$ is smooth and quasi-projective over $\Spec O_{E,(\Diamond)}$, the quotient $X_{L^\Diamond,U^\Diamond,(\Diamond)}/\Xi_{U^\Diamond}$ is again a smooth and quasi-projective scheme over $\Spec O_{E,(\Diamond)}$. It is clear that for $U^{\prime\Diamond}\subseteq U^\Diamond$ in $\Hk(T(\bA^{\infty,\Diamond}))$, the natural morphism $X_{L^\Diamond,U^{\prime\Diamond},(\Diamond)}/\Xi_{U^{\prime\Diamond}}\to X_{L^\Diamond,U^\Diamond,(\Diamond)}/\Xi_{U^\Diamond}$ is an isomorphism.

To summarize, we obtain a functor
\begin{align*}
X_{\bullet,(\Diamond)}\colon\Hk(H(\bA_F^{\infty,\Diamond}))\to\Sch^\smq_{/O_{E,(\Diamond)}}
\end{align*}
such that for every $L^\Diamond\in\Hk(H(\bA_F^{\infty,\Diamond}))$, $X_{L^\Diamond,(\Diamond)}=X_{L,U^\Diamond,(\Diamond)}/\Xi_{U^\Diamond}$ -- the quotient scheme obtained from the above process, which is canonically independent of the choice of $U^\Diamond\in\Hk(T(\bA^{\infty,\Diamond}))$. By construction, it is clear that
\[
X_{L^\Diamond,(\Diamond)}\otimes_{O_{E,(\Diamond)}}E=X_{L_\Diamond L^\Diamond}
\]
functorial in $L^\Diamond$.
\end{remark}

\section{Height pairing on motivic cohomology}
\label{ss:4}

In this chapter, we construct and study a height pairing on motivic cohomology of those automorphic motives constructed in the previous chapter. We
\begin{itemize}
  \item fix a totally positive definite incoherent hermitian space $V$ over $\bA_E/\bA_F$ of rank $m\geq 1$ satisfying that the set $\fV(V)\coloneqq\{v\in\fV_F\res V_v\text{ is anisotropic}\}$, which contains $\fV_F^{(\infty)}$, is \emph{not} a singleton (in particular, the $E$-schemes $\{X_L\}_L$ are projective);

  \item fix a hermitian $E$-weight $\lambda$ of rank $m$ (Definition \ref{de:weight}).
\end{itemize}
Put $H\coloneqq\rU(V)$. In what follows, we denote by $-_u\colon\Sch_{/E}\to\Sch_{/E_u}$ the base change functor for every place $u$ of $E$.

\subsection{Definition of height pairing}
\label{ss:height_pairing}

Take a rational prime $\ell$.

\begin{notation}\label{no:cohomological_trivial_ell}
For integers $q,r$, an object $L\in\Hk(H(\bA_F^\infty))$, and a Zariski closed subset $Z$ of $X_L$, we put
\[
\rH^q_Z(X_L,\rM_{L,\lambda}(r))^{\ang\ell}\coloneqq
\Ker\(\rH^q_Z(X_L,\rM_{L,\lambda}(r))\to\prod_{u\in\fV_E^\fin\setminus\fV_E^{(\ell)}}
\rH^q_{\et}(X_{L,u},\rM_{L,\ell,\lambda}(r))\).
\]
\end{notation}

Suppose that $\ell$ is unramified in $E$. The goal of this section is to define, for every object $L\in\Hk(H(\bA_F^\infty))$ such that $L_\ell$ is the stabilizer of a self-dual lattice of $V_\ell$, a canonical hermitian\footnote{Here, \emph{hermitian} means that the pairing below is $\dE^\lambda$-linear in the first variable and complex-conjugate symmetric.} pairing
\begin{align}\label{eq:height_pairing}
\llangle\;,\;\rrangle_{\lambda,L}^{\ell}\colon
\(\bigoplus_{r=1}^{m-1}\rH^{2r}(X_L,\rM_{L,\lambda}(r))^{\ang\ell}\)
\times\(\bigoplus_{r=1}^{m-1}\rH^{2r}(X_L,\rM_{L,\lambda}(r))^{\ang\ell}\)\to\dC\otimes\dQ_\ell
\end{align}
that vanishes on $\rH^{2r_1}(X_L,\rM_{L,\lambda}(r_1))^{\ang\ell}\times\rH^{2r_2}(X_L,\rM_{L,\lambda}(r_2))^{\ang\ell}$ unless when $r_1+r_2=m$.

Choose an object $U$ of $\Hk(T(\bA^\infty))$ such that $U_\ell$ is maximal so that we have the natural quotient map $\pi_{L,U}\colon X_{L,U}\to X_L$ \eqref{eq:quotient_morphism}. Put $\rM_{L,U,\lambda}\coloneqq\pi_{L,U}^*\rM_{L,\lambda}$ as in \S\ref{ss:motive_automorphic}. The straightening map $\st_{L,U}^{\lambda,+}$ and complementary straightening map $\st_{L,U}^{\lambda,-}$ (Definition \ref{de:straightening}) induce $\dE^\lambda$-linear maps (with the same notation)
\begin{align*}
\st_{L,U}^{\lambda,+}\colon\rH^{2r}(X_{L,U},\rM_{L,U,\lambda}(r))
&\to\rH^{2r+2|b|}(B_{L,U}^{|b|},\mathbbm{1}_{B_{L,U}^{|b|},\dE^\lambda}(r+|b|)), \\
\st_{L,U}^{\lambda,-}\colon\rH^{2r}(X_{L,U},\rM_{L,U,\lambda}(r))
&\to\rH^{2r+2(g-1)|b|}(B_{L,U}^{|b|},\mathbbm{1}_{B_{L,U}^{|b|},\dE^\lambda}(r+(g-1)|b|)),
\end{align*}
respectively. For every term $\rH$ in the above two maps, one may define the subspace $\rH^{\ang\ell}$ similarly as in Notation \ref{no:cohomological_trivial_ell}. It follows easily that we have induced $\dE^\lambda$-linear maps
\begin{align*}
\st_{L,U}^{\lambda,+}\circ\pi_{L,U}^*\colon\rH^{2r}(X_L,\rM_{L,\lambda}(r))^{\ang\ell}
&\to\rH^{2r+2|b|}(B_{L,U}^{|b|},\mathbbm{1}_{B_{L,U}^{|b|},\dE^\lambda}(r+|b|))^{\ang\ell}, \\
\st_{L,U}^{\lambda,-}\circ\pi_{L,U}^*\colon\rH^{2r}(X_L,\rM_{L,\lambda}(r))^{\ang\ell}
&\to\rH^{2r+2(g-1)|b|}(B_{L,U}^{|b|},\mathbbm{1}_{B_{L,U}^{|b|},\dE^\lambda}(r+(g-1)|b|))^{\ang\ell},
\end{align*}
respectively. Note that for every integer $s\geq 0$, under the natural isomorphism
\[
\rH^{2s}(B_{L,U}^{|b|},\mathbbm{1}_{B_{L,U}^{|b|},\dE^\lambda}(s))=\CH^s(B_{L,U}^{|b|})\otimes\dE^\lambda,
\]
the subspace $\rH^{2s}(B_{L,U}^{|b|},\mathbbm{1}_{B_{L,U}^{|b|},\dE^\lambda}(s))^{\ang\ell}$ coincides with $\CH^s(B_{L,U}^{|b|})^{\ang\ell}\otimes\dE^\lambda$, where $\CH^s(B_{L,U}^{|b|})^{\ang\ell}$ was defined in \cite{LL21}*{Notation~2.5(C3)}. In summary, we have induced $\dE^\lambda$-linear maps
\begin{align*}
\st_{L,U}^{\lambda,+}\circ\pi_{L,U}^*\colon\rH^{2r}(X_L,\rM_{L,\lambda}(r))^{\ang\ell}
&\to\CH^{r+|b|}(B_{L,U}^{|b|})^{\ang\ell}\otimes\dE^\lambda, \\
\st_{L,U}^{\lambda,-}\circ\pi_{L,U}^*\colon\rH^{2r}(X_L,\rM_{L,\lambda}(r))^{\ang\ell}
&\to\CH^{r+(g-1)|b|}(B_{L,U}^{|b|})^{\ang\ell}\otimes\dE^\lambda,
\end{align*}
respectively.

By the discussion in \S\ref{ss:integral_model}, we know that $B_{L,U}^{|b|}$ has smooth projective reduction at every $\ell$-adic place of $E$. Thus, we may apply Beilinson's construction of the height pairing in \cite{Bei87}*{\S4} (see also \cite{LL21}*{\S6}) to obtain a map
\[
\llangle\;,\;\rrangle_{B_{L,U}^{|b|}}^\ell\colon
\(\bigoplus_{s=1}^{g|b|+m-1}\CH^s(B_{L,U}^{|b|})^{\ang\ell}\otimes\dC\)
\times\(\bigoplus_{s=1}^{g|b|+m-1}\CH^s(B_{L,U}^{|b|})^{\ang\ell}\otimes\dC\)\to\dC\otimes\dQ_\ell
\]
that is hermitian and vanishes on $\CH^{s_1}(B_{L,U}^{|b|})^{\ang\ell}\otimes\dC\times\CH^{s_2}(B_{L,U}^{|b|})^{\ang\ell}\otimes\dC$ unless when $s_1+s_2=g|b|+m=\dim B_{L,U}^{|b|}+1$. Put
\begin{align*}
\llangle\;,\;\rrangle_{\lambda,L}^{\ell}&\coloneqq\frac{1}{|\Xi_U|}\cdot
\llangle\;,\;\rrangle_{B_{L,U}^{|b|},E}^\ell\circ
\((\st_{L,U}^{\lambda,+}\circ\pi_{L,U}^*)\times(\st_{L,U}^{\lambda,-}\circ\pi_{L,U}^*)\) \\
&\colon\(\bigoplus_{r=1}^{m-1}\rH^{2r}(X_L,\rM_{L,\lambda}(r))^{\ang\ell}\)
\times\(\bigoplus_{r=1}^{m-1}\rH^{2r}(X_L,\rM_{L,\lambda}(r))^{\ang\ell}\)\to\dC\otimes\dQ_\ell.
\end{align*}

\begin{lem}
The map $\llangle\;,\;\rrangle_{\lambda,L}^{\ell}$ does not depend on the choice of $U$ (hence justifying the disappearance of $U$ in the notation). Moreover,
\begin{enumerate}
  \item it is $\dE^\lambda$-linear in the first variable and complex-conjugate symmetric;

  \item it vanishes on $\rH^{2r_1}(X_L,\rM_{L,\lambda}(r_1))^{\ang\ell}\times\rH^{2r_2}(X_L,\rM_{L,\lambda}(r_2))^{\ang\ell}$ unless when $r_1+r_2=m$.
\end{enumerate}
\end{lem}

\begin{proof}
Note that for an open compact subgroup $U'\subseteq U$ with $U'_\ell=U_\ell$, the natural transition map $B_{L,U'}^{|b|}\to B_{L,U}^{|b|}$ is finite \'{e}tale of degree $|\Xi_{U'}|/|\Xi_U|$. Then the independence of $U$ of $\llangle\;,\;\rrangle_{\lambda,L}^{\ell}$ follows from the functoriality of Beilinson's height pairing \cite{Bei87}*{4.0.3}.

Part (2) follows from the similar property for $\llangle\;,\;\rrangle_{B_{L,U}^{|b|},E}^\ell$.

For part (1), the $\dE^\lambda$-linearity in the first variable is obvious. It remains to show the complex-conjugate symmetry, that is, for every pair $(c_1,c_2)\in\rH^{2r_1}(X_L,\rM_{L,\lambda}(r_1))^{\ang\ell}\times\rH^{2r_2}(X_L,\rM_{L,\lambda}(r_2))^{\ang\ell}$ with $r_1+r_2=m$, we have
\[
\llangle c_2,c_1\rrangle_{\lambda,L}^{\ell}=\ol{\llangle c_1,c_2\rrangle_{\lambda,L}^{\ell}}.
\]
For $i=1,2$, put $c_i^\pm\coloneqq\st_{L,U}^{\lambda,\pm}(\pi_{L,U}^*(c_i))$. It suffices to show that
\[
\llangle c_2^+,c_1^-\rrangle_{B_{L,U}^{|b|}}^\ell
=\ol{\llangle c_1^+,c_2^-\rrangle_{B_{L,U}^{|b|}}^\ell}.
\]
For this, note that the map
\[
\CH^{r+|b|}(B_{L,U}^{|b|})^{\ang\ell}\otimes\dQ\to\CH^{r+(g-1)|b|}(B_{L,U}^{|b|})^{\ang\ell}\otimes\dQ
\]
for every integer $r\geq 0$ induced from the map
\[
((\lambda_0\times\lambda^{-1})^{|b|})^*\circ\tF^{2(g-1)|b|}_{B_{L,U}^{|b|}/X_{L,U}}\ang{|b|}\colon
\sfh^{2|b|}(B_{L,U}^{|b|}/X_{L,U},\dQ)\ang{|b|}\to
\sfh^{2(g-1)|b|}(B_{L,U}^{|b|}/X_{L,U},\dQ)\ang{(g-1)|b|}
\]
of relative motives over $X_{L,U}$ is simply the pullback given by the correspondence $p^{\otimes 2(g-1)|b|}$, where $p\in\CH^1(B_{L,U}^{|b|}\times_{X_{L,U}}B_{L,U}^{|b|})\otimes\dQ$ is the relative Poincar\'{e} line bundle (after identifying $\CH^1(B_{L,U}^{|b|}\times_{X_{L,U}}B_{L,U}^{|b|})\otimes\dQ$ with $\CH^1(B_{L,U}^{|b|}\times_{X_{L,U}}B_{L,U}^{|b|\vee})\otimes\dQ$ via $(\lambda_0\times\lambda^{-1})^{|b|}$). In particular, $c_i^-=(p^{\otimes 2(g-1)|b|})^*c_i^+$ for $i=1,2$ by Definition \ref{de:straightening}. On the other hand, $p$ is symmetric, which implies that
\[
\llangle c_2^+,c_1^-\rrangle_{B_{L,U}^{|b|}}^\ell=\llangle (p^{\otimes 2(g-1)|b|})^*c_2^-,c_1^-\rrangle_{B_{L,U}^{|b|}}^\ell
=\llangle c_2^-,(p^{\otimes 2(g-1)|b|})^*c_1^-\rrangle_{B_{L,U}^{|b|}}^\ell
=\llangle c_2^-,c_1^+\rrangle_{B_{L,U}^{|b|}}^\ell
=\ol{\llangle c_1^+,c_2^-\rrangle_{B_{L,U}^{|b|}}^\ell},
\]
in which the second identity uses \cite{Bei87}*{4.0.3}.

The lemma is proved.
\end{proof}

We end this section by the following easy lemma.

\begin{lem}\label{le:height_pairing_level}
For another object $L'\subseteq L$ of $\Hk(H(\bA_F^\infty))$ satisfying $L'_\ell=L_\ell$, if we denote by $\tau_{L',L}\colon X_{L'}\to X_L$ the transition map, then
\[
\llangle\tau_{L',L}^*(c_1),\tau_{L',L}^*(c_2)\rrangle_{\lambda,L'}^{\ell}
=|L/L'|\cdot\llangle c_1,c_2\rrangle_{\lambda,L}^{\ell}
\]
holds for every pair $(c_1,c_2)\in\rH^{2r_1}(X_L,\rM_{L,\lambda}(r_1))^{\ang\ell}\times\rH^{2r_2}(X_L,\rM_{L,\lambda}(r_2))^{\ang\ell}$.
\end{lem}

\begin{proof}
The lemma follows from \cite{Bei87}*{4.0.3} and the fact that the natural morphism $X_{L',U}\to X_{L,U}$ is finite \'{e}tale of degree $|L/L'|$ for every object $U$ of $\Hk(T(\bA^\infty))$ such that $U_\ell$ is maximal.
\end{proof}

\subsection{Decomposition into local heights}
\label{ss:height_decomposition}

Let $\ell$ be an arbitrary rational prime.

Take an object $L\in\Hk(H(\bA_F^\infty))$. Fix a pair $(r_1,r_2)$ of positive integers satisfying $r_1+r_2=m$ and a pair $\ul{Z}=(Z_1,Z_2)$ of disjoint Zariski closed subsets of $X_L$. In this section, we define local height pairings on $\rH^{2r_1}_{Z_1}(X_L,\rM_{L,\lambda}(r_1))^{\ang\ell}\times\rH^{2r_2}_{Z_2}(X_L,\rM_{L,\lambda}(r_2))^{\ang\ell}$ indexed by elements in $(\fV_E^\fin\setminus\fV_E^{(\ell)})\cup\Hom(E,\dC)$. Write $U_i\coloneqq X_L\setminus Z_i$ for $i=1,2$.

We start from the nonarchimedean case. Take an element $u\in\fV_E^\fin\setminus\fV_E^{(\ell)}$, and we will define a pairing
\[
\llangle\;,\;\rrangle_{\lambda,L,\ul{Z},u}^{\ell}\colon
\rH^{2r_1}_{\et,Z_{1,u}}(X_{L,u},\rM_{L,\ell,\lambda}(r_1))^0
\times\rH^{2r_2}_{\et,Z_{2,u}}(X_{L,u},\rM_{L,\ell,\lambda}(r_2))^0
\to\dE^\lambda_\ell,
\]
in which
\[
\rH^{2r_i}_{\et,Z_{i,u}}(X_{L,u},\rM_{L,\ell,\lambda}(r_i))^0\coloneqq
\Ker\(\rH^{2r_i}_{\et,Z_{i,u}}(X_{L,u},\rM_{L,\ell,\lambda}(r_i))\to
\rH^{2r_i}_{\et}(X_{L,u},\rM_{L,\ell,\lambda}(r_i))\),
\]
as follows. Take a pair $(c_1,c_2)$ from the source. By the Gysin long exact sequence, we may choose an element $\gamma_i\in\rH^{2r_i-1}_{\et}(U_{i,u},\rM_{L,\ell,\lambda}(r_i))$ that maps to $c_i$ under the coboundary map
\[
\rH^{2r_i-1}_{\et}(U_{i,u},\rM_{L,\ell,\lambda}(r_i))\to\rH^{2r_i}_{\et,Z_{i,u}}(X_{L,u},\rM_{L,\ell,\lambda}(r_i))
\]
for $i=1,2$. We then have the cup product
\[
\gamma_1\cup\ol{\gamma_2}\in\rH^{2m-2}_{\et}(U_{1,u}\cap U_{2,u},\rM_{L,\ell,\lambda}(r_1)\otimes_{\dE^\lambda_\ell}\rM_{L,\ell,\ol\lambda}(r_2)).
\]
If we denote by $\gamma_1\vee\ol{\gamma_2}$ the image of $\gamma_1\cup\ol{\gamma_2}$ under the map
\[
\rH^{2m-2}_{\et}(U_{1,u}\cap U_{2,u},\rM_{L,\ell,\lambda}(r_1)\otimes_{\dE^\lambda_\ell}\rM_{L,\ell,\ol\lambda}(r_2))
\to\rH^{2m-2}_{\et}(U_{1,u}\cap U_{2,u},\dE^\lambda_\ell(m))
\]
induced by \eqref{eq:evaluation_etale_automorphic}, then $\llangle c_1,c_2\rrangle_{\lambda,L,\ul{Z},u}^{\ell}$ is defined to be the image of $\gamma_1\vee\ol{\gamma_2}$ under the composite map
\[
\rH^{2m-2}_{\et}(U_{1,u}\cap U_{2,u},\dE^\lambda_\ell(m))\to\rH^{2m-1}_{\et}(X_{L,u},\dE^\lambda_\ell(m))\to
\rH^1(\Spec E_u,\dE^\lambda_\ell(1))=\dE^\lambda_\ell,
\]
in which the first map is the coboundary map in the Mayer--Vietoris exact sequence for the Zariski open covering $X_{L,u}=U_{1,u}\cup U_{2,u}$, and the second map is the trace map for the proper smooth structural morphism $X_{L,u}\to\Spec E_u$. It is easy to check that $\llangle c_1,c_2\rrangle_{\lambda,L,\ul{Z},u}^{\ell}$ does not depend on the choices of $\gamma_1,\gamma_2$.

We then consider the archimedean case. Take an element $\bbu\in\Hom(E,\dC)$, and we will define a pairing
\[
\llangle\;,\;\rrangle_{\lambda,L,\ul{Z}}^{\bbu}\colon
\rH^{2r_1}_{\hs,Z_1^{\bbu}}(X_L^{\bbu},\rM_{L,\lambda}^{\bbu}(r_1))^0
\times\rH^{2r_2}_{\hs,Z_2^{\bbu}}(X_L^{\bbu},\rM_{L,\lambda}^{\bbu}(r_2))^0
\to\dC,
\]
in which
\[
\rH^{2r_i}_{\hs,Z_i^{\bbu}}(X_L^{\bbu},\rM_{L,\lambda}^{\bbu}(r_i))^0\coloneqq
\Ker\(\rH^{2r_i}_{\hs,Z_i^{\bbu}}(X_L^{\bbu},\rM_{L,\lambda}^{\bbu}(r_i))\to
\rH^{2r_i}_{\hs}(X_L^{\bbu},\rM_{L,\lambda}^{\bbu}(r_i))\),
\]
as follows. Take a pair $(c_1,c_2)$ from the source. By the Gysin long exact sequence, we may choose an element $\gamma_i\in\rH^{2r_i-1}_{\hs}(U_i^{\bbu},\rM_{L,\lambda}^{\bbu}(r_i))$ that maps to $c_i$ under the coboundary map
\[
\rH^{2r_i-1}_{\hs}(U_i^{\bbu},\rM_{L,\lambda}^{\bbu}(r_i))\to\rH^{2r_i}_{\hs,Z_i^{\bbu}}(X_L^{\bbu},\rM_{L,\lambda}^{\bbu}(r_i))
\]
for $i=1,2$. We then have the cup product
\[
\gamma_1\cup\ol{\gamma_2}\in\rH^{2m-2}_{\hs}(U_1^{\bbu}\cap U_2^{\bbu},\rM_{L,\lambda}^{\bbu}(r_1)\otimes_{\dE^\lambda_\infty}\rM_{L,\ol\lambda}^{\bbu}(r_2)).
\]
If we denote by $\gamma_1\vee\ol{\gamma_2}$ the image of $\gamma_1\cup\ol{\gamma_2}$ under the map
\[
\rH^{2m-2}_{\hs}(U_1^{\bbu}\cap U_2^{\bbu},\rM_{L,\lambda}^{\bbu}(r_1)\otimes_{\dE^\lambda_\infty}\rM_{L,\ol\lambda}^{\bbu}(r_2))
\to\rH^{2m-2}_{\hs}(U_1^{\bbu}\cap U_2^{\bbu},\dE^\lambda_\infty(m))
\]
induced by \eqref{eq:realization_hodge_automorphic}, then $\llangle c_1,c_2\rrangle_{\lambda,L,\ul{Z}}^{\bbu}$ is defined to be the image of $\gamma_1\vee\ol{\gamma_2}$ under the composite map
\[
\rH^{2m-2}_{\hs}(U_1^{\bbu}\cap U_2^{\bbu},\dE^\lambda_\infty(m))\to\rH^{2m-1}_{\hs}(X_{L,u},\dE^\lambda_\infty(m))\to
\rH^1(\Spec E_u,\dE^\lambda_\infty(1))=\dE^\lambda_\infty\xrightarrow{\kappa}\dC,
\]
in which the first map is the coboundary map in the Mayer--Vietoris exact sequence for the Zariski open covering $X_L^{\bbu}=U_1^{\bbu}\cup U_2^{\bbu}$, the second map is the trace map for the projective complex manifold $X_L^{\bbu}$, and the last map $\kappa$ is the contraction (Notation \ref{no:contraction}). It is easy to check that $\llangle c_1,c_2\rrangle_{\lambda,L,\ul{Z}}^{\bbu}$ does not depend on the choices of $\gamma_1,\gamma_2$.

\begin{lem}\label{le:hodge_purity}
For every integer $r\geq 0$, the natural map
\[
\rH^{2r}_{\hs}(X_L^{\bbu},\rM_{L,\lambda}^{\bbu}(r))\to\rH^{2r}(X_L^{\bbu},\ul\rM_{L,\lambda}^{\bbu}(r))
\]
is injective, where $\ul\rM_{L,\lambda}^{\bbu}$ denotes the underlying $\dE^\lambda_\infty$-linear local system of $\rM_{L,\lambda}^{\bbu}$.
\end{lem}

\begin{proof}
This lemma follows from the essentially same proofs of \cite{Bry89}*{Lemmas~2.2\&2.3}. We sketch the argument. Since the injectivity property does not depend what coefficients we use, we may forget the $\dE^\lambda_\infty$-linear structure and regard $\rM_{L,\lambda}^{\bbu}$ as a real polarized Hodge structure of weight zero. Write $X$ for $X_L^{\bbu}$ and $\rM$ for $\rM_{L,\lambda}^{\bbu}$ for short. We have an exact sequence
\[
\cdots\to\rH^{2r-1}(X,\ul\rM(r))\to\rH^{2r-1}(X,\Omega^\bullet_X(\rM)/\Fil^r\Omega^\bullet_X(\rM))
\to\rH^{2r}_{\hs}(X,\rM(r))\to\rH^{2r}(X,\ul\rM(r))\to\cdots
\]
where $\Omega^\bullet_X(\rM)\coloneqq\Omega^\bullet_X\otimes_\dR\rM$ is a filtered complex of sheaves on $X$. It is known that $\rH^{2r-1}(X,\ul\rM(r))$ admits a natural pure real Hodge structure $\rW$ of weight $2r-1$ so that $\rH^{2r-1}(X,\Omega^\bullet_X(\rM)/\Fil^r\Omega^\bullet_X(\rM))=\rW_\dC/\Fil^r\rW_\dC$ and the first map in the above sequence is given by the natural map $\rW(r)\to\rW_\dC/\Fil^r\rW_\dC$. Then the lemma follows since $\rW(r)\to\rW_\dC/\Fil^r\rW_\dC$ is a bijection.
\end{proof}

By the above lemma and the Betti-\'{e}tale comparison theorem, we know that for every $\bbu\in\Hom(E,\dC)$, the Hodge realization map
\[
\rH^{2r_i}_{Z_i}(X_L,\rM_{L,\lambda}(r_i))^{\ang\ell}\to
\rH^{2r_i}_{\hs,Z_i^{\bbu}}(X_L^{\bbu},\rM_{L,\lambda}^{\bbu}(r_i))
\]
factors through the subspace $\rH^{2r_i}_{\hs,Z_i^{\bbu}}(X_L^{\bbu},\rM_{L,\lambda}^{\bbu}(r_i))^0$. The similar statement obviously holds for every $u\in\fV_E^\fin\setminus\fV_E^{(\ell)}$. In particular, we may define a pairing
\begin{align}\label{eq:height_pairing_bis}
\llangle\;,\;\rrangle_{\lambda,L,\ul{Z}}^{\ell}\colon
\rH^{2r_1}_{Z_1}(X_L,\rM_{L,\lambda}(r_1))^{\ang\ell}\times\rH^{2r_2}_{Z_2}(X_L,\rM_{L,\lambda}(r_2))^{\ang\ell}
\to\dC\otimes\dQ_\ell
\end{align}
as
\begin{align}\label{eq:height_decomposition}
\(\sum_{u\in\fV_E^\fin\setminus\fV_E^{(\ell)}}\log\|u\|\cdot\llangle\;,\;\rrangle_{\lambda,L,\ul{Z},u}^\ell\)+
\(\sum_{\bbu\in\Hom(E,\dC)}\llangle\;,\;\rrangle_{\lambda,L,\ul{Z}}^{\bbu}\).
\end{align}

\begin{proposition}\label{pr:height_decomposition}
Let $(r_1,r_2)$ and $\ul{Z}=(Z_1,Z_2)$ be as above.
\begin{enumerate}
  \item The sum \eqref{eq:height_decomposition} is finite, uniformly for all elements in the source.

  \item Suppose that $\ell$ is unramified in $E$, $L_\ell$ is the stabilizer of a self-dual lattice of $V_\ell$, and that the Zariski closures of $Z_1$ and $Z_2$ in $X_{L,(\ell)}$ do not intersect as well. Then the pairing $\llangle\;,\;\rrangle_{\lambda,L,\ul{Z}}^{\ell}$ \eqref{eq:height_pairing_bis} coincides with the composition of the natural map
      \[
      \rH^{2r_1}_{Z_1}(X_L,\rM_{L,\lambda}(r_1))^{\ang\ell}\times\rH^{2r_2}_{Z_2}(X_L,\rM_{L,\lambda}(r_2))^{\ang\ell}\to
      \rH^{2r_1}(X_L,\rM_{L,\lambda}(r_1))^{\ang\ell}\times\rH^{2r_2}(X_L,\rM_{L,\lambda}(r_2))^{\ang\ell}
      \]
      and the pairing $\llangle\;,\;\rrangle_{\lambda,L}^{\ell}$ \eqref{eq:height_pairing}.
\end{enumerate}
\end{proposition}

\begin{proof}
Denote by $\Diamond_{L,\ul{Z}}$ the smallest set of rational primes such that for every rational prime $p\not\in\Diamond_{L,\ul{Z}}$,
\begin{itemize}
  \item $p$ is unramified in $E$;

  \item $L_p$ is the stabilizer of a self-dual lattice $V_p$;

  \item the Zariski closures of $Z_1$ and $Z_2$ in $X_{L,(p)}$ do not intersect.
\end{itemize}
Then $\Diamond_{L,\ul{Z}}$ is a finite set.

Choose an object $U\in\Hk(T(\bA^\infty))$ such that $U_p$ is maximal for every rational prime $p\not\in\Diamond_{L,\ul{Z}}$. Denote by $\ul{Z}_U\coloneqq(Z_{U,1},Z_{U,2})$ and $\ul{D}_U\coloneqq(D_{U,1},D_{U,2})$ the inverse images of $\ul{Z}=(Z_1,Z_2)$ in $X_{L,U}$ and $B_{L,U}^{|b|}$ \eqref{eq:kuga_sato}, respectively. The maps $\st_{L,U}^{\lambda,+}$ and $\st_{L,U}^{\lambda,-}$ in Definition \ref{de:straightening} induce $\dE^\lambda$-linear maps
\begin{align*}
\st_{L,U}^{\lambda,+}\circ\pi_{L,U}^*\colon\rH^{2r_1}_{Z_1}(X_L,\rM_{L,\lambda}(r_1))^{\ang\ell}
&\to\CH^{r_1+|b|}_{D_{U,1}}(B_{L,U}^{|b|})^{\ang\ell}\otimes\dE^\lambda, \\
\st_{L,U}^{\lambda,-}\circ\pi_{L,U}^*\colon\rH^{2r_2}_{Z_2}(X_L,\rM_{L,\lambda}(r_2))^{\ang\ell}
&\to\CH^{r_2+(g-1)|b|}_{D_{U,2}}(B_{L,U}^{|b|})^{\ang\ell}\otimes\dE^\lambda,
\end{align*}
respectively. Since $D_{U,1}$ and $D_{U,2}$ do not intersect, we have a (hermitian) local height
\[
\llangle\;,\;\rrangle_{B_{L,U}^{|b|},\ul{D}_{L,U},u}^\ell\colon
\CH^{r_1+|b|}_{D_{U,1}}(B_{L,U}^{|b|})^{\ang\ell}\otimes\dE^\lambda\times
\CH^{r_2+(g-1)|b|}_{D_{U,2}}(B_{L,U}^{|b|})^{\ang\ell}\otimes\dE^\lambda\to\dE^\lambda_\ell
\]
for every $u\in\fV_F^\fin$, and a (hermitian) local height
\[
\llangle\;,\;\rrangle_{B_{L,U}^{|b|},\ul{D}_{L,U}}^{\bbu}\colon
\CH^{r_1+|b|}_{D_{U,1}}(B_{L,U}^{|b|})^{\ang\ell}\otimes\dE^\lambda\times
\CH^{r_2+(g-1)|b|}_{D_{U,2}}(B_{L,U}^{|b|})^{\ang\ell}\otimes\dE^\lambda\to\dC
\]
for every $\bbu\in\Hom(E,\dC)$.

Now Lemma \ref{le:straightening} implies that for $u\in\fV_F^\fin\setminus\fV_F^{(\ell)}$, the pairing
\[
\llangle\;,\;\rrangle_{\lambda,L,\ul{Z},u}^\ell\colon
\rH^{2r_1}_{Z_1}(X_L,\rM_{L,\lambda}(r_1))^{\ang\ell}\times\rH^{2r_2}_{Z_2}(X_L,\rM_{L,\lambda}(r_2))^{\ang\ell}
\to\dE^\lambda_\ell
\]
coincides with
\[
\frac{1}{|\Xi_U|}\cdot \llangle\;,\;\rrangle_{B_{L,U}^{|b|},\ul{D}_{L,U},u}^\ell \circ
(\st_{L,U}^{\lambda,+}\circ\pi_{L,U}^*)\times(\st_{L,U}^{\lambda,-}\circ\pi_{L,U}^*),
\]
and that the similar statement holds for $\bbu\in\Hom(E,\dC)$. Note that
\begin{align}\label{eq:height_decomposition_1}
\llangle\;,\;\rrangle_{B_{L,U}^{|b|},\ul{D}_{L,U},u}^\ell \equiv 0, \quad u\in\fV_F^\fin\setminus\fV_F^{(\Diamond_{L,\ul{Z}})}
\end{align}
due to the presence of the smooth projective integral model
\[
B_{L^p,U^p,(p)}^{|b|}\coloneqq
\underbrace{B_{L^p,U^p,(p)}\times_{X_{L^p,U^p,(p)}}\cdots\times_{X_{L^p,U^p,(p)}}B_{L^p,U^p,(p)}}_{|b|}
\]
of $B_{L,U}^{|b|}$, where $B_{L^p,U^p,(p)}\coloneqq A_{0,L^p,U^p,(p)}\times_{X_{L^p,U^p,(p)}}A_{L^p,U^p,(p)}^\vee$ from \S\ref{ss:integral_model} (with $\Diamond=\{p\}$).

Now (1) follows from the above discussion and \eqref{eq:height_decomposition_1}.

For (2), the condition implies that $\ell\not\in\Diamond_{L,\ul{Z}}$. Then it follows from the above discussion, \eqref{eq:height_decomposition_1}, and the decomposition of $\llangle\;,\;\rrangle_{B_{L,U}^{|b|}}^\ell$ into local heights.

The proposition is proved.
\end{proof}

At the end of this section, we describe a way to compute the local height $\llangle\;,\;\rrangle_{\lambda,L,\ul{Z},u}^{\ell}$ for $u\in\fV_E^\fin\setminus\fV_E^{(\ell)}$ in the presence of nice regular models. Fix a prime $u\in\fV_E^\fin\setminus\fV_E^{(\ell)}$ and an object $L$ of $\Hk(H(\bA_F^\infty))$.

\begin{assumption}\label{as:integral_model}
There is a functor
\[
\cX_\bullet\colon \Hk(H(\bA_F^{\infty,p}))\to\Sch^\proj_{/O_{E_u}}
\]
that maps objects to regular schemes and arrows to finite \'{e}tale morphisms, satisfying
\[
\cX_{L^{\prime p}}\otimes_{O_{E_u}}E_u=X_{L_pL^{\prime p}}\otimes_EE_u
\]
functorial in $L^{\prime p}$. To ease notation, we will write $\cX_{L^\prime}$ for $\cX_{L^{\prime p}}$ as long as $L'_p=L_p$, and write $Y_{L'}$ for the special fiber of $\cX_{L'}$.
\end{assumption}

Under Assumption \ref{as:integral_model}, the definition of the $\dE^\lambda_\ell$-linear local system $V_{\ell,L,\lambda}$ from \S\ref{ss:realization_etale} extends to $\cX_L$ as well, so that its restriction to the generic fiber canonically identifies with $\rM_{L,\ell,\lambda}$ by \eqref{eq:realization_local_system}.

\begin{lem}\label{le:height_integral}
Suppose that we are given the functor in Assumption \ref{as:integral_model}. Let $(\cC_1,\cC_2)$ be a pair of elements in $\rH^{2r}_{\et,Z_{1,u}\cup Y_L}(\cX_L,V_{\ell,L,\lambda}(r))\times\rH^{2r}_{\et,Z_{1,u}\cup Y_L}(\cX_L,V_{\ell,L,\lambda}(r))$ satisfying that
\begin{itemize}
  \item if we denote by $c_i$ the image of $\cC_i$ in $\rH^{2r}_{\et,Z_{i,u}}(X_{L,u},\rM_{L,\ell,\lambda}(r))$, then $c_i$ belongs to $\rH^{2r}_{\et,Z_{i,u}}(X_{L,u},\rM_{L,\ell,\lambda}(r))^0$ for $i=1,2$;

  \item the image of either $\cC_1$ or $\cC_2$ in $\rH^{2r}_{\et}(\cX_L,V_{\ell,L,\lambda}(r))$ vanishes.
\end{itemize}
Denote by $\cC_1\overset{\lambda}\bullet\cC_2$ the image of $(\cC_1,\cC_2)$ under the composite map
\begin{align*}
&\quad\rH^{2r}_{\et,Z_{1,u}\cup Y_L}(\cX_L,V_{\ell,L,\lambda}(r))
\times\rH^{2r}_{\et,Z_{2,u}\cup Y_L}(\cX_L,V_{\ell,L,\lambda}(r))  \\ &\to\rH^{4r}_{\et,Y_L}(\cX_L,O_{\dE^\lambda,\ell}(2r))
\to\rH^2_{\et,\Spec\kappa_u}(\Spec O_{E_u},O_{\dE^\lambda,\ell}(1))=O_{\dE^\lambda,\ell},
\end{align*}
in which the first arrow is the cup product induced by the natural hermitian pairing on $V_{\ell,L,\lambda}$, the second arrow is the trace map for $\cX_L/O_{E_u}$, and $\kappa_u$ denotes the residue field of $E_u$. Then
\[
\llangle c_1,c_2\rrangle_{\lambda,L,\ul{Z},u}^\ell=\cC_1\overset\lambda\bullet\cC_2
\]
holds.
\end{lem}

\begin{proof}
The proof of this lemma follows from the same argument for \cite{LL21}*{Lemma~B.4}, together with the obvious analogue of \cite{LL21}*{Remark~B.5}.
\end{proof}

\subsection{Motivic arithmetic inner product formula}
\label{ss:maipf}

In this section, we assume $m=2r$ and that $\lambda$ is balanced (Definition \ref{de:weight}).

Let $V$ be as in previous sections (but of rank $m=2r$). Assume Conjecture \ref{co:modularity} for $V$ so that one can apply Definition \ref{de:theta_lift}.

\begin{definition}\label{de:ell}
Let $\ell$ be a rational prime.
\begin{enumerate}
  \item For an object $L$ of $\Hk(H(\bA_F^\infty))$, we say that $\ell$ is \emph{$L$-good} if $\ell\nmid D_E$ and $L_\ell$ is the stabilizer of a self-dual lattice of $V_\ell$.

  \item For a finite collection $\ul\phi^\infty=(\phi^\infty_1,\phi^\infty_2,\dots)$ of elements in $\sS((V^\infty)^r)$, we say that $\ell$ is \emph{$\ul\phi^\infty$-good} if there exists an object $L$ of $\Hk(H(\bA_F^\infty))$ that fixes $\phi^\infty_i$ for every $i$ such that $\ell$ is $L$-good.
\end{enumerate}
\end{definition}

Consider a holomorphic tempered cuspidal automorphic representation $(\pi,\cA_\pi)$ of $G_r(\bA_F)$ of weight $\lambda$ (Definition \ref{de:hermitian_representation}). Take a pair $\ul\varphi=(\varphi_1,\varphi_2)$ of automorphic forms in $\(\tS_{\lambda^+}(\dC^r)^\vee\otimes_\dC\cA_\pi\)\cap\cA_{\lambda^+}(G_r)_\cusp^\hol$ and a pair  $\ul\phi^\infty=(\phi^\infty_1,\phi^\infty_2)$ of Schwartz functions in $\sS((V^\infty)^r)$. Consider a rational prime $\ell$ that is $\ul\phi^\infty$-good. We define, using \eqref{eq:height_pairing}, an element
\[
\llangle\Theta_{\phi^\infty_1}^\lambda(\varphi_1),\Theta_{\phi^\infty_2}^\lambda(\varphi_2)\rrangle_\lambda^\ell
\in\dC\otimes\dQ_\ell
\]
depending only on $\ul\phi^\infty$ and $\ul\varphi$ as follows.

Choose an intermediate number field $\dE^\lambda\subseteq\dL\subseteq\dC$ over which $\pi^\infty$ is defined. Choose an object $L$ of $\Hk(H(\bA_F^\infty))$ that fixes $\phi^\infty_1$ and $\phi^\infty_2$ such that $\ell$ is $L$-good (which is possible by definition). Choose a finite set $\Diamond$ of rational primes containing $\ell$, such that for every $v\in\fV_F^\fin\setminus\fV_F^{(\Diamond)}$, $(\pi_v)^{K_{r,v}}\neq\{0\}$ (Notation \ref{no:hermitian_siegel}(1)) and $L_v$ is the stabilizer of a self-dual lattice of $V_v$. In particular, $\pi^{\infty,\Diamond}$ gives rise to a Satake character
\begin{align}\label{eq:sakata_character}
\chi_\pi\colon\dL[L^\Diamond\backslash H(\bA_F^{\infty,\Diamond})/L^\Diamond]\to\dL
\end{align}
using the canonical isomorphism $\dZ[L^\Diamond\backslash H(\bA_F^{\infty,\Diamond})/L^\Diamond]\simeq
\dZ[K_r^{\infty,\Diamond}\backslash G_r(\bA_F^{\infty,\Diamond})/K_r^{\infty,\Diamond}]$ of spherical Hecke algebras.

\begin{lem}\label{le:cohomological_annihilation}
If we denote by $\dL[L^\Diamond\backslash H(\bA_F^{\infty,\Diamond})/L^\Diamond]^{\ang\ell}_{\lambda,L}$ the subring of $\dL[L^\Diamond\backslash H(\bA_F^{\infty,\Diamond})/L^\Diamond]$ consisting of elements that annihilate
\[
\prod_{u\in\fV_E^\fin\setminus\fV_E^{(\ell)}}
\rH^{2r}_{\et}(X_{L,u},\rM_{L,\ell,\lambda}(r))\otimes_{\dE^\lambda}\dL
\]
(under pullback), then $\dL[L^\Diamond\backslash H(\bA_F^{\infty,\Diamond})/L^\Diamond]^{\ang\ell}_{\lambda,L}\setminus\Ker\chi_\pi$ is nonempty.
\end{lem}

\begin{proof}
This follows from the same proof of \cite{LL21}*{Proposition~6.9(2)}.
\end{proof}

Choose a pair $(\ts_1,\ts_2)$ of elements in $\dL[L^\Diamond\backslash H(\bA_F^{\infty,\Diamond})/L^\Diamond]^{\ang\ell}_{\lambda,L}\setminus\Ker\chi_\pi$, so that $\ts_i^*\Theta_{\phi^\infty_i}^\lambda(\varphi_i)_L$ belongs to $\rH^{2r_i}(X_L,\rM_{L,\lambda}(r_i))^{\ang\ell}$ for $i=1,2$. We set
\[
\llangle\Theta_{\phi^\infty_1}^\lambda(\varphi_1),\Theta_{\phi^\infty_2}^\lambda(\varphi_2)\rrangle_\lambda^\ell\coloneqq
\vol^\natural(L)\cdot\(\ol{\chi_\pi(\ts_1)}\chi_\pi(\ts_2)\)^{-1}\cdot
\llangle\ts_1^*\Theta_{\phi^\infty_1}^\lambda(\varphi_1)_L,
\ts_2^*\Theta_{\phi^\infty_2}^\lambda(\varphi_2)_L\rrangle_{\lambda,L}^{\ell}
\in\dC\otimes\dQ_\ell,
\]
where $\vol^\natural(L)$ is the normalized volume of $L$ \cite{LL21}*{Definition~3.8}.

\begin{lem}
The quantity $\llangle\Theta_{\phi^\infty_1}^\lambda(\varphi_1),\Theta_{\phi^\infty_2}^\lambda(\varphi_2)\rrangle_\lambda^\ell$ does not depend on the choices of $\dL$, $L$ (that fixes $\ul\phi^\infty$ with prescribed $L_\ell$ as above), $\Diamond$, and $(\ts_1,\ts_2)$.
\end{lem}

\begin{proof}
The independence of the choice of $\dL$ is obvious.

We first show the independence of $(\ts_1,\ts_2)$. Choose another pair $(\ts'_1,\ts'_2)$ with respect to a possibly different set $\Diamond'$. By Lemma \ref{le:hecke_translate}, we have
\[
\ts_i^{\prime*}\ts_i^*\Theta_{\phi^\infty_i}^\lambda(\varphi_i)_L=
\ts_i^*\Theta_{\phi^\infty_i}^\lambda(\ts'_i\varphi_i)_L
\]
for $i=1,2$. By \cite{Liu11}*{Proposition~A.5}, we have
\[
\ts_i^*\Theta_{\phi^\infty_i}^\lambda(\ts'_i\varphi_i)_L=\ts_i^*\Theta_{\ts'_i\phi^\infty_i}^\lambda(\varphi_i)_L
=\ol{\chi_\pi(\ts'_i)}\cdot\ts_i^*\Theta_{\phi^\infty_i}^\lambda(\varphi_i)_L.
\]
It follows that
\begin{align*}
&\quad\(\ol{\chi_\pi(\ts_1)}\chi_\pi(\ts_2)\)^{-1}\cdot
\llangle\ts_1^*\Theta_{\phi^\infty_1}^\lambda(\varphi_1)_L,
\ts_2^*\Theta_{\phi^\infty_2}^\lambda(\varphi_2)_L\rrangle_{\lambda,L}^{\ell} \\
&=\(\ol{\chi_\pi(\ts'_1\ts_1)}\chi_\pi(\ts'_2\ts_2)\)^{-1}\cdot
\llangle\ts_1^{\prime*}\ts_1^*\Theta_{\phi^\infty_1}^\lambda(\varphi_1)_L,
\ts_2^{\prime*}\ts_2^*\Theta_{\phi^\infty_2}^\lambda(\varphi_2)_L\rrangle_{\lambda,L}^{\ell} \\
&=\(\ol{\chi_\pi(\ts_1\ts'_1)}\chi_\pi(\ts_2\ts'_2)\)^{-1}\cdot
\llangle\ts_1^*\ts_1^{\prime*}\Theta_{\phi^\infty_1}^\lambda(\varphi_1)_L,
\ts_2^*\ts_2^{\prime*}\Theta_{\phi^\infty_2}^\lambda(\varphi_2)_L\rrangle_{\lambda,L}^{\ell} \\
&=\(\ol{\chi_\pi(\ts'_1)}\chi_\pi(\ts'_2)\)^{-1}\cdot
\llangle\ts_1^{\prime*}\Theta_{\phi^\infty_1}^\lambda(\varphi_1)_L,
\ts_2^{\prime*}\Theta_{\phi^\infty_2}^\lambda(\varphi_2)_L\rrangle_{\lambda,L}^{\ell}.
\end{align*}
Thus, we obtain the independence of $(\ts_1,\ts_2)$. The independence of $\Diamond$ follows from the independence of $(\ts_1,\ts_2)$. Finally, the independence of $L$ follows from Lemma \ref{le:height_pairing_level} the independence of $\Diamond$.

The lemma is proved.
\end{proof}

\begin{conjecture}[Motivic arithmetic inner product formula]\label{co:maipf}
Suppose that the hermitian $E$-weight $\lambda$ is of rank $2r$ and is balanced. Let $(\pi,\cA_\pi)$ be a holomorphic tempered cuspidal automorphic representation of $G_r(\bA_F)$ of weight $\lambda$ (Definition \ref{de:hermitian_representation}). For every pair $\ul\varphi=(\varphi_1,\varphi_2)$ of automorphic forms in $\(\tS_{\lambda^+}(\dC^r)^\vee\otimes_\dC\cA_\pi\)\cap\cA_{\lambda^+}(G_r)_\cusp^\hol$, and very pair $\ul\phi^\infty=(\phi^\infty_1,\phi^\infty_2)$ of Schwartz functions in $\sS((V^\infty)^r)$,
\[
\llangle\Theta_{\phi^\infty_1}^\lambda(\varphi_1),\Theta_{\phi^\infty_2}^\lambda(\varphi_2)\rrangle_\lambda^\ell
=\frac{L'(\tfrac{1}{2},\pi)}{b_{2r}(0)}\cdot C(\lambda)\cdot
\fZ^\natural_\fin(\ul\varphi,\phi^\infty_1\otimes\ol{\phi^\infty_2})
\]
holds for all rational primes $\ell$ that are $\ul\phi^\infty$-good, in which $C(\lambda)$ is the real constant from Proposition \ref{pr:ripf}, and $\fZ^\natural_\fin(\ul\varphi,\phi^\infty_1\otimes\ol{\phi^\infty_2})$ is introduced in Definition \ref{de:doubling_zeta}. In particular, the left-hand side belongs to $\dC$ and is independent of $\ell$.
\end{conjecture}

\begin{remark}\label{re:dichotomy}
We have the following remarks concerning the above conjecture.
\begin{enumerate}
  \item By Theorem \ref{th:independence} below, we know a priori that, when $r=1$,
      \[
      \llangle\Theta_{\phi^\infty_1}^\lambda(\varphi_1),\Theta_{\phi^\infty_2}^\lambda(\varphi_2)\rrangle_\lambda^\ell
      \]
      takes value in $\dC$ and is independent of $\ell$.

  \item By the local theta dichotomy, for every $v\in\fV_F$, there exists a unique up to isomorphism hermitian space $V(\pi_v)$ over $E_v$ of rank $m=2r$, such that the local theta lift of $\pi_v$ to $\rU(V(\pi_v))$ is nontrivial. Put $V(\pi)\coloneqq\bigotimes_v V(\pi_v)$ as a hermitian space over $\bA_E/\bA_F$. Then the condition that $\lambda$ is balanced implies that $V(\pi)$ is totally positive definite. Moreover, $V(\pi)$ is coherent/incoherent (Definition \ref{de:space_coherence}) if and only if $\epsilon(\pi)=+1/{-1}$ (Definition \ref{de:hermitian_representation}).

  \item If $V$ is not isomorphic to $V(\pi)$, then both terms
      \[
      \llangle\Theta_{\phi^\infty_1}^\lambda(\varphi_1),\Theta_{\phi^\infty_2}^\lambda(\varphi_2)\rrangle_\lambda^\ell,\quad
      \fZ^\natural_\fin(\ul\varphi,\phi^\infty_1\otimes\ol{\phi^\infty_2})
      \]
      always vanish by the theta dichotomy. In particular, when $\epsilon(\pi)=1$, both sides in the conjecture are zero.

  \item When $\epsilon(\pi)=-1$ (so that $L(\tfrac{1}{2},\pi)=0$) and $V=V(\pi)$,\footnote{Due to our technical assumption on $V$, this forces $|\fV(V(\pi))|>1$, which always holds if $F\neq \dQ$.} one can always find $\ul\varphi$ and $\ul\phi^\infty$ such that $\fZ^\natural_\fin(\ul\varphi,\phi^\infty_1\otimes\ol{\phi^\infty_2})\neq 0$.\footnote{By \cite{LL21}*{Proposition~3.13}, it is even true that for every finite set $\Diamondblack$ of nonarchimedean places of $F$, there exists a pair $\ul\phi^\infty$ of decomposable Schwartz functions such that $(\phi^\infty_1\otimes\ol{\phi^\infty_2})_v$ is regularly supported (Definition \ref{de:regular_support}) for $v\in\Diamondblack$ and that the functional $\ul\varphi\mapsto\fZ^\natural_\fin(\ul\varphi,\phi^\infty_1\otimes\ol{\phi^\infty_2})$ is nonzero.} Thus, the conjecture indeed provides a formula for the central $L$-derivative $L'(\tfrac{1}{2},\pi)$.

  \item When ``everything is unramified'', that is, $E/F$ is everywhere unramified, $(\phi^\infty_0\coloneqq)\phi^\infty_1=\phi^\infty_2$ is the characteristic function of $\Lambda^r$ for a self-dual lattice $\Lambda$ of $V^\infty$,\footnote{The existence of $\Lambda$ implies that $r[F:\dQ]$ is odd, and $F\neq\dQ$ when $r>1$.} $(\varphi_0\coloneqq)\varphi_1=\varphi_2$ is fixed by $\prod_{v\in\fV_F^\fin}K_{r,v}$, the conjecture particularly implies that
      \[
      \llangle\Theta_{\phi^\infty_0}^\lambda(\varphi_0),\Theta_{\phi^\infty_0}^\lambda(\varphi_0)\rrangle_\lambda^\ell
      =\frac{L'(\tfrac{1}{2},\pi)}{b_{2r}(0)}\cdot C(\lambda)\cdot(\varphi_0,\varphi_0)_{G_r},
      \]
      where $(\;,\;)_{G_r}$ is the inner product from Definition \ref{de:hermitian_siegel}(2).
\end{enumerate}
\end{remark}

To end this section, we formula an analogue of the Beilinson--Bloch conjecture for the motivic cohomology
\[
\rH^{2r}(X,\rM_{\lambda}(r))\coloneqq\varinjlim_L\rH^{2r}(X_L,\rM_{L,\lambda}(r))
\]
as a module over $\dE^\lambda[H(\bA_F^\infty)]$.

\begin{remark}\label{re:endoscopic}
Consider a holomorphic tempered cuspidal automorphic representation $(\pi,\cA_\pi)$ of $G_r(\bA_F)$ of weight $\lambda$ with $\Pi\coloneqq\BC(\pi)$ its base change to $E$, which is a hermitian isobaric automorphic representation of $\GL_{2r}(\bA_E)$ in the sense of \cite{BPLZZ21}*{Definition~1.5}. For every irreducible admissible representation $\tilde\pi^\infty$ of $H(\bA_F^\infty)$ that is nearly equivalent to $\pi^\infty$, one has an associated isobaric factor $\Pi_{\tilde\pi^\infty}$ of $\Pi$, determined via the way in \cite{LL21}*{Lemma~3.5}.
\end{remark}

\begin{conjecture}[Beilinson--Bloch]\label{co:beilinson}
Suppose that the hermitian $E$-weight $\lambda$ is of rank $2r$ and is balanced. Let $V$ be a totally positive definite incoherent hermitian space over $\bA_E/\bA_F$ of rank $2r$ satisfying $|\fV(V)|>1$.
\begin{enumerate}
  \item The $\dE^\lambda[H(\bA_F^\infty)]$-module $\rH^{2r}(X,\rM_{\lambda}(r))$ is admissible and semisimple.

  \item For every holomorphic tempered cuspidal automorphic representation $(\pi,\cA_\pi)$ of $G_r(\bA_F)$ of weight $\lambda$, the relation
      \[
      \dim_\dC\Hom_{H(\bA_F^\infty)}\(\tilde\pi^\infty,\rH^{2r}(X,\rM_{\lambda}(r))\otimes_{\dE^\lambda}\dC\)=
      \ord_{s=1/2}L(s,\Pi_{\tilde\pi^\infty})
      \]
      holds for every irreducible admissible representation $\tilde\pi^\infty$ of $H(\bA_F^\infty)$ that is nearly equivalent to $\pi^\infty$.
\end{enumerate}
\end{conjecture}

For $(\pi,\cA_\pi)$ in Conjecture \ref{co:beilinson}(2) satisfying $\epsilon(\pi)=-1$ and $|\fV(V(\pi))|>1$ (and we take $V=V(\pi)$), the assignment $(\phi^\infty,\varphi)\mapsto\Theta_{\phi^\infty}^\lambda(\varphi)$ indeed induces an element in $\Hom_{H(\bA_F^\infty)}\(\theta(\pi^\infty),\rH^{2r}(X,\rM_{\lambda}(r))\otimes_{\dE^\lambda}\dC\)$. Thus, Conjecture \ref{co:maipf} provides with the following evidence:
\[
\ord_{s=1/2}L(s,\pi)=1 \quad\Longrightarrow\quad
\dim_\dC\Hom_{H(\bA_F^\infty)}\(\theta(\pi^\infty),\rH^{2r}(X,\rM_{\lambda}(r))\otimes_{\dE^\lambda}\dC\)\geq 1
\]
toward Conjecture \ref{co:beilinson}(2). Indeed, in this case, $\Pi_{\theta(\pi^\infty)}$ is the unique isobaric factor of $\Pi$ with global root number $-1$.

\subsection{Independence-of-$\ell$ for height pairing over unitary Shimura curves}
\label{ss:independence}

In this section, we assume $m=2$.

\begin{theorem}\label{th:independence}
For every object $L$ of $\Hk(H(\bA_F^\infty))$, the height pairing $\llangle\;,\;\rrangle_{\lambda,L}^{\ell}$ \eqref{eq:height_pairing} (for $L$-good $\ell$ in the sense of Definition \ref{de:ell}(1)) takes values in $\dC$ and is independent of $\ell$.
\end{theorem}

When $\lambda=\b{0}$, the above theorem is well-known as the height pairing coincides with the N\'{e}ron--Tate height pairing. It remains to consider the case where $\lambda\neq\b{0}$. In this case,
\[
\rH^2(X_L,\rM_{L,\lambda}(1))^{\ang\ell}=\rH^2(X_L,\rM_{L,\lambda}(1))
\]
since $\rH^2_{\et}(X_{L,u},r_\ell\rM_{L,\lambda}(1))$ vanishes for every $u\in\fV_E^\fin\setminus\fV_E^{(\ell)}$, so at least the source of $\llangle\;,\;\rrangle_{\lambda,L}^{\ell}$ is independent of $\ell$.

\begin{definition}\label{re:optimal}
For every $v\in\fV_F^\fin$, we define an integer $t_v$ (depending also on $V$) as follows:
\center
\begin{tabular}{|l|c|c|}
  \hline
   &  $V_v$ split &  $V_v$ nonsplit \\ \hline
   $v$ split in $E$ & $0$ & n/a \\ \hline
   $v$ inert in $E$  & $0$ & $1$ \\ \hline
   $v$ ramified in $E$ & $0$ & $2$ \\
  \hline
\end{tabular}
\center
\end{definition}

Take a finite set $\Diamond$ of rational primes. The definition below extends Definition \ref{de:moduli_integral}.

\begin{definition}\label{de:moduli_integral_curve}
We define a functor
\begin{align*}
X^\flat_{\bullet,\bullet,(\Diamond)}\colon\Hk(H(\bA_F^{\infty,\Diamond}))\times\Hk(T(\bA^{\infty,\Diamond}))
\to\sfP\Sch'_{/O_{E,(\Diamond)}}
\end{align*}
such that for every objects $(L^\Diamond,U^\Diamond)\in\Hk(H(\bA_F^{\infty,\Diamond}))\times\Hk(T(\bA^{\infty,\Diamond}))$ and $S\in\Sch'_{/O_{E,(\Diamond)}}$, $X^\flat_{L^\Diamond,U^\Diamond,(\Diamond)}(S)$ is the set of equivalence classes of sextuples $(A_0,\lambda_0,\bar\eta_0^\Diamond;A,\lambda,\bar\eta^\Diamond)$, where
\begin{itemize}
  \item $(A_0,\lambda_0,\bar\eta_0^\Diamond)$ is an object of $Y_{U^\Diamond,(\Diamond)}(S)$ (Definition \ref{de:moduli_auxiliary_integral});

  \item $A$ is an object of $\AV(S,O_{E,(\Diamond)})$ (Definition \ref{de:category_av}) satisfying the Eisenstein condition and the signature $(1,1)$ condition \cite{LMZ}*{Definition~6.2} for
      \[
      A=\{\phi\in\Hom(E,E^\ac)\res\phi\res_{E_0}=\id_{E_0}\}\setminus\{\id_E\},\quad B=\{\phi\in\Hom(E,E^\ac)\res\phi\res_{E_0}=\id_{E_0}\}\cup\{\tc\};
      \footnote{Here, we are using the abstract algebraic closure $E^\ac$ instead of $\ol\dQ$ in \cite{LMZ}.}
      \]

  \item $\lambda\colon A\to A^\vee$ is a polarization that turns the complex conjugation of $O_{E,(\Diamond)}$ into the Rosati involution, such that for every $v\in\fV_F^{(\Diamond)}$, $\Ker\lambda[\fp_v^\infty]$ is a semisimple $O_{E_v}$-module of length $t_v$;

  \item $\bar\eta^\Diamond$ is a \emph{hermitian $L^\Diamond$-level structure}, that is, for a chosen geometric point $s$ on every connected component of $S$, a $\pi_1(S,s)$-invariant $L^\Diamond$-orbit of isometries
      \[
      \eta^\Diamond\colon V^{\infty,\Diamond}\xrightarrow\sim\Hom_{\bA_E^{\infty,\Diamond}}^{\lambda_0,\lambda}
      \(\rH_1^{\et}(A_{0s},\bA^{\infty,\Diamond})\otimes_{\bA^{\infty,\Diamond}}\bA_F^{\infty,\Diamond},
      \rH_1^{\et}(A_s,\bA^{\infty,\Diamond})\)
      \]
      of hermitian spaces over $\bA_E^{\infty,\Diamond}/\bA_F^{\infty,\Diamond}$.
\end{itemize}
Two objects $(A_0,\lambda_0,\bar\eta_0^\Diamond;A,\lambda,\bar\eta^\Diamond)$ and $(A'_0,\lambda'_0,\bar\eta_0^{\Diamond\prime};A',\lambda',\bar\eta^{\Diamond\prime})$ are equivalent if there exist isomorphisms $\varphi_0\colon A_0\to A'_0$ in the category $\AV(S,O_{E_0,(\Diamond)})$ and $\varphi\colon A\to A'$ in the category $\AV(S,O_{E,(\Diamond)})$ satisfying that
\begin{itemize}
  \item $\varphi_0^\vee\circ\lambda'_0\circ\varphi_0$ and $\varphi^\vee\circ\lambda'\circ\varphi$ are the \emph{same} prime-to-$\Diamond$ multiple of $\lambda_0$ and $\lambda$, respectively;

  \item the $U^\Diamond$-orbits $\varphi_{0*}\circ\bar\eta_0^\Diamond$ and $\bar\eta_0^{\Diamond\prime}$ coincide;

  \item the $L^\Diamond$-orbit of maps $v\mapsto\varphi_*\circ\bar\eta^\Diamond(v)\circ(\varphi_{0*})^{-1}$ for $v\in V^{\infty,\Diamond}$ coincides with $\bar\eta^{\Diamond\prime}$.
\end{itemize}
The values on morphisms are defined similarly as in Definition \ref{de:moduli_integral}.
\end{definition}

Note that over $S\otimes_{O_{E,(\Diamond)}}E$, the Eisenstein condition plus the signature $(1,1)$ condition is equivalent to that $\tr_A(e)=2\Tr_{E/E_0}(e)-e+\ol{e}$ for $e\in O_{E,(\Diamond)}$ by \cite{LMZ}*{Example~6.3}. It follows that if we choose an open compact subgroup $L^\flat_\Diamond$ of $H(F_\Diamond)$ that is optimal (see Definition \ref{de:principal_level} below) and denote by $U_\Diamond$ the maximal open compact subgroup of $T(\dQ_\Diamond)$, then it is clear from the construction that
\[
X^\flat_{L^\Diamond,U^\Diamond,(\Diamond)}\otimes_{O_{E,(\Diamond)}}E=X_{L^\flat_\Diamond L^\Diamond,U_\Diamond U^\Diamond}
\]
functorial in $(L^\Diamond,U^\Diamond)$.

\begin{definition}\label{de:principal_level}
We say that an open compact subgroup $L^\flat_\Diamond$ of $H(F_\Diamond)$ is \emph{optimal} if $(L^\flat_\Diamond)_v$ is the stabilizer of a (necessarily unique) vertex lattice of $V_v$ of index $t_v$ for every $v\in\fV_F^{(\Diamond)}$.
\end{definition}

Similar to Definition \ref{de:moduli_integral}, $X^\flat_{\bullet,\bullet,(\Diamond)}$ factors through the full subcategory $\Sch^\proj_{/O_{E,(\Diamond)}}$ of $\sfP\Sch'_{/O_{E,(\Diamond)}}$, and we have universal objects
\begin{align*}
A_{0,\bullet,\bullet,(\Diamond)}/X^\flat_{\bullet,\bullet,(\Diamond)}&\colon
\Hk(H(\bA_F^{\infty,\Diamond}))\times\Hk(T(\bA^{\infty,\Diamond}))
\to\AV(\Sch^\proj_{/O_{E,(\Diamond)}},O_{E_0,(\Diamond)}),\\
A_{\bullet,\bullet,(\Diamond)}/X^\flat_{\bullet,\bullet,(\Diamond)}&\colon
\Hk(H(\bA_F^{\infty,\Diamond}))\times\Hk(T(\bA^{\infty,\Diamond}))
\to\AV(\Sch^\proj_{/O_{E,(\Diamond)}},O_{E,(\Diamond)})
\end{align*}
over $X^\flat_{\bullet,\bullet,(\Diamond)}$ that are Cartesian (Definition \ref{de:category_section}), and the functor
\begin{align*}
X^\flat_{\bullet,(\Diamond)}\colon\Hk(H(\bA_F^{\infty,\Diamond}))\to\Sch^\proj_{/O_{E,(\Diamond)}}
\end{align*}
such that for every $L^\Diamond\in\Hk(H(\bA_F^{\infty,\Diamond}))$, $X^\flat_{L^\Diamond,(\Diamond)}=X^\flat_{L^\Diamond,U^\Diamond,(\Diamond)}/\Xi_{U^\Diamond}$. By construction, it is clear that
\[
X^\flat_{L^\Diamond,(\Diamond)}\otimes_{O_{E,(\Diamond)}}E=X_{L^\flat_\Diamond L^\Diamond}
\]
functorial in $L^\Diamond$.

\begin{remark}\label{re:uniformization}
Choose an object $U^\Diamond\in\Hk(T(\bA^{\infty,\Diamond}))$ and an optimal open compact subgroup $L^\flat_\Diamond$ of $H(F_\Diamond)$. In this remark, we recall the basic uniformization of $X^\flat_{\bullet,U^\Diamond,(\Diamond)}$ at a place $u\in\fV_E^{(\Diamond)}$ with $p=p(u)$ the underlying rational prime.

Denote by $\breve{E}_u$ a complete maximal unramified extension of $E_u$ and $\kappa_u$ its residue field, which is an algebraic closure of $\dF_p$. Choose a \emph{framing}, that is, an element
\[
P^\star\coloneqq (A_0^\star,\lambda_0^\star,\bar\eta_0^{\Diamond\star};A^\star,\lambda^\star,\eta^{\Diamond\star})\in
\varprojlim_{L^\Diamond}X^\flat_{L^\Diamond,U^\Diamond,(\Diamond)}(O_{\breve{E}_u})
\]
whose reduction is basic, where the inverse limit is taken over all $L^\Diamond\in\Hk(H(\bA_F^{\infty,\Diamond}))$ with respect to inclusions. The framing $P^\star$ gives rise to
\begin{description}
  \item[(U1)] a totally positive definite hermitian space
      \[
      V^\star\coloneqq
      \Hom_{O_{E,(\Diamond)}}(A^\star_0\otimes_{O_{\breve{E}_u}}\kappa_u,A^\star\otimes_{O_{\breve{E}_u}}\kappa_u)\otimes\dQ
      \]
      over $E/F$ (and put $H^\star\coloneqq\rU(V^\star)$);

  \item[(U2)] an isometry $(V^\star)^{\infty,\Diamond}\simeq V^{\infty,\Diamond}$ given by $\eta^{\Diamond\star}$, and an $(L^\flat_\Diamond)_v$-class of isometries $V^\star\otimes_FF_v\simeq V\otimes_{\bA_F}F_v$ for every $v\in\fV_F^{(\Diamond)}\setminus\{v(u)\}$;

  \item[(U3)] a relative Rapoport--Zink space, that is, a formal scheme $\fN^\flat$ locally of finite type over $O_{\breve{E}_u}$ classifying quasi-isogenies of $(A^\star[\fp_{v(u)}],\lambda^\star[\fp_{v(u)}])\otimes_{O_{\breve{E}_u}}\kappa_u$ of signature $(1,1)$ and type $t_{v(u)}$ \cite{LMZ}*{(5.20)}.
\end{description}

Denote by $\fX^\flat_{L^\Diamond,U^\Diamond}$ the formal completion of $X^\flat_{L^\Diamond,U^\Diamond,(\Diamond)}\otimes_{O_{E,(\Diamond)}}O_{\breve{E}_u}$ along the basic locus. By \cite{LMZ}*{Theorem~6.8}, we have the following isomorphism of basic uniformization
\begin{align}\label{eq:uniformization_1}
\fX^\flat_{L^\Diamond,U^\Diamond}\simeq (H^\star(F)\times O_{E_0,(\Diamond)}^\times)\backslash\(\fN^\flat\times H(\bA_F^{\infty,v(u)})/(L^\flat_\Diamond L^\Diamond)^{v(u)}\times T(\bA^{\infty,\Diamond})/U^\Diamond\)
\end{align}
functorial in $L^\Diamond\in\Hk(H(\bA_F^{\infty,\Diamond}))$, together with a quasi-isogeny
\begin{align}\label{eq:uniformization_2}
A^\star\otimes_{O_{\breve{E}_u}}(\fN^\flat_{L^\Diamond,U^\Diamond}\otimes_{O_{\breve{E}_u}}\kappa_u)\dashrightarrow
A_{L^\Diamond,U^\Diamond,(\Diamond)}\otimes_{X^\flat_{L^\Diamond,U^\Diamond,(\Diamond)}}
(\fN^\flat_{L^\Diamond,U^\Diamond}\otimes_{O_{\breve{E}_u}}\kappa_u)
\end{align}
preserving polarizations, where
\[
\fN^\flat_{L^\Diamond,U^\Diamond}\coloneqq
\fN^\flat\times H(\bA_F^{\infty,v(u)})/(L^\flat_\Diamond L^\Diamond)^{v(u)}\times T(\bA^{\infty,\Diamond})/U^\Diamond.
\]
Note that when $V_{v(u)}$ is nonsplit, the entire special fiber $X^\flat_{L^\Diamond,U^\Diamond,(\Diamond)}\otimes_{O_{E,(\Diamond)}}\kappa_u$ is basic.
\end{remark}

\begin{definition}\label{de:regular_closure}
Let $R$ be a Dedekind domain with the generic point $\eta$ and $\cX$ a scheme over $R$ locally of finite type.
\begin{enumerate}
  \item We denote by $\cX^{\r{fl}}$ the flat closure of $\cX$ in its generic fiber $\cX_\eta$.

  \item For a finite \'{e}tale cover $X'\to\cX_\eta$ over $\eta$ in which $X'$ is smooth of dimension one, we define the \emph{regular closure} of $\cX$ in $X'$ to be the canonical resolution of (isolated) singularities of the normalization of $\cX^{\r{fl}}$ in $X'$.
\end{enumerate}
Note that the regular closure is regular and admits a natural map to $\cX$.

We have the similar definition when $R$ is a complete local ring, $\cX$ is a formal scheme over $R$ locally of finite type, and $X'\to\cX_\eta$ is a finite \'{e}tale cover (of rigid analytic spaces) over $\eta$ in which $X'$ is a smooth rigid analytic space of dimension one.
\end{definition}

Now we consider an open compact subgroup $L_\Diamond$ of $H(F_\Diamond)$ together with an optimal one $L_\Diamond^\flat$ containing $L_\Diamond$.

\begin{definition}\label{de:moduli_integral_regular}
We define a functor
\begin{align*}
X_{\bullet,\bullet,(\Diamond)}\colon\Hk(H(\bA_F^{\infty,\Diamond}))\times\Hk(T(\bA^{\infty,\Diamond}))
\to\Sch^\proj_{/O_{E,(\Diamond)}}
\end{align*}
such that for every objects $(L^\Diamond,U^\Diamond)\in\Hk(H(\bA_F^{\infty,\Diamond}))\times\Hk(T(\bA^{\infty,\Diamond}))$, $X_{L^\Diamond,U^\Diamond,(\Diamond)}$ is the regular closure of $X_{L^\Diamond,U^\Diamond,(\Diamond)}^\flat$ in $X_{L_\Diamond L^\Diamond, U_\Diamond U^\Diamond}$.\footnote{The functor $X_{\bullet,\bullet,(\Diamond)}$ depends on both $L_\Diamond$ and $L^\flat_\Diamond$. Moreover, even when $L_\Diamond=L^\flat_\Diamond$, the natural natural transformation $X_{\bullet,\bullet,(\Diamond)}\to X^\flat_{\bullet,\bullet,(\Diamond)}$ may not be an isomorphism.}
\end{definition}

Choose an object $U^\Diamond\in\Hk(T(\bA^{\infty,\Diamond}))$ and take a place $u\in\fV_E^{(p)}$ as above. Denote by $\fN_\eta\to\fN^\flat_\eta$ the classifying space of $L_{v(u)}$-level structure over $\fN^\flat_\eta$ -- the rigid analytic space of $\fN^\flat$, which is a finite \'{e}tale cover, and define $\fN$ to be the regular closure of $\fN^\flat$ in $\fN_\eta$. If we denote by $\fX_{L^\Diamond, U^\Diamond}$ the formal completion of $X_{L^\Diamond,U^\Diamond,(\Diamond)}\otimes_{O_{E,(\Diamond)}}O_{\breve{E}_u}$ along the basic locus, then since taking regular closure is functorial and commutes with taking formalization, \eqref{eq:uniformization_1} implies the following isomorphism of basic uniformization
\begin{align}\label{eq:uniformization_3}
\fX_{L^\Diamond,U^\Diamond}\simeq (H^\star(F)\times O_{E_0,(\Diamond)}^\times)\backslash\(\fN\times H(\bA_F^{\infty,v(u)})/(L_\Diamond L^\Diamond)^{v(u)}\times T(\bA^{\infty,\Diamond})/U^\Diamond\)
\end{align}
functorial in $L^\Diamond\in\Hk(H(\bA_F^{\infty,\Diamond}))$.

\begin{proof}[Proof of Theorem \ref{th:independence}]
As we have explained, it suffices to consider the case where $\lambda\neq\b{0}$.

Take $\Diamond$ to be the smallest set of rational primes that contains those ramified in $E$ such that $L^\Diamond$ is the stabilizer of a self-dual lattice of $V^{\infty,\Diamond}$. We may shrink $L_\Diamond$ so that it is contained in an optimal open compact subgroup $L^\flat_\Diamond$ of $H(F_\Diamond)$. We still denote by $A_{0,\bullet,\bullet,(\Diamond)}$ and $A_{\bullet,\bullet,(\Diamond)}$ their restrictions to $X_{\bullet,\bullet,(\Diamond)}$.

Choose an object $U^\Diamond\in\Hk(T(\bA^{\infty,\Diamond}))$ and put $U=U_\Diamond U^\Diamond$ with $U_\Diamond$ the maximal open compact subgroup of $T(\dQ_\Diamond)$. Recall that we have the map
\[
\st_{L,U}^{\lambda,+}\colon\rH^2(X_{L,U},\rM_{L,U,\lambda}(1))
\to\CH^{1+|b|}(B_{L,U}^{|b|})\otimes\dE^\lambda
\]
which is canonically a direct summand. By Remark \ref{re:abstract_motive}, we may extend the relative motive $\rM_{L^\Diamond,U^\Diamond,\lambda}$ to an object in $\Mot(X_{L^\Diamond,U^\Diamond,(\Diamond)},\dE^\lambda)$ (with the same notation) as the abstract automorphic motive over $X_{L^\Diamond,U^\Diamond,(\Diamond)}$ of weight $\lambda$ associated with $(A_{0,L^\Diamond,U^\Diamond,(\Diamond)},\lambda_0)$ and $(A_{L^\Diamond,U^\Diamond,(\Diamond)},\lambda)$, together with a canonical direct summand map
\begin{align}\label{eq:independence}
\st_{L^\Diamond,U^\Diamond,(\Diamond)}^{\lambda,+}\colon\rH^2(X_{L^\Diamond,U^\Diamond,(\Diamond)},\rM_{L,U,\lambda}(1))
\to\CH^{1+|b|}(B_{L^\Diamond,U^\Diamond,(\Diamond)}^{|b|})\otimes\dE^\lambda
\end{align}
by Remark \ref{re:relative_motivic}(1), where $B_{L^\Diamond,U^\Diamond,(\Diamond)}\coloneqq A_{0,L^\Diamond,U^\Diamond,(\Diamond)}\times_{X_{L^\Diamond,U^\Diamond,(\Diamond)}}A_{L^\Diamond,U^\Diamond,(\Diamond)}^\vee$.

Now take a pair $(c_1,c_2)$ of elements in $\rH^2(X_{L,U},\rM_{L,\lambda}(1))$. For $i=1,2$, put $c_i^\pm\coloneqq\st_{L,U}^{\lambda,\pm}(\pi_{L,U}^*(c_i))$. By the construction of \eqref{eq:height_pairing}, the goal is to show that
\[
\llangle c_1^+,c_2^-\rrangle_{B_{L,U}^{|b|}}^\ell
\]
takes value in $\dC$ and is independent of $\ell$ for every $\ell\not\in\Diamond$. Now, using the direct summand map $\st_{L^\Diamond,U^\Diamond,(\Diamond)}^{\lambda,+}$ above, we may assume that $c_1^+$ extends to a class $c_{1,(\Diamond)}\in\CH^{1+|b|}(B_{L^\Diamond,U^\Diamond,(\Diamond)}^{|b|})\otimes\dE^\lambda$ that is contained in the image of $\st_{L^\Diamond,U^\Diamond,(\Diamond)}^{\lambda,+}$. Choose an arbitrary extension $c_{2,(\Diamond)}\in\CH^{1+(g-1)|b|}(B_{L^\Diamond,U^\Diamond,(\Diamond)}^{|b|})\otimes\dE^\lambda$ of $c_2^-$. By the moving lemma, we may assume that $c_{1,(\Diamond)}$ and $c_{2,(\Diamond)}$ are represented by algebraic cycles $C_{1,(\Diamond)}$ and $C_{2,(\Diamond)}$ on $B_{L^\Diamond,U^\Diamond,(\Diamond)}^{|b|}$ (with coefficients in $\dE^\lambda$) with disjoint support on $B_{L,U}^{|b|}$. For $i=1,2$, put $C_i\coloneqq C_{i,(\Diamond)}\res_{B_{L,U}^{|b|}}$ whose support is denoted by $D_i$. By the decomposition formula for the usual Beilinson--Bloch height pairing, it suffices to show that for every $u\in\fV_E^\fin$,
\[
\llangle C_1,C_2\rrangle_{B_{L,U}^{|b|},\ul{D},u}^\ell
\]
takes value in $\dE^\lambda$ and is independent of $\ell$ for every $\ell\not\in\Diamond$, where $\ul{D}=(D_1,D_2)$. When $u\not\in\fV_E^{(\Diamond)}$, this holds since $B_{L,U}^{|b|}$ has smooth projective reduction at $u$. It remains to consider places in $\fV_E^{(\Diamond)}$.

Take an element $u\in\fV_E^{(\Diamond)}$ with $v=v(u)$. Write
\[
\cX_{L^\Diamond,U^\Diamond}\coloneqq X_{L^\Diamond,U^\Diamond,(\Diamond)}\otimes_{O_{E,(\Diamond)}}O_{\breve{E}_u},\quad
\cB_{L^\Diamond,U^\Diamond}\coloneqq B_{L^\Diamond,U^\Diamond,(\Diamond)}\otimes_{O_{E,(\Diamond)}}O_{\breve{E}_u},
\]
and $\cC_i\coloneqq C_{i,(\Diamond)}\otimes_{O_{E,(\Diamond)}}O_{\breve{E}_u}$ for $i=1,2$. Below we will regard $\rH^2(\cX_{L^\Diamond,U^\Diamond},\rM_{L,U,\lambda}(1))$ as a direct summand of $\CH^{1+|b|}(\cB_{L^\Diamond,U^\Diamond}^{|b|})\otimes\dE^\lambda$ via \eqref{eq:independence}. For every $\ell\not\in\Diamond$, put
\[
\rH^2(\cX_{L^\Diamond,U^\Diamond},\rM_{L,U,\lambda}(1))^{\ang{\ell}}\coloneqq
\Ker\(\rH^2(\cX_{L^\Diamond,U^\Diamond},\rM_{L,U,\lambda}(1))
\to\rH^2_{\et}(\cX_{L^\Diamond,U^\Diamond},r_\ell\rM_{L,U,\lambda}(1))\).
\]
Similar to Lemma \ref{le:realization_etale} and \eqref{eq:realization_local_system}, we have a canonical isomorphism $r_\ell\rM_{L,U,\lambda}\simeq V_{\ell,L,\lambda}$ over $\cX_{L^\Diamond,U^\Diamond}$.

Our strategy is the following: Find an algebraic cycle $\cC_1^\dag$ on $\cB_{L^\Diamond,U^\Diamond}^{|b|}$ with coefficients in $\dE^\lambda$ whose restriction to $B_{L,U}^{|b|}\otimes_{E}\breve{E}_u$ coincides with that of $\cC_1$, such that its image in $\CH^{1+|b|}(\cB_{L^\Diamond,U^\Diamond}^{|b|})\otimes\dE^\lambda$ belongs to the subspace $\rH^2(\cX_{L^\Diamond,U^\Diamond},\rM_{L,U,\lambda}(1))^{\ang{\ell}}$ for every $\ell\not\in\Diamond$. Then by \cite{LL21}*{Proposition~B.10},
\[
\llangle C_1,C_2\rrangle_{B_{L,U}^{|b|},\ul{D},u}^\ell=\cC_1^\dag.\ol{\cC_2},
\]
hence is in $\dE^\lambda$ and is independent of $\ell$ for every $\ell\not\in\Diamond$. By Lemma \ref{le:vanishing_ordinary}(2) below, it suffices to verify that the cycle class of $\cC_1^\dag$ in $\rH^2_{\et}(C,V_{\ell,L,\lambda}(1))$ vanishes for every irreducible component $C$ of $\cX_{L^\Diamond,U^\Diamond}\otimes_{O_{\breve{E}_u}}\kappa_u$ contained in the basic locus and every $\ell\not\in\Diamond$.\footnote{Once again, when $V_{v(u)}$ is nonsplit, every such irreducible component is contained in the basic locus.}

Denote by $\rM^\star_\lambda$ the abstract automorphic motive (Remark \ref{re:abstract_motive}) over $\kappa_u$ of weight $\lambda$ associated with $(A^\star_0,\lambda^\star_0)$ and $(A^\star,\lambda^\star)$. The two maps $\st^{\lambda,+}$ and $\st^{\lambda,-}$ from the same remark induce a hermitian pairing
\[
\rH^0(\Spec\kappa_u,\rM^\star_\lambda)\times\rH^0(\Spec\kappa_u,\rM^\star_\lambda)\to\dE^\lambda
\]
via the Poincar\'{e} duality on $(B^\star)^{|b|}$, where $B^\star\coloneqq A_0^\star\times_{\kappa_u}(A^\star)^\vee$. Put $V^\star_\lambda\coloneqq\rH^0_\num(\Spec\kappa_u,\rM^\star_\lambda)$ -- the $\dE^\lambda$-linear space of numerical equivalence classes of $\rH^0(\Spec\kappa_u,\rM^\star_\lambda)$ under the above pairing.\footnote{Indeed, the proof of Lemma \ref{le:abstract_motive} below implies that $V^\star_\lambda$ is canonically isomorphic to $(V^\star)_\lambda$ via Construction \ref{co:harmonic_function}(2), hence justifying the notation.} By Lemma \ref{le:abstract_motive} below, the natural map $V^\star_\lambda\otimes\dQ_\ell\to r_\ell\rM^\star_\lambda$ is an isomorphism for every $\ell\not\in\Diamond$. On the other hand, the level structure $\eta^{\Diamond\star}$ (at $\ell$ alone) induces an isomorphism $r_\ell\rM^\star_\lambda\simeq V_{\ell,\lambda}\coloneqq(V_\ell)_\lambda$. In summary, the $\ell$-adic realization map induces an isomorphism
\begin{align}\label{eq:independence_5}
V^\star_\lambda\otimes\dQ_\ell\simeq V_{\ell,\lambda}
\end{align}
for every $\ell\not\in\Diamond$, depending only on the framing $P^\star$.

Denote by $N_{L^\Diamond,U^\Diamond}$ the normalization of the induced reduced subscheme of
\[
\fN_{L^\Diamond,U^\Diamond}\coloneqq\fN\times H(\bA_F^{\infty,v(u)})/(L_\Diamond L^\Diamond)^{v(u)}\times T(\bA^{\infty,\Diamond})/U^\Diamond,
\]
and by $|N_{L^\Diamond,U^\Diamond}|$ the set of its irreducible components. Consider an element $C\in|N_{L^\Diamond,U^\Diamond}|$ with $\pi_C\colon C\to\Spec\kappa_u$ the (smooth projective) structure map. Using the maps $\st^{\lambda,+}$ and $\st^{\lambda,-}$ in Remark \ref{re:abstract_motive} (with $S=C$), we obtain a hermitian pairing
\[
\rH^2(C,\rM_{L,U,\lambda}(1))\times\rH^2(C,\rM_{L,U,\lambda}(1))\to\dE^\lambda
\]
via the Poincar\'{e} duality on $\cB_{L^\Diamond,U^\Diamond}^{|b|}\times_{\cX_{L^\Diamond,U^\Diamond}}C$. Denote by $\rH^2_\num(C,\rM_{L,U,\lambda}(1))$ the $\dE^\lambda$-linear space of numerical equivalence classes of $\rH^2(C,\rM_{L,U,\lambda}(1))$ under the above pairing. By \cite{Kle68}*{Theorem~3.5}, the map
\begin{align}\label{eq:independence_2}
\rH^2_\num(C,\rM_{L,U,\lambda}(1))\otimes{\dQ_\ell}\to\rH^2_{\et}(C,r_\ell\rM_{L,U,\lambda}(1))
\end{align}
induced by the $\ell$-adic realization is injective for every $\ell\not\in\Diamond$. The universal quasi-isogeny \eqref{eq:uniformization_2} gives rise to an isomorphism $\rM_{L,U,\lambda}\res_{C}\simeq\pi_{C}^*\rM^\star_\lambda$, which in particular gives rise to an isomorphism $r_\ell\rM_{L,U,\lambda}\res_{C}\simeq\pi_{C}^* V_{\ell,\lambda}$. Then the map \eqref{eq:independence_2} can be written as
\[
\rH^2_\num(C,\pi_{C}^*\rM^\star_\lambda(1))\otimes{\dQ_\ell}\to
\rH^2_{\et}(C,V_{\ell,\lambda}(1))
=\rH^2_{\et}(C,\dQ_\ell(1))\otimes_{\dQ_\ell}V_{\ell,\lambda}\simeq V_{\ell,\lambda}.
\]
Consider the Gysin map
\begin{align}\label{eq:independence_3}
\rH^0_\num(\Spec\kappa_u,\rM^\star_\lambda)\to\rH^2_\num(C,\pi_{C}^*\rM^\star_\lambda(1))
\end{align}
at an arbitrary closed point of $C$ (which does not depend on the choice of such a point). By the injectivity of \eqref{eq:independence_2} and the isomorphism \eqref{eq:independence_5}, the consecutive maps
\[
\rH^0_\num(\Spec\kappa_u,\rM^\star_\lambda)\otimes{\dQ_\ell}
\to\rH^2_\num(C,\pi_{C}^*\rM^\star_\lambda(1))\otimes{\dQ_\ell}
\to\rH^2_{\et}(C,\dQ_\ell(1))\otimes_{\dQ_\ell}V_{\ell,\lambda}
\]
are both isomorphisms for every $\ell\not\in\Diamond$, which implies that \eqref{eq:independence_3} is also an isomorphism.

By \eqref{eq:uniformization_3} and the isomorphism \eqref{eq:independence_3}, the image of the natural composite map
\begin{align}\label{eq:independence_4}
\rH^2(\cX_{L^\Diamond,U^\Diamond},\rM_{L,U,\lambda}(1))
&\to\prod_{C\in|N_{L^\Diamond,U^\Diamond}|}\rH^2(C,\rM_{L,U,\lambda}(1)) \\
&\to\prod_{C\in|N_{L^\Diamond,U^\Diamond}|}\rH^2_\num(C,\rM_{L,U,\lambda}(1))\xrightarrow\sim
\Fun(|N_{L^\Diamond,U^\Diamond}|,V^\star_\lambda) \notag
\end{align}
is contained in $\Fun(|N_{L^\Diamond,U^\Diamond}|,V^\star_\lambda)^\Gamma$, where $\Gamma\coloneqq H^\star(F)\times O_{E_0,(\Diamond)}^\times$ acts on $V^\star_\lambda$ via the projection $H^\star(F)\times O_{E_0,(\Diamond)}^\times\to H^\star(F)$. On the other hand, if we denote by $Y_{L^\Diamond,U^\Diamond}$ the special fiber of $\cX_{L^\Diamond,U^\Diamond}$, then we have the Gysin map
\[
\zeta\colon\rH^2_{Y_{L^\Diamond,U^\Diamond}}(\cX_{L^\Diamond,U^\Diamond},\rM_{L,U,\lambda}(1))\to
\rH^2(\cX_{L^\Diamond,U^\Diamond},\rM_{L,U,\lambda}(1)).
\]
Composing $\zeta$ with \eqref{eq:independence_4}, we obtain the map
\begin{align}\label{eq:independence_1}
\rH^2_{Y_{L^\Diamond,U^\Diamond}}(\cX_{L^\Diamond,U^\Diamond},\rM_{L,U,\lambda}(1))
\to\Fun(|N_{L^\Diamond,U^\Diamond}|,V^\star_\lambda)^\Gamma.
\end{align}
We claim that \eqref{eq:independence_1} is surjective. Assuming the claim, we may find an element
\[
\cC'_1\in\rH^2_{Y_{L^\Diamond,U^\Diamond}}(\cX_{L^\Diamond,U^\Diamond},\rM_{L,U,\lambda}(1))
\]
whose image under \eqref{eq:independence_1} coincides with the image of $\cC_1$ under \eqref{eq:independence_4}. Then $\cC^\dag_1\coloneqq\cC_1-\zeta\cC'_1$ works for our strategy.

It remains to show the claim. Composing the natural composite map
\begin{align*}
\bigoplus_{C\in|N_{L^\Diamond,U^\Diamond}|}\rH^0(\Spec\kappa_u,\rM^\star_\lambda)
\xrightarrow{\oplus\pi_C^*}
\bigoplus_{C\in|N_{L^\Diamond,U^\Diamond}|}\rH^0(C,\rM_{L,U,\lambda})\to
\rH^2_{Y_{L^\Diamond,U^\Diamond}}(\cX_{L^\Diamond,U^\Diamond},\rM_{L,U,\lambda}(1))
\end{align*}
with \eqref{eq:independence_1}, we obtain the map
\begin{align}\label{eq:independence_6}
\bigoplus_{C\in|N_{L^\Diamond,U^\Diamond}|}\rH^0(\Spec\kappa_u,\rM^\star_\lambda)
\to\Fun(|N_{L^\Diamond,U^\Diamond}|,V^\star_\lambda).
\end{align}
It suffices to show that \eqref{eq:independence_6} is surjective. For every abelian group $A$, denote by $\Fun_c(|N_{L^\Diamond,U^\Diamond}|,A)\subseteq\Fun(|N_{L^\Diamond,U^\Diamond}|,A)$ the subgroup of finitely supported functions (valued in $A$), for which we have the map
\[
\Delta\colon\dZ[|N_{L^\Diamond,U^\Diamond}|]\to\Fun_c(|N_{L^\Diamond,U^\Diamond}|,\dZ)
\]
given by the intersection in (the regular formal scheme) $\fN_{L^\Diamond,U^\Diamond}$, and the map
\[
\r{av}^\Gamma\colon\Fun_c(|N_{L^\Diamond,U^\Diamond}|,V^\star_\lambda)\to
\Fun(|N_{L^\Diamond,U^\Diamond}|,V^\star_\lambda)^\Gamma
\]
sending $f$ to $\sum_{\gamma\in\Gamma}\gamma f$. Then \eqref{eq:independence_6} is given by the composition
\begin{align*}
\bigoplus_{C\in|N_{L^\Diamond,U^\Diamond}|}\rH^0(\Spec\kappa_u,\rM^\star_\lambda)
&\to\bigoplus_{C\in|N_{L^\Diamond,U^\Diamond}|}\rH^0_\num(\Spec\kappa_u,\rM^\star_\lambda)
=\dZ[|N_{L^\Diamond,U^\Diamond}|]\otimes V^\star_\lambda \\
&\xrightarrow{\Delta\otimes V^\star_\lambda}
\Fun_c(|N_{L^\Diamond,U^\Diamond}|,V^\star_\lambda)\xrightarrow{\r{av}^\Gamma}
\Fun(|N_{L^\Diamond,U^\Diamond}|,V^\star_\lambda)^\Gamma.
\end{align*}
Since the first map is surjective, it suffices to show that $\r{av}^\Gamma\circ(\Delta\otimes V^\star_\lambda)$ is surjective, or rather the induced map
\[
\Delta_\lambda\colon\(\dZ[|N_{L^\Diamond,U^\Diamond}|]\otimes V^\star_\lambda\)_\Gamma
\to\Fun(|N_{L^\Diamond,U^\Diamond}|,V^\star_\lambda)^\Gamma
\]
is surjective. However, by a similar argument, $\Delta_\lambda\otimes\dQ_\ell$ (for an arbitrary rational prime $\ell\not\in\Diamond$) can be canonically identified with the $\ell$-adic Gysin map
\[
\zeta_\ell\colon\rH^2_{\et,Y_{L^\Diamond,U^\Diamond}}(\cX_{L^\Diamond,U^\Diamond},r_\ell\rM_{L,U,\lambda}(1))\to
\rH^2_{\et}(\cX_{L^\Diamond,U^\Diamond},r_\ell\rM_{L,U,\lambda}(1)).
\]
Note the the cokernel of $\zeta_\ell$ is contained in the $\Gal(\breve{E}_u/E_u)$-finite subspace of $\rH^2_{\et}(X_{L,U}\otimes_E\breve{E}_u,r_\ell\rM_{L,U,\lambda}(1))$, which vanishes when $\lambda\neq\b{0}$. It follows that $\r{av}^\Gamma\circ(\Delta\otimes V^\star_\lambda)$ hence \eqref{eq:independence_1} are surjective.

The theorem is finally proved!
\end{proof}

\begin{lem}\label{le:vanishing_ordinary}
Let $u\in\fV_E^{(\Diamond)}$ be an element such that $V_{v(u)}$ is split.
\begin{enumerate}
  \item The base change $(X^\flat_{L^\Diamond,U^\Diamond,(\Diamond)})^{\r{fl}}\otimes_{O_{E,(\Diamond)}}O_{E_u}$ is smooth over $O_{E_u}$.

  \item For every irreducible component $C$ of $X_{L^\Diamond,U^\Diamond,(\Diamond)}\otimes_{O_{E,(\Diamond)}}\kappa_u$ that is not contained in the basic locus, we have
      \[
      \rH^2_{\et}(C,V_{\ell,L,\lambda}(1))=0
      \]
      for every $\ell\not\in\Diamond$ when $\lambda\neq\b{0}$.
\end{enumerate}
\end{lem}

\begin{proof}
For (1), it is straightforward to see that $X^\flat_{L^\Diamond,U^\Diamond,(\Diamond)}\otimes_{O_{E,(\Diamond)}}O_{E_u}$ is already smooth away from the basic locus. Thus, by \eqref{eq:uniformization_1}, it suffices to show that $(\fN^\flat)^{\r{fl}}$ is a formal open unit disc of relative dimension one. This is known when $v(u)$ splits in $E$. When $v(u)$ is nonsplit in $E$, there is a natural map from the formal moduli $\fM$ of Lubin--Tate formal $O_{E_u}$-divisible groups to $\fN^\flat$ sending $Y$ to $Y\oplus\ol{Y}$ with the natural polarization $\lambda$, which factors through $(\fN^\flat)^{\r{fl}}$ as $\fM$ is flat and induces an isomorphism between generic fibers -- this implies that $(\fN^\flat)^{\r{fl}}\simeq\fM$ which is a formal open unit disc of relative dimension one.

For (2), denote by $C^\flat$ the image of $C$ in $X^\flat_{L^\Diamond,U^\Diamond,(\Diamond)}\otimes_{O_{E,(\Diamond)}}\kappa_u$, which is a smooth irreducible curve by (1). Denote by $H^\der$ denotes the derived subgroup of $H$. It follows by (1) that the monodromy group of $V_{\ell,L,\lambda}$ over $C^\flat$ is an open compact subgroup of $H^\der(F_\ell)$, acting via the representation $(V_\ell)_\lambda$ which remains irreducible when restricting to $H^\der(F_\ell)$. It follows that $\rH^2_{\et}(C,V_{\ell,L,\lambda}(1))=0$.
\end{proof}

\begin{lem}\label{le:abstract_motive}
In Remark \ref{re:abstract_motive} suppose that $S=\Spec\kappa$ for an algebraically closed field $\kappa$. If the $E$-vector space $\Hom_{E_0}(A_0,A)$ is of dimension $m$, then the natural map
\[
\rH^0(\Spec\kappa,\rM_\lambda)\otimes\dQ_\ell\to\rH^0_\et(\Spec\kappa,r_\ell\rM_\lambda)=r_\ell\rM_\lambda
\]
is surjective for every rational prime $\ell$ invertible on $\kappa$.
\end{lem}

\begin{proof}
Take a rational prime $\ell$ that is invertible on $\kappa$. It is clear that
\[
r_\ell\rM=\Hom_{E_{0,\ell}}\(\sfh^1_{\et}(A_0,\dQ_\ell),\sfh^1_{\et}(A,\dQ_\ell)\)
\]
in which $\sfh^1_{\et}(-,\dQ_\ell)=\rV_\ell(-)[-1]$ is nothing but the rational $\ell$-adic Tate module placed in degree $1$. It is well-known that the natural map
\begin{align}\label{eq:abstract_motive}
\Hom_{E_0}(A_0,A)\otimes\dQ_\ell\to\Hom_{E_{0,\ell}}\(\rV_\ell(A_0),\rV_\ell(A)\)
\end{align}
is injective. Since the target is a free $E_\ell$-module of rank $m$, it follows that \eqref{eq:abstract_motive} is an isomorphism. On the other hand, \eqref{eq:abstract_motive} factors as
\[
\Hom_{E_0}(A_0,A)\otimes\dQ_\ell\to\Hom_{E_0}(\sfh_1(A_0),\sfh_1(A))\otimes\dQ_\ell=\rH^0(\Spec\kappa,\rM)\otimes\dQ_\ell
\to\Hom_{E_{0,\ell}}\(\rV_\ell(A_0),\rV_\ell(A)\),
\]
which implies that the map
\[
\rH^0(\Spec\kappa,\rM)\otimes\dQ_\ell\to r_\ell\rM
\]
of $E_\ell$-modules is surjective. Applying Construction \ref{co:harmonic_function}(2) to $\rH^0(\Spec\kappa,\rM)$ (resp.\ $r_\ell\rM$), we obtain an $\dE^\lambda$-module (resp.\ $\dE^\lambda_\ell$-module) $\rH^0(\Spec\kappa,\rM)_\lambda$ (resp.\ $(r_\ell\rM)_\lambda$), such that the induced map
\[
\rH^0(\Spec\kappa,\rM)_\lambda\otimes\dQ_\ell\to (r_\ell\rM)_\lambda
\]
is surjective. However, $(r_\ell\rM)_\lambda$ is canonically isomorphic to $r_\ell\rM_\lambda$ and the above map factors naturally as
\[
\rH^0(\Spec\kappa,\rM)_\lambda\otimes\dQ_\ell\to \rH^0(\Spec\kappa,\rM_\lambda)\otimes\dQ_\ell\to r_\ell\rM_\lambda.
\]
The lemma follows.
\end{proof}

\section{Computation of local heights of motivic special classes}
\label{ss:5}

In this chapter, we compute local heights of motivic special classes in various situations.

We fix a totally positive definite incoherent hermitian space $V$ over $\bA_E/\bA_F$ of rank $2r$ satisfying that the set $\fV(V)$ is not a singleton (and put $H\coloneqq\rU(V)$ as always), and a balanced hermitian $E$-weight $\lambda$ of rank $2r$ (Definition \ref{de:weight}).

To simplify notation, write $\Psi^\lambda_\infty$ and $\Phi^\lambda_\infty$ for $\Psi^\lambda_{V_\infty^{2r}}$ and $\Phi^\lambda_{V_\infty^{2r}}$ from Construction \ref{co:harmonic_function}(4) and Definition \ref{de:doubling_schwartz}, respectively. It is clear that as functions on $V_\infty^{2r}$ valued in $\tS_{\lambda^+}(\dC^r)^\vee\otimes_\dC\ol{\tS_{\lambda^+}(\dC^r)^\vee}$, both $\Psi^\lambda_\infty$ and $\Phi^\lambda_\infty$ factor through the moment map $T\colon V_\infty^{2r}\to\Herm_{2r}(F_\infty)^+$.

\subsection{Explicit local heights at inert places}
\label{ss:height_inert}

In this section, we take a nonarchimedean place $u$ of $E$ such that $v\coloneqq v(u)$ is inert in $E$ and $p\coloneqq p(u)$ is odd.

Consider an object $L\in\Hk(H(\bA_F^\infty))$ such that
\begin{itemize}
  \item $L_v$ is the stabilizer of a self-dual or almost self-dual lattice (Definition \ref{de:lattice}), say $\Lambda_v$, of $V_v$, depending on whether $V_v$ is split or nonsplit;

  \item $L_{v'}$ is the stabilizer of a vertex lattice (Definition \ref{de:lattice}), say $\Lambda_{v'}$, for every place $v'\in\fV_F^{(p)}$ that is nonsplit in $E$.
\end{itemize}

Let
\[
\cX_\bullet\colon \Hk(H(\bA_F^{\infty,p}))\to\Sch^\proj_{/O_{E_u}}
\]
be the functor (with respect to $L$) constructed in \cite{LZ22}*{\S11.2~or~\S11.3} (see Remark \ref{re:moduli_integral} below), satisfying Assumption \ref{as:integral_model} with the further property that $\cX_{L'}$ is smooth or strictly semistable over $O_{E_u}$, depending on whether $V_v$ is split or not.

\begin{remark}\label{re:moduli_integral}
The integral model $\cX_{L'}=\cX_{L^{\prime p}}$ is constructed in the following way: First, choose an object $U^p$ of $\Hk(T(\bA^{\infty,p}))$. In \cite{LZ22}*{\S11.2~or~\S11.3},\footnote{In fact, in \cite{LZ22} the authors assumed that $F_v/\dQ_p$ is unramified; however, this assumption can be removed with all the properties unchanged by the recent work \cite{LMZ}.} the authors constructed a moduli problem $\cX_{L^{\prime p},U^p}$ classifying similar data in Definition \ref{de:moduli_integral} (together with the Eisenstein condition), hence admits a natural free action of $\Xi_{U^p}$. We then put $\cX_{L^{\prime p}}\coloneqq\cX_{L^{\prime p},U^p}/\Xi_{U^p}$, which is canonically independent of the choice of $U^p$.
\end{remark}

\begin{proposition}\label{pr:height_inert}
Let $L\in\Hk(H(\bA_F^\infty))$ be as above. Consider
\begin{itemize}
  \item an intermediate number field $\dE^\lambda\subseteq\dL\subseteq\dC$,

  \item a finite set $\Diamond$ of rational primes containing $p$ and a nonempty subset $\Diamondblack\subseteq\fV_F^{(\Diamond)}\setminus\{v\}$,

  \item a pair $\ul\phi^\infty=(\phi^\infty_1,\phi^\infty_2)$ of decomposable Schwartz functions in $\sS((V^\infty)^r,\dL)^L$ satisfying that $(\phi^\infty_1\otimes\ol{\phi^\infty_2})_{v'}$ is regularly supported (Definition \ref{de:regular_support}) for $v'\in\Diamondblack$; and $(\phi^\infty_i)_{v'}$ is the characteristic function of $\Lambda_{v'}^r$ for $i=1,2$ for $v'\in\fV_F^{(p)}$ that is nonsplit in $E$,

  \item a rational prime $\ell\neq p$, and

  \item a pair $(\ts_1,\ts_2)$ of elements in $\dL[L^\Diamond\backslash H(\bA_F^{\infty,\Diamond})/L^\Diamond]$ that annihilate both $\rH^{2r}_{\et}(X_{L,u},V_{\ell,L,\lambda}(r))\otimes_{\dE^\lambda}\dL$ and $\rH^{2r}_{\et}(\cX_L,V_{\ell,L,\lambda}(r))\otimes_{\dE^\lambda}\dL$.
\end{itemize}
\begin{enumerate}
  \item Suppose that $V_v$ is split. Then for every pair $(T_1,T_2)$ of elements in $\Herm_r^\circ(F)^+$ (so that $\ul{Z}=(Z_1\coloneqq Z_{T_1}[\ts_1\phi^\infty_1]_L,Z_2\coloneqq Z_{T_2}[\ts_2\phi^\infty_2]_L)$ (Definition \ref{de:special_cycle}) is a disjoint pair) and every pair $(g_1,g_2)\in G_r(F_\infty)^2$,
      \begin{align*}
      \log\|u\|\cdot\vol^\natural(L)\cdot
      \llangle\ts_1^*\bbZ^\lambda_{T_1}(\phi^\infty_1)_L(g_1),
      \ts_2^*\bbZ^\lambda_{T_2}(\phi^\infty_2)_L(g_2)\rrangle_{\lambda,L,\ul{Z},u}^\ell
      =\fE_{T_1,T_2}((g_1,g_2),\Phi^\lambda_\infty\cdot(\ts_1\phi^\infty_1\otimes\ol{\ts_2\phi^\infty_2}))_v
      \end{align*}
      holds in $\tS_{\lambda^+}(\dC^r)^\vee\otimes_\dC\ol{\tS_{\lambda^+}(\dC^r)^\vee}\otimes\dQ_\ell$, where $\bbZ^\lambda_{T_i}(\phi^\infty_i)_L$ and $\fE_{T_1,T_2}$ are introduced in Definition \ref{de:special_cycle}(1) and Definition \ref{de:analytic_kernel}(2), respectively.

  \item Suppose that $V_v$ is nonsplit. Then for every pair $(T_1,T_2)$ of elements in $\Herm_r^\circ(F)^+$ (so that $\ul{Z}=(Z_1\coloneqq Z_{T_1}[\ts_1\phi^\infty_1]_L,Z_2\coloneqq Z_{T_2}[\ts_2\phi^\infty_2]_L)$ is a disjoint pair) and every pair $(g_1,g_2)\in G_r(F_\infty)^2$,
      \begin{align*}
      &\quad\log\|u\|\cdot\vol^\natural(L)\cdot
      \llangle\ts_1^*\bbZ^\lambda_{T_1}(\phi^\infty_1)_L(g_1),
      \ts_2^*\bbZ^\lambda_{T_2}(\phi^\infty_2)_L(g_2)\rrangle_{\lambda,L,\ul{Z},u}^\ell\\
      &=\fE_{T_1,T_2}((g_1,g_2),\Phi^\lambda_\infty\cdot(\ts_1\phi^\infty_1\otimes\ol{\ts_2\phi^\infty_2}))_v
      -\frac{\log\|u\|}{\|u\|^r-1}
      E_{T_1,T_2}((g_1,g_2),\Phi^\lambda_\infty\cdot(\ts_1\phi^\infty_1\otimes\ol{\ts_2\phi^\infty_2})^v\cdot\Phi^\star_v)
      \end{align*}
      holds in $\tS_{\lambda^+}(\dC^r)^\vee\otimes_\dC\ol{\tS_{\lambda^+}(\dC^r)^\vee}\otimes\dQ_\ell$, where $E_{T_1,T_2}$ is introduced in Definition \ref{de:analytic_kernel}(1) and $\Phi_v^\star$ is the characteristic function of $(\Lambda_v^\star)^{2r}$ in which $\Lambda_v^\star$ denotes the (unique up to isomorphism) hermitian space over $O_{E_v}/O_{F_v}$ of rank $2r$.
\end{enumerate}
\end{proposition}

\begin{proof}
We show (2) and leave the similar and less involved case (1) to readers.

By Lemma \ref{le:doubling_schwartz} and Definition \ref{de:special_cycle}(1), it suffices to consider the case where $g_i=m(a_i)$ with $a_i\in\GL_r(E_\infty)$ for $i=1,2$. Furthermore, since for every $T^\Box\in\Herm_{2r}^\circ(F)^+$,
\[
\(\rho_{\lambda^+}(m(a_1))^\vee\otimes\ol{\rho_{\lambda^+}(m(a_2))^\vee}\)\Phi^\lambda_\infty(T^\Box)
=\Phi^\lambda_\infty(\pres{t}a T^\Box a)
\]
where $a=\(\begin{smallmatrix} a_1 &  \\  & a_2 \end{smallmatrix}\)$, we may just assume $g_1=g_2=1_{2r}$.

Write $\Phi_v\coloneqq(\phi^\infty_1)_v\otimes\ol{(\phi^\infty_2)_v}$ (which is the characteristic function of $\Lambda_v^{2r}$) and $\Phi^{\infty,v}\coloneqq\ts_1(\phi^\infty_1)^v\otimes\ol{\ts_2(\phi^\infty_2)^v}$ for short.

For every element $T^\Box\in\Herm_{2r}^\circ(F_v)$ whose determinant has odd valuation, put
\[
\chi(T^\Box)\coloneqq\frac{b_{2r,v}(0)}{\log\|u\|}\(\frac{\|u\|^{2r}-1}{\|u\|+1}W'_{T^\Box}(0,1_{4r},\Phi_v)
-\frac{\log\|u\|}{\|u\|+1}W_{T^\Box}(0,1_{4r},\Phi_v^\star)\),
\]
which is indeed an integer. By the definition of $\vol^\natural(L)$ in \cite{LL21}*{Definition~3.8}, (2) reduces to the identity
\begin{align}\label{eq:height_inert}
\llangle\ts_1^*\bbZ^\lambda_{T_1}(\phi^\infty_1)_L,\ts_2^*\bbZ^\lambda_{T_2}(\phi^\infty_2)_L\rrangle_{\lambda,L,\ul{Z},u}^\ell
=\sum_{\substack{T^\Box\in\Herm_{2r}^\circ(F)^+\\ \partial T^\Box=(T_1,T_2)}}
\kappa\Psi^\lambda_\infty(T^\Box)\cdot\chi(T^\Box)\cdot\(\sum_{x\in L^v\backslash (V^{\infty,v})^{2r}_{T^\Box}}
\Phi^{\infty,v}(x)\),
\end{align}
where we recall that $\Psi^\lambda_\infty$ is from Construction \ref{co:harmonic_function}(4), which factors through the moment map. Note that on the right-hand side of \eqref{eq:height_inert}, if $T^\Box$ has even valuation at $v$, then $(V^{\infty,v})^{2r}_{T^\Box}$ is empty.

For every integer $1\leq s\leq 2r$ and every element $x^v\in(V^{\infty,v})^s$ with $T=T(x^v)\in\Herm_s^\circ(F)^+$, we have a natural finite unramified morphism $\varsigma_T(x^v)_L\colon\cZ_T(x^v)_L\to\cX_L$ whose generic fiber is
\[
\coprod_{\substack{x=x_v x^v \in (V^\infty)^s_T \\ x_v\in L_v\backslash\Lambda_v^r}} Z(x)_{L,u} \to X_{L,u}
\]
via moduli interpretation of $\cX_{L,U}$ (after choosing $U$ as in Remark \ref{re:moduli_integral}), and denote by $\cZ_T[x^v]_L$ its image as a Zariski closed subset of $\cX_L$. It also carries a map
\begin{align}\label{eq:special_morphism_4}
\mu_T(x^v)_L\colon\varsigma_T(x^v)_L^*V_{\ell,L}^\vee\to(\mathbbm{1}^\et_{\cZ_T(x^v)_L,E_\ell})^{\oplus s}
\end{align}
similar to \eqref{eq:special_morphism_3}.

Now take $x^v=(x^v_1,\dots,x^v_r)\in(V^{\infty,v})^r$ with $T=T(x^v)\in\Herm_r^\circ(F)^+$ and $t_1,\dots,t_r$ diagonal entries of $T$. For $1\leq i\leq r$, write $\cC_{t_i}(x^v_i)_L\in\rH^2_{\et,\cZ_{t_i}[x^v_i]_L}(\cX_L,\dQ_\ell(1))$ for the refined absolute $\ell$-adic cycle class of $\cZ_{t_i}(x^v_i)_L$, so that we have
\[
\cC_{t_1}(x^v_1)_L\cup\cdots\cup\cC_{t_r}(x^v_r)_L\in
\bigoplus_{\substack{T'\in\Herm_r^\circ(F)^+ \\ T'\asymp T}}\rH^{2r}_{\et,\cZ_{T'}[x^v]_L}(\cX_L,\dQ_\ell(r)),
\]
where $T'\asymp T$ means that they have the same diagonal entries. Define $\cC_T(x^v)_L$ to be its component in $\rH^{2r}_{\et,\cZ_T[x^v]_L}(\cX_L,\dQ_\ell(r))$. Similar to \eqref{eq:special_map_motive_0}, the map \eqref{eq:special_morphism_4} induces maps
\begin{align*}
\mu^\lambda_T(x^v)_L\colon\tS_{\lambda^+}((\dE^\lambda)^r)\otimes\varsigma_T(x^v)_L^!\mathbbm{1}^\et_{\cX_L,\dQ_\ell}
\to\varsigma_T(x^v)_L^!V_{\ell,L,\lambda},
\end{align*}
hence
\begin{align}\label{eq:special_map_etale}
\tS_{\lambda^+}((\dE^\lambda)^r)\otimes\rH^{2r}_{\et,\cZ_T[x^v]_L}(\cX_L,\dQ_\ell(r))
\to\rH^{2r}_{\et,\cZ_T[x^v]_L}(\cX_L,V_{\ell,L,\lambda}(r)).
\end{align}
Denote by
\[
\cC^\lambda_T(x^v)_L\colon\tS_{\lambda^+}((\dE^\lambda)^r)\to\rH^{2r}_{\et,\cZ_T[x^v]_L}(\cX_L,V_{\ell,L,\lambda}(r))
\]
the image of $\id\otimes\cC_T(x^v)_L$ under \eqref{eq:special_map_etale}.

For $i=1,2$, put
\[
\cZ_{T_i}[\ts_i\phi^\infty_i]_L\coloneqq\coprod_{\substack{x^v_i\in (V^{\infty,v})^r_{T_i}\\ (\ts_i\phi^{\infty,v}_i)(x^v_i)\neq0}}
\cZ_{T_i}[x^v_i]_L
\]
whose generic fiber coincides with $Z_{i,u}$, and put
\[
\cC^\lambda_{T_i}(\ts_i\phi^\infty_i)_L\coloneqq
\sum_{x^v_i\in L^v\backslash (V^{\infty,v})^r_{T_i}}(\ts_i\phi^{\infty,v}_i)(x^v_i)\cdot \cC^\lambda_{T_i}(x^v_i)_L
\in\tS_{\lambda^+}(\dL^r)^\vee\otimes_{\dE^\lambda}\rH^{2r}_{\et,\cZ_{T_i}[\ts_i\phi^\infty_i]_L}(\cX_L,V_{\ell,L,\lambda}(r))
\]
whose natural image in $\rH^{2r}_{\et,Z_{i,u}}(X_{L,u},V_{\ell,L,\lambda}(r))$ coincides with that of $\ts_i^*\bbZ^\lambda_{T_i}(\phi^\infty_i)_L$.

Since $\ts_1$ annihilates $\rH^{2r}_{\et}(\cX_L,V_{\ell,L,\lambda}(r))\otimes_{\dE^\lambda}\dL$, by Lemma \ref{le:height_integral}, we have
\[
\llangle\ts_1^*\bbZ^\lambda_{T_1}(\phi^\infty_1)_L,\ts_2^*\bbZ^\lambda_{T_2}(\phi^\infty_2)_L\rrangle_{\lambda,L,\ul{Z},u}^\ell
=\cC^\lambda_{T_1}(\ts_1\phi^\infty_1)_L\overset{\lambda}\bullet\cC^\lambda_{T_2}(\ts_2\phi^\infty_2)_L.
\]
Note that the intersection of $\cZ_{T_1}[\ts_1\phi^\infty_1]_L$ and $\cZ_{T_2}[\ts_2\phi^\infty_2]_L$ is contained in the basic locus of $\cX_L$.

Denote by $V^\star$ the unique totally positive definite hermitian space over $E/F$ such that $V^\star\otimes_F\bA_F^{\infty,v}$ is isometric to $V\otimes_{\bA_F}\bA_F^{\infty,v}$. Put $H^\star\coloneqq\rU(V^\star)$, and fix an isometry $V^\star\otimes_F\bA_F^{\infty,v}\simeq V\otimes_{\bA_F}\bA_F^{\infty,v}$ under which we identify $H(\bA_F^{\infty,v})$ with $H^\star(\bA_F^{\infty,v})$. For every $L^{\prime p}\in\Hk(H(\bA_F^{\infty,p}))$, denote by $\fX_{L'}$ the formal completion of $\cX_{L'}\otimes_{O_{E_u}}O_{\breve{E}_u}$ along its basic locus, where $\breve{E}_u$ denotes a complete maximal unramified extension of $E_u$. Then we have the uniformization isomorphism
\[
\fX_{L'}\simeq H^\star(F)\backslash\(\fN\times H(\bA_F^{\infty,v})/{L'}^v\)
\]
functorial in $L'$, where $\fN$ is a certain formal Rapoport--Zink space over $\Spf O_{\breve{E}_u}$ (see \cite{LZ22}*{\S13.1} for more details). The $\dE^\lambda_\ell$-linear local system $V_{\ell,L',\lambda}$ over $\cX_{L'}$ has a compatible uniformization:
\[
V_{\ell,L',\lambda}\res_{\fX_{L'}}=H^\star(F)\backslash\(\fN\times H(\bA_F^{\infty,v})/{L'}^v\times V_{\ell,\lambda}\),
\]
in which $H^\star(F)$ acts on $V_{\ell,\lambda}\coloneqq(V_\ell)_\lambda$ via the natural map $H^\star(F)\to H(F_\ell)$.

For $i=1,2$, if we denote by $\fZ_{T_i}[\ts_i\phi^\infty_i]_L$ the inverse image of $\cZ_{T_i}[\ts_i\phi^\infty_i]_L$ in $\fN\times H(\bA_F^{\infty,v})/L^v$, then $\cC^\lambda_{T_i}(\ts_i\phi^\infty_i)_L\res_{\fX_L}$ is represented by the $H^\star(F)$-invariant class
\begin{align*}
&\quad\sum_{x^\star_i\in(V^\star)^r_{T_i}}\(\sum_{h\in H(\bA_F^{\infty,v})/L^v}\ts_i\phi^{\infty,v}_i(h^{-1}x^\star_i)\cdot\fC(x^\star_{i,v})_h\)\otimes\mu^\lambda(x^\star_i) \\
&\in\rH^{2r}_{\et,\fZ_{T_i}[\ts_i\phi^\infty_i]_L}(\fN\times H(\bA_F^{\infty,v})/L^v,\dQ_\ell(r))
\otimes\(\tS_{\lambda^+}((\dE^\lambda)^r)^\vee\otimes_{\dE^\lambda}V_\lambda\) \\
&=\tS_{\lambda^+}((\dE^\lambda)^r)^\vee\otimes_{\dE^\lambda}\(\rH^{2r}_{\et,\fZ_{T_i}[\ts_i\phi^\infty_i]_L}(\fN\times H(\bA_F^{\infty,v})/L^v,\dQ_\ell(r))\otimes_{\dQ_\ell}V_{\ell,\lambda}\),
\end{align*}
where $\fC(x^\star_{i,v})_h\in\rH^{2r}_{\et,\fZ_{T_i}[x^v_i]_L}(\fN\times H(\bA_F^{\infty,v})/L^v,\dQ_\ell(r))$ denotes the $\ell$-adic cycle class of the derived special cycle of $x^\star_{i,v}$ on $\fN$ defined in \cite{LZ22}*{Definition~2.4.1}  but put in the component indexed by $h$; and $\mu^\lambda$ is defined in Construction \ref{co:harmonic_function}(3) (with $R=\dQ$). Combining this with \cite{LZ22}*{Theorem~10.5.1~\&~Remark~10.5.4} and using the argument for \cite{LL21}*{Proposition~8.1}, we obtain
\[
\cC^\lambda_{T_1}(\ts_1\phi^\infty_1)_L\overset{\lambda}\bullet\cC^\lambda_{T_2}(\ts_2\phi^\infty_2)_L
=\sum_{\substack{T^\Box\in\Herm_{2r}^\circ(F)^+\\ \partial T^\Box=(T_1,T_2)}}
\kappa\Psi^\lambda_\infty(T^\Box)\cdot\chi(T^\Box)\cdot\(\sum_{x\in L^v\backslash (V^{\infty,v})^{2r}_{T^\Box}}
\Phi^{\infty,v}(x)\).
\]
The proposition is proved.
\end{proof}

\begin{remark}\label{re:height_inert}
By the similar argument for \cite{LL21}*{Lemma~8.2}, Proposition \ref{pr:height_inert} holds as well if we replace $(\phi^\infty_1,\phi^\infty_2)$ by $(\omega_r^\infty(g_1)\phi^\infty_1,\omega_r^\infty(g_2)\phi^\infty_2)$ for every pair $(g_1,g_2)\in G_r(\bA_F^{\infty,\Diamondblack})^2$.
\end{remark}

\begin{remark}
Let $(\pi,\cA_\pi)$ be a holomorphic tempered cuspidal automorphic representation of $G_r(\bA_F)$ of weight $\lambda$ (Definition \ref{de:hermitian_representation}). Suppose that $\pi^\infty$ is defined over $\dL$ and that $(\pi_v)^{K_{r,v}}\neq\{0\}$ (Notation \ref{no:hermitian_siegel}(1)) for every $v\in\fV_F^\fin\setminus\fV_F^{(\Diamond)}$.
\begin{enumerate}
  \item Suppose that $V_v$ is split. By Lemma \ref{le:cohomological_annihilation} and the proper base change, there exist elements in $\dL[L^\Diamond\backslash H(\bA_F^{\infty,\Diamond})/L^\Diamond]\setminus\Ker\chi_\pi$ \eqref{eq:sakata_character} that annihilate $\rH^{2r}_{\et}(\cX_L,V_{\ell,L,\lambda}(r))\otimes_{\dE^\lambda}\dL$ and support only at places split in $E$.

  \item Suppose that $V_v$ is nonsplit. By the similar argument for \cite{LL21}*{Lemma~9.2}, there exist elements in $\dL[L^\Diamond\backslash H(\bA_F^{\infty,\Diamond})/L^\Diamond]\setminus\Ker\chi_\pi$ that annihilate $\rH^{2r}_{\et}(\cX_L,V_{\ell,L,\lambda}(r))\otimes_{\dE^\lambda}\dL$ and support only at places split in $E$.
\end{enumerate}
\end{remark}

\subsection{Vanishing local heights at split places}
\label{ss:height_split}

In this section, we take a nonarchimedean place $u$ of $E$ such that $v\coloneqq v(u)$ is split in $E$.

Consider an object $L\in\Hk(H(\bA_F^\infty))$ such that
\begin{itemize}
  \item $L_v=\Ker\(\rU(\Lambda_v)\to\rU(\Lambda_v/v^\delta)\)$ for some self-dual lattice $\Lambda_v$ of $V_v$ and some integer $\delta\geq 0$.

  \item when $r>1$, $L_{v'}$ is the stabilizer of a vertex lattice (Definition \ref{de:lattice}) for every place $v'\in\fV_F^{(p)}$ that is nonsplit in $E$, where $p\coloneqq p(u)$.
\end{itemize}

Let
\[
\cX_\bullet\colon \Hk(H(\bA_F^{\infty,p}))\to\Sch^\proj_{/O_{E_u}}
\]
be the functor constructed in \cite{LZ22}*{\S11.5},\footnote{We have a remark similar to Remark \ref{re:moduli_integral}. Also, when $r=1$, the restriction on $L_{v'}$ for $v'\in\fV_F^{(p)}$ nonsplit in $E$ is not needed.} satisfying Assumption \ref{as:integral_model}.

\begin{proposition}\label{pr:height_split}
Let $L\in\Hk(H(\bA_F^\infty))$ be as above. Consider
\begin{itemize}
  \item an intermediate number field $\dE^\lambda\subseteq\dL\subseteq\dC$,

  \item a finite set $\Diamond$ of rational primes containing $p$ and a nonempty subset $\Diamondblack\subseteq\fV_F^{(\Diamond)}\setminus\{v\}$,

  \item a pair $\ul\phi^\infty=(\phi^\infty_1,\phi^\infty_2)$ of decomposable Schwartz functions in $\sS((V^\infty)^r,\dL)^L$ satisfying that $(\phi^\infty_1\otimes\ol{\phi^\infty_2})_{v'}$ is regularly supported (Definition \ref{de:regular_support}) for $v'\in\Diamondblack$,

  \item a rational prime $\ell\neq p$, and

  \item a pair $(\ts_1,\ts_2)$ of elements in $\dL[L^\Diamond\backslash H(\bA_F^{\infty,\Diamond})/L^\Diamond]$ that annihilate both $\rH^{2r}_{\et}(X_{L,u},V_{\ell,L,\lambda}(r))\otimes_{\dE^\lambda}\dL$ and  $\rH^{2r}_{\et}(\cX_L,V_{\ell,L,\lambda}(r))\otimes_{\dE^\lambda}\dL$.
\end{itemize}
Then for every pair $(T_1,T_2)$ of elements in $\Herm_r^\circ(F)^+$ (so that $\ul{Z}=(Z_1\coloneqq Z_{T_1}[\ts_1\phi^\infty_1]_L,Z_2\coloneqq Z_{T_2}[\ts_2\phi^\infty_2]_L)$ (Definition \ref{de:special_cycle}) is a disjoint pair),
\[
\llangle\ts_1^*\bbZ^\lambda_{T_1}(\phi^\infty_1)_L,\ts_2^*\bbZ^\lambda_{T_2}(\phi^\infty_2)_L\rrangle_{\lambda,L,\ul{Z},u}^\ell=0
\]
holds.
\end{proposition}

\begin{proof}
Denote by $\cZ_i$ the Zariski closure of $Z_{i,u}$ in $\cX_L$ for $i=1,2$. Then $\ts_1\cZ_1$ and $\ts_2\cZ_2$ have empty intersection by (the proof of) \cite{LL21}*{Lemma~7.2}. For $i=1,2$, we may find a class
\[
\cC_i\in\rH^{2r}_{\et,\cZ_i}(\cX_L,V_{\ell,L,\lambda}(r))\otimes_{\dE^\lambda}\dL
\]
whose image in $\rH^{2r}_{\et,Z_{i,u}}(X_{L,u},V_{\ell,L,\lambda}(r))\otimes_{\dE^\lambda}\dL$ coincides with that of $\bbZ^\lambda_{T_i}(\phi^\infty_i)_L$. Since the image of $\ts_1^*\cC_1$ in $\rH^{2r}_{\et}(\cX_L,V_{\ell,L,\lambda}(r))\otimes_{\dE^\lambda}\dL$ vanishes, the proposition follows from (the $\dL$-sesquilinear extension of) Lemma \ref{le:height_integral}.
\end{proof}

\begin{remark}
Let $(\pi,\cA_\pi)$ be a holomorphic tempered cuspidal automorphic representation of $G_r(\bA_F)$ of weight $\lambda$ (Definition \ref{de:hermitian_representation}). Suppose that $\pi^\infty$ is defined over $\dL$ and $(\pi_v)^{K_{r,v}}\neq\{0\}$ (Notation \ref{no:hermitian_siegel}(1)) for every $v\in\fV_F^\fin\setminus\fV_F^{(\Diamond)}$. By the similar argument for \cite{LL22}*{Theorem~4.21}, there exist elements in $\dL[L^\Diamond\backslash H(\bA_F^{\infty,\Diamond})/L^\Diamond]\setminus\Ker\chi_\pi$ \eqref{eq:sakata_character} that annihilate $\rH^{2r}_{\et}(\cX_L,V_{\ell,L,\lambda}(r))\otimes_{\dE^\lambda}\dL$.
\end{remark}

\subsection{Nonarchimedean local heights over unitary moduli curves}
\label{ss:height_curve}

In this section, we assume $r=1$ and take a nonarchimedean place $u$ of $E$ such that $v\coloneqq v(u)$ is nonsplit in $E$. Put $p\coloneqq p(u)$.

For every intermediate field $\dE^\lambda\subseteq\dL\subseteq\dC$, we canonically identify both $\tS_{\lambda^+}(\dL)$ and $\ol{\tS_{\lambda^+}(\dL)}$ with $\dL$ as $\dL$-vector spaces.

\begin{definition}\label{de:schwartz_polynomial}
For every $\sfv\in V_\lambda\otimes_{\dE^\lambda}\dC$, we have a polynomial function $P[\sfv]$ on $V_\infty$ given by the formula
\[
P[\sfv](x)=(\mu^\lambda(x),\sfv)_{V_\lambda}\in\dC,\qquad x\in V_\infty,
\]
where $\mu^\lambda$ is defined in Construction \ref{co:harmonic_function}(3) (with $R=\dR$).
\end{definition}

\begin{definition}\label{de:eisenstein_ideal}
For a field $\dL$ containing $\dE^\lambda$, a finite subset $\Diamondblack$ of $\fV_F^\fin$, and an object $L\in\Hk(H(\bA_F^\infty))$, we define the \emph{special Eisenstein ideal} of $\dL[L^\Diamondblack\backslash H(\bA_F^{\infty,\Diamondblack})/L^\Diamondblack]$ with respect to the weight $\lambda$, denoted by $\dL[L^\Diamondblack\backslash H(\bA_F^{\infty,\Diamondblack})/L^\Diamondblack]^\lambda$, to be
\begin{itemize}
  \item when $\lambda=\b{0}$, the (two-sided) ideal consisting of elements that are of degree zero and supported on double cosets with unit determinants;

  \item when $\lambda\neq\b{0}$, the entire ring $\dL[L^\Diamondblack\backslash H(\bA_F^{\infty,\Diamondblack})/L^\Diamondblack]$.
\end{itemize}
In both cases, elements in $\dL[L^\Diamondblack\backslash H(\bA_F^{\infty,\Diamondblack})/L^\Diamondblack]^\lambda$ annihilate $\rH^2_{\et}(X_L\otimes_EE^\ac,V_{\ell,L,\lambda}(1))\otimes_{\dE^\lambda}\dL$.
\end{definition}

Denote by $V(v)$ the unique totally positive definite hermitian space over $E/F$ such that $V(v)\otimes_F\bA_F^{\infty,v}$ is isometric to $V\otimes_{\bA_F}\bA_F^{\infty,v}$. Put $H(v)\coloneqq\rU(V(v))$, and fix an isometry $V(v)\otimes_F\bA_F^{\infty,v}\simeq V\otimes_{\bA_F}\bA_F^{\infty,v}$ under which we identify $H(\bA_F^{\infty,v})$ with $H(v)(\bA_F^{\infty,v})$.

\begin{proposition}\label{pr:height_curve}
Consider
\begin{itemize}
  \item an object $L$ of $\Hk(H(\bA_F^\infty))$,

  \item an intermediate number field $\dE^\lambda\subseteq\dL\subseteq\dC$,

  \item a finite set $\Diamond$ of rational primes containing $p$ and a nonempty subset $\Diamondblack\subseteq\fV_F^{(\Diamond)}\setminus\{v\}$,

  \item a pair $\ul\phi^\infty=(\phi^\infty_1,\phi^\infty_2)$ of decomposable Schwartz functions in $\sS(V^\infty,\dL)^L$ satisfying that $(\phi^\infty_1\otimes\ol{\phi^\infty_2})_{v'}$ is regularly supported (Definition \ref{de:regular_support}) for $v'\in\Diamondblack\cup\{v\}$; and $\phi^\infty_{i,v}$ is regularly supported for $i=1,2$.\footnote{When $V_v$ is nonsplit, that $(\phi^\infty_1\otimes\ol{\phi^\infty_2})_v$ is regularly supported implies that $\phi^\infty_{i,v}$ is regularly supported for $i=1,2$.}
\end{itemize}
\begin{enumerate}
  \item Suppose that $V_v$ is split. There exists a Schwartz function $\Phi^\star_v\in\sS(V(v)_v^2,\dL)$ such that for every pair $(\ts_1,\ts_2)$ in $\dL[L^\Diamond\backslash H(\bA_F^{\infty,\Diamond})/L^\Diamond]^\lambda$ and every pair $(t_1,t_2)$ in $(F^\times)^+$ (so that $\ul{Z}=(Z_1\coloneqq Z_{t_1}[\ts_1\phi^\infty_1]_L,Z_2\coloneqq Z_{t_2}[\ts_2\phi^\infty_2]_L)$ (Definition \ref{de:special_cycle}) is a disjoint pair),
      \begin{align*}
      \vol^\natural(L)\cdot
      \llangle\ts_1^*\bbZ^\lambda_{t_1}(\phi^\infty_1)_L(g_1),
      \ts_2^*\bbZ^\lambda_{t_2}(\phi^\infty_2)_L(g_2)\rrangle_{\lambda,L,\ul{Z},u}^\ell
      =E_{t_1,t_2}((g_1,g_2),\Phi^\lambda_\infty\cdot(\ts_1\phi^\infty_1\otimes\ol{\ts_2\phi^\infty_2})^v\cdot\Phi^\star_v)
      \end{align*}
      holds (in $\dC\otimes\dQ_\ell$) for every pair $(g_1,g_2)\in G_1(F_\infty)^2$ and every prime $\ell\neq p$, where $E_{t_1,t_2}$ is introduced in Definition \ref{de:analytic_kernel}(1).

  \item Suppose that $V_v$ is nonsplit. Consider further a pair $(\ts_1,\ts_2)$ in $\dL[L^\Diamond\backslash H(\bA_F^{\infty,\Diamond})/L^\Diamond]^\lambda$. For every element $t_2\in(F^\times)^+$, there exists a $V_\lambda\otimes_{\dE^\lambda}\dL$-valued Schwartz function $\phi^\infty_{t_2}$ on $V(v)^\infty$ (depending on the data so far) such that
      \begin{align*}
      \vol^\natural(L)\cdot
      \llangle\ts_1^*\bbZ^\lambda_{t_1}(\phi^\infty_1)_L(g_1),
      \ts_2^*\bbZ^\lambda_{t_2}(\phi^\infty_2)_L\rrangle_{\lambda,L,\ul{Z},u}^\ell
      =\sum_{x^\star_1\in V(v)_{t_1}}\omega_{1,\infty}(g_1)\(P[\phi^\infty_{t_2}(x^\star_1)]\cdot\phi^0_{V_\infty}\)(x^\star_1)
      \end{align*}
      holds (in $\dC\otimes\dQ_\ell$) for every $t_1\in(F^\times)^+$, every element $g_1\in G_1(F_\infty)$ and every prime $\ell\neq p$.
\end{enumerate}
\end{proposition}

We will prove the above proposition in the next two sections (as Proposition \ref{pr:height_curve_split} and Proposition \ref{pr:height_curve_nonsplit}), according to whether $V_v$ is split or nonsplit.

\begin{remark}
The conclusion in Proposition \ref{pr:height_curve}(1) is stronger than the conclusion in Proposition \ref{pr:height_curve}(2). The reason we are not able to obtain the similar statement as in (1) when $V_v$ is nonsplit is that, in this case, it is difficult to control flat extensions of special classes in the integral model $\cX_L$ we use below.
\end{remark}

In the rest of this section, we make some common preparation for both cases. Fix a complete maximal unramified extension $\breve{E}_u$ of $E_u$, with $\kappa$ its residue field. Put $R_0\coloneqq O_{E_u}$ and $R\coloneqq O_{\breve{E}_u}$ for short.

By the projection formula for local heights, we may assume, without lost of generality, that the following holds.

\begin{assumption}\label{as:height_curve}
There exist integers $a<b$ such that
\begin{itemize}
  \item $L_p$ is contained in an optimal open compact subgroup $L_p^\flat$ (Definition \ref{de:principal_level} with $\Diamond=\{p\}$) such that $L_v$ is the subgroup of $L_v^\flat$ consisting of elements that act trivially on $\Lambda_v/\fp_v^{b-a}\Lambda_v$, where $\Lambda_v$ is the lattice of $V_v$ that $L_v^\flat$ stabilizes;

  \item $(\phi^\infty_i)_v$ is a function on $\fp_v^a\Lambda_v/\fp_v^b\Lambda_v$ (that is, supported on $\fp_v^a\Lambda_v$ and invariant under $\fp_v^b\Lambda_v$) for $i=1,2$.
\end{itemize}
\end{assumption}

Applying Definition \ref{de:moduli_integral_regular} to $L_p\subseteq L^\flat_p$, we obtain the functor
\begin{align*}
X_{\bullet,\bullet,(p)}\colon\Hk(H(\bA_F^{\infty,p}))\times\Hk(T(\bA^{\infty,p}))
\to\Sch^\proj_{/O_{E,(p)}},
\end{align*}
and put $\cX_{L',U}\coloneqq X_{L^{\prime p},U^p,(p)}\otimes_{O_{E,(p)}}R_0$ for every $L'\in\Hk(H(\bA_F^\infty))$ with $L'_p=L_p$ and $U_p$ maximal.

Choose an object $U^p$ of $\Hk(T(\bA^{\infty,p}))$ and a framing $P^\star$ as in Remark \ref{re:uniformization}. Put $\cX_{L'}\coloneqq\cX_{L',U}/\Xi_{U^p}$ (which is canonically independent of the choice of $U^p$). Denote by $\fX_{L'}$ the formal completion of $\cX_{L'}\otimes_{R_0}R$ along the basic locus. Then \eqref{eq:uniformization_3} implies the following isomorphism of basic uniformization
\begin{align}\label{eq:uniformization_4}
\fX_{L'}\simeq H^\star(F)\backslash\(\fN^{b-a}\times H(\bA_F^{\infty,v})/L^{\prime v}\)
\end{align}
functorial in $L^{\prime\Diamond}\in\Hk(H(\bA_F^{\infty,p}))$, where $\fN^{b-a}$ is described in Definition \ref{de:rapoport_zink} below. Note that when $V_v$ is nonsplit, the entire special fiber $\cX_{L'}\otimes_{R_0}\kappa$ is basic.

\begin{definition}\label{de:rapoport_zink}
Put $(X^\star,\lambda^\star)\coloneqq(A^\star[\fp_v^\infty],\lambda^\star[\fp_v^\infty])$ by abuse of notation on $\lambda^\star$ and similarly $(X_0^\star,\lambda_0^\star)\coloneqq(A_0^\star[p^\infty],\lambda_0^\star[p^\infty])$. Over $\fN^\flat$ (see Remark \ref{re:uniformization} (U3)), we have the universal quasi-isogeny
\[
\rho\colon (X,\lambda)\otimes_{R}\kappa \dashrightarrow (X^\star,\lambda^\star)\otimes_R(\fN^\flat\otimes_R\kappa)
\]
of polarized $O_{F_v}$-divisible groups with $R_0$-actions. Write $\fN^\flat_\eta$ for the generic fiber of $\fN^\flat$ as a smooth rigid analytic space over $R[p^{-1}]$ of dimension one.

For every integer $\delta\geq 0$, denote by $\fN^\delta_\eta\to\fN^\flat_\eta$ the finite \'{e}tale cover classifying isometries
\[
\eta_v\colon\Lambda_v/\fp_v^\delta\Lambda_v\xrightarrow\sim\Hom_{O_{E_0,p}}\(X^\star_0[\fp_v^\delta](S), X[\fp_v^\delta](S)\)
\]
over any rigid analytic space $S$ over $\fN^\flat_\eta$ of finite type, where the hermitian form on the latter module is induced from the polarizations $\lambda^\star_0$ and $\lambda$, and define $\fN^\delta$ to be the regular closure of $\fN^\flat$ in $\fN^\delta_\eta$ (Definition \ref{de:regular_closure}). Then $\fN^\delta$ is a regular formal scheme over $R$ of relative dimension one.
\end{definition}

\begin{definition}\label{de:special_cycle_formal}
Take a nonzero element $x^\star\in V^\star_v$ and an integer $\delta\geq 0$.
\begin{enumerate}
  \item We denote by $\fZ_a^\flat(x^\star)$ the maximal closed subscheme of $\fN^\flat$ over which the composite quasi-map
      \[
      \varpi^{-a}\cdot\rho^{-1}\circ x^\star\colon X^\star_0\otimes_R\kappa \to A^\star[\fp_v^\infty]\otimes_R\kappa \to  X\times_{\fN^\flat}\Spec\kappa
      \]
      lifts to a map of $O_{F_v}$-divisible groups for some hence every uniformizer $\varpi$ of $F_v$, where the last base change is with respect to the closed point $\Spec\kappa\to\fN^\flat$ given by the framing $P^\star$. Write $\fZ_a^\delta(x^\star)$ for the inverse image of $\fZ_a^\flat(x^\star)$ in $\fN^\delta$.

  \item For every element $x\in\fp_v^a\Lambda_v/\fp_v^{a+\delta}\Lambda_v$, define $\fZ_a^\delta(x^\star,x)$ to be the maximal vertical divisor contained in $\fZ_a^\delta(x^\star)$ such that over the generic fiber,
      \[
      (\varpi^{-a}\cdot\rho^{-1}\circ x^\star)[\fp_v^\delta]=\eta_v(\varpi^{-a}\cdot x)
      \]
      holds for some hence every uniformizer $\varpi$ of $F_v$.

  \item For every $\dL$-valued function $\phi$ on $\fp_v^a\Lambda_v/\fp_v^{a+\delta}\Lambda_v$, put
      \[
      \fZ_a^\delta(x^\star,\phi)\coloneqq\sum_{x\in\fp_v^a\Lambda_v/\fp_v^{a+\delta}\Lambda_v}\phi(x)\cdot\fZ_a^\delta(x^\star,x)
      \in\Div(\fN^\delta)\otimes\dL.
      \]
\end{enumerate}
\end{definition}

\begin{lem}\label{le:special_divisor_formal}
Consider an integer $\delta\geq0$ and an element $x^\star\in V^\star_v$ with $T(x^\star)\neq 0$. For every integer $N\geq 1$, there exists an open neighbourhood $U_N$ of $x^\star$ in $V^\star_v$, such that for every $y^\star\in U_N$, $\fZ_a^\delta(x^\star,x)$ and $\fZ_a^\delta(y^\star,x)$ have the same image in $\Div(\fN^\delta\otimes_RR/p^N)$ for every $x\in\fp_v^a\Lambda_v/\fp_v^{a+\delta}\Lambda_v$, where $\Div(\fN^\delta\otimes_RR/p^N)$ denotes the abelian group of Cartier divisors on $\fN^\delta\otimes_RR/p^N$.
\end{lem}

\begin{proof}
Without lost of generality, we may assume $a=0$ and suppress it from the subscripts.

We first consider the case where $\delta=0$. Denote by $\rH$ the (topological) subgroup of $\GU(V^\star_v)(F_v)$ consisting of elements whose similitudes are in $O_{F_v}^\times$. Then $\rH$ acts on $(A^\star[p^\infty],\lambda^\star[p^\infty])$ by quasi-automorphisms of height zero, which induces a continuous action on $\fN^\flat$ by composition hence on $\fN^0$ by the functoriality of taking regular closure. Take an element $x^\star\in V^\star_v$ with $T(x^\star)\neq 0$. Then for $h\in\rH$, $h\fZ^0(x^\star)=\fZ^0(h x^\star)$ by Definition \ref{de:special_cycle_formal}(1). For every $N\geq 1$, there exists an open compact subgroup $\rI_N$ of $\rH$ such that the image of $h\fZ^0(x^\star)$ in $\Div(\fN^\delta\otimes_RR/p^N)$ does not depend on $h\in\rI_N$. Since $\{\rI x^\star\res \rI\text{ open compact subgroups of }\rH\}$ is a topological basis of $x^\star$, the lemma follows in this case.

For general $\delta$, the action of $\rH$ on $\fN^0_\eta$ lifts naturally to a (continuous) action on $\fN^\delta_\eta$. By the functoriality of taking regular closure, we obtain a continuous action of $\rH$ on $\fN^\delta$ lifting that on $\fN^0$. Take an element $x^\star\in V^\star_v$ with $T(x^\star)\neq 0$. Then for $h\in\rH$ and every $x\in\Lambda_v/\fp_v^\delta\Lambda_v$, $h\fZ^\delta(x^\star,x)=\fZ^\delta(h x^\star,x)$ by Definition \ref{de:special_cycle_formal}(2). The lemma then follows in general by the same reason as above.
\end{proof}

\subsection{Proof of Proposition \ref{pr:height_curve}(1)}
\label{ss:height_curve_split}

In this section, we assume $V_v$ split. We keep the notation from \S\ref{ss:height_curve} and take $\delta\coloneqq b-a$.

\begin{lem}\label{le:height_curve_split_1}
For every nonzero element $x^\star\in V^\star_v$ and every element $x\in\fp_v^a\Lambda_v/\fp_v^b\Lambda_v$, there is a unique element $\fZ_a^\delta(x^\star,x)^\dag\in\Div(\fN^\delta)\otimes\dQ$ satisfying that
\begin{enumerate}[label=(\alph*)]
  \item $\fZ_a^\delta(x^\star,x)^\dag-\fZ_a^\delta(x^\star,x)$ is a ($\dQ$-)linear combination of irreducible components of the induced reduced subscheme of $\fN^\delta$;

  \item for every irreducible component $C$ of the induced reduced subscheme of $\fN^\delta$, $\fZ_a^\delta(x^\star,x)^\dag\res_C$ has degree zero.
\end{enumerate}
\end{lem}

\begin{proof}
By Lemma \ref{le:vanishing_ordinary}(1), $(\fN^\flat)^{\r{fl}}$ is formally smooth over $R$ (of relative dimension one). Since $\fN^\delta$ is the canonical resolution of (isolated) singularities of the normalization of $(\fN^\flat)^{\r{fl}}$ in $\fN^\delta_\eta$, it is well-known that the intersection matrix of the reduced fiber of $\fN^\delta$ is negative definite. The lemma then follows.
\end{proof}

For every $\dL$-valued function $\phi$ on $\fp_v^a\Lambda_v/\fp_v^b\Lambda_v$, put
\[
\fZ_a^\delta(x^\star,\phi)^\dag\coloneqq\sum_{x\in\fp_v^a\Lambda_v/\fp_v^b\Lambda_v}\phi(x)\cdot\fZ_a^\delta(x^\star,x)^\dag
\in\Div(\fN^\delta)\otimes\dL.
\]

\begin{lem}\label{le:height_curve_split_2}
There exists a unique $\dL$-valued Schwartz function $\Phi^{[a,b]}_v$ on $(V^\star_v)^2$ such that for every pair $(x^\star_1,x^\star_2)\in(V^\star_v\setminus\{0\})^2$,
\begin{align*}
\Phi^{[a,b]}_v(x^\star_1,x^\star_2)
&=
\fZ_a^\delta(x^\star_1,\phi^\infty_{1,v})^\dag\bullet\fZ_a^\delta(x^\star_2,\phi^\infty_{2,v})^\dag \\
&=\sum_{x_1\in\fp_v^a\Lambda_v/\fp_v^b\Lambda_v}\sum_{x_2\in\fp_v^a\Lambda_v/\fp_v^b\Lambda_v}
\phi^\infty_{1,v}(x_1)\ol{\phi^\infty_{2,v}(x_2)}\cdot
\fZ_a^\delta(x^\star_1,x_1)^\dag\bullet\fZ_a^\delta(x^\star_2,x_2)^\dag,
\end{align*}
in which $\bullet$ denotes the intersection pairing on $\fN^\delta$ as an $\dL$-hermitian form. Moreover, $\Phi^{[a,b]}_v$ is invariant under the action of $H^\star(F_v)$.
\end{lem}

Here, note that since $(\phi^\infty_1\otimes\ol{\phi^\infty_2})_v$ is regularly supported, the intersection $\fZ_a^\delta(x^\star_1,x_1)\cap\fZ_a^\delta(x^\star_2,x_2)$ is finite if $\phi^\infty_{1,v}(x_1)\ol{\phi^\infty_{2,v}(x_2)}\neq 0$.

\begin{proof}
It is clear that if $\Phi^{[a,b]}_v$ exists, then it is unique, $\dL$-valued, invariant under $H^\star(F_v)$. The continuity of $\Phi^{[a,b]}_v$ on $(V^\star_v\setminus\{0\})^2$ follows from Lemma \ref{le:special_divisor_formal}.

The lemma follows if we can show that $\Phi^{[a,b]}_v$ is compactly supported in $(V^\star_v\setminus\{0\})^2$. Indeed, $\fZ_a^\delta(x^\star,x)$ hence $\fZ_a^\delta(x^\star,x)^\dag$ vanish unless $(x^\star,x^\star)\equiv(x,x)\mod\fp_v^{2b}$. Now as both $\phi^\infty_{1,v}$ and $\phi^\infty_{2,v}$ are regularly supported, the support of $\Phi^{[a,b]}_v$ is contained in the subset $D_{(2b,2a]}^2$, where
\[
D_{(2b,2a]}\coloneqq\{x^\star\in V^\star_v\res (x^\star,x^\star)\in\fp_v^{2a}\setminus\fp_v^{2b}\}.
\]
Since $V^\star_v$ is anisotropic, $D_{(2b,2a]}$ is compact in $V^\star_v\setminus\{0\}$.

The lemma is proved.
\end{proof}

\begin{proposition}\label{pr:height_curve_split}
Suppose that we are in the situation of Proposition \ref{pr:height_curve}(1) together with Assumption \ref{as:height_curve}.
If we put
\[
\Phi^\star_v\coloneqq\frac{\vol(L_v)}{\vol(H^\star(F_v))}\Phi^{[a,b]}_v,
\]
then for every pair $(\ts_1,\ts_2)$ in $\dL[L^\Diamond\backslash H(\bA_F^{\infty,\Diamond})/L^\Diamond]^\lambda$ and every pair $(t_1,t_2)$ in $(F^\times)^+$,
\begin{align*}
\vol^\natural(L)\cdot
\llangle\ts_1^*\bbZ^\lambda_{t_1}(\phi^\infty_1)_L(g_1),
\ts_2^*\bbZ^\lambda_{t_2}(\phi^\infty_2)_L(g_2)\rrangle_{\lambda,L,\ul{Z},u}^\ell
=E_{t_1,t_2}((g_1,g_2),\Phi^\lambda_\infty\cdot(\ts_1\phi^\infty_1\otimes\ol{\ts_2\phi^\infty_2})^v\cdot\Phi^\star_v)
\end{align*}
holds (in $\dC\otimes\dQ_\ell$) for every pair $(g_1,g_2)\in G_1(F_\infty)^2$ and every prime $\ell\neq p$.
\end{proposition}

Take a rational prime $\ell\neq p$. Note that the $\dE^\lambda_\ell$-linear local system $V_{\ell,L',\lambda}$ over $\cX_{L'}$ has a uniformization:
\[
V_{\ell,L',\lambda}\res_{\fX_{L'}}=H^\star(F)\backslash\(\fN^\delta\times H(\bA_F^{\infty,v})/L^{\prime v}\times V_{\ell,\lambda}\)
\]
compatible with \eqref{eq:uniformization_4}, in which $H^\star(F)$ acts on $V_{\ell,\lambda}\coloneqq(V_\ell)_\lambda$ via the natural map $H^\star(F)\to H(F_\ell)$.

\begin{construction}\label{co:special_divisor_split}
Consider a Schwartz function $\phi^\infty\in\sS(V^\infty,\dL)$ of the form $\phi_v\otimes\phi^{\infty,v}$ in which $\phi_v$ is a regularly supported $\dL$-valued function on $\fp_v^a\Lambda_v/\fp_v^b\Lambda_v$.
\begin{enumerate}
  \item For $x^\star_v\in V^\star_v\setminus\{0\}$, put
      \[
      \fZ(x^\star_v,\phi_v)_h\coloneqq
      \fZ_a^\delta(x^\star_v,\phi_v)^\dag\times\{h\}\in\Div(\fN^\delta\times H(\bA_F^{\infty,v})/L^v)\otimes\dL,
      \]
      and denote by
      \[
      \fC(x^\star_v,\phi_v)_h\in\rH^2_{\et,|\fZ(x^\star_v,\phi_v)_h|}(\fN^\delta\times H(\bA_F^{\infty,v})/L^v,\dL_\ell(1))
      \]
      its refined $\ell$-adic cycle class.

  \item For every $t\in(F^\times)^+$, define $\cZ_t(\phi^\infty)_L$ to be the unique element in $\Div(\cX_L\otimes_{R_0}R)\otimes\dL$ supported on the basic locus whose pullback to $\fN^\delta\times H(\bA_F^{\infty,v})/L^v$ is given by
      \[
      \fZ_t(\phi^\infty)_L\coloneqq\sum_{x^\star\in V^\star_t}\sum_{h\in H(\bA_F^{\infty,v})/L^v}
      \phi^{\infty,v}(h^{-1}x^\star)\cdot\fZ(x^\star_v,\phi_v)_h.
      \]

  \item For every $t\in(F^\times)^+$, denote by
      \[
      \cC^\lambda_t(\phi^\infty)_L\in
      \rH^2_{\et,|\cZ_t(\phi^\infty)_L|}(\cX_L\otimes_{R_0}R,V_{\ell,L,\lambda}(1))\otimes_{\dE^\lambda}\dL
      \]
      the unique element whose restriction to $\fN^\delta\times H(\bA_F^{\infty,v})/L^v$ is given by the $H^\star(F)$-invariant class
      \begin{align*}
      &\quad\sum_{x^\star\in V^\star_t}\(\sum_{h\in H(\bA_F^{\infty,v})/L^v}
      \phi^{\infty,v}(h^{-1}x^\star)\cdot\fC(x^\star_v,\phi_v)_h\)\otimes\mu^\lambda(x^\star) \\
      &\in\rH^2_{\et,|\fZ_t(\phi^\infty)_L|}(\fN^\delta\times H(\bA_F^{\infty,v})/L^v,\dL_\ell(1))
      \otimes_{\dE^\lambda}V^\star_\lambda
      =\rH^2_{\et,|\fZ_t(\phi^\infty)_L|}(\fN^\delta\times H(\bA_F^{\infty,v})/L^v,\dL_\ell(1))
      \otimes_{\dE^\lambda_\ell}V_{\ell,\lambda},
      \end{align*}
      where $\mu^\lambda$ is defined in Construction \ref{co:harmonic_function}(3) (with $R=\dQ$).
\end{enumerate}
\end{construction}

\begin{lem}\label{le:height_curve_split_3}
Let $\phi^\infty$ be as in Construction \ref{co:special_divisor_split}. For every $\ts\in\dL[L^\Diamond\backslash H(\bA_F^{\infty,\Diamond})/L^\Diamond]^\lambda$ (Definition \ref{de:eisenstein_ideal}) and every $t\in (F^\times)^+$, the image of $\cC^\lambda_t(\ts\phi^\infty)_L$ in $\rH^2_{\et}(\cX_L\otimes_{R_0}R,V_{\ell,L,\lambda}(1))\otimes_{\dE^\lambda}\dL$ vanishes.
\end{lem}

\begin{proof}
It suffices to show that for every (reduced) irreducible component $C$ of $\cX_L\otimes_{R_0}\kappa$, the image in $\rH^2_{\et}(C,V_{\ell,L,\lambda}(1))\otimes_{\dE^\lambda}\dL$ vanishes. There are two cases.

First suppose that $C$ is contained in the basic locus. Then the vanishing follows from the property (b) of Lemma \ref{le:height_curve_split_1}.

Second suppose that $C$ is not contained in the basic locus. If $\lambda=\b{0}$, then $\ts$ preserves and annihilates $\rH^2_{\et}(C,\dL_\ell(1))$ by Lemma \ref{le:vanishing_ordinary}(1) (with $\Diamond=\{p\}$). If $\lambda\neq\b{0}$, we have $\rH^2_{\et}(C,V_{\ell,L,\lambda}(1))=0$ by Lemma \ref{le:vanishing_ordinary}(2) (with $\Diamond=\{p\}$).

The lemma is proved.
\end{proof}

\begin{proof}[Proof of Proposition \ref{pr:height_curve_split}]
Similar to the proof of Proposition \ref{pr:height_inert}, we may assume $g_1=g_2=1_2$.

By construction, for $i=1,2$, the natural image of the class $\cC^\lambda_{t_i}(\ts_i\phi^\infty_i)_L$ from Construction \ref{co:special_divisor_split}(3) in $\rH^2_{\et,Z_{i,u}\otimes_{R_0}R}(X_{L,u}\otimes_{R_0}R,V_{\ell,L,\lambda}(1))$ coincides with that of $\ts_i^*\bbZ^\lambda_{t_i}(\phi^\infty_i)_L$. Then by Lemma \ref{le:height_integral} (or rather its analogue after base change to $R$) and Lemma \ref{le:height_curve_split_3}, we have
\[
\llangle\ts_1^*\bbZ^\lambda_{t_1}(\phi^\infty_1)_L,\ts_2^*\bbZ^\lambda_{t_2}(\phi^\infty_2)_L\rrangle_{\lambda,L,\ul{Z},u}^\ell
=\cC^\lambda_{t_1}(\ts_1\phi^\infty_1)_L\overset{\lambda}\bullet\cC^\lambda_{t_2}(\ts_2\phi^\infty_2)_L.
\]
Applying the argument for \cite{LL21}*{Proposition~8.1} to Construction \ref{co:special_divisor_split}, we obtain
\[
\cC^\lambda_{t_1}(\ts_1\phi^\infty_1)_L\overset{\lambda}\bullet\cC^\lambda_{t_2}(\ts_2\phi^\infty_2)_L
=\sum_{\substack{T^\Box\in\Herm_2^\circ(F)^+\\ \partial T^\Box=(t_1,t_2)}}
\kappa\Psi^\lambda_\infty(T^\Box)\cdot\Phi^{[a,b]}_v(T^\Box)\cdot\(\sum_{x\in L^v\backslash (V^{\infty,v})^2_{T^\Box}}
\Phi^{\infty,v}(x)\).
\]
Here, we note that $\Phi^{[a,b]}_v$ factors through the moment map $T\colon(V^\star_v)^2\to\Herm_2(F_v)$ by Lemma \ref{le:height_curve_split_2}. Then the proposition follows from the local Siegel--Weil formula.
\end{proof}

\subsection{Proof of Proposition \ref{pr:height_curve}(2)}
\label{ss:height_curve_nonsplit}

In this section, we assume $V_v$ nonsplit. We keep the notation from \S\ref{ss:height_curve} and take $\delta\coloneqq b-a$.

Fix an element $x^\star_2\in V^\star$ with $T(x^\star_2)=t_2$, and denote by $(H^\star)^{x^\star_2}$ the stabilizer of $x^\star_2$ in $H^\star$. Consider the set
\[
S(x^\star_2)\coloneqq\left\{h\in H(\bA_F^{\infty,v})\left| \ts_2\phi^{\infty,v}_2(h^{-1}x^\star_2)\neq 0\right.\right\}.
\]
Then the double coset
\[
(H^\star)^{x^\star_2}(F)\backslash S(x^\star_2)/L^v
\]
is finite, for which we fix a set $\{h_1,\dots,h_n\}$ of representatives.

Take an index $1\leq j\leq n$. We
\begin{itemize}
  \item put $\Gamma_j\coloneqq H^\star(F)\cap h_j L^v h_j^{-1}$, regarded as a discrete subgroup of $H^\star(F_v)$;

  \item denote by
      \[
      \pi_j\colon\fN^\delta\times\{h_j\}\to\Gamma_j\backslash\(\fN^\delta\times h_jL^v/L^v\)
      \]
      the natural projection map;

  \item denote by $\cX_{L,j}$ the unique connected component of $\cX_L\otimes_{R_0}R$ such that the formal completion along its special fiber $Y_{L,j}\coloneqq\cX_{L,j}\otimes_R\kappa$ coincides with $\Gamma_j\backslash\(\fN^\delta\times h_jL^v/L^v\)$;

  \item denote by $N^\delta_j$ the normalization of the induced reduced subscheme of $\fN^\delta\times\{h_j\}=\fN^\delta\times h_jL^v/L^v$, and by $|N^\delta_j|$ the set of its irreducible components.
\end{itemize}
For every abelian group $A$, denote by $\Fun_c(|N^\delta_j|,A)\subseteq\Fun(|N^\delta_j|,A)$ the subgroup of finitely supported functions (valued in $A$), for which we have the map
\[
\Delta\colon\dZ[|N^\delta_j|]\to\Fun_c(|N^\delta_j|,\dZ)
\]
given by the intersection in (the regular formal scheme) $\fN^\delta\times\{h_j\}$, and the map
\[
\r{av}^{\Gamma_j}\colon\Fun_c(|N^\delta_j|,V^\star_\lambda)\to\Fun(|N^\delta_j|,V^\star_\lambda)^{\Gamma_j}
\]
sending $f$ to $\sum_{\gamma\in\Gamma_j}\gamma f$. The composition $\r{av}^{\Gamma_j}\circ(\Delta\otimes V^\star_\lambda)$ factors through uniquely a map
\[
\Delta_{\lambda,j}\colon\(\dZ[|N^\delta_j|]\otimes V^\star_\lambda\)_{\Gamma_j}\to\Fun(|N^\delta_j|,V^\star_\lambda)^{\Gamma_j}.
\]
We also have the map
\[
\nabla_j\colon\Div(\fN^\delta\times\{h_j\})\to\Fun_c(|N^\delta_j|,\dZ)
\]
given by the intersection in $\fN^\delta\times\{h_j\}$.

\begin{lem}\label{le:height_curve_nonsplit_1}
For every $\dL$-valued function $\phi$ on $\fp_v^a\Lambda_v/\fp_v^b\Lambda_v$ and every element $x^\star\in V^\star$, there exists an element
\[
\cF^\lambda(x^\star,\phi)_j\in\(\dL[|N^\delta_j|]\otimes V^\star_\lambda\)_{\Gamma_j}
\]
such that
\[
\Delta_{\lambda,j}\cF^\lambda(x^\star,\phi)_j-\r{av}^{\Gamma_j}\(C(x^\star_v,\phi)_j\otimes \mu^\lambda(x^\star)\)
\]
is a constant function on $|N^\delta_j|$ (which must be the zero function when $\lambda\neq\b{0}$). Here,
\begin{itemize}
  \item $C(x^\star_v,\phi)_j\coloneqq\nabla_j\(\fZ_a^\delta(x^\star_v,\phi)\times\{h_j\}\)\in\Fun_c(|N^\delta_j|,\dL)$, where $\fZ_a^\delta(x^\star_v,\phi)\times\{h_j\}\in\Div(\fN^\delta\times\{h_j\})\otimes\dL$ is from Definition \ref{de:special_cycle_formal}(3);

  \item $\mu^\lambda$ is defined in Construction \ref{co:harmonic_function}(3) (with $R=\dQ$).
\end{itemize}
\end{lem}

\begin{proof}
Choose an auxiliary rational prime $\ell\neq p$. By Remark \ref{re:height_curve_nonsplit} below, the map $\Delta_{\lambda,j}\otimes\dQ_\ell$ is identified with the $\ell$-adic Gysin map \eqref{eq:height_curve_nonsplit_1}. Note that the cokernel of \eqref{eq:height_curve_nonsplit_1}, which is same as the image of the natural map
\[
\rH^2_{\et}(\cX_{L,j},V_{\ell,L,\lambda}(1))\to\rH^2_{\et}(\cX_{L,j}[1/p],V_{\ell,L,\lambda}(1)),
\]
is the $\Gal(\breve{E}_u/E_u)$-finite subspace of $\rH^2_{\et}(\cX_{L,j}[1/p],V_{\ell,L,\lambda}(1))$. There are two cases.

First suppose that $\lambda\neq\b{0}$. Then the $\Gal(\breve{E}_u/E_u)$-finite subspace of $\rH^2_{\et}(\cX_{L,j}[1/p],V_{\ell,L,\lambda}(1))$ vanishes, which implies that $\Delta_{\lambda,j}\otimes\dQ_\ell$ hence $\Delta_{\lambda,j}$ are surjective. The lemma then follows.

Second suppose that $\lambda=\b{0}$. Then the $\Gal(\breve{E}_u/E_u)$-finite subspace of $\rH^2_{\et}(\cX_{L,j}[1/p],V_{\ell,L,\lambda}(1))$ is generated by the image of constant functions on $|N^\delta_j|$ under the identification $\rH^2_{\et}(\cX_{L,j},\dQ_\ell(1))\simeq\Fun(|N^\delta_j|,\dQ_\ell)^{\Gamma_j}$ from Remark \ref{re:height_curve_nonsplit} below. The lemma then follows.
\end{proof}

Denote by $N^\delta$ the normalization of the induced reduced subscheme of $\fN^\delta$, and by $|N^\delta|$ the set of its irreducible components. For every $1\leq j\leq n$, we fix an element
\[
\fF^\lambda(x^\star_2,\phi^\infty_{2,v})_j\in\dL[|N^\delta|]\otimes V^\star_\lambda
\]
such that $\fF^\lambda(x^\star_2,\phi^\infty_{2,v})_j\times\{h_j\}$ lifts $\cF^\lambda(x^\star_2,\phi^\infty_{2,v})_j$ from Lemma \ref{le:height_curve_nonsplit_1}.

\begin{lem}\label{le:height_curve_nonsplit_2}
For every $1\leq j\leq n$, there exists a unique $V_\lambda\otimes_{\dE^\lambda}\dL$-valued Schwartz function $\phi^{[a,b]}_{v,j}$ on $V^\star_v$ such that for every $x^\star_1\in V^\star_v\setminus\{0\}$,
\begin{align*}
\phi^{[a,b]}_{v,j}(x^\star_1)=
\(\fZ_a^\delta(x^\star_1,\phi^\infty_{1,v})\bullet\fZ_a^\delta(x^\star_{2,v},\phi^\infty_{2,v})\)
\cdot\mu^\lambda(x^\star_2)-\fZ_a^\delta(x^\star_1,\phi^\infty_{1,v})\bullet\fF^\lambda(x^\star_2,\phi^\infty_{2,v})_j
\end{align*}
in which $\bullet$ denotes the intersection pairing on $\fN^\delta$ as an $\dL$-hermitian form.
\end{lem}

\begin{proof}
The continuity of $\phi^{[a,b]}_{v,j}$ on $V^\star_v\setminus\{0\}$ follows from Lemma \ref{le:special_divisor_formal}. Since $\phi^\infty_{1,v}$ is regularly supported, $\fZ_a^\delta(x^\star_1,\phi^\infty_{1,v})=0$ unless
$x^\star_1\in D_{(2b,2a]}\coloneqq\{x^\star\in V^\star_v\res (x^\star,x^\star)\in\fp_v^{2a}\setminus\fp_v^{2b}\}$. Thus, $\phi^{[a,b]}_{v,j}$ extends uniquely to a $V_\lambda\otimes_{\dE^\lambda}\dL$-valued continuous function on $V^\star_v$. It remains to show that $\phi^{[a,b]}_{v,j}$ is compactly supported. Indeed, this follows from the fact that for every $t\in F_v^\star$, the subset
\[
\{x^\star\in(V^\star_v)_t\res \fZ_a^0(x^\star,0)\cap\Sigma_j\neq\emptyset\}
\]
is compact, where $\Sigma_j$ is the image of the support of $\fF^\lambda(x^\star_2,\phi^\infty_{2,v})_j$ in $\fN^0$. The lemma is proved.
\end{proof}

For $1\leq j\leq n$, we also define a function $\phi^{\infty,v}_j$ on $V^{\infty,v}$ by the formula
\[
\phi^{\infty,v}_j(x)=\ol{\ts_1\phi^{\infty,v}_1}(h_j^{-1}x)\cdot\ts_2\phi^{\infty,v}_2(h_j^{-1}x^\star_2),
\qquad x\in V^{\infty,v},
\]
which is an element in $\sS(V^{\infty,v},\dL)$.

\begin{proposition}\label{pr:height_curve_nonsplit}
Suppose that we are in the situation of Proposition \ref{pr:height_curve}(2) together with Assumption \ref{as:height_curve}. If we put
\[
\phi^\infty_{t_2}\coloneqq\vol^\natural(L)\cdot\sum_{j=1}^n\phi^{\infty,v}_j\otimes\phi^{[a,b]}_{v,j}
\]
as a $V_\lambda\otimes_{\dE^\lambda}\dL$-valued Schwartz function on $V^\star\otimes_F\bA_F^\infty$, then for every $t_1\in(F^\times)^+$,
\begin{align*}
\vol^\natural(L)\cdot
\llangle\ts_1^*\bbZ^\lambda_{t_1}(\phi^\infty_1)_L(g_1),
\ts_2^*\bbZ^\lambda_{t_2}(\phi^\infty_2)_L\rrangle_{\lambda,L,\ul{Z},u}^\ell
=\sum_{x^\star_1\in V(v)_{t_1}}\omega_{1,\infty}(g_1)\(P[\phi^\infty_{t_2}(x^\star_1)]\cdot\phi^0_{V_\infty}\)(x^\star_1)
\end{align*}
holds (in $\dC\otimes\dQ_\ell$) for every element $g_1\in G_1(F_\infty)$ and every prime $\ell\neq p$.
\end{proposition}

Take a rational prime $\ell\neq p$. Note that the $\dE^\lambda_\ell$-linear local system $V_{\ell,L',\lambda}$ over $\cX_{L'}$ has a uniformization:
\begin{align}\label{eq:height_curve_nonsplit_5}
V_{\ell,L',\lambda}\res_{\fX_{L'}}=H^\star(F)\backslash\(\fN^\delta\times H(\bA_F^{\infty,v})/{L'}^v\times V_{\ell,\lambda}\)
\end{align}
compatible with \eqref{eq:uniformization_4}, in which $H^\star(F)$ acts on $V_{\ell,\lambda}\coloneqq(V_\ell)_\lambda$ via the natural map $H^\star(F)\to H(F_\ell)$.

\begin{remark}\label{re:height_curve_nonsplit}
By \eqref{eq:height_curve_nonsplit_5}, the pullback map
\[
\rH^2_{\et}(\cX_{L,j},V_{\ell,L,\lambda}(1))\to\prod_{C\in|N^\delta_j|}\rH^2_{\et}(C,V_{\ell,L,\lambda}(1))=
\Fun(|N^\delta_j|,V_{\ell,L,\lambda})
\]
induces an isomorphism
\[
\rH^2_{\et}(\cX_{L,j},V_{\ell,L,\lambda}(1))\simeq\Fun(|N^\delta_j|,V_{\ell,L,\lambda})^{\Gamma_j}.
\]
On the other hand, the Gysin map
\[
\dZ[|N^\delta_j|]\otimes V_{\ell,L,\lambda}=
\bigoplus_{C\in|N^\delta_j|}\rH^0_{\et}(C,V_{\ell,L,\lambda})\to\rH^2_{\et,Y_{L,j}}(\cX_{L,j},V_{\ell,L,\lambda}(1))
\]
induces an isomorphism
\[
\(\dZ[|N^\delta_j|]\otimes V^\star_\lambda\)_{\Gamma_j}\simeq\rH^2_{\et,Y_{L,j}}(\cX_{L,j},V_{\ell,L,\lambda}(1)).
\]
Then it is clear that under the above two isomorphisms, the Gysin map
\begin{align}\label{eq:height_curve_nonsplit_1}
\rH^2_{\et,Y_{L,j}}(\cX_{L,j},V_{\ell,L,\lambda}(1))\to\rH^2_{\et}(\cX_{L,j},V_{\ell,L,\lambda}(1))
\end{align}
coincides with $\Delta_{\lambda,j}\otimes\dQ_\ell$.
\end{remark}

The construction below is parallel to Construction \ref{co:special_divisor_split}, except that we use the na\"{\i}ve horizontal divisor $\fZ_a^\delta(x^\star_v,\phi_v)$.

\begin{construction}\label{co:special_divisor_nonsplit}
Consider a Schwartz function $\phi^\infty\in\sS(V^\infty,\dL)$ of the form $\phi_v\otimes\phi^{\infty,v}$ in which $\phi_v$ is a regularly supported $\dL$-valued function on $\fp_v^a\Lambda_v/\fp_v^b\Lambda_v$.
\begin{enumerate}
  \item For $x^\star_v\in V^\star_v\setminus\{0\}$, put
      \[
      \fZ(x^\star_v,\phi_v)_h\coloneqq
      \fZ_a^\delta(x^\star_v,\phi_v)\times\{h\}\in\Div(\fN^\delta\times H(\bA_F^{\infty,v})/L^v)\otimes\dL
      \]
      (Definition \ref{de:special_cycle_formal}(3)), and denote by
      \[
      \fC(x^\star_v,\phi_v)_h\in\rH^2_{\et,|\fZ(x^\star_v,\phi_v)_h|}(\fN^\delta\times H(\bA_F^{\infty,v})/L^v,\dL_\ell(1))
      \]
      its refined $\ell$-adic cycle class.

  \item For every $t\in(F^\times)^+$, define $\cZ_t(\phi^\infty)_L$ to be the unique element in $\Div(\cX_L\otimes_{R_0}R)\otimes\dL$ supported on the basic locus whose pullback to $\fN^\delta\times H(\bA_F^{\infty,v})/L^v$ is given by
      \[
      \fZ_t(\phi^\infty)_L\coloneqq\sum_{x^\star\in V^\star_t}\sum_{h\in H(\bA_F^{\infty,v})/L^v}
      \phi^{\infty,v}(h^{-1}x^\star)\cdot\fZ(x^\star_v,\phi_v)_h.
      \]

  \item For every $t\in(F^\times)^+$, denote by
      \[
      \cC^\lambda_t(\phi^\infty)_L\in
      \rH^2_{\et,|\cZ_t(\phi^\infty)_L|}(\cX_L\otimes_{R_0}R,V_{\ell,L,\lambda}(1))\otimes_{\dE^\lambda}\dL
      \]
      the unique element whose restriction to $\fN^\delta\times H(\bA_F^{\infty,v})/L^v$ is given by the $H^\star(F)$-invariant class
      \begin{align*}
      &\quad\sum_{x^\star\in V^\star_t}\(\sum_{h\in H(\bA_F^{\infty,v})/L^v}
      \phi^{\infty,v}(h^{-1}x^\star)\cdot\fC(x^\star_v,\phi_v)_h\)\otimes\mu^\lambda(x^\star) \\
      &\in\rH^2_{\et,|\fZ_t(\phi^\infty)_L|}(\fN^\delta\times H(\bA_F^{\infty,v})/L^v,\dL_\ell(1))
      \otimes_{\dE^\lambda}V^\star_\lambda
      =\rH^2_{\et,|\fZ_t(\phi^\infty)_L|}(\fN^\delta\times H(\bA_F^{\infty,v})/L^v,\dL_\ell(1))
      \otimes_{\dE^\lambda_\ell}V_{\ell,\lambda},
      \end{align*}
      where $\mu^\lambda$ is defined in Construction \ref{co:harmonic_function}(3) (with $R=\dQ$).
\end{enumerate}
\end{construction}

For every $1\leq j\leq n$, denote by
\[
\cC^\lambda(x^\star_2,\phi^\infty_{2,v})_j\in
\rH^2_{\et,\pi_j|\fZ(x^\star_{2,v},\phi^\infty_{2,v})_{h_j}|}(\cX_j,V_{\ell,\lambda}(1))\otimes_{\dE^\lambda}\dL
\]
the pushforward of the class
\begin{align*}
&\quad\fC(x^\star_{2,v},\phi^\infty_{2,v})_{h_j}\otimes\mu^\lambda(x^\star_2) \\
&\in\rH^2_{\et,|\fZ(x^\star_{2,v},\phi^\infty_{2,v})_{h_j}|}(\fN^\delta\times\{h_j\},\dL_\ell(1))
\otimes_{\dE^\lambda}V^\star_\lambda
=\rH^2_{\et,|\fZ(x^\star_{2,v},\phi^\infty_{2,v})_{h_j}|}(\fN^\delta\times\{h_j\},\dL_\ell(1))
\otimes_{\dE^\lambda_\ell}V_{\ell,\lambda}
\end{align*}
(Construction \ref{co:special_divisor_nonsplit}(1)) along $\pi_j$.

\begin{lem}\label{le:height_curve_nonsplit_3}
The natural image of
\begin{align}\label{eq:height_curve_nonsplit_2}
\sum_{j=1}^n\ts_2\phi^{\infty,v}_2(h_j^{-1}x^\star_2)
\(\cC^\lambda(x^\star_2,\phi^\infty_{2,v})_j-\cF^\lambda(x^\star_2,\phi^\infty_{2,v})_j\)
\end{align}
in $\rH^2_{\et}(\cX_L\otimes_{R_0}R,V_{\ell,L,\lambda}(1))\otimes_{\dE^\lambda}\dL$ vanishes. Here, for the term $\cF^\lambda(x^\star_2,\phi^\infty_{2,v})_j$ from Lemma \ref{le:height_curve_nonsplit_1}, we have used the map $\(\dZ[|N^\delta_j|]\otimes V^\star_\lambda\)_{\Gamma_j}\to\rH^2_{\et}(\cX_{L,j},V_{\ell,L,\lambda}(1))$ from Remark \ref{re:height_curve_nonsplit}.
\end{lem}

\begin{proof}
When $\lambda\neq\b{0}$, this follows immediately from Lemma \ref{le:height_curve_nonsplit_1} and Remark \ref{re:height_curve_nonsplit}.

When $\lambda=\b{0}$, the image of \eqref{eq:height_curve_nonsplit_2} in $\rH^2_{\et}(X_{L,u}\otimes_{R_0}R,\dL_\ell(1))$ vanishes, which implies that the total degree of \eqref{eq:height_curve_nonsplit_2} as an $\dL_\ell$-valued function on the set of vertical irreducible components on each connected component of $\cX_L\otimes_{R_0}R$ is zero. Combining this with Lemma \ref{le:height_curve_nonsplit_1}, the lemma follows.
\end{proof}

\begin{proof}[Proof of Proposition \ref{pr:height_curve_nonsplit}]
Similar to the proof of Proposition \ref{pr:height_inert}, we may assume $g_1=1_2$.

Put
\begin{align}\label{eq:height_curve_nonsplit_3}
\cC^\lambda_{t_2}(\ts_2\phi^\infty_2)_L\coloneqq\sum_{j=1}^n\ts_2\phi^{\infty,v}_2(h_j^{-1}x^\star_2)
\(\cC^\lambda(x^\star_2,\phi^\infty_{2,v})_j-\cF^\lambda(x^\star_2,\phi^\infty_{2,v})_j\),
\end{align}
regarded as an element in $\rH^2_{\et,(Z_{2,u}\cup Y_L)\otimes_{R_0}R}(\cX_L\otimes_{R_0}R,V_{\ell,L,\lambda}(1))$, so that, by construction, its natural image in $\rH^2_{\et,Z_{2,u}\otimes_{R_0}R}(X_{L,u}\otimes_{R_0}R,V_{\ell,L,\lambda}(1))$ coincides with that of $\ts_2^*\bbZ^\lambda_{t_2}(\phi^\infty_2)_L$. On the other hand, the natural image of the class $\cC^\lambda_{t_1}(\ts_1\phi^\infty_1)_L$ from Construction \ref{co:special_divisor_nonsplit}(3) in $\rH^2_{\et,Z_{1,u}\otimes_{R_0}R}(X_{L,u}\otimes_{R_0}R,V_{\ell,L,\lambda}(1))$ coincides with that of $\ts_1^*\bbZ^\lambda_{t_1}(\phi^\infty_1)_L$.

Then by Lemma \ref{le:height_integral} (or rather its analogue after base change to $R$) and Lemma \ref{le:height_curve_nonsplit_3}, we have
\[
\llangle\ts_1^*\bbZ^\lambda_{t_1}(\phi^\infty_1)_L,\ts_2^*\bbZ^\lambda_{t_2}(\phi^\infty_2)_L\rrangle_{\lambda,L,\ul{Z},u}^\ell
=\cC^\lambda_{t_1}(\ts_1\phi^\infty_1)_L\overset{\lambda}\bullet\cC^\lambda_{t_2}(\ts_2\phi^\infty_2)_L.
\]
By \eqref{eq:height_curve_nonsplit_3},
\begin{align*}
\cC^\lambda_{t_1}(\ts_1\phi^\infty_1)_L\overset{\lambda}\bullet\cC^\lambda_{t_2}(\ts_2\phi^\infty_2)_L=
\sum_{j=1}^n\ol{\ts_2\phi^{\infty,v}_2}(h_j^{-1}x^\star_2)\cdot
\cC^\lambda_{t_1}(\ts_1\phi^\infty_1)_L\overset{\lambda}\bullet
\(\cC^\lambda(x^\star_2,\phi^\infty_{2,v})_j-\cF^\lambda(x^\star_2,\phi^\infty_{2,v})_j\),
\end{align*}
Note that
\[
\cC^\lambda(x^\star_2,\phi^\infty_{2,v})_j=
(\pi_j)_*\(\fC(x^\star_{2,v},\phi^\infty_{2,v})_{h_j}\otimes\mu^\lambda(x^\star_2)\),
\]
and
\[
\cF^\lambda(x^\star_2,\phi^\infty_{2,v})_j=(\pi_j)_*\(\fF^\lambda(x^\star_2,\phi^\infty_{2,v})_j\times\{h_j\}\).
\]
It follows by the projection formula that
\begin{align}\label{eq:height_curve_nonsplit_4}
&\quad\cC^\lambda_{t_1}(\ts_1\phi^\infty_1)_L\overset{\lambda}\bullet
\(\cC^\lambda(x^\star_2,\phi^\infty_{2,v})_j-\cF^\lambda(x^\star_2,\phi^\infty_{2,v})_j\) \\
&=\pi_j^*\cC^\lambda_{t_1}(\ts_1\phi^\infty_1)_L\overset{\lambda}\bullet
\(\fC(x^\star_{2,v},\phi^\infty_{2,v})_{h_j}\otimes\mu^\lambda(x^\star_2)
-\fF^\lambda(x^\star_2,\phi^\infty_{2,v})_j\times\{h_j\}\), \notag
\end{align}
where the second $\overset{\lambda}\bullet$ is the analogous pairing over $\fN^\delta\times\{h_j\}$. By Construction \ref{co:special_divisor_nonsplit}(3),
\[
\pi_j^*\cC^\lambda_{t_1}(\ts_1\phi^\infty_1)_L=\sum_{x^\star_1\in V^\star_t}
\ts_1\phi^{\infty,v}_1(h_j^{-1}x^\star_1)\cdot\fC(x^\star_{1,v},\phi^\infty_{1,v})_{h_j}\otimes\mu^\lambda(x^\star_1).
\]
It follows by Lemma \ref{le:height_curve_nonsplit_2} that
\[
\eqref{eq:height_curve_nonsplit_4}=\sum_{x^\star_1\in V^\star_t}
\ts_1\phi^{\infty,v}_1(h_j^{-1}x^\star_1)\cdot(\mu^\lambda(x^\star_1),\phi^{[a,b]}_{v,j}(x^\star_{1,v}))_{V^\star_\lambda}.
\]
Therefore,
\begin{align*}
\cC^\lambda_{t_1}(\ts_1\phi^\infty_1)_L\overset{\lambda}\bullet\cC^\lambda_{t_2}(\ts_2\phi^\infty_2)_L
&=\sum_{j=1}^n\ol{\ts_2\phi^{\infty,v}_2}(h_j^{-1}x^\star_2)\sum_{x^\star_1\in V^\star_t}
\ts_1\phi^{\infty,v}_1(h_j^{-1}x^\star_1)
\cdot(\mu^\lambda(x^\star_1),\phi^{[a,b]}_{v,j}(x^\star_{1,v}))_{V^\star_\lambda} \\
&=\sum_{x^\star_1\in V^\star_t}\sum_{j=1}^n\ts_1\phi^{\infty,v}_1(h_j^{-1}x^\star_1)\ol{\ts_2\phi^{\infty,v}_2}(h_j^{-1}x^\star_2)
\cdot(\mu^\lambda(x^\star_1),\phi^{[a,b]}_{v,j}(x^\star_{1,v}))_{V^\star_\lambda} \\
&=\vol^\natural(L)^{-1}\cdot\sum_{x^\star_1\in V^\star_t}(\mu^\lambda(x^\star_1),\phi^\infty_{t_2}(x^\star_1))_{V^\star_\lambda},
\end{align*}
so that
\[
\vol^\natural(L)\cdot
\llangle\ts_1^*\bbZ^\lambda_{t_1}(\phi^\infty_1)_L,
\ts_2^*\bbZ^\lambda_{t_2}(\phi^\infty_2)_L\rrangle_{\lambda,L,\ul{Z},u}^\ell=
\sum_{x^\star_1\in V^\star_t}(\mu^\lambda(x^\star_1),\phi^\infty_{t_2}(x^\star_1))_{V^\star_\lambda}.
\]
The proposition is proved.
\end{proof}

\subsection{Archimedean local heights over unitary moduli curves}
\label{ss:height_curve_archimedean}

In this section, we assume $r=1$ and take an element $\bbu\in\Hom(E,\dC)$. Put $\bbv\coloneqq v(\bbu)\in\fV_F^{(\infty)}$. As in Remark \ref{re:moduli_shimura}, denote by $V(\bbv)$ the unique hermitian space over $E/F$ that has signature $(1,1)$ at $\bbv$ and such that $V(\bbv)\otimes_F\bA_F^{\bbv}$ is isometric to $V\otimes_{\bA_F}\bA_F^{\bbv}$. Put $H(\bbv)\coloneqq\rU(V(\bbv))$, and fix an isometry $V(\bbv)\otimes_F\bA_F^\infty\simeq V^\infty$ under which we identify $H(\bA_F^\infty)$ with $H(\bbv)(\bA_F^\infty)$.

As before, for every intermediate field $\dE^\lambda\subseteq\dL\subseteq\dC$, we canonically identify both $\tS_{\lambda^+}(\dL)$ and $\ol{\tS_{\lambda^+}(\dL)}$ with $\dL$ as $\dL$-vector spaces.

\begin{remark}\label{re:hermitian_siegel}
Since we have identified both $\tS_{\lambda^+}(\dC)$ and $\ol{\tS_{\lambda^+}(\dC)}$ with $\dC$, now $\cA_{\lambda^+}(G_1)$ (Definition \ref{de:hermitian_siegel}) is simply the subspace of $\cA(G_1)$ consisting of $(\fg_{1,\infty},K_{1,\infty})$-finite automorphic forms $\varphi$ on $G_1(\bA_F)$ satisfying
\[
R_k(\varphi)=\rho_{\lambda^+}(k)\delta(k)\cdot\varphi
\]
for every $k\in K_{1,\infty}$, where $\rho_{\lambda^+}$ is the character of $K_{1,\infty}$ determined by $\lambda^+$ given in Definition \ref{de:hermitian_siegel}.
\end{remark}

\begin{proposition}\label{pr:height_archimedean}
Consider
\begin{itemize}
  \item an object $L$ of $\Hk(H(\bA_F^\infty))$,

  \item a pair $\ul\phi^\infty=(\phi^\infty_1,\phi^\infty_2)$ of decomposable Schwartz functions in $\sS(V^\infty)^L$ satisfying that $(\phi^\infty_1\otimes\ol{\phi^\infty_2})_v$ is regularly supported (Definition \ref{de:regular_support}) for some $v\in\fV_F^\fin$,

  \item a pair $(\ts_1,\ts_2)$ of elements in $\dC[L^v\backslash H(\bA_F^{\infty,v})/L^v]^\lambda$,

  \item an open compact subgroup $K^\infty$ of $G_1(\bA_F^\infty)$ that fixes $\phi^\infty_1$.
\end{itemize}
For every element $t_2\in(F^\times)^+$, there exists an element $\varphi^\KM_{t_2}\in\cA_{\lambda^+}(G_1)^{K^\infty}$ (depending on the data so far) such that
\begin{enumerate}[label=(\alph*)]
  \item for every $t_1\in(F^\times)^+$ (so that $\ul{Z}=(Z_1\coloneqq Z_{t_1}[\ts_1\phi^\infty_1]_L,Z_2\coloneqq Z_{t_2}[\ts_2\phi^\infty_2]_L)$ (Definition \ref{de:special_cycle}) is a disjoint pair),
      \begin{align*}
      \vol^\natural(L)\cdot
      \llangle\ts_1^*\bbZ^\lambda_{t_1}(\phi^\infty_1)_L(g_1),
      \ts_2^*\bbZ^\lambda_{t_2}(\phi^\infty_2)_L\rrangle_{\lambda,L,\ul{Z}}^{\bbu}
      =\frac{1}{2}\fE_{t_1,t_2}((g_1,1_2),\Phi^\lambda_\infty\cdot(\ts_1\phi^\infty_1\otimes\ol{\ts_2\phi^\infty_2}))_{\bbv}
      +(\varphi^\KM_{t_2})_{t_1}(g_1),
      \end{align*}

  \item for every $t_1\in F\setminus(F^\times)^+$,
      \[
      (\varphi^\KM_{t_2})_{t_1}(g_1)=0
      \]
\end{enumerate}
hold (in $\dC$) for every $g_1\in G_1(F_\infty)$, where $(\varphi^\KM_{t_2})_{t_1}$ denotes the $t_1$-th Fourier coefficient of $\varphi^\KM_{t_2}$ (as a function on $P_1(F)\backslash G_1(\bA_F)$).
\end{proposition}

\begin{proof}
For $i=1,2$ and $t_i\in(F^\times)^+$, we have a function $\rF^\lambda_{t_i}(\phi^\infty_i)$ on $G_1(F_\infty)$ valued in $\sA^{1,1}(X_L^{\bbu},\tS_{[\lambda]}(V(\bbv)_\infty))$ defined by the formula
\[
\rF^\lambda_{t_i}(\phi^\infty_i)(g)=\rho_{\lambda^+}(k)\delta(k)
\sum_{(x,h)\in H(\bbv)(F)\backslash V(\bbv)_{t_1}\times H(\bA_F^\infty)/L}\phi^\infty_i(h^{-1}x)\cdot \varphi_\lambda(xa,h)_L\cdot\prod_{v\in\Hom(F,\dR)}\te^{2\pi\ti(t_i b_v)}
\]
if we write $g=m(a)n(b)k$ with $a\in E_\infty^\times$, $b\in F_\infty$, and $k\in K_{1,\infty}$. Here, $\varphi_\lambda(xa,h)_L\in\sA^{1,1}(X_L^{\bbu},\tS_{[\lambda]}(V(\bbv)_\infty))$ is the $\tS_{[\lambda]}(V(\bbv)_\infty)$-valued Funke--Millson (Kudla--Millson if $\lambda=\b{0}$) form constructed above \cite{Liu26}*{Definition~3.9}. Furthermore, we have another function $\rG^\lambda_{t_i}(\phi^\infty_i)$ on $G_1(F_\infty)$ valued in the space of $\tS_{[\lambda]}(V(\bbv)_\infty)$-valued Green functions on $X_L^{\bbu}$ defined by the formula
\[
\rG^\lambda_{t_i}(\phi^\infty_i)(g)=\rho_{\lambda^+}(k)\delta(k)
\sum_{(x,h)\in H(\bbv)(F)\backslash V(\bbv)_{t_1}\times H(\bA_F^\infty)/L}\phi^\infty_i(h^{-1}x)\cdot g_\lambda(xa,h)_L\cdot\prod_{v\in\Hom(F,\dR)}\te^{2\pi\ti(t_i b_v)},
\]
where $g_\lambda(xa,h)_L$ is defined above \cite{Liu26}*{Lemma~8.4}, so that $\rG^\lambda_{t_i}(\phi^\infty_i)(g)$ is a Green function for the $\tS_{[\lambda]}(V(\bbv)_\infty)$-valued divisor $\bbZ^\lambda_{t_i}(\phi^\infty_i)(g)$ with tail form $\rF^\lambda_{t_i}(\phi^\infty_i)(g)$ \cite{Liu26}*{Proposition~8.6}.

For $i=1,2$ and $t_i\in(F^\times)^+$, fix a harmonic Green function $\rG_{t_i}^\lambda(\ts_i\phi^\infty_i)^\heartsuit$ (that is, Green function with vanishing tail form \cite{Liu26}*{Definition~5.9}) for $\bbZ^\lambda_{t_i}(\ts_i\phi^\infty_i)$, so that $\rG_{t_i}^\lambda(\ts_i\phi^\infty_i)^\heartsuit-\rG_{t_i}^\lambda(\ts_i\phi^\infty_i)(1_2)$ is a smooth function on $X_L^{\bbu}$.

By \cite{Liu26}*{Proposition~3.10}, there exists a unique element
\[
\rF^\lambda(\phi^\infty_1)\in\cA_{\lambda^+}(G_1)^{K^\infty}\otimes_\dC\sA^{1,1}(X_L^{\bbu},\tS_{[\lambda]}(V(\bbv)_\infty))
\]
satisfying
\begin{itemize}
  \item $\rF^\lambda(\phi^\infty_1)$ is supported on $G_1(F)G_1(F_\infty)K^\infty$;

  \item for every $t_1\in(F^\times)^+$, $\rF^\lambda(\phi^\infty_1)_{t_1}\res_{G_1(F_\infty)}=\rF^\lambda_{t_1}(\phi^\infty_1)$;

  \item for every $t_1\in F\setminus (F^\times)^+$, $\rF^\lambda(\phi^\infty_1)_{t_1}\res_{G_1(F_\infty)}=0$.
\end{itemize}
Put
\[
\varphi^\KM_{t_2}\coloneqq\frac{\vol^\natural(L)}{2}\int_{X_L^{\bbu}}\(\rF^\lambda(\ts_1\phi^\infty_1),
\rG_{t_2}^\lambda(\ts_2\phi^\infty_2)^\heartsuit-\rG_{t_2}^\lambda(\ts_2\phi^\infty_2)(1_2)\)_{\tS_{[\lambda]}(V(\bbv)_\infty)}
\in \cA_{\lambda^+}(G_1)^{K^\infty}.
\]
We claim that the proposition holds with this function $\varphi^\KM_{t_2}$. Property (2) is obvious. It remains to show (1). It is clear that we may assume $g_1=m(a)$ for $a\in E_\infty^\times$.

By \cite{Liu26}*{Definition~5.12~\&~Remark~5.13}, we have
\[
\llangle\ts_1^*\bbZ^\lambda_{t_1}(\phi^\infty_1)_L(g_1),
\ts_2^*\bbZ^\lambda_{t_2}(\phi^\infty_2)_L\rrangle_{\lambda,L,\ul{Z}}^{\bbu}
=\frac{1}{2}\int_{X_L^{\bbu}} \rG_{t_1}^\lambda(\ts_1\phi^\infty_1)^\heartsuit
\ast \rG_{t_2}^\lambda(\ts_2\phi^\infty_2)^\heartsuit,
\]
where $\ast$ is defined above \cite{Liu26}*{Definition~5.12}. We compute the difference
\begin{align*}
&\quad\int_{X_L^{\bbu}} \rG_{t_1}^\lambda(\ts_1\phi^\infty_1)^\heartsuit
\ast \rG_{t_2}^\lambda(\ts_2\phi^\infty_2)^\heartsuit
-\int_{X_L^{\bbu}} \rG_{t_1}^\lambda(\ts_1\phi^\infty_1)(m(a))
\ast \rG_{t_2}^\lambda(\ts_2\phi^\infty_2)(1_2) \\
&=\int_{X_L^{\bbu}}
\(\rG_{t_1}^\lambda(\ts_1\phi^\infty_1)^\heartsuit-\rG_{t_1}^\lambda(\ts_1\phi^\infty_1)(m(a)),
\delta_{\bbZ^\lambda_{t_2}(\ts_2\phi^\infty_2)}\)_{\tS_{[\lambda]}(V(\bbv)_\infty)}
-\int_{X_L^{\bbu}}\(\rF_{t_1}^\lambda(\ts_1\phi^\infty_1)(m(a)),
\rG_{t_2}^\lambda(\ts_2\phi^\infty_2)(1_2)\)_{\tS_{[\lambda]}(V(\bbv)_\infty)}
\end{align*}
which, by Stokes' formula, equals
\begin{align*}
&\quad\int_{X_L^{\bbu}}
\(\rF_{t_1}^\lambda(\ts_1\phi^\infty_1)(m(a)),
\rG_{t_2}^\lambda(\ts_2\phi^\infty_2)^\heartsuit\)_{\tS_{[\lambda]}(V(\bbv)_\infty)}
-\int_{X_L^{\bbu}}\(\rF_{t_1}^\lambda(\ts_1\phi^\infty_1)(m(a)),
\rG_{t_2}^\lambda(\ts_2\phi^\infty_2)(1_2)\)_{\tS_{[\lambda]}(V(\bbv)_\infty)} \\
&=\int_{X_L^{\bbu}}
\(\rF_{t_1}^\lambda(\ts_1\phi^\infty_1)(m(a)),
\rG_{t_2}^\lambda(\ts_2\phi^\infty_2)^\heartsuit-\rG_{t_2}^\lambda(\ts_2\phi^\infty_2)(1_2)\)_{\tS_{[\lambda]}(V(\bbv)_\infty)} \\
&=2\vol^\natural(L)^{-1}\varphi^\KM_{t_2}(m(a)).
\end{align*}
Thus, (1) is reduced to the identity
\[
\vol^\natural(L)\int_{X_L^{\bbu}} \rG_{t_1}^\lambda(\ts_1\phi^\infty_1)(m(a))
\ast \rG_{t_2}^\lambda(\ts_2\phi^\infty_2)(1_2)=
\fE_{t_1,t_2}((m(a),1_2),\Phi^\lambda_\infty\cdot(\ts_1\phi^\infty_1\otimes\ol{\ts_2\phi^\infty_2}))_{\bbv}.
\]
However, this is exactly \cite{Liu26}*{Theorem~8.10}.
\end{proof}

\begin{proposition}\label{pr:height_archimedean_bis}
Consider
\begin{itemize}
  \item an object $L$ of $\Hk(H(\bA_F^\infty))$,

  \item a pair $\ul\phi^\infty=(\phi^\infty_1,\phi^\infty_2)$ of decomposable Schwartz functions in $\sS(V^\infty)^L$ satisfying that $(\phi^\infty_1\otimes\ol{\phi^\infty_2})_v$ is regularly supported (Definition \ref{de:regular_support}) for some $v\in\fV_F^\fin$,

  \item a pair $(\ts_1,\ts_2)$ in $\dC[L^v\backslash H(\bA_F^{\infty,v})/L^v]^\lambda$,

  \item an open compact subgroup $K^\infty$ of $G_1(\bA_F^\infty)$ that fixes $\phi^\infty_1$.
\end{itemize}
For every element $\varphi\in\cA_{\lambda^+}(G_1)_\cusp^\hol$ invariant under $K^\infty$, the identity
\begin{align*}
&\quad\vol^\natural(L)
\int_{\Omega}\ol{\varphi(g_1)}
\sum_{t_1\in (F^\times)^+}\llangle\ts_1^*\bbZ^\lambda_{t_1}(\phi^\infty_1)_L(g_1),
\ts_2^*\bbZ^\lambda_{t_2}(\phi^\infty_2)_L\rrangle_{\lambda,L,\ul{Z}}^{\bbu} \rd g_1\\
&=\frac{1}{2}\int_{\Omega}\ol{\varphi(g_1)}\cdot
\fE_{t_1,t_2}((g_1,1_2),\Phi^\lambda_\infty\cdot(\ts_1\phi^\infty_1\otimes\ol{\ts_2\phi^\infty_2}))_{\bbv}  \rd g_1
\end{align*}
holds in which both sides are absolutely convergent, where $\Omega\subseteq G_1(F_\infty)$ is an arbitrary fundamental region for the arithmetic subgroup $G_1(F)\cap K^\infty$.
\end{proposition}

\begin{proof}
This is \cite{Liu26}*{Theorem~8.11}.
\end{proof}

\subsection{Abundance of nonarchimedean Schwartz functions}
\label{ss:abundance}

We continue to assume $r=1$. Consider
\begin{itemize}
  \item a finite set $\Diamond$ of rational primes containing all those ramified in $E$ and all those underlying $\fV(V)$;

  \item an open compact subgroup $L^\Diamond$ of $H(\bA^{\infty,\Diamond}_F)$ that is the stabilizer of a self-dual lattice of $V^{\infty,\Diamond}$;

  \item a subset $\Diamondblack\subseteq\fV_F^{(\Diamond)}$ containing all those nonsplit in $E$ and satisfying $d(\Diamondblack)\geq 2$ (Notation \ref{no:cm}(5)).
\end{itemize}

\begin{definition}\label{de:schwartz_special}
For every field $\dL$ of characteristic zero, we define $\sS((V^\infty)^2,\dL)^{L^\Diamond}_\Diamondblack$ to be the subspace of $\sS((V^\infty)^2,\dL)^{L^\Diamond}$ spanned by decomposable functions $\Phi^\infty=\otimes_v\Phi^\infty_v$ satisfying
\begin{enumerate}
  \item for $v\in\Diamondblack$, $\Phi^\infty_v$ is regularly supported;

  \item for $v\in\Diamondblack$ nonsplit in $E$, $\Phi^\infty_v$ is deeply regularly supported;

  \item for $v\in\Diamondblack$ nonsplit in $E$, $\Phi^\infty_v$ is regularly supported along both variables of $V_v^2$.\footnote{Note that when $V_v$ is nonsplit, this is implied by (2).}
\end{enumerate}
See Definition \ref{de:regular_support} for the above terminologies. It is clear that $\sS((V^\infty)^2,\dL)^{L^\Diamond}_\Diamondblack$ is stable under the action of $H(F_\Diamond)$. Furthermore, we denote by
\[
\sS((V^\infty)^2,\dL)^{L^\Diamond}_{(\Diamondblack)}\subseteq\sS((V^\infty)^2,\dL)^{L^\Diamond}
\]
the $G_1(\bA_F^\infty)^2$-stable subspace generated by $\sS((V^\infty)^2,\dL)^{L^\Diamond}_\Diamondblack$.
\end{definition}

\begin{lem}\label{le:abundance}
Let $(\pi,\cA_\pi)$ be a holomorphic tempered cuspidal automorphic representation of $G_1(\bA_F)$ of weight $\lambda$ (Definition \ref{de:hermitian_representation}) and $\dL\subseteq\dC$ a subfield over which $\pi$ is defined. Suppose that $V=V(\pi)$ as in Remark \ref{de:hermitian_representation}(1) and $(\pi_v)^{K_{1,v}}\neq\{0\}$ (Notation \ref{no:hermitian_siegel}(1)) for every $v\in\fV_F^{\fin}\setminus\fV_F^{(\Diamond)}$. Then one can find a pair $\ul\varphi=(\varphi_1,\varphi_2)$ of automorphic forms in $\cA_\pi\cap\cA_{\lambda^+}(G_1)_\cusp^\hol$ and a pair $\ul\phi^\infty=(\phi^\infty_1,\phi^\infty_2)$ of Schwartz functions in $\sS(V^\infty,\dL)$ satisfying $\phi^\infty_1\otimes\ol{\phi^\infty_2}\in\sS((V^\infty)^2,\dL)^{L^\Diamond}_\Diamondblack$, such that
\[
\fZ^\natural_\fin(\ul\varphi,\phi^\infty_1\otimes\ol{\phi^\infty_2})\neq 0
\]
(Definition \ref{de:doubling_zeta}).
\end{lem}

\begin{proof}
By \cite{LL21}*{Proposition~3.13}, we may find $\ul\varphi$ and $\ul\phi^\infty$ as in the statement satisfying Definition \ref{de:schwartz_special}(1) such that $\fZ^\natural_\fin(\ul\varphi,\phi^\infty_1\otimes\ol{\phi^\infty_2})\neq 0$.

For each $v\in\Diamondblack$ nonsplit in $E$, we may replace $\phi^\infty_i$ by $\omega_{1,v}(m(a_i))\phi^\infty_i$ (and same for $\varphi_i$) with $a_i\in E_v^\times$ sufficiently close to zero to further meet \ref{de:schwartz_special}(2) without changing the value of $\fZ^\natural_\fin$.

Finally, for each $v\in\Diamondblack$ nonsplit in $E$, we may choose elements $b_1,b_2\in F_v$ such that $\varphi'_i\coloneqq(\pi_v(n(b_i))-1)\varphi_i\neq 0$ for $i=1,2$. Put $\phi^{\infty\prime}_i\coloneqq\omega_{1,v}(n(b_i)-1)\phi^\infty_i$ for $i=1,2$. Then $\phi^{\infty\prime}_1\otimes\ol{\phi^{\infty\prime}_2}$ belongs to $\sS((V^\infty)^2,\dL)^{L^\Diamond}_\Diamondblack$, satisfying $\fZ^\natural_\fin(\ul\varphi',\phi^{\infty\prime}_1\otimes\ol{\phi^{\infty\prime}_2})\neq 0$, where $\ul\varphi'=(\varphi'_1,\varphi'_2)$.

The lemma is proved.
\end{proof}

At the end of this section, we explain that the ingredients we obtained so far are enough to give a proof of Conjecture \ref{co:maipf} when $\lambda=\b{0}$ (and $r=1$), simplifying to some degree the one given in \cite{Liu12}; see Remark \ref{re:aipf} below.

For $(\pi,\cA_\pi)$ in Conjecture \ref{co:maipf}, we choose $\Diamond$ so that $(\pi_v)^{K_{1,v}}\neq\{0\}$ (Notation \ref{no:hermitian_siegel}(1)) for every $v\in\fV_F^{\fin}\setminus\fV_F^{(\Diamond)}$. We may assume $V=V(\pi)$ as in Remark \ref{re:dichotomy}(2), since otherwise both sides of the identity in Conjecture \ref{co:maipf} are zero. In particular, the global root number $\epsilon(\pi)=-1$.

Note that Conjecture \ref{co:modularity} is known when $\lambda=\b{0}$ by \cite{Liu11}*{Theorem~3.5}. It follows that for every pair $\ul\phi^\infty=(\phi^\infty_1,\phi^\infty_2)$ of Schwartz functions in $\sS(V^\infty)^L$ (with $L^\Diamond$ the stabilizer of a self-dual lattice) and every pair $\ul{\ts}=(\ts_1,\ts_2)$ of elements in $\dC[L^\Diamond\backslash H(\bA_F^{\infty,\Diamond})/L^\Diamond]^{\b{0}}$, the function $\fI_{\ul\phi^\infty}^{\ul{\ts}}$ on $G_1(\bA_F)\times G_1(\bA_F)$ defined by the formula
\[
\fI_{\ul\phi^\infty}^{\ul{\ts}}(g_1,g_2)\coloneqq\vol^\natural(L)
\llangle\ts_1^*Z^{\b{0}}(\phi^\infty_1)_L(g_1),\ts_2^*Z^{\b{0}}(\phi^\infty_2)_L(g_2)\rrangle_{\b{0},L}
\]
belongs to the space $\cA_{\b{0}}(G_1)^\hol\otimes_\dC\ol{\cA_{\b{0}}(G_1)^\hol}$.

Now suppose that
\begin{enumerate}[label=(\alph*)]
  \item $\ul\phi^\infty$ is decomposable such that $\phi^\infty_1\otimes\ol{\phi^\infty_2}\in\sS((V^\infty)^2,\dL)^{L^\Diamond}_\Diamondblack$ for some $\Diamondblack$ as above, and that $\phi^{\infty,\Diamond}_1=\phi^{\infty,\Diamond}_2$ is the characteristic function of the self-dual lattice $L^\Diamond$ stabilizes;

  \item both $\ts_1$ and $\ts_2$ are supported at places split in $E$.
\end{enumerate}
Using the height decomposition formula (Proposition \ref{pr:height_decomposition}) and combining Proposition \ref{pr:height_inert}(1) (plus Remark \ref{re:height_inert}), Proposition \ref{pr:height_split}, Proposition \ref{pr:height_curve} and Proposition \ref{pr:height_archimedean_bis} for various local heights, we conclude that for every $g_2\in G_1(\bA_F^{\Diamondblack})$, the difference
\[
\fI_{\ul\phi^\infty}^{\ul{\ts}}(-,g_2)-E'(0,(-,g_2),\Phi^0_\infty\otimes(\ts_1\phi^\infty_1\otimes\ol{\ts_2\phi^\infty_2})),
\]
as an element in $\cA_{\b{0}}(G_1)^\hol$, is perpendicular to $\cA_\pi\cap\cA_{\b{0}}(G_1)_\cusp^\hol$ under the inner product in Definition \ref{de:hermitian_siegel}(2).\footnote{Here, we need to use the fact that since $\epsilon(\pi)=-1$, all global theta lifts of $\cV_\pi$ to $\rU(V)$ vanish.} Then by Proposition \ref{pr:ripf} and Remark \ref{re:constant}(2), we obtain the identity
\[
\vol^\natural(L)
\llangle\Theta_{\ts_1\phi^\infty_1}^\lambda(\varphi_1),\Theta_{\ts_2\phi^\infty_2}^\lambda(\varphi_2)\rrangle_{\b{0},L}
=\frac{L'(\tfrac{1}{2},\pi)}{b_2(0)}\cdot\(-\frac{\pi}{4}\)^{[F:\dQ]}\cdot
\fZ^\natural_\fin(\ul\varphi,\ts_1\phi^\infty_1\otimes\ol{\ts_2\phi^\infty_2})
\]
for $\ul\phi^\infty$ satisfying (a), $\ul\ts$ satisfying (b), and every pair $\ul\varphi=(\varphi_1,\varphi_2)$ in $\cA_\pi\cap\cA_{\b{0}}(G_1)_\cusp^\hol$ satisfying
\begin{enumerate}
  \item[(c)] both $\varphi_1$ and $\varphi_2$ are invariant under $K_{1,v}$ (Notation \ref{no:hermitian_siegel}(1)) for every $v\in\fV_F^\fin\setminus\fV_F^{(\Diamond)}$.
\end{enumerate}

Now we can choose $\ul\ts$ satisfying (b) such that $\ts_i$ fixes $\varphi_i$ for $i=1,2$. Then we arrive at the identity
\begin{align}\label{eq:aipf}
\vol^\natural(L)\llangle\Theta_{\phi^\infty_1}^\lambda(\varphi_1),\Theta_{\phi^\infty_2}^\lambda(\varphi_2)\rrangle_{\b{0},L}
=\frac{L'(\tfrac{1}{2},\pi)}{b_2(0)}\cdot\(-\frac{\pi}{4}\)^{[F:\dQ]}\cdot
\fZ^\natural_\fin(\ul\varphi,\phi^\infty_1\otimes\ol{\phi^\infty_2}).
\end{align}
By Lemma \ref{le:abundance}, we may find $\ul\varphi$ satisfying (a) and $\ul\phi^\infty$ satisfying (c) such that $\fZ^\natural_\fin(\ul\varphi,\phi^\infty_1\otimes\ol{\phi^\infty_2})\neq 0$. Finally, by the multiplicity one property \cite{Liu12}*{Proposition~6.7}, \eqref{eq:aipf} holds for arbitrary pairs $\ul\varphi$ and $\ul\phi^\infty$.

\begin{remark}\label{re:aipf}
The proof given above simplifies the one given in \cite{Liu12} in the following three aspects:
\begin{enumerate}
  \item We have simplified the argument for nonarchimedean local heights at bad reductions without using the ``change of Shimura data'' argument in \cite{Liu12}; this is possible partially due to the recent advances on integral models of unitary Shimura varieties.

  \item We have used the notion of special Eisenstein ideal to annihilate certain degenerate terms that appeared in the computation of nonarchimedean local heights at bad reductions in \cite{Liu12}.

  \item We have avoided the step of holomorphic projection for the archimedean computation by using Proposition \ref{pr:height_archimedean_bis}, an idea originally adopted in \cite{LL21}.
\end{enumerate}
\end{remark}

\section{Proof of main results}
\label{ss:6}

In this chapter, we prove our main results. Fix a totally positive definite incoherent hermitian space $V$ over $\bA_E/\bA_F$ of rank $2$ satisfying that the set $\fV(V)$ is not a singleton (and put $H\coloneqq\rU(V)$ as always), and a \emph{nontrivial} balanced hermitian $E$-weight $\lambda$ of rank $2$ (Definition \ref{de:weight}).

\subsection{Modularity of motivic doubling kernel functions}
\label{ss:modularity}

By Theorem \ref{th:independence}, for every object $L$ of $\Hk(H(\bA_F^\infty))$, we have a well-defined height pairing
\[
\llangle\;,\;\rrangle_{\lambda,L}\colon\rH^2(X_L,\rM_{L,\lambda}(1))\times\rH^2(X_L,\rM_{L,\lambda}(1))\to\dC.
\]
As in \S\ref{ss:abundance}, we consider
\begin{itemize}
  \item a finite set $\Diamond$ of rational primes containing all those ramified in $E$ and all those underlying $\fV(V)$;

  \item an open compact subgroup $L^\Diamond$ of $H(\bA^{\infty,\Diamond}_F)$ that is the stabilizer of a self-dual lattice of $V^{\infty,\Diamond}$;

  \item a subset $\Diamondblack\subseteq\fV_F^{(\Diamond)}$ containing all those nonsplit in $E$ and satisfying $d(\Diamondblack)\geq 2$ (Notation \ref{no:cm}(5)).
\end{itemize}
Recall that we have the $G_1(\bA_F^\infty)^2$-stable subspace $\sS((V^\infty)^2,\dC)^{L^\Diamond}_{(\Diamondblack)}\subseteq\sS((V^\infty)^2,\dC)^{L^\Diamond}$ from Definition \ref{de:schwartz_special}.

\begin{theorem}\label{th:modularity}
There is a unique $G_1(\bA_F^\infty)^2$-equivariant complex linear map
\[
\fI\colon\sS((V^\infty)^2,\dC)^{L^\Diamond}_{(\Diamondblack)}\to
\cA_{\lambda^+}(G_1)_\cusp^\hol\otimes_\dC\ol{\cA_{\lambda^+}(G_1)_\cusp^\hol}
\]
(Remark \ref{re:hermitian_siegel}) such that for every pair $\ul\phi^\infty=(\phi^\infty_1,\phi^\infty_2)$ of Schwartz functions in $\sS(V^\infty,\dC)^{L^\Diamond}$ satisfying $\phi^\infty_1\otimes\ol{\phi^\infty_2}\in\sS((V^\infty)^2,\dC)^{L^\Diamond}_\Diamondblack$, the identity
\[
\fI(\phi^\infty_1\otimes\ol{\phi^\infty_2})(g_1,g_2)=\vol^\natural(L)\cdot\sum_{t_1\in(F^\times)^+}\sum_{t_2\in(F^\times)^+}
\llangle Z^\lambda_{t_1}(\phi^\infty_1)_L(g_1),Z^\lambda_{t_2}(\phi^\infty_2)_L(g_2)\rrangle_{\lambda,L}
\]
holds with absolute convergence on the right-hand side for every $(g_1,g_2)\in G_1(\bA_F^\Diamondblack)^2$ and every $L\in\Hk(H(\bA_F^\infty))$ fixing both $\phi^\infty_1$ and $\phi^\infty_2$.
\end{theorem}

\begin{proof}
The uniqueness is clear since $G_1(F)G_1(\bA_F^\Diamondblack)$ is dense in $G_1(\bA_F)$, and the functions of the form $\phi^\infty_1\otimes\ol{\phi^\infty_2}$ in the statement generate $\sS((V^\infty)^2,\dC)^{L^\Diamond}_{(\Diamondblack)}$ as a $G_1(\bA_F^\infty)^2$-module.

For every $\Phi^\infty\in\sS((V^\infty)^2,\dC)^{L^\Diamond}_{(\Diamondblack)}$, we may express it as
\begin{align}\label{eq:modularity_1}
\Phi^\infty=\sum_{\{\ul\phi^\infty\}}\omega_1^\infty(g'_1)\phi^\infty_1\otimes\ol{\omega_1^\infty(g'_2)\phi^\infty_2},
\end{align}
in which the sum is taken over finitely many pairs $\ul\phi^\infty=(\phi^\infty_1,\phi^\infty_2)$ in $\sS(V^\infty,\dC)^{L^\Diamond}$ satisfying $\phi^\infty_1\otimes\ol{\phi^\infty_2}\in\sS((V^\infty)^2,\dC)^{L^\Diamond}_\Diamondblack$ with $(g'_1,g'_2)\in G_1(\bA_F^\infty)^2$. For $(g_1,g_2)\in G_1(\bA_F^\Diamondblack)^2$, put
\[
\fI(\Phi^\infty)(g_1,g_2)\coloneqq\vol^\natural(L)\sum_{\{\ul\phi^\infty\}}\sum_{t_1\in(F^\times)^+}\sum_{t_2\in(F^\times)^+}
\llangle Z^\lambda_{t_1}(\phi^\infty_1)_L(g_1 g'_1),Z^\lambda_{t_2}(\phi^\infty_2)_L(g_2 g'_2)\rrangle_{\lambda,L},
\]
where $L$ is some object of $\Hk(H(\bA_F^\infty))$ fixing both $\phi^\infty_1$ and $\phi^\infty_2$ for every $\ul\phi^\infty$ appearing in the above sum. Assume the absolute convergence on the right-hand side. Then by Lemma \ref{le:height_pairing_level}, $\fI(\Phi^\infty)(g_1,g_2)$ does not depend on the choice of $L$; from this, it is straightforward to see that $\fI(\Phi^\infty)(g_1,g_2)$ also does not depend on the choice of the expression of \eqref{eq:modularity_1}.

Thus, it remains to show that for $\Phi^\infty=\phi^\infty_1\otimes\ol{\phi^\infty_2}\in\sS((V^\infty)^2,\dC)^{L^\Diamond}_\Diamondblack$ with $\phi^\infty_1$ and $\phi^\infty_2$ in $\sS(V^\infty,\dC)^{L^\Diamond}$, the sum
\[
\vol^\natural(L)\sum_{t_1\in(F^\times)^+}\sum_{t_2\in(F^\times)^+}
\llangle Z^\lambda_{t_1}(\phi^\infty_1)_L(g_1),Z^\lambda_{t_2}(\phi^\infty_2)_L(g_2)\rrangle_{\lambda,L}
\]
is absolutely convergent for $(g_1,g_2)\in G_1(\bA_F^\Diamondblack)^2$ and extends to a function in $\cA_{\lambda^+}(G_1)_\cusp^\hol\otimes_\dC\ol{\cA_{\lambda^+}(G_1)_\cusp^\hol}$. However, this is equivalent to the following claim: If $K^\infty\subset G_1(\bA_F^\infty)$ is an open compact subgroup that fixes both $\phi^\infty_1$ and $\phi^\infty_2$, then
\begin{enumerate}
  \item for every $t_2\in(F^\times)^+$ and $g_2\in G_1(\bA_F^\Diamondblack)$, the sum
      \[
      \vol^\natural(L)\sum_{t_1\in(F^\times)^+}
      \llangle Z^\lambda_{t_1}(\phi^\infty_1)_L(g_1),Z^\lambda_{t_2}(\phi^\infty_2)_L(g_2)\rrangle_{\lambda,L}
      \]
      is absolutely convergent for $g_1\in G_1(\bA_F^\Diamondblack)$ and extends to a function in $(\cA_{\lambda^+}(G_1)_\cusp^\hol)^{K^\infty}$;

  \item for every $t_1\in(F^\times)^+$ and $g_1\in G_1(\bA_F^\Diamondblack)$, the sum
      \[
      \vol^\natural(L)\sum_{t_2\in(F^\times)^+}
      \ol{\llangle Z^\lambda_{t_1}(\phi^\infty_1)_L(g_1),Z^\lambda_{t_2}(\phi^\infty_2)_L(g_2)\rrangle_{\lambda,L}}
      \]
      is absolutely convergent for $g_2\in G_1(\bA_F^\Diamondblack)$ and extends to a function in $(\cA_{\lambda^+}(G_1)_\cusp^\hol)^{K^\infty}$.
\end{enumerate}
We only show the property (1) by symmetry, that is,
\begin{align}\label{eq:modularity_2}
\sum_{t_1\in(F^\times)^+}\fI_{t_1}(g_1)
\end{align}
is absolutely convergent and extends to a function in $(\cA_{\lambda^+}(G_1)_\cusp^\hol)^{K^\infty}$, where
\[
\fI_{t_1}(g_1)\coloneqq\vol^\natural(L)\cdot
\llangle Z^\lambda_{t_1}(\phi^\infty_1)_L(g_1),Z^\lambda_{t_2}(\phi^\infty_2)_L(g_2)\rrangle_{\lambda,L},
\]
which extends uniquely to a smooth function on $G_1(\bA_F^\Diamondblack)K^\infty$. It is clear from its expression that it suffices to show that \eqref{eq:modularity_2} is absolutely convergent and extends to a function on $G_1(\bA_F)$ that is left invariant under $G_1(F)$. Without loss of generality, we may assume that $L$ is of the form $L_\Diamond L^\Diamond$ (with $L^\Diamond$ the given one at the beginning of this section).

Choose an auxiliary rational prime $\ell\not\in\Diamond$ that splits in $E_0$. For $g_1\in G_1(\bA_F^\Diamondblack)K^\infty$, the supports of $\bbZ^\lambda_{t_1}(\phi^\infty_1)_L(g_1)$ and $\bbZ^\lambda_{t_2}(\phi^\infty_2)_L(g_2)$ have empty intersection. By Proposition \ref{pr:height_decomposition} and \cite{LL21}*{Lemma~7.2}, we have for every $t_1\in(F^\times)^+$ that
\begin{align*}
\fI_{t_1}(g_1)
&=\vol^\natural(L)\(\sum_{u\in\fV_E^\fin\setminus\fV_E^{(\ell)}}\log\|u\|\cdot\llangle \bbZ^\lambda_{t_1}(\phi^\infty_1)_L(g_1),\bbZ^\lambda_{t_2}(\phi^\infty_2)_L(g_2)\rrangle_{\lambda,L,u}^\ell\) \\
&\quad +\vol^\natural(L)\(\sum_{\bbu\in\Hom(E,\dC)}\llangle \bbZ^\lambda_{t_1}(\phi^\infty_1)_L(g_1),\bbZ^\lambda_{t_2}(\phi^\infty_2)_L(g_2)\rrangle_{\lambda,L}^{\bbu}\),
\end{align*}
where we have suppress the pair $\ul{Z}$ of supports in the subscripts as it is clear from the context.

Denote by $\Diamondblack^{(1)}$ and $\Diamondblack^{(2)}$ the subsets of $\Diamondblack$ consisting of $v$ that is nonsplit in $E$ such that $V_v$ is split and nonsplit, respectively. For every $t_1\in(F^\times)^+$, put
\begin{align*}
\fK_{t_1}(g_1)&\coloneqq \sum_{v\in\fV_F\setminus(\Diamondblack^{(1)}\cup\Diamondblack^{(2)})}
\fE_{t_1,t_2}((g_1,g_2),\Phi^\lambda_\infty\cdot\Phi^\infty)_v \\
&\quad +\sum_{v\in\Diamondblack^{(1)}}E_{t_1,t_2}((g_1,g_2),\Phi^\lambda_\infty\cdot\Phi^{\infty,v}\cdot\Phi^\star_v) \\
&\quad +\sum_{v\in\Diamondblack^{(2)}}\sum_{x^\star_1\in V(v)_{t_1}}\omega_{1,\infty}(g_1)\(P[\phi^\infty_{t_2}(x^\star_1)]\cdot\phi^0_{V_\infty}\)(x^\star_1) \\
&\quad +\sum_{\bbu\in\Hom(E,\dC)}(\varphi^\KM_{t_2})_{t_1}(g_1),
\end{align*}
in which the terms $\Phi^\star_{v}$, $P[\phi^\infty_{t_2}(x^\star_1)]$, and $\varphi^\KM_{t_2}$ are from Proposition \ref{pr:height_curve}(1), Proposition \ref{pr:height_curve}(2), and Proposition \ref{pr:height_archimedean}, respectively. By Proposition \ref{pr:height_inert}(1) (plus Remark \ref{re:height_inert}), Proposition \ref{pr:height_split}, Proposition \ref{pr:height_curve} and Proposition \ref{pr:height_archimedean}, we have
\[
\fI_{t_1}(g_1)=\fK_{t_1}(g_1)
\]
for every $t_1\in(F^\times)^+$ and $g_1\in G_1(\bA_F^\Diamondblack)$. Thus, it is equivalent to showing that
\begin{align}\label{eq:modularity_3}
\sum_{t_1\in(F^\times)^+}\fK_{t_1}(g_1)
\end{align}
is absolutely convergent and extends to a function on $G_1(\bA_F)$ that is left invariant under $G_1(F)$. We analyze term by term.
\begin{enumerate}
  \item By Proposition \ref{pr:eisenstein_series}(2,3), we know that
      \[
      \sum_{t_1\in(F^\times)^+}\sum_{v\in\fV_F}\fE_{t_1,t_2}((g_1,g_2),\Phi^\lambda_\infty\cdot\Phi^\infty)_v
      \]
      is the (full) Fourier expansion of $E'(0,(g_1,g_2),\Phi^\lambda_\infty\cdot\Phi^\infty)_{t_2}$ (regarded as an automorphic form on $G_1(\bA_F)$ in the $g_1$-variable), hence is absolutely convergent and extends to a function on $G_1(\bA_F)$ that is left invariant under $G_1(F)$. On the other hand, by Lemma \ref{le:eisenstein_series}(2), we know that
      \[
      \sum_{t_1\in(F^\times)^+}\sum_{v\in\Diamondblack^{(1)}\cup\Diamondblack^{(2)}}
      \fE_{t_1,t_2}((g_1,g_2),\Phi^\lambda_\infty\cdot\Phi^\infty)_v
      \]
      is the (full) Fourier expansion of a doubling Eisenstein series (regarded as an automorphic form on $G_1(\bA_F)$ in the $g_1$-variable), hence is absolutely convergent and extends to a function on $G_1(\bA_F)$ that is left invariant under $G_1(F)$.

  \item For every $v\in\Diamondblack^{(1)}$,
      \[
      \sum_{t_1\in(F^\times)^+}E_{t_1,t_2}((g_1,g_2),\Phi^\lambda_\infty\cdot\Phi^{\infty,v}\cdot\Phi^\star_v)
      \]
      is the (full) Fourier expansion of $E(0,(g_1,g_2),\Phi^\lambda_\infty\cdot\Phi^{\infty,v}\cdot\Phi^\star_v)_{t_2}$ (regarded as an automorphic form on $G_1(\bA_F)$ in the $g_1$-variable), hence is absolutely convergent and extends to a function on $G_1(\bA_F)$ that is left invariant under $G_1(F)$.

  \item For every $v\in\Diamondblack^{(2)}$,
      \[
      \sum_{t_1\in(F^\times)^+}\sum_{x^\star_1\in V(v)_{t_1}}\omega_{1,\infty}(g_1)\(P[\phi^\infty_{t_2}(x^\star_1)]\cdot\phi^0_{V_\infty}\)(x^\star_1)
      \]
      is the (full) Fourier expansion of the theta series defined by the Schwartz function 
      \[
      x\mapsto \(P[\phi^\infty_{t_2}(x^\infty)]\cdot\phi^0_{V_\infty}\)(x_\infty),
      \]
      hence is absolutely convergent and extends to a function on $G_1(\bA_F)$ that is left invariant under $G_1(F)$.

  \item For every $\bbu\in\Hom(E,\dC)$,
      \[
      \sum_{t_1\in(F^\times)^+}(\varphi^\KM_{t_2})_{t_1}(g_1)
      \]
      is the (full) Fourier expansion of $\varphi^\KM_{t_2}$ by Proposition \ref{pr:height_archimedean}(b), hence is absolutely convergent and extends to a function on $G_1(\bA_F)$ that is left invariant under $G_1(F)$.
\end{enumerate}
Combining the four parts together, it follows that \eqref{eq:modularity_3} is absolutely convergent and extends to a function on $G_1(\bA_F)$ that is left invariant under $G_1(F)$.

The theorem is proved.
\end{proof}

\subsection{Central derivatives of $L$-functions}
\label{ss:derivative}

In this section, we consider a holomorphic tempered cuspidal automorphic representation $(\pi,\cA_\pi)$ of $G_1(\bA_F)$ of weight $\lambda$ (Definition \ref{de:hermitian_representation}, also Remark \ref{re:hermitian_siegel}).

Let $(\Diamond,L^\Diamond,\Diamondblack)$ be as in \S\ref{ss:modularity} such that $(\pi_v)^{K_{1,v}}\neq\{0\}$ (Notation \ref{no:hermitian_siegel}(1)) for every $v\in\fV_F^\fin\setminus\fV_F^{(\Diamond)}$. We define a functional
\[
\fZ\colon (\cA_\pi\cap\cA_{\lambda^+}(G_1)_\cusp^\hol)^2\times\sS((V^\infty)^2,\dC)^{L^\Diamond}_{(\Diamondblack)}
\to\dC
\]
(Definition \ref{de:schwartz_special}) by the formula
\[
\fZ(\ul\varphi,\Phi^\infty)\coloneqq
\int_{G_1(F)\backslash G_1(\bA_F)}\int_{G_1(F)\backslash G_1(\bA_F)}
\fI(\Phi^\infty)(g_1,g_2)\cdot\ol{\varphi_1(g_1)}\varphi_2(g_2)\rd g_1\rd g_2
\]
for $\ul\varphi=(\varphi_1,\varphi_2)\in(\cA_\pi\cap\cA_{\lambda^+}(G_1)_\cusp^\hol)^2$ and $\Phi^\infty\in\sS((V^\infty)^2,\dC)^{L^\Diamond}_{(\Diamondblack)}$, where $\fI$ is from Theorem \ref{th:modularity}. It is clear that $\fZ$ is an element in
\[
\Hom_{G_1(\bA_F^\infty)^2}\(\sS((V^\infty)^2,\dC)^{L^\Diamond}_{(\Diamondblack)},\pi^\infty\boxtimes(\pi^\infty)^\vee\).
\]

\begin{lem}\label{le:dichotomy}
We have
\[
\dim_\dC\Hom_{G_1(\bA_F^\infty)^2}\(\sS((V^\infty)^2,\dC)^{L^\Diamond}_{(\Diamondblack)},\pi^\infty\boxtimes(\pi^\infty)^\vee\)
\leq 1.
\]
Moreover, it equals one if and only if $V=V(\pi)$ (Remark \ref{re:dichotomy}(2)).
\end{lem}

\begin{proof}
Since the space $\sS((V^\infty)^2,\dC)^{L^\Diamond}$ admits a $G_1(\bA_F^\infty)^2$-invariant inner product, $\sS((V^\infty)^2,\dC)^{L^\Diamond}_{(\Diamondblack)}$ is a $G_1(\bA_F^\infty)^2$-equivariant direct summand of $\sS((V^\infty)^2,\dC)^{L^\Diamond}$. It follows from \cite{Liu12}*{Proposition~6.7} or \cite{LL21}*{Proposition~3.6} that
\[
\dim_\dC\Hom_{G_1(\bA_F^\infty)^2}\(\sS((V^\infty)^2,\dC)^{L^\Diamond}_{(\Diamondblack)},\pi^\infty\boxtimes(\pi^\infty)^\vee\)
\leq\dim_\dC\Hom_{G_1(\bA_F^\infty)^2}\(\sS((V^\infty)^2,\dC)^{L^\Diamond},\pi^\infty\boxtimes(\pi^\infty)^\vee\)\leq 1,
\]
and it equals zero if $V\neq V(\pi)$. When $V=V(\pi)$, the above Hom space is nonzero by Lemma \ref{le:abundance}.

The lemma is proved.
\end{proof}

\begin{theorem}\label{th:maipf}
Let $(\pi,\cA_\pi)$ be a holomorphic tempered cuspidal automorphic representation of $G_1(\bA_F)$ of weight $\lambda$. For every pair $\ul\varphi=(\varphi_1,\varphi_2)$ of automorphic forms in $\cA_\pi\cap\cA_{\lambda^+}(G_1)_\cusp^\hol$ and every $\Phi^\infty\in\sS((V^\infty)^2,\dC)^{L^\Diamond}_{(\Diamondblack)}$,
\[
\fZ(\ul\varphi,\Phi^\infty)=\frac{L'(\tfrac{1}{2},\pi)}{b_2(0)}\cdot \(-\frac{\pi}{4}\)^{[F:\dQ]}\cdot
\fZ^\natural_\fin(\ul\varphi,\Phi^\infty)
\]
(Definition \ref{de:doubling_zeta}) holds.
\end{theorem}

\begin{proof}
By Lemma \ref{le:dichotomy}, it suffices to consider the case where $V=V(\pi)$ since otherwise both sides vanish. In particular, we may assume $\epsilon(\pi)=-1$ (Remark \ref{re:dichotomy}(2)).

Similar to the proof of Theorem \ref{th:modularity}, it suffices to consider the case for $\Phi^\infty=\phi^\infty_1\otimes\ol{\phi^\infty_2}\in\sS((V^\infty)^2,\dC)^{L^\Diamond}_\Diamondblack$ with $\phi^\infty_1$ and $\phi^\infty_2$ in $\sS(V^\infty,\dC)^{L^\Diamond}$.

Take elements $t_2\in (F^\times)^+$ and $g_2\in G_1(\bA_F^\Diamondblack)$. Put $\fI_{t_2,g_2}(\Phi^\infty)(g_1)\coloneqq\fI(\Phi^\infty)_{t_2}(g_1,g_2)$ as a function in $g_1\in G_1(\bA_F)$, so that we have
\[
\fI_{t_2,g_2}(\Phi^\infty)(g_1)=\sum_{t_1\in (F^\times)^+}\fI_{t_2,g_2}(\Phi^\infty)_{t_1}(g_1).
\]
As in the proof of Theorem \ref{th:modularity}, we denote by $\Diamondblack^{(1)}$ and $\Diamondblack^{(2)}$ the subsets of $\Diamondblack$ consisting of $v$ that is nonsplit in $E$ such that $V_v$ is split and nonsplit, respectively. Take an element $v\in\Diamondblack^{(2)}$, by the theory of theta series, there exists a unique element $\theta^{(v)}_{t_2,g_2}(\Phi^\infty)\in\cA_{\lambda^+}(G_1)_\cusp^\hol$ (which is a theta series) such that for $g_1\in G_1(\bA_F^\Diamondblack)$,
\[
\theta^{(v)}_{t_2,g_2}(\Phi^\infty)(g_1)=
\sum_{x^\star_1\in V(v)}\omega_{1,\infty}(g_{1,\infty})
\(P[\omega_1^\infty(g_1^\infty)\phi^\infty_{t_2}(x^\star_1)]\cdot\phi^0_{V_\infty}\)(x^\star_1)
\]
from Proposition \ref{pr:height_curve_nonsplit}.\footnote{It is a cusp form since $\phi^\infty_{t_2}(0)=0$ when $\lambda$ is nontrivial.}

Then by Proposition \ref{pr:height_inert}(1) (plus Remark \ref{re:height_inert}), Proposition \ref{pr:height_split}, Proposition \ref{pr:height_curve}(1), Proposition \ref{pr:height_curve_nonsplit} and Proposition \ref{pr:height_archimedean_bis}, we have
\begin{align*}
&\quad\int_{G_1(F)\backslash G_1(\bA_F)}\ol{\varphi_1(g_1)}\fI_{t_2,g_2}(\Phi^\infty)(g_1)\rd g_1 \\
&=\int_{G_1(F)\backslash G_1(\bA_F)}\ol{\varphi_1(g_1)}
\sum_{t_1\in(F^\times)^+}\sum_{v\in\fV_F\setminus(\Diamondblack^{(1)}\cup\Diamondblack^{(2)})}
\fE_{t_1,t_2}((g_1,g_2),\Phi^\lambda_\infty\cdot\Phi^\infty)_v \rd g_1\\
&\quad +\int_{G_1(F)\backslash G_1(\bA_F)}\ol{\varphi_1(g_1)}
\sum_{t_1\in(F^\times)^+}
\sum_{v\in\Diamondblack^{(1)}}E_{t_1,t_2}((g_1,g_2),\Phi^\lambda_\infty\cdot\Phi^{\infty,v}\cdot\Phi^\star_v)\rd g_1 \\
&\quad +\int_{G_1(F)\backslash G_1(\bA_F)}\ol{\varphi_1(g_1)}\sum_{v\in\Diamondblack^{(2)}}\theta^{(v)}_{t_2,g_2}(\Phi^\infty)(g_1)\rd g_1,
\end{align*}
in which the terms $\Phi^\star_v$ is from Proposition \ref{pr:height_curve}(1).

We have three sorts of vanishing:
\begin{enumerate}
  \item first, by Lemma \ref{le:eisenstein_series}(2), for every $v\in \Diamondblack^{(1)}\cup\Diamondblack^{(2)}$, there exists a regularly supported Schwartz function $\Phi^\star_v\in\sS((V^\star_v)^2)$ such that
      \[
      \sum_{t_1\in(F^\times)^+}\fE_{t_1,t_2}((g_1,g_2),\Phi^\lambda_\infty\cdot\Phi^\infty)_v=
      \sum_{t_1\in(F^\times)^+}E_{t_1,t_2}((g_1,g_2),\Phi^\lambda_\infty\cdot\Phi^{\infty,v}\cdot\Phi^\star_v),
      \]
      which implies that
      \begin{align*}
      &\quad\int_{G_1(F)\backslash G_1(\bA_F)}\ol{\varphi_1(g_1)}\sum_{t_1\in(F^\times)^+}
      \sum_{v\in\Diamondblack^{(1)}\cup\Diamondblack^{(2)}}
      \fE_{t_1,t_2}((g_1,g_2),\Phi^\lambda_\infty\cdot\Phi^\infty)_v \rd g_1\\
      &=\sum_{v\in\Diamondblack^{(1)}\cup\Diamondblack^{(2)}}
      \int_{G_1(F)\backslash G_1(\bA_F)}\ol{\varphi_1(g_1)}
      \sum_{t_1\in(F^\times)^+}\fE_{t_1,t_2}((g_1,g_2),\Phi^\lambda_\infty\cdot\Phi^\infty)_v \rd g_1\\
      &=\sum_{v\in\Diamondblack^{(1)}\cup\Diamondblack^{(2)}}
      \int_{G_1(F)\backslash G_1(\bA_F)}\ol{\varphi_1(g_1)}\cdot
      E(0,(g_1,g_2),\Phi^\lambda_\infty\cdot\Phi^{\infty,v}\cdot\Phi^\star_v)_{t_2}\rd g_1,
      \end{align*}
      in which every integration in the last summation vanishes as $\epsilon(\pi)=-1$;

  \item second,
      \begin{align*}
      &\quad\int_{G_1(F)\backslash G_1(\bA_F)}\ol{\varphi_1(g_1)}
      \sum_{t_1\in(F^\times)^+}
      \sum_{v\in\Diamondblack^{(1)}}E_{t_1,t_2}((g_1,g_2),\Phi^\lambda_\infty\cdot\Phi^{\infty,v}\cdot\Phi^\star_v)\rd g_1 \\
      &=\int_{G_1(F)\backslash G_1(\bA_F)}\ol{\varphi_1(g_1)}
      \sum_{v\in\Diamondblack^{(1)}}\sum_{t_1\in(F^\times)^+}
      E_{t_1,t_2}((g_1,g_2),\Phi^\lambda_\infty\cdot\Phi^{\infty,v}\cdot\Phi^\star_v)\rd g_1 \\
      &=\sum_{v\in\Diamondblack^{(1)}}\int_{G_1(F)\backslash G_1(\bA_F)}\ol{\varphi_1(g_1)}\cdot
      E(0,(g_1,g_2),\Phi^\lambda_\infty\cdot\Phi^{\infty,v}\cdot\Phi^\star_v)_{t_2}\rd g_1,
      \end{align*}
      in which every integration in the last summation vanishes as $\epsilon(\pi)=-1$;

  \item third,
      \begin{align*}
      \int_{G_1(F)\backslash G_1(\bA_F)}\ol{\varphi_1(g_1)}\sum_{v\in\Diamondblack^{(2)}}\theta^{(v)}_{t_2,g_2}(\Phi^\infty)(g_1)\rd g_1
      =\sum_{v\in\Diamondblack^{(2)}}\int_{G_1(F)\backslash G_1(\bA_F)}
      \ol{\varphi_1(g_1)}\cdot\theta^{(v)}_{t_2,g_2}(\Phi^\infty)(g_1)\rd g_1,
      \end{align*}
      in which every integration in the last summation vanishes as $\epsilon(\pi)=-1$.
\end{enumerate}
It follows that
\begin{align*}
&\quad\int_{G_1(F)\backslash G_1(\bA_F)}\ol{\varphi_1(g_1)}\fI_{t_2,g_2}(\Phi^\infty)(g_1)\rd g_1 \\
&=\int_{G_1(F)\backslash G_1(\bA_F)}\ol{\varphi_1(g_1)}
\sum_{t_1\in(F^\times)^+}\sum_{v\in\fV_F}\fE_{t_1,t_2}((g_1,g_2),\Phi^\lambda_\infty\cdot\Phi^\infty)_v \rd g_1.
\end{align*}
If we sum over all $t_2\in(F^+)^\times$ and take the convolution against $\varphi_2$, we obtain that
\begin{align*}
\fZ(\ul\varphi,\Phi^\infty)=\iint_{(G_1(F)\backslash G_1(\bA_F))^2}
\sum_{(t_1,t_2)\in((F^\times)^+)^2}\sum_{v\in\fV_F}\fE_{t_1,t_2}((g_1,g_2),\Phi^\lambda_\infty\cdot\Phi^\infty)_v
\cdot\ol{\varphi_1(g_1)}\varphi_2(g_2)\rd g_1\rd g_2,
\end{align*}
which by Proposition \ref{pr:eisenstein_series}(2,3) equals
\[
\iint_{(G_1(F)\backslash G_1(\bA_F))^2}E'(0,(g_1,g_2),\Phi^\lambda_\infty\cdot\Phi^\infty)
\cdot\ol{\varphi_1(g_1)}\varphi_2(g_2)\rd g_1\rd g_2.
\]
Then the theorem follows from Proposition \ref{pr:ripf} and Remark \ref{re:constant}(2).
\end{proof}

\begin{definition}\label{de:special_numerical}
For every intermediate field $\dE^\lambda\subseteq\dL\subseteq\dC$ and every object $L$ of $\Hk(H(\bA_F^\infty))$, we
\begin{enumerate}
  \item define $\SC_\lambda(X_L,\dL)$ to be the subspace of $\rH^2(X_L,\rM_{L,\lambda})\otimes_{\dE^\lambda}\dL$ spanned by motivic special classes $z^\lambda(x)_L$ for $x\in V^\infty$ satisfying $T(x)\in F^+$;

  \item denote by $\SC_\lambda^\circ(X_L,\dL)\subseteq\SC_\lambda(X_L,\dL)$ the subspace contained in the radical of the height pairing
      \[
      \llangle\;,\;\rrangle_{\lambda,L}\colon
      \rH^2(X_L,\rM_{L,\lambda})\otimes_{\dE^\lambda}\dL\times\rH^2(X_L,\rM_{L,\lambda})\otimes_{\dE^\lambda}\dL\to\dC
      \]
      extended to a hermitian form with coefficients in $\dL$;

  \item put
      \[
      \ol\SC_\lambda(X_L,\dL)\coloneqq\frac{\SC_\lambda(X_L,\dL)}{\SC_\lambda^\circ(X_L,\dL)},
      \]
      which receives a natural action by $\dL[L\backslash H(\bA_F^\infty)/L]$.
\end{enumerate}
\end{definition}

\begin{corollary}\label{co:derivative}
Let $(\pi,\cA_\pi)$ be a holomorphic tempered cuspidal automorphic representation of $G_1(\bA_F)$ of weight $\lambda$ such that $V(\pi)$ (Remark \ref{re:dichotomy}(2)) is incoherent (Definition \ref{de:space_coherence}) satisfying $|\fV(V(\pi))|>1$. Let $\dE^\lambda\subseteq\dL_\pi\subseteq\dC$ be the field of definition of $\pi^\infty$. If $L'(\tfrac{1}{2},\pi)\neq 0$, then for $V=V(\pi)$, every finite set $\Diamond$ of rational primes containing all those ramified in $E$ such that $(\pi_v)^{K_{1,v}}\neq\{0\}$ (Notation \ref{no:hermitian_siegel}(1)) for every $v\in\fV_F^\fin\setminus\fV_F^{(\Diamond)}$, and every open compact subgroup $L^\Diamond$ of $H(\bA^{\infty,\Diamond}_F)$ that is the stabilizer of a self-dual lattice of $V^{\infty,\Diamond}$,
\[
\varinjlim_{L_\Diamond}\ol\SC_\lambda(X_{L_\Diamond L^\Diamond},\dL_{\pi})_{\fm_\pi}\neq 0
\]
holds, where $\fm_\pi\coloneqq\Ker\chi_\pi\subseteq\dL_\pi[L^\Diamond\backslash H(\bA_F^{\infty,\Diamond})/L^\Diamond]$ from \eqref{eq:sakata_character}.
\end{corollary}

\begin{proof}
This follows from Theorem \ref{th:maipf} and Lemma \ref{le:abundance}.
\end{proof}

The corollary below confirms the (conjectural) motivic arithmetic inner product formula (Conjecture \ref{co:maipf}) for $r=1$, assuming the modularity of motivic generating functions (Conjecture \ref{co:modularity}).

\begin{corollary}\label{co:maipf_bis}
Let $(\pi,\cA_\pi)$ be a holomorphic tempered cuspidal automorphic representation of $G_1(\bA_F)$ of weight $\lambda$. Assume Conjecture \ref{co:modularity} for $V$ so that one can apply Definition \ref{de:theta_lift}. For every pair $\ul\varphi=(\varphi_1,\varphi_2)$ of automorphic forms in $\cA_\pi\cap\cA_{\lambda^+}(G_1)_\cusp^\hol$, and very pair $\ul\phi^\infty=(\phi^\infty_1,\phi^\infty_2)$ of Schwartz functions in $\sS(V^\infty)$,
\[
\llangle\Theta_{\phi^\infty_1}^\lambda(\varphi_1),\Theta_{\phi^\infty_2}^\lambda(\varphi_2)\rrangle_\lambda
=\frac{L'(\tfrac{1}{2},\pi)}{b_2(0)}\cdot \(-\frac{\pi}{4}\)^{[F:\dQ]}\cdot
\fZ^\natural_\fin(\ul\varphi,\phi^\infty_1\otimes\ol{\phi^\infty_2})
\]
holds.
\end{corollary}

\begin{proof}
By \cite{Liu12}*{Proposition~6.7} or \cite{LL21}*{Proposition~3.6}, it suffices to consider the case where $V=V(\pi)$ since otherwise both sides vanish; and in this case, it suffices to show the identity in the statement for a single choice of data with $\fZ^\natural_\fin(\ul\varphi,\phi^\infty_1\otimes\ol{\phi^\infty_2})\neq 0$. Then the corollary follows from Theorem \ref{th:maipf} and Lemma \ref{le:abundance}.
\end{proof}

\begin{remark}\label{re:general_cm}
In this remark, we explain how to modify our theory when $E$ is a general CM field.

We fix an \emph{abstract CM data} of $E$, which is a pair $(\widetilde{E}/E,\tr)$ in which
\begin{itemize}
  \item $\widetilde{E}/E$ is a finite extension with $\widetilde{E}$ also a CM field; and

  \item $\tr\colon E\to\widetilde{E}$ is a function such that for some hence every embedding $\tilde{u}\colon\widetilde{E}\to\dC$, $\tilde{u}\circ\tr$ is the trace function of a CM type of $E$.
\end{itemize}
Note that the case in the main text corresponds to the abstract CM data in which $\widetilde{E}=E$ and $\tr$ is the composition of $\tr_{E/E_0}$ and the inclusion $E_0\hookrightarrow E$.

We modify the functor $Y_\bullet$ in Definition \ref{de:moduli_auxiliary} to classify triples $(A_0,\lambda_0,\bar\eta_0)$ over $\widetilde{E}$ in which $A_0$ is an object of $\AV(S,E)$ (for $S\in\Sch'_{/\widetilde{E}}$) such that $\tr_{A_0}(e)=\tr(e)$ for $e\in E$; and the others are similar.\footnote{Indeed, we also need to replace $Y_\bullet$ by an $\bA_E^\infty$-orbit of it in case the Hasse principle fails for the group of rational similitudes of $E/F$; see \cite{LTXZZ}*{\S3.5} for more details.} We then modify the functor $X_{\bullet,\bullet}$ in Definition \ref{de:moduli_rational} to classify sextuples $(A_0,\lambda_0,\bar\eta_0;A,\lambda,\bar\eta)$ over $\widetilde{E}$ (rather than $E$) in which $A$ is an object of $\AV(S,E)$ (for $S\in\Sch'_{/\widetilde{E}}$) such that $\tr_A(e)=m\tr(e)-e+\ol{e}$ for $e\in E$; and the others are similar. By the similar descent process, we obtain a functor
\[
X_\bullet\colon\Hk(H(\bA_F^\infty))\to\Sch^\smq_{/\widetilde{E}}
\]
such that for every $\tilde{u}\in\Hom(\widetilde{E},\dC)$, $X_\bullet\otimes_{\widetilde{E}}\tilde{u}(\widetilde{E})\simeq \Sh(H(v))_\bullet^u\otimes_{u(E)}\tilde{u}(\widetilde{E})$ where $v\coloneqq\tilde{u}\res_F$ and $u\coloneqq\tilde{u}\res_E$, together with a Cartesian functor
\[
\rM_\bullet\colon\Hk(H(\bA_F^\infty))\to\Mot(\Sch^\smq_{/\widetilde{E}},E)
\]
above $X_\bullet$ as the \emph{standard motive}.

The remaining part is essentially same to the case in the main text, except that in the definition of the height pairing $\llangle\;,\;\rrangle_{\lambda,L}^{\ell}$ \eqref{eq:height_pairing} (or rather $\llangle\;,\;\rrangle_{\lambda,L}$ by Theorem \ref{th:independence}) we divide by $[\widetilde{E}:E]$ at the end.

Our theory a priori depends on the choice of the abstract CM data. For a general CM field $E$, our approach may not produce an optimal theory over the base field $E$ since it is not true that we can always find an abstract CM data with $\widetilde{E}=E$.
\end{remark}

\appendix

\section{Relative motives}
\label{ss:a}

In the appendix, we collect some known facts concerning relative motives.

\subsection{Relative motives and their realizations}
\label{ss:relative_motive}

Let $L$ be a $\dQ$-ring that is finite-dimensional as a $\dQ$-vector space. For every place $\diamond\in\{\infty,2,3,5,\dots\}$ of $\dQ$, put $L_\diamond\coloneqq L\otimes\dQ_\diamond$.

\begin{notation}\label{no:relative_motivic}
For every Noetherian scheme $S$, we have the following.
\begin{enumerate}
  \item Denote by $\Mot(S,\dQ)$ the category of relative (Chow) motives over $S$ with coefficients in $\dQ$ (see, for example, \cite{MNP13}*{Definition~8.1.4}\footnote{The setting that $S$ is quasi-projective over the complex field is not needed in this definition.}), which is a closed symmetric monoidal $\dQ$-linear pseudo-abelian category.\footnote{Recall that a pseudo-abelian category is an additive category in which every projector has a kernel.} Denote by $\Mot(S,L)$ the category of $L$-linear objects in $\Mot(S,\dQ)$, which is a closed symmetric monoidal $L$-linear pseudo-abelian category.

  \item Denote by $\sfD(S,L)$ the category of constructible \'{e}tale motivic sheaves over $S$ with coefficients in $L$ (see \cite{CD16} for the theory), which is a closed symmetric monoidal $L$-linear triangulated category whose unit element we denote by $\mathbbm{1}_{S,L}$.\footnote{By \cite{CD16}*{Corollaries~5.5.5~\&~5.5.7}, $\sfD(S,L)$ remains canonically the same whether we defined it using \'{e}tale or Nisnevich topology, with or without transfers.}

  \item For every object $\rM$ in $\sfD(S,L)$ and integers $i,j$, we have the \emph{motivic cohomology}
      \[
      \rH^i(S,\rM(j))\coloneqq\Hom_{\sfD(S,L)}(\mathbbm{1}_{S,L},\rM[i](j))
      \]
      as an $L$-module, and for every closed subset $T\subseteq S$ the \emph{motivic cohomology with support}
      \[
      \rH^i_T(S,\rM(j))\coloneqq\rH^{i+1}(S,\r{cone}(\rM\to\rR a_* a^*\rM)(j)),
      \]
      where $a\colon S\setminus T\to S$ denotes the open complement.
\end{enumerate}
\end{notation}

\begin{remark}\label{re:relative_motivic}
Suppose that $S$ is a regular Noetherian scheme.
\begin{enumerate}
  \item For every smooth projective morphism $f\colon X\to S$, we have a canonical isomorphism
      \[
      \rH^{2r}(S,\rR f_*\mathbbm{1}_{X,L}(r))\simeq\CH^r(X)\otimes L
      \]
      by \cite{Deg18}*{Remark~1.3.5(4)} (and the fact that the motivic absolutely purity is known for coefficients in $\dQ$ hence $L$).\footnote{When $S$ is the spectrum of a field, this is a well-known result proved in \cite{Voe02}.}

  \item Part (1) implies that the assignment sending $(X,\id)$ for every smooth projective morphism $f\colon X\to S$ to $\rR f_*\mathbbm{1}_{X,L}$ induces a fully faithful $L$-linear functor $\Mot(S,L)\to\sfD(S,L)$.
\end{enumerate}
\end{remark}

\begin{notation}\label{no:relative_motive_etale}
For every Noetherian scheme $S$ and every rational prime $\ell$ that is invertible on $S$, we have the following.
\begin{enumerate}
  \item Put $L_\ell\coloneqq L\otimes\dQ_\ell$ and denote by $\sfD_{\et}(S,L_\ell)$ the derived category of $L_\ell$-linear \'{e}tale sheaves over $S$, which is a closed symmetric monoidal $L_\ell$-linear triangulated category whose unit element we denote by $\mathbbm{1}^\et_{S,L_\ell}$. Denote by $\sfD_{\et}^\heartsuit(S,L_\ell)$ the heart of $\sfD_{\et}(S,L_\ell)$ with respect to the usual t-structure.

  \item For every object $\rM$ in $\sfD_{\et}(S,L_\ell)$ and integers $i,j$, we have the usual \'{e}tale cohomology
      \[
      \rH^i_{\et}(S,\rM(j))\coloneqq\Hom_{\sfD_{\et}(S,L_\ell)}(\mathbbm{1}^\et_{S,L_\ell},\rM[i](j))
      \]
      as an $L_\ell$-module, and for every closed subset $T\subseteq S$ the \'{e}tale cohomology with support
      \[
      \rH^i_{\et,T}(S,\rM(j))\coloneqq\rH^{i+1}_{\et}(S,\r{cone}(\rM\to\rR a_* a^*\rM)(j)),
      \]
      where $a\colon S\setminus T\to S$ denotes the open complement.

  \item We have the \emph{\'{e}tale realization functor}
      \[
      r_\ell\colon\sfD(S,L)\to\sfD_{\et}(S,L_\ell),
      \]
      compatible with the six operations \cite{Ayo14}. In particular, for every morphism $f\colon X\to S$ of finite type, $r_\ell\rR f_*\mathbbm{1}_{X,L}=\rR f_*\mathbbm{1}^\et_{X,L_\ell}$.
\end{enumerate}
\end{notation}

\begin{notation}\label{no:relative_motive_hodge}
For every complex manifold $S$, we have the following.
\begin{enumerate}
  \item Put $L_\infty\coloneqq L\otimes\dR$ and denote by $\sfD_\hs(S,L_\infty)$
      the derived category of $L_\infty$-linear objects in the category of (real) mixed Hodge modules over $S$ \cite{Sai90}, which is a closed symmetric monoidal $L_\infty$-linear triangulated category whose unit element we denote by $\mathbbm{1}^\hs_{S,L_\infty}$. Denote by $\sfD_\hs^\heartsuit(S,L_\infty)$ the heart of $\sfD_{\et}(S,L_\infty)$ with respect to the usual t-structure.

  \item For every object $\rM$ in $\sfD_\hs(S,L_\infty)$ and integers $i,j$, we have the absolute Hodge cohomology \footnote{This can be computed in the way described in \cite{Bei86}.}
      \[
      \rH^i_\hs(S,\rM(j))\coloneqq\Hom_{\sfD_\hs(S,L_\infty)}(\mathbbm{1}^\hs_{S,L_\infty},\rM[i](j))
      \]
      as an $L_\infty$-module, and for every closed subset $T\subseteq S$ the absolute Hodge cohomology with support
      \[
      \rH^i_{\hs,T}(S,\rM(j))\coloneqq\rH^{i+1}_{\hs}(S,\r{cone}(\rM\to\rR a_* a^*\rM)(j)),
      \]
      where $a\colon S\setminus T\to S$ denotes the open complement.

  \item Consider a field $k$ together with an embedding $u\colon k\to \dC$. For every smooth quasi-projective scheme $S$ over $k$, we have the \emph{Hodge realization functor}
      \[
      r^u\colon\sfD(S,L)\to\sfD_\hs(S^u,L_\infty),
      \]
      compatible with the six operations \cites{Ivo16,Tub25}. In particular, for every smooth quasi-projective morphism $f\colon X\to S$, $r^u\rR f_*\mathbbm{1}_{X,L}=\rR f^u_*\mathbbm{1}^\hs_{X^u,L_\infty}$.
\end{enumerate}
\end{notation}

\subsection{Motives of abelian schemes}
\label{ss:abelian_motive}

\begin{definition}\label{de:category_av}
Consider a ring $R$ and an $R$-scheme $T$.
\begin{itemize}
  \item We have a natural functor $\AV(-,R)\colon(\Sch_{/T})^\op\to\Cat$ such that for every $T$-scheme $S$, $\AV(S,R)$ is the category of abelian schemes $A$ over $T$ together with a unit-preserving algebra homomorphism
      \[
      R\to\End_S(A)\otimes\dZ[p^{-1}\res p\in R^\times].
      \]

  \item For $S\in\Sch_{/T}$ and $A\in \AV(S,R)$, write $\tr_A(r)\coloneqq\tr(r\res\Lie_SA)\in\sO(S)$ for $r\in R$.
\end{itemize}
\end{definition}

Let $S$ be a regular Noetherian scheme. Consider an abelian scheme $A/S$ of relative dimension $g$. We have the corresponding cohomological and homological relative motives
\[
\sfh^*(A/S,\dQ),\quad\sfh_*(A/S,\dQ)
\]
in $\Mot(S,\dQ)$, together with the evaluation map (Definition \ref{de:dualizable})
\begin{align}\label{eq:abelian_duality}
\sfh^*(A/S,\dQ)\otimes\sfh_*(A/S,\dQ)\to\mathbbm{1}_{S,\dQ}.
\end{align}
By \cite{DM91}*{Corollary~3.2} (see Remark \ref{re:motivie_abelian} below), there is a canonical decomposition\footnote{Our motive $\sfh^*(A/S,\dQ)$ is denoted as $R(A/S)$ in \cite{DM91}.}
\[
\sfh^*(A/S,\dQ)=\bigoplus_{i=0}^{2g}\sfh^i(A/S,\dQ)
\]
in the category $\Mot(S,\dQ)$ such that $n^*$ acts on $\sfh^i(A/S,\dQ)$ by multiplication with $n^i$ for every $n\in\dZ$. It follows by \eqref{eq:abelian_duality} that there is another canonical decomposition
\[
\sfh_*(A/S,\dQ)=\bigoplus_{i=0}^{2g}\sfh_i(A/S,\dQ)
\]
in the category $\Mot(S,\dQ)$ such that $n_*$ acts on $\sfh_i(A/S,\dQ)$ by multiplication with $n^i$ for every $n\in\dZ$. By \cite{DM91}*{Proposition~2.11} (and the following remark), the (relative) Fourier transform $\tF_{A/S}$ induces an isomorphism
\begin{align}\label{eq:fourier}
\tF_{A/S}^i\colon\sfh^{2g-i}(A/S,\dQ)\xrightarrow\sim\sfh^{i}(A^\vee/S,\dQ)\ang{i-g}
\end{align}
in the category $\Mot(S,\dQ)$ for every $0\leq i\leq 2g$.

\begin{remark}\label{re:motivie_abelian}
The article \cite{DM91} only considered the case where $S$ is smooth quasi-projective over a base field. However, as pointed out in \cite{Kun94}*{Remark~1.1}, all results in \cite{DM91} remain true as long as the Chow ring (with coefficients in $\dQ$) is defined for every smooth projective scheme over $S$, which is case when $S$ is a regular Noetherian scheme by Remark \ref{re:relative_motivic}(1).
\end{remark}

We define a map
\begin{align}\label{eq:poincare}
\sfh^1(A/S,\dQ)\otimes\sfh^1(A^\vee/S,\dQ)\to\mathbbm{1}_{S,\dQ}\ang{-1}
\end{align}
to be the composition of the isomorphism
\[
1\otimes(\tF_{A/S}^1)^{-1}\colon\sfh^1(A/S,\dQ)\otimes\sfh^1(A^\vee/S,\dQ)
\xrightarrow\sim\sfh^1(A/S,\dQ)\otimes\sfh^{2g-1}(A/S,\dQ)\ang{g-1}
\]
and the Poincar\'{e} duality map
\[
\sfh^1(A/S,\dQ)\otimes\sfh^{2g-1}(A/S,\dQ)\to\mathbbm{1}_{S,\dQ}\ang{-g}.
\]

\begin{lem}\label{le:poincare}
The map \eqref{eq:poincare} induces an isomorphism $\sfh^1(A^\vee/S,\dQ)\xrightarrow\sim\sfh^1(A/S,\dQ)^\vee\ang{-1}$.
\end{lem}

\begin{proof}
It suffices to show that the induced map
\[
\sfh^1(A/S,\dQ)\otimes\sfh^1(A^\vee/S,\dQ)\ang{1}\to\mathbbm{1}_{S,\dQ}
\]
is an evaluation map. By construction, this map is isomorphic to the map
\[
\sfh^1(A/S,\dQ)\otimes\sfh^{2g-1}(A/S,\dQ)\ang{g}\to\mathbbm{1}_{S,\dQ}
\]
induced from the Poincar\'{e} duality map, which is an evaluation map.
\end{proof}

\begin{remark}\label{re:abelian_realization}
In this remark, we identify the \'{e}tale and Hodge realizations of $\sfh^*(A/S,\dQ)$ and $\sfh_*(A/S,\dQ)$. Write $f\colon A\to S$ for the structure morphism.
\begin{enumerate}
  \item Take a rational prime $\ell$ that is invertible on $S$. By the compatibility of six operations for the \'{e}tale realization (Notation \ref{no:relative_motive_etale}), we have
      \[
      r_\ell\sfh^*(A/S,\dQ)=\rR f_*\mathbbm{1}^\et_{A,\dQ_\ell}.
      \]
      For the right-hand side, the corresponding decomposition is
      \[
      \rR f_*\mathbbm{1}^\et_{A,\dQ_\ell}=\bigoplus_{i=0}^{2g}\rR^i f_*\mathbbm{1}^\et_{A,\dQ_\ell}[-i].
      \]
      It follows that
      \[
      r_\ell\sfh^i(A/S,\dQ)=\sfh^i_{\et}(A/S,\dQ_\ell)\coloneqq\rR^i f_*\mathbbm{1}^\et_{A,\dQ_\ell}[-i]
      \]
      for every $0\leq  i\leq 2g$. See the remark after \cite{DM91}*{Corollary~3.2}. Dually, we have
      \[
      r_\ell\sfh_i(A/S,\dQ)=\sfh_i^{\et}(A/S,\dQ_\ell)\coloneqq(\rR^i f_*\mathbbm{1}^\et_{A,\dQ_\ell})^\vee[i]
      \]
      for every $0\leq  i\leq 2g$.

  \item Suppose that $S$ is smooth quasi-projective over a field $k$ with an embedding $u\colon k\to \dC$. By the compatibility of six operations for the Hodge realization (Notation \ref{no:relative_motive_hodge}), we have
      \[
      r_\ell\sfh^*(A/S,\dQ)=\rR f^u_*\mathbbm{1}^\hs_{A^u,\dR}.
      \]
      For the right-hand side, the corresponding decomposition is
      \[
      \rR f^u_*\mathbbm{1}^\hs_{A^u,\dR}=\bigoplus_{i=0}^{2g}\rR^i f^u_*\mathbbm{1}^\hs_{A^u,\dR}[-i].
      \]
      It follows that
      \[
      r_\ell\sfh^i(A/S,\dQ)=\sfh^i_{\hs}(A/S,\dR)\coloneqq\rR^i f^u_*\mathbbm{1}^\hs_{A^u,\dR}[-i]
      \]
      for every $0\leq  i\leq 2g$. Dually, we have
      \[
      r_\ell\sfh_i(A/S,\dQ)=\sfh_i^{\hs}(A/S,\dR)\coloneqq(\rR^i f^u_*\mathbbm{1}^\hs_{A^u,\dR})^\vee[i]
      \]
      for every $0\leq  i\leq 2g$.
\end{enumerate}

\end{remark}

\end{document}